\documentclass{book}

\input{diagrams}
\diagramstyle[PostScript=dvips,nohug]
\usepackage{epsfig}

\newtheorem{thm}{Theorem}[section]
\newtheorem{defn}[thm]{Definition}
\newtheorem{propn}[thm]{Proposition}

\newtheorem{cor}[thm]{Corollary}
\newtheorem{constn}[thm]{Construction}

\newtheorem{lotsofremarks}[thm]{Remarks}
\newenvironment{remarks}[1]
	{\begin{lotsofremarks}\label{#1}\end{lotsofremarks}\begin{enumerate}}
	{\end{enumerate}}

\newtheorem{lotsofegs}[thm]{Examples}
\newenvironment{eg}[1]
	{\begin{lotsofegs}\label{#1}\end{lotsofegs}\begin{enumerate}}
	{\end{enumerate}}

\newcommand{\mcm}[3]{\newcommand{#1}[#2]{{\ensuremath{#3}}}}

\usepackage{latexsym}
\usepackage{amssymb}

\mcm{\blank}{0}{(\emptybk)}
\mcm{\dashbk}{0}{\mbox{---}}
\mcm{\emptybk}{0}{\:\:}
\mcm{\hyph}{0}{\mbox{-}}

\mcm{\diagspace}{0}{\mbox{\hspace{2em}}}

\newcommand{\done}{\hfill\ensuremath{\Box}}

\newcommand{\bref}[1]{(\ref{#1})}
\newcommand{\ucontents}[2]{\addcontentsline{toc}{#1}{\numberline{}{#2}}}

\mcm{\cat}{1}{\mc{#1}}
\mcm{\fcat}{1}{\mb{#1}}
\mcm{\mc}{1}{\mathcal{#1}}
\mcm{\mr}{1}{\mathrm{#1}}
\mcm{\mi}{1}{\mathit{#1}}
\mcm{\mb}{1}{\mathbf{#1}}
\mcm{\scat}{1}{\Bbb{#1}}
\mcm{\twid}{1}{\widetilde{#1}}
\newcommand{\url}[1]{\mbox{\tt #1 }}

\mcm{\elt}{0}{\in}
\mcm{\sub}{0}{\,\subseteq\,}
\mcm{\such}{0}{\:|\:}
\mcm{\without}{0}{\setminus}

\mcm{\atsr}{0}{\Box}
\mcm{\eqv}{0}{\,\simeq\,}
\mcm{\iso}{0}{\,\cong\,}
\mcm{\of}{0}{\raisebox{0.2mm}{\ensuremath{\scriptstyle\circ}}}

\mcm{\bdry}{0}{\partial}

\mcm{\Bee}{0}{\cat{B}}
\mcm{\Beep}{0}{\cat{B'}}
\mcm{\Eee}{0}{\cat{E}}
\mcm{\Eeep}{0}{\cat{E'}}
\newcommand{\epsln}{\varepsilon}
\mcm{\Ess}{0}{\cat{S}}
\mcm{\Tee}{0}{\cat{T}}
\mcm{\Teep}{0}{\cat{T'}}
\mcm{\Stee}{0}{\scat{T}}
\mcm{\Steep}{0}{\scat{T'}}

\mcm{\blbk}{0}{\blank^{\blob}}
\mcm{\blob}{0}{\scriptscriptstyle{\bullet}}
\mcm{\stbk}{0}{\blank^{*}}
\mcm{\ubl}{0}{{}^{\blob}}
\mcm{\ust}{0}{{}^{*}}

\mcm{\Cartpr}{0}{\pr{\Eee}{T}}
\mcm{\Cartprp}{0}{\pr{\Eeep}{T'}}
\mcm{\Mnd}{0}{\triple{T}{\eta}{\mu}}
\mcm{\Zeropr}{0}{\pr{\Set}{\id}}

\mcm{\dopset}{0}{\ftrcat{\Delta^{\op}}{\Set}}
\mcm{\tropset}{0}{\ftrcat{\fcat{TR}^{\op}}{\Set}}

\mcm{\cod}{0}{\mr{cod}}
\mcm{\dom}{0}{\mr{dom}}
\mcm{\End}{0}{\mr{End}}
\mcm{\Hom}{0}{\mr{Hom}}
\mcm{\ob}{0}{\mr{ob}\,}
\mcm{\op}{0}{\mr{op}}

\mcm{\comp}{0}{\mi{comp}}
\mcm{\id}{0}{\mi{id}}
\mcm{\ids}{0}{\mi{ids}}
\mcm{\mult}{0}{\mi{mult}}
\mcm{\unit}{0}{\mi{unit}}

\mcm{\Ab}{0}{\fcat{Ab}}
\mcm{\Alg}{0}{\fcat{Alg}}
\mcm{\Bim}{1}{\fcat{Bim}(#1)}
\mcm{\Cat}{0}{\fcat{Cat}}
\mcm{\Cay}{0}{\fcat{Cay}}
\mcm{\Cpn}{1}{\pr{\Set/S_{#1}}{T_{#1}}}
\mcm{\fc}{0}{\fcat{fc}}
\mcm{\fm}{0}{\fcat{fm}}
\mcm{\Graph}{0}{\fcat{Graph}}
\mcm{\Gy}{0}{\fcat{Gy}}
\mcm{\Hpn}{1}{\pr{\Eee_{#1}}{P_{#1}}}
\mcm{\Mon}{0}{\mb{Mon}}
\mcm{\Multicat}{0}{\fcat{Multicat}}
\mcm{\One}{0}{\fcat{1}}
\mcm{\PD}{1}{\fcat{PD}_{#1}}
\mcm{\Prof}{0}{\fcat{Prof}}
\mcm{\Set}{0}{\fcat{Set}}
\mcm{\Span}{0}{\fcat{Span}}
\mcm{\Ssq}{0}{\fcat{Ssq}}
\mcm{\Struc}{0}{\fcat{Struc}}
\mcm{\Sym}{0}{\fcat{Sym}}
\mcm{\TR}{1}{\fcat{TR}(#1)}
\mcm{\Tr}{0}{\fcat{Tr}}
\mcm{\Twocat}{0}{\fcat{2\hyph\Cat}}

\mcm{\integers}{0}{\mathbb{Z}}

\mcm{\range}{2}{#1,\,\ldots\,,#2}
\mcm{\bftuple}{2}{\tuplebts{\range{#1}{#2}}}
\mcm{\tuple}{3}{\tuplebts{\range{#1,#2}{#3}}}
\mcm{\rttuple}{1}{\tuplebts{\,\ldots\,,#1}}
\mcm{\abftuple}{2}{\atuplebts{\range{#1}{#2}}}
\mcm{\atuple}{3}{\atuplebts{\range{#1,#2}{#3}}}
\mcm{\arttuple}{1}{\atuplebts{\,\ldots\,,#1}}
\mcm{\sqbftuple}{2}{\obt\range{#1}{#2}\cbt}
\mcm{\pr}{2}{\tuplebts{#1,#2}}
\mcm{\triple}{3}{\tuplebts{#1,#2,#3}}

\mcm{\eend}{2}{#1[#2]}
\mcm{\ehom}{3}{#1[#2,#3]}
\mcm{\ftrcat}{2}{[#1,#2]}
\mcm{\homset}{3}{#1(#2,#3)}
\mcm{\multihom}{3}{#1(#2;#3)} 
\mcm{\relhom}{5}{#1_{#2}(\range{#3}{#4};#5)}

\mcm{\go}{0}{\rTo}
\mcm{\goby}{1}{\rTo^{#1}}
\mcm{\goesto}{0}{\,\longmapsto\,}
\mcm{\goiso}{0}{\goby{\diso}}
\mcm{\monic}{0}{\rMonic}
\mcm{\og}{0}{\lTo}
\mcm{\ogby}{1}{\lTo^{#1}}

\mcm{\gph}{2}{\spn{#1}{T #2}{#2}}
\mcm{\graph}{4}{\spaan{#1}{T #2}{#2}{#3}{#4}}
\mcm{\oppair}{2}{\stackrel{\rTo^{#1}}{\lTo_{#2}}}
\mcm{\parpair}{2}{\stackrel{\rTo^{#1}}{\rTo_{#2}}}
\mcm{\spn}{3}{#2 \og #1 \go #3}
\mcm{\spaan}{5}{#2 \ogby{#4} #1 \goby{#5} #3}

\mcm{\bktdvslob}{3}
	{\left(
	\begin{diagram}[height=1.5em]
	#1		\\
	\dTo>{\,#2}	\\
	#3		\\
	\end{diagram}
	\right)}
\mcm{\slob}{3}{(#1 \goby{#2} #3)}
\mcm{\vslob}{3}
	{\left.
	\begin{diagram}[height=1.5em]
	#1		\\
	\dTo>{\,#2}	\\
	#3		\\
	\end{diagram}
	\right.}

\newarrow{Edge}---->		% for opetopes
\newarrow{Equals}=====
\newarrow{Get}....>
\newarrow{Goesto}|---{->}
\newarrow{Mod}--+->
\newarrow{Monic}{vee}---{vee}
\newenvironment{opetope}
	{\begin{diagram}[size=1em,abut,tight,noPS]}
	{\end{diagram}}
\newenvironment{slopeydiag}
	{\begin{diagram}[size=2em]}
	{\end{diagram}}
\newenvironment{tree}
	{\begin{diagram}[height=1em,width=.75em,abut,noPS,tight]}	
	{\end{diagram}}
\newenvironment{triangdiag}
	{\begin{diagram}[width=1em,height=1.5em]}
	{\end{diagram}}

\newcommand{\dn}{\dLine}
\mcm{\enode}{0}{\circ}
\newcommand{\lt}[1]{\ldLine(#1,2)}
\mcm{\nl}{1}{\stackrel{\textstyle #1}{\node}}
\mcm{\node}{0}{\bullet}
\newcommand{\rt}[1]{\rdLine(#1,2)}
\mcm{\utree}{0}{\node}

\newcommand{\cnr}{}	% for opetopes
\mcm{\diso}{0}{\sim}
\newcommand{\pullshape}
	{\setlength{\unitlength}{1em}
	\begin{picture}(2,5)(-1,-5)
	\put(0,-5){\line(1,1){1}}
	\put(0,-5){\line(-1,1){1}}
	\end{picture}}
\newcommand{\Spbk}{\overprint{\raisebox{-2.5em}{\pullshape}}}
\mcm{\vdiso}{0}{\wr}

\mcm{\nat}{0}{\mathbb{N}}	

\newcommand{\piccy}[1]{\epsfig{file=#1}}

\mcm{\Onepr}{0}{\pr{\Graph}{\fc}}
\newlength{\nllwidth}
\newlength{\nllheight}
\newcommand{\stackbelow}[2]{%
\settowidth{\nllwidth}{\ensuremath{#1}\ensuremath{#2}}%
\settoheight{\nllheight}{\ensuremath{#2}}%
\addtolength{\nllheight}{.3ex}%
\mbox{%
\ensuremath{#1}%
\hspace{-.5\nllwidth}%
\raisebox{-1\nllheight}{\ensuremath{#2}}}}
\mcm{\nlal}{2}{\stackbelow{\nl{#1}}{#2}}
\mcm{\nll}{1}{\stackbelow{\node}{#1}}
\mcm{\wun}{0}{\fcat{1}}
\mcm{\atuplebts}{1}{\langle #1 \rangle}
\mcm{\tuplebts}{1}{(#1)}
\mcm{\bo}{0}{(}
\mcm{\bc}{0}{)}
\mcm{\UBilax}{0}{\fcat{UBicat}_\mr{lax}}
\mcm{\UBiwk}{0}{\fcat{UBicat}_\mr{wk}}
\mcm{\UBistr}{0}{\fcat{UBicat}_\mr{str}}
\mcm{\Bilax}{0}{\fcat{Bicat}_\mr{lax}}
\mcm{\Biwk}{0}{\fcat{Bicat}_\mr{wk}}
\mcm{\Bistr}{0}{\fcat{Bicat}_\mr{str}}
\mcm{\rotsub}{0}{\cup \raisebox{0.1em}{$\scriptstyle{|}$}}
\mcm{\pd}{0}{\fcat{pd}}
\newenvironment{proof}{\paragraph*{Proof}}{\paragraph*{}}

\mcm{\rep}{1}{\widehat{#1}}
\mcm{\ovln}{1}{\overline{#1}}
\newarrow{Incl}C--->
\mcm{\Gph}{0}{\fcat{Gph}}
\mcm{\tr}{0}{\fcat{tr}}
\newcommand{\UB}{\textbf{(UB)}}
\newcommand{\UF}{\textbf{(UF)}}
\newcommand{\CB}{\textbf{(CB)}}
\newcommand{\CF}{\textbf{(CF)}}
\mcm{\ladj}{0}{\,\dashv\,}
\mcm{\zeropd}{0}{\node}
\newenvironment{ntdiag}
	{\begin{diagram}[size=1.5em,noPS]}	% added in noPS - looks nicer
	{\end{diagram}}
\mcm{\END}{0}{\fcat{End}}
\mcm{\HOM}{0}{\fcat{Hom}}

\newlength{\gwidth}	% the width of a glob
\newlength{\gvert}	% the overall vertical measurement
\newlength{\gdrop}	% the distance a labelled glob protrudes below the
\newlength{\gbaredrop}	% the distance from the textline to the bottom of the
\newlength{\goffset}	% the distance from the centre of the glob to the
\newlength{\gtemp}	% temporary register

\newcommand{\present}[1]{%
\makebox[1\gwidth]{%
\rule[-1\gdrop]{0ex}{1\gvert}%
\raisebox{-1\gbaredrop}{#1}}}

\newcommand{\presentl}[1]{%
\makebox[1\gwidth][l]{%
\rule[-1\gdrop]{0ex}{1\gvert}%
\raisebox{-1\gbaredrop}{#1}}}

\newcommand{\presentr}[1]{%
\makebox[1\gwidth][r]{%
\rule[-1\gdrop]{0ex}{1\gvert}%
\raisebox{-1\gbaredrop}{#1}}}

\newcommand{\ginitdims}[2]{%		% GLOBULAR VERSION
\setlength{\unitlength}{1em}%		% unitlength = 1em
\setlength{\goffset}{.25\unitlength}%	% globular offset = .25em
\setlength{\gwidth}{#1\unitlength}%	% width as specified
\setlength{\gvert}{#2\unitlength}%	% vert = #2
\setlength{\gdrop}{.5\gvert}%		% 
\addtolength{\gdrop}{-1\goffset}%	% 
\setlength{\gbaredrop}{1\gdrop}%	% gdrop = drop = half(vert) - offset
\addtolength{\gvert}{.6\unitlength}%	% total extra clearance of .6em...
\addtolength{\gdrop}{.3\unitlength}}	% ...half of which is at bottom

\newcommand{\cinitdims}[2]{%		% CELLULAR VERSION
\setlength{\unitlength}{1em}%		% unitlength = 1em
\setlength{\goffset}{.35\unitlength}%	% cellular offset = .35em
\setlength{\gwidth}{#1\unitlength}%	% width as specified
\setlength{\gvert}{#2\unitlength}%	% vert = #2
\setlength{\gdrop}{.5\gvert}%		% 
\addtolength{\gdrop}{-1\goffset}%	% 
\setlength{\gbaredrop}{1\gdrop}%	% gdrop = drop = half(vert) - offset
\addtolength{\gvert}{.6\unitlength}%	% total extra clearance of .6em...
\addtolength{\gdrop}{.3\unitlength}}	% ...half of which is at bottom

\newcommand{\gsinitdims}[2]{%		% SMALL GLOBULAR VERSION
\setlength{\unitlength}{0.5em}%		% unitlength = 0.5em
\setlength{\goffset}{.25\unitlength}%	% globular offset = .25em
\setlength{\gwidth}{#1\unitlength}%	% width as specified
\setlength{\gvert}{#2\unitlength}%	% vert = #2
\setlength{\gdrop}{.5\gvert}%		% 
\addtolength{\gdrop}{-1\goffset}%	% 
\setlength{\gbaredrop}{1\gdrop}%	% gdrop = drop = half(vert) - offset
\addtolength{\gvert}{.6\unitlength}%	% total extra clearance of .6em...
\addtolength{\gdrop}{.3\unitlength}}	% ...half of which is at bottom

\newcommand{\sidespic}[1]{%
\settowidth{\gtemp}{\ensuremath{#1}}%
\addtolength{\gwidth}{1\gtemp}}

\newcommand{\abovepic}[1]{%
\settoheight{\gtemp}{\ensuremath{#1}}%
\addtolength{\gvert}{1\gtemp}%
\settodepth{\gtemp}{\ensuremath{#1}}%
\addtolength{\gvert}{1\gtemp}}

\newcommand{\belowpic}[1]{%
\settoheight{\gtemp}{\ensuremath{#1}}%
\addtolength{\gvert}{1\gtemp}%
\addtolength{\gdrop}{1\gtemp}%
\settodepth{\gtemp}{\ensuremath{#1}}%
\addtolength{\gvert}{1\gtemp}%
\addtolength{\gdrop}{1\gtemp}}

\newcommand{\cell}[4]{\put(#1,#2){\makebox(0,0)[#3]{\ensuremath{#4}}}}
\mcm{\zmark}{0}{\scriptstyle{\bullet}}

\newcommand{\pregfst}[1]{%
\begin{picture}(0.5,0.2)(-0.5,-0.2)%
\cell{-0.1}{-0.2}{tr}{#1}%
\cell{0}{0}{c}{\zmark}%
\end{picture}}

\mcm{\gfst}{1}{%
\ginitdims{0.5}{0.4}%
\sidespic{#1}%
\belowpic{#1}%
\presentr{\pregfst{#1}}}

\newcommand{\preglst}[1]{%
\begin{picture}(0.5,0.2)(0,-0.2)%
\cell{0.1}{-0.2}{tl}{#1}%
\cell{0.05}{0}{c}{\zmark}%
\end{picture}}

\mcm{\glst}{1}{%
\ginitdims{.5}{.4}%
\sidespic{#1}%
\belowpic{#1}%
\presentl{\preglst{#1}}}

\newcommand{\preglft}[1]{%
\begin{picture}(0,0.2)(0,-0.2)%
\cell{-0.1}{-0.2}{tr}{#1}%
\cell{0.05}{0}{c}{\zmark}%
\end{picture}}

\mcm{\glft}{1}{%
\ginitdims{0}{.4}%
\belowpic{#1}%
\present{\preglft{#1}}}

\newcommand{\pregrgt}[1]{%
\begin{picture}(0,0.2)(0,-0.2)%
\cell{0.1}{-0.2}{tl}{#1}%
\cell{0.05}{0}{c}{\zmark}%
\end{picture}}

\mcm{\grgt}{1}{%
\ginitdims{0}{.4}%
\belowpic{#1}%
\present{\pregrgt{#1}}}

\newcommand{\pregblw}[1]{%
\begin{picture}(0,0.3)(0,-0.3)
\cell{0}{-0.3}{t}{#1}%
\cell{0.05}{0}{c}{\zmark}%
\end{picture}}

\mcm{\gblw}{1}{%
\ginitdims{0}{.6}%
\belowpic{#1}%
\present{\pregblw{#1}}}

\newcommand{\pregfbw}[1]{%
\begin{picture}(0,0.65)(0,-0.65)
\cell{0}{-0.65}{t}{#1}%
\cell{0.05}{0}{c}{\zmark}%
\end{picture}}

\mcm{\gfbw}{1}{%
\ginitdims{0}{1.3}%
\belowpic{#1}%
\present{\pregfbw{#1}}}

\newcommand{\pregzero}[1]{%
\begin{picture}(0.8,0.4)(-0.4,-0.4)
\cell{0}{-0.4}{t}{#1}%
\cell{0}{0}{c}{\zmark}%
\end{picture}}

\mcm{\gzero}{1}{%
\ginitdims{0.8}{.6}%
\belowpic{#1}%
\sidespic{#1}%	
\present{\pregzero{#1}}}

\newcommand{\pregone}[1]{%
\begin{picture}(5,0.4)(0,-0.2)%
\cell{2.5}{0.2}{b}{#1}%
\put(0,0){\vector(1,0){5}}%
\end{picture}}

\mcm{\gone}{1}{%
\ginitdims{5}{0.4}%
\abovepic{#1}%
\present{\pregone{#1}}}

\newcommand{\pregtwo}[3]{%
\begin{picture}(5,3.4)(0,-0.2)%
\cell{2.5}{3.2}{b}{#1}%
\cell{2.5}{-.2}{t}{#2}%
\cell{2.7}{1.5}{l}{#3}%
\qbezier(0,1.5)(2.5,4.5)(5,1.5)%
\qbezier(0,1.5)(2.5,-1.5)(5,1.5)%
\put(5,1.5){\vector(1,-1){0}}%
\put(5,1.5){\vector(1,1){0}}%
\put(2.5,2.5){\vector(0,-1){2}}%
\end{picture}}

\mcm{\gtwo}{3}{%
\ginitdims{5}{3.4}%
\abovepic{#1}%
\belowpic{#2}%
\present{\pregtwo{#1}{#2}{#3}}}

\newcommand{\pregthree}[5]{%
\begin{picture}(5,5.4)(0,-1.2)%
\cell{2.5}{4.2}{b}{#1}%
\cell{1.5}{1.7}{b}{#2}%
\cell{2.5}{-1.2}{t}{#3}%
\cell{2.7}{2.75}{l}{#4}%
\cell{2.7}{0.25}{l}{#5}%
\qbezier(0,1.5)(2.5,6.5)(5,1.5)%
\qbezier(0,1.5)(2.5,-3.5)(5,1.5)%
\put(0,1.5){\vector(1,0){5}}%
\put(2.5,3.5){\vector(0,-1){1.5}}%
\put(2.5,1){\vector(0,-1){1.5}}%
\put(5,1.5){\vector(1,-3){0}}%
\put(5,1.5){\vector(1,3){0}}%
\end{picture}}

\mcm{\gthree}{5}{%
\ginitdims{5}{5.4}%
\abovepic{#1}%
\belowpic{#3}%
\present{\pregthree{#1}{#2}{#3}{#4}{#5}}}

\newcommand{\pregfour}[7]{%
\begin{picture}(5,8.4)(0,-2.7)%
\cell{2.5}{5.7}{b}{#1}%
\cell{1.5}{2.8}{b}{#2}%
\cell{1.5}{0.2}{t}{#3}%
\cell{2.5}{-2.7}{t}{#4}%
\cell{2.7}{4.25}{l}{#5}%
\cell{2.7}{1.5}{l}{#6}%
\cell{2.7}{-1.25}{l}{#7}%
\qbezier(0,1.5)(2.5,9.5)(5,1.5)%
\qbezier(0,1.5)(2.5,4)(5,1.5)%
\qbezier(0,1.5)(2.5,-1)(5,1.5)%
\qbezier(0,1.5)(2.5,-6.5)(5,1.5)%
\put(2.5,5.25){\vector(0,-1){2}}%
\put(2.5,2.5){\vector(0,-1){2}}%
\put(2.5,-0.25){\vector(0,-1){2}}%
\put(5,1.5){\vector(1,-4){0}}%
\put(5,1.5){\vector(4,-3){0}}%
\put(5,1.5){\vector(4,3){0}}%
\put(5,1.5){\vector(1,4){0}}%
\end{picture}}

\mcm{\gfour}{7}{%
\ginitdims{5}{8.4}%
\abovepic{#1}%
\belowpic{#4}%
\present{\pregfour{#1}{#2}{#3}{#4}{#5}{#6}{#7}}}

\newcommand{\pregthreecell}[5]{%
\begin{picture}(8,5)(-4,-2.5)%
\cell{0}{2.5}{b}{#1}%
\cell{0}{-2.5}{t}{#2}%
\cell{-1.7}{0}{r}{#3}%
\cell{1.7}{0}{l}{#4}%
\cell{0}{0.2}{b}{#5}%
\qbezier(-4,0)(0,4.2)(4,0)%
\qbezier(-4,0)(0,-4.2)(4,0)%
\qbezier(-0.5,1.8)(-2.5,0)(-0.5,-1.8)%
\qbezier(0.5,1.8)(2.5,0)(0.5,-1.8)%
\put(-1,0){\vector(1,0){2}}%
\put(4,0){\vector(1,-1){0}}%
\put(4,0){\vector(1,1){0}}%
\put(-0.5,-1.8){\vector(1,-1){0}}%
\put(0.5,-1.8){\vector(-1,-1){0}}%
\end{picture}}

\mcm{\gthreecell}{5}{%
\ginitdims{8}{5}%
\abovepic{#1}%
\belowpic{#2}%
\present{\pregthreecell{#1}{#2}{#3}{#4}{#5}}}

\newcommand{\pregthreecellu}{%
\begin{picture}(5,3.4)(-0.5,-0.2)%
\qbezier(-.5,1.5)(2,4.5)(4.5,1.5)%
\qbezier(-.5,1.5)(2,-1.5)(4.5,1.5)%
\qbezier(1.5,2.7)(0.5,1.5)(1.5,0.3)%
\qbezier(2.5,2.7)(3.5,1.5)(2.5,0.3)%
\put(1.3,1.5){\vector(1,0){1.4}}%
\put(4.5,1.5){\vector(1,-1){0}}%
\put(4.5,1.5){\vector(1,1){0}}%
\put(1.5,0.3){\vector(2,-3){0}}%
\put(2.5,0.3){\vector(-2,-3){0}}%
\end{picture}}

\mcm{\gthreecellu}{0}{%
\ginitdims{5}{3.4}%
\present{\pregthreecellu}}

\newcommand{\pregtwocentre}[3]{%
\begin{picture}(5,3.4)(0,-0.2)%
\cell{2.5}{3.2}{b}{#1}%
\cell{2.5}{-.2}{t}{#2}%
\cell{2.5}{1.5}{c}{#3}%
\qbezier(0,1.5)(2.5,4.5)(5,1.5)%
\qbezier(0,1.5)(2.5,-1.5)(5,1.5)%
\put(5,1.5){\vector(1,-1){0}}%
\put(5,1.5){\vector(1,1){0}}%
\put(2.5,2.5){\vector(0,-1){2}}%
\end{picture}}

\mcm{\gtwocentre}{3}{%
\ginitdims{5}{3.4}%
\abovepic{#1}%
\belowpic{#2}%
\present{\pregtwocentre{#1}{#2}{#3}}}

\newcommand{\pregspecialone}[9]{%
\begin{picture}(8,8)(-4,-4)%
\cell{0}{3.9}{b}{#1}%
\cell{-2}{-0.2}{t}{#2}%
\cell{0}{-3.9}{t}{#3}%
\cell{-1.5}{1.1}{r}{#4}%
\cell{0.2}{1.5}{l}{#5}%
\cell{1.5}{1.1}{l}{#6}%
\cell{0.2}{-2}{l}{#7}%
\cell{-0.9}{2.3}{b}{#8}%
\cell{0.9}{2.3}{b}{#9}%
\qbezier(-4,0)(0,8)(4,0)%
\qbezier(-4,0)(0,-8)(4,0)%
\qbezier(-0.5,3.4)(-3.5,2)(-0.5,0.6)%
\qbezier(0.5,3.4)(3.5,2)(0.5,0.6)%
\put(-4,0){\vector(1,0){8}}%
\put(0,3.4){\vector(0,-1){2.8}}%
\put(0,-0.8){\vector(0,-1){2.4}}%
\put(-1.5,2.2){\vector(1,0){1.2}}%
\put(0.3,2.2){\vector(1,0){1.2}}%
\put(4,0){\vector(1,-2){0}}%
\put(4,0){\vector(1,2){0}}%
\put(-0.5,0.6){\vector(2,-1){0}}%
\put(0.5,0.6){\vector(-2,-1){0}}%
\end{picture}}

\mcm{\gspecialone}{9}{%
\ginitdims{8}{8}%
\abovepic{#1}%
\belowpic{#3}%
\present{\pregspecialone{#1}{#2}{#3}{#4}{#5}{#6}{#7}{#8}{#9}}}

\newcommand{\pregspecialtwo}{%
\begin{picture}(5,3.4)(0,-0.2)%
\qbezier(0,1.5)(2.5,4.5)(5,1.5)%
\qbezier(0,1.5)(2.5,-1.5)(5,1.5)%
\qbezier(1.7,2.5)(0,1.5)(1.7,0.5)%
\qbezier(3.3,2.5)(5,1.5)(3.3,0.5)%
\put(5,1.5){\vector(1,-1){0}}%
\put(5,1.5){\vector(1,1){0}}%
\put(1.7,0.5){\vector(3,-2){0}}%
\put(3.3,0.5){\vector(-3,-2){0}}%
\put(2.5,2.5){\vector(0,-1){2}}%
\put(1.2,1.5){\vector(1,0){1}}%
\put(2.8,1.5){\vector(1,0){1}}%
\end{picture}}

\mcm{\gspecialtwo}{0}{%
\ginitdims{5}{3.4}%
\present{\pregspecialtwo}}

\newcommand{\pregspecialthree}{%
\begin{picture}(5,5.4)(0,-1.2)%
\qbezier(0,1.5)(2.5,6.5)(5,1.5)%
\qbezier(0,1.5)(2.5,-3.5)(5,1.5)%
\qbezier(2,3.5)(1,2.75)(2,2)%
\qbezier(3,3.5)(4,2.75)(3,2)%
\qbezier(2,1)(1,0.25)(2,-0.5)%
\qbezier(3,1)(4,0.25)(3,-0.5)%
\put(0,1.5){\vector(1,0){5}}%
\put(1.5,2.75){\vector(1,0){2}}%
\put(1.5,0.25){\vector(1,0){2}}%
\put(5,1.5){\vector(1,-3){0}}%
\put(5,1.5){\vector(1,3){0}}%
\put(2,2){\vector(1,-1){0}}%
\put(3,2){\vector(-1,-1){0}}%
\put(2,-0.5){\vector(1,-1){0}}%
\put(3,-0.5){\vector(-1,-1){0}}%
\end{picture}}

\mcm{\gspecialthree}{0}{%
\ginitdims{5}{5.4}%
\present{\pregspecialthree}}

\newcommand{\pregonew}[1]{%
\begin{picture}(8,0.4)(0,-0.2)%
\cell{4}{0.2}{b}{#1}%
\put(0,0){\vector(1,0){8}}%
\end{picture}}

\mcm{\gonew}{1}{%
\ginitdims{8}{0.4}%
\abovepic{#1}%
\present{\pregonew{#1}}}

\mcm{\gzersu}{0}{%
\gsinitdims{0}{.6}%
\present{\pregblw{}}}

\mcm{\gonesu}{0}{%
\gsinitdims{5}{0.4}%
\present{\pregone{}}}

\mcm{\gtwosu}{0}{%
\gsinitdims{5}{3.4}%
\present{\pregtwo{}{}{}}}

\mcm{\gthreesu}{0}{%
\gsinitdims{5}{5.4}%
\present{\pregthree{}{}{}{}{}}}

\mcm{\gfoursu}{0}{%
\gsinitdims{5}{8.4}%
\present{\pregfour{}{}{}{}{}{}{}}}

\newcommand{\prectwo}[3]%
{\begin{picture}(4.2,3.4)(-0.1,-0.2)%
\cell{2}{3.2}{b}{#1}%
\cell{2}{-0.2}{t}{#2}%
\cell{2.2}{1.5}{l}{#3}%
\qbezier(0,2)(2,4)(4,2)%
\qbezier(0,1)(2,-1)(4,1)%
\put(4,2){\vector(1,-1){0}}%
\put(4,1){\vector(1,1){0}}%
\put(2,2.5){\vector(0,-1){2}}%
\end{picture}}

\mcm{\ctwo}{3}{%
\cinitdims{4.2}{3.4}%
\abovepic{#1}%
\belowpic{#2}%
\present{\prectwo{#1}{#2}{#3}}}

\newcommand{\precthree}[5]{%
\begin{picture}(4.2,5.4)(-0.1,-0.2)%
\cell{2}{5.2}{b}{#1}%
\cell{1}{2.7}{b}{#2}%
\cell{2}{-.2}{t}{#3}%
\cell{2.2}{3.75}{l}{#4}%
\cell{2.2}{1.25}{l}{#5}%
\qbezier(0,3)(2,7)(4,3)%
\qbezier(0,2)(2,-2)(4,2)%
\put(0,2.5){\vector(1,0){4}}%
\put(2,4.5){\vector(0,-1){1.5}}%
\put(2,2){\vector(0,-1){1.5}}%
\put(4,3){\vector(1,-3){0}}%
\put(4,2){\vector(1,3){0}}%
\end{picture}}

\mcm{\cthree}{5}{%
\cinitdims{4.2}{5.4}%
\abovepic{#1}%
\belowpic{#3}%
\present{\precthree{#1}{#2}{#3}{#4}{#5}}}

\newcommand{\prectwoop}[3]%
{\begin{picture}(4.2,3.4)(-0.1,-0.2)%
\cell{2}{3.2}{b}{#1}%
\cell{2}{-0.2}{t}{#2}%
\cell{2.2}{1.5}{l}{#3}%
\qbezier(0,2)(2,4)(4,2)%
\qbezier(0,1)(2,-1)(4,1)%
\put(0,2){\vector(-1,-1){0}}%
\put(0,1){\vector(-1,1){0}}%
\put(2,2.5){\vector(0,-1){2}}%
\end{picture}}

\mcm{\ctwoop}{3}{%
\cinitdims{4.2}{3.4}%
\abovepic{#1}%
\belowpic{#2}%
\present{\prectwoop{#1}{#2}{#3}}}

\newcommand{\prectwopar}[4]{%
\begin{picture}(4.2,3.4)(-0.1,-0.2)%
\cell{2}{3.2}{b}{#1}%
\cell{2}{-0.2}{t}{#2}%
\cell{1.6}{1.5}{r}{#3}%
\cell{2.4}{1.5}{l}{#4}%
\qbezier(0,2)(2,4)(4,2)%
\qbezier(0,1)(2,-1)(4,1)%
\put(4,2){\vector(1,-1){0}}%
\put(4,1){\vector(1,1){0}}%
\put(1.8,2.5){\vector(0,-1){2}}%
\put(2.2,2.5){\vector(0,-1){2}}%
\end{picture}}

\mcm{\ctwopar}{4}{%
\cinitdims{4.2}{3.4}%
\abovepic{#1}%
\belowpic{#2}%
\present{\prectwopar{#1}{#2}{#3}{#4}}}

\newcommand{\precthreein}[5]{%
\begin{picture}(4.2,5.4)(-0.1,-0.2)%
\cell{2}{5.2}{b}{#1}%
\cell{1}{2.7}{b}{#2}%
\cell{2}{-.2}{t}{#3}%
\cell{2.2}{3.75}{l}{#4}%
\cell{2.2}{1.25}{l}{#5}%
\qbezier(0,3)(2,7)(4,3)%
\qbezier(0,2)(2,-2)(4,2)%
\put(0,2.5){\vector(1,0){4}}%
\put(2,4.5){\vector(0,-1){1.5}}%
\put(2,0.5){\vector(0,1){1.5}}%
\put(4,3){\vector(1,-3){0}}%
\put(4,2){\vector(1,3){0}}%
\end{picture}}

\mcm{\cthreein}{5}{%
\cinitdims{4.2}{5.4}%
\abovepic{#1}%
\belowpic{#3}%
\present{\precthreein{#1}{#2}{#3}{#4}{#5}}}

\newcommand{\precthreecell}[5]{%
\begin{picture}(8.2,5)(-4.1,-2.5)%
\cell{0}{2.5}{b}{#1}%
\cell{0}{-2.5}{t}{#2}%
\cell{-1.7}{0}{r}{#3}%
\cell{1.7}{0}{l}{#4}%
\cell{0}{0.2}{b}{#5}%
\qbezier(-4,0.5)(0,4)(4,0.5)%
\qbezier(-4,-0.5)(0,-4)(4,-0.5)%
\qbezier(-0.5,2)(-2.5,0)(-0.5,-2)%
\qbezier(0.5,2)(2.5,0)(0.5,-2)%
\put(-1,0){\vector(1,0){2}}%
\put(4,0.5){\vector(1,-1){0}}%
\put(4,-0.5){\vector(1,1){0}}%
\put(-0.5,-2){\vector(1,-1){0}}%
\put(0.5,-2){\vector(-1,-1){0}}%
\end{picture}}

\mcm{\cthreecell}{5}{%
\cinitdims{8.2}{5}%
\abovepic{#1}%
\belowpic{#2}%
\present{\precthreecell{#1}{#2}{#3}{#4}{#5}}}

\newcommand{\precthreecellpar}[6]{%
\begin{picture}(8.2,5)(-4.1,-2.5)%
\cell{0}{2.5}{b}{#1}%
\cell{0}{-2.5}{t}{#2}%
\cell{-1.7}{0}{r}{#3}%
\cell{1.7}{0}{l}{#4}%
\cell{0}{0.4}{b}{#5}%
\cell{0}{-0.4}{t}{#6}%
\qbezier(-4,0.5)(0,4)(4,0.5)%
\qbezier(-4,-0.5)(0,-4)(4,-0.5)%
\qbezier(-0.5,2)(-2.5,0)(-0.5,-2)%
\qbezier(0.5,2)(2.5,0)(0.5,-2)%
\put(-1,0.2){\vector(1,0){2}}%
\put(-1,-0.2){\vector(1,0){2}}%
\put(4,0.5){\vector(1,-1){0}}%
\put(4,-0.5){\vector(1,1){0}}%
\put(-0.5,-2){\vector(1,-1){0}}%
\put(0.5,-2){\vector(-1,-1){0}}%
\end{picture}}

\mcm{\cthreecellpar}{6}{%
\cinitdims{8.2}{5}%
\abovepic{#1}%
\belowpic{#2}%
\present{\precthreecellpar{#1}{#2}{#3}{#4}{#5}{#6}}}

\newcommand{\prectwov}[5]{%
\begin{picture}(3.4,4.2)(0.8,0.9)%
\cell{2.5}{5.1}{b}{#1}%
\cell{2.5}{0.9}{t}{#2}%
\cell{0.8}{3}{r}{#3}%
\cell{4.2}{3}{l}{#4}%
\cell{2.5}{3.2}{b}{#5}%
\qbezier(2,5)(0,3)(2,1)%
\qbezier(3,5)(5,3)(3,1)%
\put(2,1){\vector(1,-1){0}}%
\put(3,1){\vector(-1,-1){0}}%
\put(1.5,3){\vector(1,0){2}}%
\end{picture}}

\mcm{\ctwov}{5}{%
\cinitdims{3.4}{4.2}%
\abovepic{#1}%
\belowpic{#2}%
\sidespic{#3}%
\sidespic{#4}%
\present{\prectwov{#1}{#2}{#3}{#4}{#5}}}

\newcommand{\precthreecellv}[7]{%
\begin{picture}(5,8.2)(0.5,-1.6)%
\cell{3}{6.6}{b}{#1}%
\cell{3}{-1.6}{t}{#2}%
\cell{0.5}{2.5}{r}{#3}%
\cell{5.5}{2.5}{l}{#4}%
\cell{3}{4.2}{b}{#5}%
\cell{3}{0.8}{t}{#6}%
\cell{3.2}{2.5}{l}{#7}%
\qbezier(3.5,6.5)(7,2.5)(3.5,-1.5)%
\qbezier(2.5,6.5)(-1,2.5)(2.5,-1.5)%
\put(2.5,-1.5){\vector(1,-1){0}}%
\put(3.5,-1.5){\vector(-1,-1){0}}%
\qbezier(1,3)(3,5)(5,3)%
\qbezier(1,2)(3,0)(5,2)%
\put(5,3){\vector(1,-1){0}}%
\put(5,2){\vector(1,1){0}}%
\put(3,3.5){\vector(0,-1){2}}%
\end{picture}}

\mcm{\cthreecellv}{7}{%
\cinitdims{5}{8.2}%
\abovepic{#1}%
\belowpic{#2}%
\sidespic{#3}%
\sidespic{#4}%
\present{\precthreecellv{#1}{#2}{#3}{#4}{#5}{#6}{#7}}}

\newcommand{\pretopez}[2]{%
\begin{picture}(2.6,2.3)(-1.3,-2.2)% 
\cell{0}{-2.2}{t}{#1}%
\cell{0}{-1.2}{c}{#2}%
\qbezier(0,0)(-2,-2)(0,-2)%
\qbezier(0,0)(2,-2)(0,-2)%
\put(0,0){\vector(-1,1){0}}%
\end{picture}}

\mcm{\topez}{2}{%
\ginitdims{2.6}{2.3}% 
\belowpic{#1}%
\present{\pretopez{#1}{#2}}}

\newcommand{\pretopea}[3]{%
\begin{picture}(4,1.9)(-2,-0,2)%
\cell{0}{1.7}{b}{#1}%
\cell{0}{-0.2}{t}{#2}%
\cell{0}{0.7}{c}{#3}%
\qbezier(-2,0)(0,3)(2,0)%
\put(-2,0){\vector(1,0){4}}%
\put(2,0){\vector(2,-3){0}}%
\end{picture}}

\mcm{\topea}{3}{%
\ginitdims{4}{1.9}%
\abovepic{#1}%
\belowpic{#2}%
\present{\pretopea{#1}{#2}{#3}}}

\newcommand{\pretopeb}[4]{%
\begin{picture}(4,2.2)(-2,-0.2)%
\cell{-1.1}{1}{br}{#1}%
\cell{1.1}{1}{bl}{#2}%
\cell{0}{-0.2}{t}{#3}%
\cell{0}{0.8}{c}{#4}%
\put(-2,0){\vector(1,1){2}}%
\put(0,2){\vector(1,-1){2}}%
\put(-2,0){\vector(1,0){4}}%
\end{picture}}

\mcm{\topeb}{4}{%
\ginitdims{4}{2.2}%
\belowpic{#3}%
\present{\pretopeb{#1}{#2}{#3}{#4}}}

\newcommand{\pretopec}[5]{%
\begin{picture}(4,2.2)(-2,-0.2)%
\cell{-1.8}{1}{br}{#1}%
\cell{0}{2.2}{b}{#2}%
\cell{1.8}{1}{bl}{#3}%
\cell{0}{-0.2}{t}{#4}%
\cell{0}{0.8}{c}{#5}%
\put(-2,0){\vector(1,2){1}}%
\put(-1,2){\vector(1,0){2}}%
\put(1,2){\vector(1,-2){1}}%
\put(-2,0){\vector(1,0){4}}%
\end{picture}}

\mcm{\topec}{5}{%
\ginitdims{4}{2.2}%
\sidespic{#1}%
\abovepic{#2}%
\sidespic{#3}%
\belowpic{#4}%
\present{\pretopec{#1}{#2}{#3}{#4}{#5}}}

\newcommand{\pretoped}[6]{%
\begin{picture}(4,2.5)(-2,-0.2)%
\cell{-2}{0.6}{br}{#1}%
\cell{-0.7}{2.2}{br}{#2}%
\cell{0.7}{2.2}{bl}{#3}%
\cell{2}{0.6}{bl}{#4}%
\cell{0}{-0.2}{t}{#5}%
\cell{0}{0.8}{c}{#6}%
\put(-2,0){\vector(1,3){0.5}}%
\put(-1.5,1.5){\vector(3,2){1.5}}%
\put(0,2.5){\vector(3,-2){1.5}}%
\put(1.5,1.5){\vector(1,-3){0.5}}%
\put(-2,0){\vector(1,0){4}}%
\end{picture}}

\mcm{\toped}{6}{%
\ginitdims{4}{2.5}%
\sidespic{#1}%
\abovepic{#2}%
\abovepic{#3}%
\sidespic{#4}%
\belowpic{#5}%
\present{\pretoped{#1}{#2}{#3}{#4}{#5}{#6}}}

\newcommand{\pretopeq}[5]{%
\begin{picture}(4,2.5)(-2,-0.2)%
\cell{-2}{0.6}{br}{#1}%
\cell{-1}{2.2}{br}{#2}%
\cell{2}{0.6}{bl}{#3}%
\cell{0}{-0.2}{t}{#4}%
\cell{0}{0.8}{c}{#5}%
\put(-2,0){\vector(1,3){0.5}}%
\put(-1.5,1.5){\vector(1,1){1}}%
\cell{0.9}{2.3}{c}{\ddots}
\put(1.5,1.5){\vector(1,-3){0.5}}%
\put(-2,0){\vector(1,0){4}}%
\end{picture}}

\mcm{\topeq}{5}{%
\ginitdims{4}{2.5}%
\sidespic{#1}%
\abovepic{#2}%
\sidespic{#3}%
\belowpic{#4}%
\present{\pretopeq{#1}{#2}{#3}{#4}{#5}}}

\newcommand{\pretopebase}[1]{%
\begin{picture}(4,0.4)(0,-0.2)%
\cell{2}{0.2}{b}{#1}%
\put(0,0){\vector(1,0){4}}%
\end{picture}}

\mcm{\topebase}{1}{%
\ginitdims{4}{0.4}%
\abovepic{#1}%
\present{\pretopebase{#1}}}

\newcommand{\pretopezs}[2]{%
\begin{picture}(2.6,2.3)(-1.3,-2.2)% 
\cell{0}{-2.2}{t}{#1}%
\cell{0}{-1.2}{c}{#2}%
\qbezier(0,0)(-2,-2)(0,-2)%
\qbezier(0,0)(2,-2)(0,-2)%
\end{picture}}

\mcm{\topezs}{2}{%
\ginitdims{2.6}{2.3}% 
\belowpic{#1}%
\present{\pretopezs{#1}{#2}}}

\newcommand{\pretopeas}[3]{%
\begin{picture}(4,1.9)(-2,-0,2)%
\cell{0}{1.7}{b}{#1}%
\cell{0}{-0.2}{t}{#2}%
\cell{0}{0.7}{c}{#3}%
\qbezier(-2,0)(0,3)(2,0)%
\put(-2,0){\line(1,0){4}}%
\end{picture}}

\mcm{\topeas}{3}{%
\ginitdims{4}{1.9}%
\abovepic{#1}%
\belowpic{#2}%
\present{\pretopeas{#1}{#2}{#3}}}

\newcommand{\pretopebs}[4]{%
\begin{picture}(4,2.2)(-2,-0.2)%
\cell{-1.1}{1}{br}{#1}%
\cell{1.1}{1}{bl}{#2}%
\cell{0}{-0.2}{t}{#3}%
\cell{0}{0.8}{c}{#4}%
\put(-2,0){\line(1,1){2}}%
\put(0,2){\line(1,-1){2}}%
\put(-2,0){\line(1,0){4}}%
\end{picture}}

\mcm{\topebs}{4}{%
\ginitdims{4}{2.2}%
\belowpic{#3}%
\present{\pretopebs{#1}{#2}{#3}{#4}}}

\newcommand{\pretopecs}[5]{%
\begin{picture}(4,2.2)(-2,-0.2)%
\cell{-1.8}{1}{br}{#1}%
\cell{0}{2.2}{b}{#2}%
\cell{1.8}{1}{bl}{#3}%
\cell{0}{-0.2}{t}{#4}%
\cell{0}{0.8}{c}{#5}%
\put(-2,0){\line(1,2){1}}%
\put(-1,2){\line(1,0){2}}%
\put(1,2){\line(1,-2){1}}%
\put(-2,0){\line(1,0){4}}%
\end{picture}}

\mcm{\topecs}{5}{%
\ginitdims{4}{2.2}%
\sidespic{#1}%
\abovepic{#2}%
\sidespic{#3}%
\belowpic{#4}%
\present{\pretopecs{#1}{#2}{#3}{#4}{#5}}}

\newcommand{\pretopeds}[6]{%
\begin{picture}(4,2.5)(-2,-0.2)%
\cell{-2}{0.6}{br}{#1}%
\cell{-0.7}{2.2}{br}{#2}%
\cell{0.7}{2.2}{bl}{#3}%
\cell{2}{0.6}{bl}{#4}%
\cell{0}{-0.2}{t}{#5}%
\cell{0}{0.8}{c}{#6}%
\put(-2,0){\line(1,3){0.5}}%
\put(-1.5,1.5){\line(3,2){1.5}}%
\put(0,2.5){\line(3,-2){1.5}}%
\put(1.5,1.5){\line(1,-3){0.5}}%
\put(-2,0){\line(1,0){4}}%
\end{picture}}

\mcm{\topeds}{6}{%
\ginitdims{4}{2.5}%
\sidespic{#1}%
\abovepic{#2}%
\abovepic{#3}%
\sidespic{#4}%
\belowpic{#5}%
\present{\pretopeds{#1}{#2}{#3}{#4}{#5}{#6}}}

\newcommand{\pretopeqs}[5]{%
\begin{picture}(4,2.5)(-2,-0.2)%
\cell{-2}{0.6}{br}{#1}%
\cell{-1}{2.2}{br}{#2}%
\cell{2}{0.6}{bl}{#3}%
\cell{0}{-0.2}{t}{#4}%
\cell{0}{0.8}{c}{#5}%
\put(-2,0){\line(1,3){0.5}}%
\put(-1.5,1.5){\line(1,1){1}}%
\cell{0.9}{2.3}{c}{\ddots}
\put(1.5,1.5){\line(1,-3){0.5}}%
\put(-2,0){\line(1,0){4}}%
\end{picture}}

\mcm{\topeqs}{5}{%
\ginitdims{4}{2.5}%
\sidespic{#1}%
\abovepic{#2}%
\sidespic{#3}%
\belowpic{#4}%
\present{\pretopeqs{#1}{#2}{#3}{#4}{#5}}}

\newcommand{\pretopebases}[1]{%
\begin{picture}(4,0.4)(0,-0.2)%
\cell{2}{0.2}{b}{#1}%
\put(0,0){\line(1,0){4}}%
\end{picture}}

\mcm{\topebases}{1}{%
\ginitdims{4}{0.4}%
\abovepic{#1}%
\present{\pretopebases{#1}}}

\newcommand{\pregdots}[6]{%
\begin{picture}(5,8.4)(0,-2.7)%
\cell{2.5}{5.7}{b}{#1}%
\cell{1.5}{2.8}{b}{#2}%
\cell{1.5}{0.2}{t}{#3}%
\cell{2.5}{-2.7}{t}{#4}%
\cell{2.7}{4.25}{l}{#5}%
\cell{2.7}{-1.25}{l}{#6}%
\qbezier(0,1.5)(2.5,9.5)(5,1.5)%
\qbezier(0,1.5)(2.5,4)(5,1.5)%
\qbezier(0,1.5)(2.5,-1)(5,1.5)%
\qbezier(0,1.5)(2.5,-6.5)(5,1.5)%
\put(2.5,5.25){\vector(0,-1){2}}%
\put(2.5,-0.25){\vector(0,-1){2}}%
\cell{2.5}{1.7}{c}{\vdots}%
\put(5,1.5){\vector(1,-4){0}}%
\put(5,1.5){\vector(4,-3){0}}%
\put(5,1.5){\vector(4,3){0}}%
\put(5,1.5){\vector(1,4){0}}%
\end{picture}}

\mcm{\gdots}{6}{%
\ginitdims{5}{8.4}%
\abovepic{#1}%
\belowpic{#4}%
\present{\pregdots{#1}{#2}{#3}{#4}{#5}{#6}}}

\newcommand{\EmptyOne}
{\begin{tree}
\enode	\\
\dn	\\
\node	\\
\end{tree}}

\newcommand{\Oak}[5]
{\begin{tree}
 & & & &\enode& & & & & & & & \\
 & & & &\dn& & & & & & & & \\
\nl{#1}& &\nl{#2}& &\node& & & &\nl{#4}& & & &\nl{#5}\\
 &\rt{2}&\dn&\lt{2}& & & & & &\rt{2}& &\lt{2}& \\
 & &\node& & & &\nl{#3}& & & &\node& & \\
 & & &\rt{4}& & &\dn& & &\lt{4}& & & \\
 & & & & & &\node& & & & & & \\
\end{tree}}

\newcommand{\Pear}[3]
{\begin{tree}
	&	&	&	&\enode	&&	&	&	\\
	&	&	&	&\dn	&&	&	&	\\
\nl{#1}&	&	&	&\node	&&	&	&	\\
	&\rt{2}	&	&\lt{2}	&	&&	&	&	\\
	&	&\node	&	&	&&	&	&\nl{#2}\\
	&	&	&\rt{4}	&	&&	&\lt{2}	&	\\
	&	&	&	&	&&\nll{#3}&	&	\\
\end{tree}}

\newcommand{\Orange}[3]
{\begin{tree}
	&	&	&	&\nl{#2}\\
	&	&	&	&\dn	\\
\nl{#1}&	&	&	&\node	\\
	&\rt{2}	&	&\lt{2}	&	\\
	&	&\nll{#3}&	&	\\
\end{tree}}

\newcommand{\Apple}[3]
{\begin{tree}
	&	&\enode	&	&	\\
	&	&\dn	&	&	\\
\nl{#1}	&	&\node	&	&\nl{#2}\\
	&\rt{2}	&\dn	&\lt{2}	&	\\
	&	&\node	&	&	\\
	&	&\dn	&	&	\\
	&	&\nll{#3}&	&	\\
\end{tree}}

\newcommand{\MixedFruit}[5]
{\begin{tree}
 & & & &\nl{#2}& & & & & & &\enode& & \\
 & & & &\dn    & & & & & & &\dn   & & \\
\nl{#1}& & & &\node& & &\enode& &\nl{#3}& &\node& &\nl{#4}\\
 &\rt{2}& &\lt{2}& & & &\dn& & &\rt{2}&\dn&\lt{2}& \\
 & &\node& & & & &\node& & & &\node& & \\
 & & &\rt{3}& & &\lt{2}& & & & &\dn& & \\
 & & & & &\node& & & & & &\node& & \\
 & & & & & &\rt{3}& & & &\lt{3}& & & \\
 & & & & & & & &\nll{#5}& & & & & \\
\end{tree}}

\newcommand{\treedc}{
\begin{tree}
\node &\node            &\node &      &      &\node \\
      &\rt{1} \dn \lt{1}&      &      &      & \dn  \\
      &\node            &      &\node &      &\node \\
      &                 &\rt{2}&\dn   &\lt{2}&      \\
      &                 &      &\node &      &      \\
\end{tree}}

\newcommand{\treec}{
\begin{tree}
\node	&	&\node	&	&\node	\\
	&\rt{2}	&\dn	&\lt{2}	&	\\
	&	&\node	&	&	\\
\end{tree}}

\newlength{\volt}
\newcommand{\transistor}[5]
{\setlength{\unitlength}{1\volt}
\begin{picture}(18,12)(-5,-6)
\put(0,6){\line(0,-1){12}}
\put(0,-6){\line(3,2){9}}
\put(0,6){\line(3,-2){9}}
\put(9,0){\line(1,0){2}}
\put(-2,4){\line(1,0){2}}
\put(-2,2){\line(1,0){2}}	
\put(-2,-4){\line(1,0){2}}
\put(12,-0.5){\ensuremath{#5}}
\put(-5,3.5){\ensuremath{#2}}
\put(-5,1.5){\ensuremath{#3}}		
\put(-5,-4.5){\ensuremath{#4}}
\thicklines
\put(-1.5,0){\line(1,0){.1}}	
\put(-1.5,-1){\line(1,0){.1}}		
\put(-1.5,-2){\line(1,0){.1}}		
\thinlines
\put(2,-0.5){\ensuremath{#1}}
\end{picture}}
\newcommand{\ctransistor}[5]
	{\raisebox{-6\volt}{\transistor{#1}{#2}{#3}{#4}{#5}}}

\newcommand{\bftransistor}[4]
{\setlength{\unitlength}{1\volt}
\begin{picture}(18,12)(-5,-6)
\put(0,6){\line(0,-1){12}}
\put(0,-6){\line(3,2){9}}
\put(0,6){\line(3,-2){9}}
\put(9,0){\line(1,0){2}}
\put(-2,4){\line(1,0){2}}
\put(-2,-4){\line(1,0){2}}
\put(12,-0.5){\ensuremath{#4}}
\put(-5,3.5){\ensuremath{#2}}
\put(-5,-4.5){\ensuremath{#3}}
\thicklines
\put(-1.5,1){\line(1,0){.1}}		
\put(-1.5,0){\line(1,0){.1}}		
\put(-1.5,-1){\line(1,0){.1}}
\thinlines		
\put(2,-0.5){\ensuremath{#1}}
\end{picture}}
\newcommand{\comptrans}[4]
{\setlength{\volt}{.5ex}
\setlength{\unitlength}{1\volt}
\begin{picture}(80,72)(0,-36)
\put(0,23){\bftransistor{#1}{}{}{}}
\put(0,6){\bftransistor{#2}{}{}{}}
\put(0,-35){\bftransistor{#3}{}{}{}}
\put(63,-6){\transistor{#4}{}{}{}}
\put(17,29){\line(2,-1){50}}
\put(17,12){\line(5,-1){50}}
\put(17,-29){\line(2,1){50}}
\thicklines
\put(24,-5){\line(1,0){.4}}		
\put(24,-9){\line(1,0){.4}}		
\put(24,-13){\line(1,0){.4}}
\thinlines		
\setlength{\volt}{1ex}
\end{picture}}

\newcommand{\threetransistor}[5]
{\setlength{\unitlength}{1\volt}
\begin{picture}(18,12)(-5,-6)
\put(0,6){\line(0,-1){12}}
\put(0,-6){\line(3,2){9}}
\put(0,6){\line(3,-2){9}}
\put(9,0){\line(1,0){2}}
\put(-2,4){\line(1,0){2}}
\put(-2,0){\line(1,0){2}}	
\put(-2,-4){\line(1,0){2}}
\put(12,-0.5){\ensuremath{#5}}
\put(-5,3.5){\ensuremath{#2}}
\put(-5,-0.5){\ensuremath{#3}}		
\put(-5,-4.5){\ensuremath{#4}}
\put(2,-0.5){\ensuremath{#1}}
\end{picture}}
\newcommand{\twotransistor}[4]
{\setlength{\unitlength}{1\volt}
\begin{picture}(18,12)(-5,-6)
\put(0,6){\line(0,-1){12}}
\put(0,-6){\line(3,2){9}}
\put(0,6){\line(3,-2){9}}
\put(9,0){\line(1,0){2}}
\put(-2,4){\line(1,0){2}}
\put(-2,-4){\line(1,0){2}}
\put(12,-0.5){\ensuremath{#4}}
\put(-5,3.5){\ensuremath{#2}}
\put(-5,-4.5){\ensuremath{#3}}
\put(2,-0.5){\ensuremath{#1}}
\end{picture}}
\newcommand{\notransistor}[2]
{\setlength{\unitlength}{1\volt}
\begin{picture}(18,12)(-5,-6)
\put(0,6){\line(0,-1){12}}
\put(0,-6){\line(3,2){9}}
\put(0,6){\line(3,-2){9}}
\put(9,0){\line(1,0){2}}
\put(12,-0.5){\ensuremath{#2}}
\put(2,-0.5){\ensuremath{#1}}
\end{picture}}
\newcommand{\freemoncatpic}
{\setlength{\volt}{.6ex}
\setlength{\unitlength}{1\volt}
\begin{picture}(18,48)
\put(0,3){\twotransistor{\theta_3}{a_4}{a_5}{a'_3}}
\put(0,18){\notransistor{\theta_2}{a'_2}}
\put(0,33){\threetransistor{\theta_1}{a_1}{a_2}{a_3}{a'_1}}
\put(5.7,0){\framebox(10.3,48){}}
\end{picture}}
\newcommand{\dotty}{\raisebox{-.4ex}{\ensuremath{\!\scriptstyle\bullet}}}
\newcommand{\discfibpic}
{\setlength{\unitlength}{1em}		% if you change this, also
\begin{picture}(16,9)(-6,0)
\put(4,2){\oval(8,4)}
\put(4,8){\dotty}
\put(2,2){\dotty}
\put(4,1){\dotty}
\put(6,1){\dotty}
\put(5,3){\dotty}
\put(4,8){\vector(-1,-3){2}}
\put(2,2){\vector(2,-1){2}}
\put(4,8){\vector(0,-1){7}}
\put(2,2){\vector(3,1){3}}
\put(4,1){\vector(1,0){2}}
\put(4.2,8){0}
\put(4.2,5.5){\ensuremath{\theta\of y}}
\put(2.5,5.5){\ensuremath{y}}
\put(1.2,1.7){\ensuremath{a}}
\put(3.6,.2){\ensuremath{a'}}
\put(2.1,.9){\ensuremath{\theta}}
\put(8.5,1.5){\ensuremath{C_0}}
\end{picture}}
\newcommand{\cdiscfibpic}{\raisebox{-4.5em}{\discfibpic}}

\newcommand{\clearemptydoublepage}{\newpage{\pagestyle{empty}\cleardoublepage}}
\usepackage{fancyhdgs}
\begin{document}
\frontmatter

\title{Operads in Higher-Dimensional Category Theory}
\author{Thomas S. H. Leinster\\ \\
        \normalsize{Trinity College and St John's College, Cambridge}}
\date{\normalsize
	\mbox{}\vspace*{28mm}\\
	\piccy{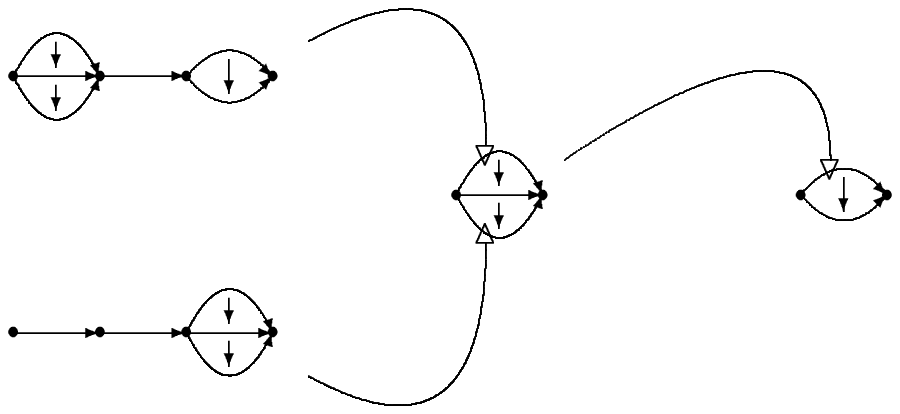}\\
	\mbox{}\vspace*{27mm}\\
	A dissertation submitted for the degree of Doctor of Philosophy at
	the University of Cambridge\\
	\mbox{}\vspace*{1mm}\\
	July 2000}

\maketitle
\clearemptydoublepage

\begin{titlepage}
\mbox{}\vspace*{60mm}\\
\noindent
This dissertation is the result of my own work and includes nothing which is
the outcome of work done in collaboration. The statements made in the
`Related Work' section of the Introduction, concerning which ideas are
original or novel, are to the best of my knowledge correct.
\vspace*{3mm}\\
\noindent
This dissertation is not substantially the same as any that I have submitted
for a degree or diploma or any other qualification at any other university.
\vspace*{30mm}\\
\begin{flushright}
Thomas S. H. Leinster\\
10 July 2000
\end{flushright}
\end{titlepage}
\clearemptydoublepage

\tableofcontents

\chapter{Introduction}

The purpose of this dissertation is to set up a theory of generalized operads
and multicategories, and to use it as a language in which to propose a
definition of weak $\omega$-category. This theory of operads and
multicategories has various other applications too: for instance, to the
opetopic approach to $n$-categories expounded by Baez, Dolan and others, and
to the theory of enrichment of higher-dimensional categorical structures. We
sketch some of these further developments, without exploring them in full.

We start with a look at bicategories~(Chapter~\ref{ch:bicats}). Having
reviewed the basics of the classical definition, we define `unbiased
bicategories', in which $n$-fold composites of 1-cells are specified for all
natural $n$ (rather than the usual nullary and binary presentation). We go on
to show that the theories of (classical) bicategories and of unbiased
bicategories are equivalent, in a strong sense.

The heart of this work is the theory of generalized operads and
multicategories. More exactly, given a monad $T$ on a category \Eee,
satisfying simple conditions, there is a theory of $T$-operads and
$T$-multicategories.  (As explained in `Terminology' below, a $T$-operad is a
special kind of $T$-multicategory.) In Chapter~\ref{ch:mtione} we set up the
basic concepts of the theory, including the important definition of an
algebra for a $T$-multicategory. In Chapter~\ref{ch:mtitwo} we cover an
assortment of further operadic topics, some of which are used in later parts
of the thesis, and some of which pertain to the applications mentioned in the
first paragraph.

Chapter~\ref{ch:defn} is a definition of weak $\omega$-category. (That is, it
is a \emph{proposed} definition; there are many such proposals out there, and
no attempt at a comparison is made.) As discussed at more length under
`Related Work', it is a modification of Batanin's
definition~\cite{Bat}. Having given the definition formally, we take a long
look at why it is a \emph{reasonable} definition. We then explore weak
$n$-categories (for finite $n$), and show that weak 2-categories are exactly
unbiased bicategories.

The four appendices take care of various details which would have been
distracting in the main text. Appendix~\ref{app:unbiased} contains the proof
that unbiased bicategories are essentially the same as classical
bicategories. Appendix~\ref{app:free} describes how to form the free
$T$-multicategory on a given $T$-graph. In Appendix~\ref{app:strict-omega} we
discuss various facts about strict $\omega$-categories, including a proof
that the category they form is monadic over an appropriate category of
graphs. Finally, Appendix~\ref{app:initial} is a proof of the existence of an
initial object in a certain category, as required in Chapter~\ref{ch:defn}.

\subsection*{Terminology}

The terminology for `strength' in higher-dimensional category theory is
rather in disarray. For example, when something works up to coherent
isomorphism, it is variously described as `pseudo', `weak' and `strong', or
not given a qualifier at all. In the context of maps between bicategories
another word altogether is often used (`homomorphism'---see \cite{Ben}).  Not
quite as severe a problem is the terminology for $n$-categories themselves:
the version where things hold up to coherent isomorphism or equivalence is
(almost) invariably called weak, and the version where everything holds up to
equality is always called strict, but `$n$-category' on its own is sometimes
used to mean the weak one and sometimes the strict one. The tradition has
been for `$n$-category' to mean `strict $n$-category'. However, Baez has
argued (convincingly) that the terminology should reflect the fact that the
weak version is much more abundant in nature; so in his work `$n$-category'
means `weak $n$-category'.

I have tried to bring some unity to the situation. When an entity is
characterized by things holding on the nose (i.e.\ up to equality), it will
be called strict. When they hold up to coherent isomorphism or equivalence it
will be called weak. When they hold up to a not-necessarily-invertible
connecting map (which does not happen often here), it will be called lax. The
term `$n$-category' will not (I hope) be used in isolation, but will always
be qualified by either `strict' or `weak', except in informal discussion
where both possibilities are intended. However, in deference to tradition,
`2-category' will always mean `strict 2-category', and `bicategory' will be
used for the notion of weak 2-category proposed by B\'enabou in \cite{Ben}.

We will, of course, be talking about operads and multicategories. Again the
terminology has been a bit messy: topologists, who by and large do not seem
to be aware of Lambek's (late 1960s) definition of multicategory, call
multicategories `coloured operads'; whereas amongst category theorists, the
notion of multicategory seems much more widely known than that of
operad. Basically, an operad is a one-object multicategory. This is also the
way the terminology will work when we are dealing with generalized operads
and multicategories, from Chapter~\ref{ch:mtione} onwards: a $T$-operad will
be a one-object $T$-multicategory, in a sense made precise just after the
definition~(\ref{defn:multicat}) of $T$-multicategory. So a $T$-operad is a
special kind of $T$-multicategory. This means that in the title of this work,
the word `operads' would more accurately be `multicategories': but, of
course, euphony is paramount.

I have not been very conscientious about the distinction between small and
large (sets and classes), and hope that the reader will find the issue no
more disturbing than usual.

$0$ is a member of the natural numbers,~\nat.

\subsection*{Related Work}

This section only describes those pieces of work which are very directly
related to the material of this thesis. I apologise for any omissions.

I have written up most parts of this thesis before, in papers available
electronically. The main references are \cite{GOM} and Chapters I and II of
\cite{SHDCT}, and to a lesser extent \cite{GECM}. In many places I have added
detail and rigour; indeed, much of the new writing is in the appendices.

The first chapter, \textbf{Bicategories}, is also largely new
writing. However, the results it contains are unlikely to surprise anyone:
they have certainly been in the air for a while, even if they have not been
written up in full detail before. See \cite[9.1]{HerRM}, \cite[p.~8]{SHDCT},
\cite[4.4]{GECM} and \cite[4.3]{HAO} for more or less explicit references to
the idea.  Closely related issues have been considered in the study of
2-monads made by the (largely) Australian school: see, for instance,
\cite{BKP}, \cite{KelCD} and \cite{Pow}. The virtues of the main proof of
this chapter (which is actually in Appendix~\ref{app:unbiased}) are its
directness, and that it uses an operad where a 2-monad might be used instead,
which is more in the spirit of this work. Similar methods to those used here
also provide a way of answering more general questions concerning possible
ways of defining `bicategory', as explained in~\cite{WMC}.

I first wrote up the material of Chapter~\ref{ch:mtione} , \textbf{Operads
and Multicategories}, in \cite{GOM} (and another account appears in
\cite{SHDCT}). At that time the ideas were new to me, but subsequently I
discovered that the definition of $T$-multicategory had appeared in Burroni's
1971 paper \cite{Bur}. Very similar ideas were also being developed, again in
ignorance of Burroni, by Hermida: \cite{HerRM}. However, one important part
of Chapter~\ref{ch:mtione} which does not seem to be anywhere else
is~\ref{sec:alg}, on algebras for a multicategory.

Burroni's paper is in French, which I do not read well. This has had two
effects: firstly, that I have not used it as a source at all, and secondly,
that I cannot accurately tell what is in it and what is not. I have attempted
to make correct attributions, but I may not entirely have succeeded here.

Chapter~\ref{ch:mtitwo}, \textbf{More on Operads and Multicategories}, is a
selection of further topics concerning multicategories. Sections
\ref{sec:struc}--\ref{sec:fib} all appear, more or less, in both \cite{GOM}
and \cite{SHDCT}. Other work related to~\ref{sec:free-multicats} (Free
Multicategories) is described in the paragraph on Appendix~\ref{app:free}
below. A shorter version of section~\ref{sec:endo} is in
\cite{GOM}. Section~\ref{sec:fcm} (on \fc-multicategories) is covered in each
of \cite{FCM}, \cite{GECM} and \cite{GEC}. \fc-multicategories are another of
those ideas that seem to have been in the air; they also seem to be in
\cite{Bur} (p.~280), and appear in \cite[10.2]{HerCSU}. Moreover, Burroni's
section~IV.3 is entitled `$T$-profunctors and $T$-natural transformations'
(in French), and these entities presumably bear some resemblance to the
profunctors and natural transformations discussed in~\ref{sec:bim}.

Chapter~\ref{ch:defn} is \textbf{A Definition of Weak $\omega$-category},
based on the definition given by Batanin in \cite{Bat} (and summarized by
Street in \cite{StrRMB}). I first wrote a version of this chapter in
\cite{SHDCT}. At the time I thought I was writing an account of Batanin's
definition, reshaped and very much simplified but with the same end result
mathematically. In fact, in trying to understand the meaning of a difficult
part of \cite{Bat}, I had made a guess which turned out to be inaccurate (as
Batanin informed me), but still provided a reasonable definition of weak
$\omega$-category. 

As far as originality and novelty go, the upshot for Chapter~\ref{ch:defn} is
this. The chapter contains two main ideas: globular operads and
contractions. Globular operads were proposed in~\cite{Bat}, but in a rather
complicated way; here, we are able to give a one-line definition (`operads
for the free strict $\omega$-category monad'). Contractions were the concept
in \cite{Bat} of which I had made a creative and inaccurate interpretation,
so our two definitions of contraction differ; the definition given here seems
more economical than that in~\cite{Bat}. There is a comparison of the two
strategies at the end of~\ref{sec:the-defn}. Overall, the present definition
of weak $\omega$-category is very economical conceptually, and short too:
given the basic language of general multicategories, it only takes a page or
two~(\pageref{sec:formal}--\pageref{p:end-of-defn}).

Appendix~\ref{app:unbiased}, \textbf{Biased \textit{vs.}\ Unbiased
Bicategories}, is commented on with Chapter~\ref{ch:bicats} above.

Appendix~\ref{app:free}, \textbf{The Free Multicategory Construction}, is
almost exactly the same as the appendix of \cite{GECM}. It is very like the
free monoid construction in Appendix B of \cite{BJT}, although I did not see
this until after writing \cite{GECM}. This is a more subtle free monoid
construction than most: it does not require the tensor (with respect to which
we are taking monoids) either to be symmetric or to preserve sums on each
side. In our context, the latter condition translates to saying that the
functor $T$ preserves sums, where we are trying to form free
$T$-multicategories. This is often not the case: for instance, if $T$ is the
free monoid functor on \Set. There is a version of the free multicategory
construction in Burroni's paper \cite{Bur}~(III.III), but he does insist that
$T$ preserves sums.

Most of Appendix~\ref{app:strict-omega}, \textbf{Strict $\omega$-Categories},
sets out results which are widely assumed (e.g. \cite[\S 10.1]{HerCSU} or
\cite[Ch.~II]{SHDCT}). However, I do not know of another place where the main
result, that strict $\omega$-categories are monadic over globular sets and
the induced monad is cartesian and finitary, is actually proved. The material
in the last section~(\ref{sec:pasting-rep}) is not so widely known, but is a
reworking of results in~\cite{Bat}.

\begin{sloppypar}
Appendix~\ref{app:initial} proves the \textbf{Existence of an Initial
Operad-with-Contraction}. This is new material, and fills a gap left in
\cite{SHDCT}~(II.5). Experts in these matters will probably be able to wave
their hands and say with conviction that the initial object exists, on the
general principle of there being free models for finitary essentially
algebraic theories.
\end{sloppypar}

\subsection*{Acknowledgements}

This work was supported by a PhD grant from EPSRC, by a Research Scholarship
at Trinity College, Cambridge, and subsequently by the Laurence Goddard
Fellowship at St John's College, Cambridge.

The document was prepared in \LaTeX, using Paul Taylor's diagrams package for
some of the diagrams. Audrey Tan helped me prettify it.

I am very grateful to the many people who have helped me in this project, and
would especially like to thank Martin Hyland, who has been a wonderful
supervisor.

\mainmatter

\chapter{Bicategories}
\label{ch:bicats}

The main purpose of this chapter is to provide an alternative definition of
bicategory in which, instead of having a specified identity 1-cell on each
object and a specified binary composite of any pair of adjacent 1-cells, one
has a specified composite of any string of $n$ 1-cells
\[
\gfst{}\gone{}\gblw{}\gone{}\diagspace\cdots\diagspace\gone{}\glst{} 
\]
for each $n\in\nat$. We then prove that this definition is equivalent, in a
strong sense, to the classical definition. The details of the proof are
relegated to Appendix~\ref{app:unbiased}.

This alternative definition of bicategory---which we call an \emph{unbiased
bicategory}---is very natural, and in many ways more natural than the
classical definition. But this is not why it appears in this work: the reason
is that we will need it in Chapter~\ref{ch:defn}, where we show that for
$n=2$, our weak $n$-categories are just unbiased bicategories.

More information on the pedigree of these ideas is contained in the `Related
Work' section of the Introduction.

\section{Review of Classical Material}

Here we review the basic properties of bicategories and state our
terminology. The original definition of bicategory was made in B\'enabou's
paper \cite{Ben}, along with the definition of lax functor (called `morphism'
there). Other references for these definitions are \cite{BB} and
\cite{StrCS}, which also include definitions of transformation and
modification; but we will not need these further concepts here. 

We will typically denote 0-cells (or `objects') of a bicategory \Bee\ by $A$,
$B$, \ldots, 1-cells by $f$, $g$, \ldots\ and 2-cells by $\alpha$, $\beta$,
\ldots, e.g.
\[
A\ctwo{f}{g}{\alpha}B.
\]
The `vertical' composite of 2-cells
\[
\cdot\cthree{}{}{}{\alpha}{\beta}\cdot
\]
is written $\beta\of\alpha$ or $\beta\alpha$, and the `horizontal' composite
of 2-cells
\[
\cdot\ctwo{}{}{\alpha}\cdot\ctwo{}{}{\alpha'}\cdot
\]
is written $\alpha'*\alpha$. We will not need names for the associativity and
unit isomorphisms; when they are all identities, the bicategory is called a
2-category. 

A \emph{lax functor}\/ $\pr{F}{\phi}: \Bee\go\Beep$ (between bicategories
\Bee\ and \Beep) consists of a function $F_0: \Bee_0 \go \Bee'_0$ on objects,
a functor 
\[
F_{A,B}: \homset{\Bee}{A}{B} \go \homset{\Beep}{F_0 A}{F_0 B}
\]
for each pair $A,B$ of objects of \Bee, and `coherence' 2-cells
\[
\phi_{f,g}: Fg\of Ff \go F(g\of f), \diagspace 
\phi_{A}: 1_{FA}\go F1_A
\]
satisfying some axioms. If these 2-cells are all invertible then $F$ is
called a \emph{weak functor} (B\'enabou: `homomorphism'). If they are
identities (so that $Fg\of Ff = F(g\of f)$ and $F1=1$) then $F$ is called a
\emph{strict functor}.

Lax functors can be composed, and this composition obeys strict associativity
and identity laws, so that we obtain a category \Bilax. Moreover,
the class of weak functors is closed under composition, and the same goes for
strict functors, and the identity functor on a bicategory is strict; thus we
have categories
\[
\Bistr \sub \Biwk \sub \Bilax,
\]
all with the same objects. (A more categorical way of putting it is that
there are faithful functors
\[
\Bistr \go \Biwk \go \Bilax
\]
which are the identity on objects, but I will continue to use the \sub\
notation for brevity.)

A \emph{monad} in a bicategory \Bee\ is a lax functor from the terminal
bicategory \One\ to \Bee. Explicitly, this consists of a 0-cell $A$, a 1-cell
$A \goby{t} A$, and 2-cells
\[
A\cthreein{1}{t}{t\of t}{\eta}{\mu}A,
\]
such that the diagrams
\[\left.
\begin{diagram}
t\of 1	&\rTo^{t\eta}	&t\of t		&\lTo^{\eta t}	&1\of t	\\
   	&\rdTo<{\diso}	&\dTo<{\mu}	&\ldTo>{\diso}	&	\\
   	&     		&t  		&     		&   	\\
\end{diagram}
\right.
\ \ \ 
\left.
\begin{diagram}
   		&     		&t\of(t\of t)	&\rTo^{t\mu}	&t\of t\\
   		&\ldLine<{\diso}&   		&     		&   \\
(t\of t)\of t	&     		&   		& 		&\dTo>{\mu}\\
\dTo<{\mu t}	&     		&   		&     		&   \\
t\of t		&     		&\rTo_{\mu}	&     		&t  \\
\end{diagram}
\right.\]
commute.

There is a one-to-one correspondence between bicategories with precisely one
0-cell and monoidal categories. Given such a bicategory, \Bee, there is a
monoidal category whose objects are the 1-cells of \Bee\ and whose morphisms
are the 2-cells, and with $p\otimes q = p\of q$ and $\alpha\otimes\beta =
\alpha*\beta$, where $p$, $q$ are 1-cells of \Bee\ and $\alpha$, $\beta$ are
2-cells. Lax, weak and strict functors between the bicategories then
correspond to lax monoidal functors, (weak) monoidal functors and strict
monoidal functors.

We could equally well have chosen the opposite orientation, so that $p\otimes
q = q\of p$ and $\alpha\otimes\beta = \beta*\alpha$. However, we stick with
our choice. The consequence is that `$\otimes$ and $\of$ go in the same
direction'. (This accounts for the apparently odd reversal of $R$
and $R'$ in Example \ref{egs:fcms}\bref{eg:weak-double}.)

\section{Unbiased Bicategories}

The traditional definition of a bicategory is `biased' towards binary and
nullary compositions, in that only these are given explicit mention. For
instance, there is no specified ternary composite of 1-cells,
$\triple{h}{g}{f} \goesto hgf$, only the derived ones like $h(gf)$ and
$((h1)g)(f1)$. It is necessary to be biased in order to achieve a finite
axiomatization. However, it is useful in this work (and elsewhere) to have a
notion of `unbiased bicategory', in which all arities are treated
even-handedly. In this section we define unbiased bicategory and unbiased
weak functor, and in the next we compare this approach to the classical one.

\begin{defn}
An \emph{unbiased bicategory} \Bee\ consists of 
\begin{itemize}
\item
a class $\Bee_0$, whose elements are called \emph{objects} or
\emph{0-cells}
\item
for each pair $A, B$ of objects, a category \homset{\Bee}{A}{B}, whose
objects are called \emph{1-cells} and whose morphisms are called \emph{2-cells}
\item
for each sequence \range{A_0}{A_n} of objects ($n\geq 0$), a `composition'
functor 
\[
\!\!\!\!\!\!\!\!
\!\!\!\!\!\!\!
\begin{array}{rrcl}
\comp_{\bftuple{A_0}{A_n}}:	&
\homset{\Bee}{A_{n-1}}{A_n} \times\cdots\times \homset{\Bee}{A_0}{A_1}	&
\go	&
\homset{\Bee}{A_0}{A_n},	\\
&\bftuple{f_n}{f_1}	&\goesto	&\bo f_n \of\cdots\of f_1 \bc,	\\
&\bftuple{\alpha_n}{\alpha_1}&\goesto	&\bo \alpha_n *\cdots* \alpha_1 \bc,
\end{array}
\]
where the $f_i$'s are 1-cells and the $\alpha_i$'s are 2-cells
\item
for each double sequence
\bftuple{\bftuple{f_1^1}{f_1^{k_1}}}{\bftuple{f_n^1}{f_n^{k_n}}} of 1-cells
such that the composite $\bo f_n^{k_n} \of\cdots\of f_n^1 \of \cdots \of
f_1^{k_1} \of\cdots\of f_1^1 \bc$ makes sense, an invertible 2-cell
\begin{eqnarray*}
\lefteqn{\gamma_{\bftuple{\bftuple{f_1^1}{f_1^{k_1}}}{\bftuple{f_n^1}{f_n^{k_n}}}}:}\\
&\bo \bo f_n^{k_n} \of\cdots\of f_n^1 \bc \of\cdots\of 
\bo f_1^{k_1} \of\cdots\of f_1^1 \bc \bc
\goiso
\bo f_n^{k_n} \of\cdots\of f_n^1 \of \cdots \of f_1^{k_1} \of\cdots\of f_1^1
\bc
\end{eqnarray*}
\item
for each 1-cell $f$, an invertible 2-cell
\[
\iota_f: f \goiso \bo f \bc
\]
\end{itemize}
with the following properties:
\begin{itemize}
\item
$\gamma_{\bftuple{\bftuple{f_1^1}{f_1^{k_1}}}{\bftuple{f_n^1}{f_n^{k_n}}}}$
is natural in each of the $f_i^j$'s, and $\iota_f$ is natural in $f$
\item
associativity: for any triple sequence $((( f_{p,q,r}
)_{r=1}^{k_p^q} )_{q=1}^{m_p} )_{p=1}^{n}$ of 1-cells
such that the following composites make sense, the diagram
\[
\!\!\!\!\!\!\!\!
\!\!\!\!\!\!\!\!
\!\!\!\!\!\!\!\!
\!\!\!\!\!\!\!\!
\!\!\!\!\!\!\!\!
\begin{diagram}[width=5em,tight,scriptlabels,leftflush]
&&\scriptstyle\bo 
      \bo 
          \bo f_{n, m_n, k_n^{m_n}} \of\cdots\of f_{n, m_n, 1} \bc
          \of\cdots\of
          \bo f_{n, 1, k_n^1} \of\cdots\of f_{n, 1, 1} \bc
      \bc
      \of\cdots\of
      \bo 
          \bo f_{1, m_1, k_1^{m_1}} \of\cdots\of f_{1, m_1, 1} \bc
          \of\cdots\of
          \bo f_{1, 1, k_1^1} \of\cdots\of f_{1, 1, 1} \bc 
      \bc 
  \bc
& & \\
&\ldTo<{\bo \gamma_{D_n} *\cdots * \gamma_{D_1} \bc}
& &\rdTo>{\gamma_{D'}} & \\
\scriptstyle\bo 
    \bo f_{n, m_n, k_n^{m_n}} \of\cdots\of f_{n, 1, 1} \bc
    \of\cdots\of
    \bo f_{1, m_1, k_1^{m_1}} \of\cdots\of f_{1, 1, 1} \bc
\bc
& & & & 
\scriptstyle\bo
    \bo f_{n, m_n, k_n^{m_n}} \of\cdots\of f_{n, m_n, 1} \bc
    \of\cdots\of
    \bo f_{1, 1, k_1^1} \of\cdots\of f_{1, 1, 1} \bc
\bc \\
&\rdTo<{\gamma_D} & &\ldTo>{\gamma_{D''}} & \\
& &\scriptstyle\bo f_{n, m_n, k_n^{m_n}} \of\cdots\of f_{1, 1, 1} \bc & & \\
\end{diagram}
% \end{flushleft}
\]
commutes, where the double sequences $D_p, D, D', D''$ are
\begin{eqnarray*}
D_p	&=& 
\bftuple{\bftuple{f_{p, 1, 1}}{f_{p, 1, k_p^1}}}{\bftuple{f_{p, m_p,
1}}{f_{p, m_p, k_p^{m_p}}}},	\\
D	&=&	
\bftuple{\bftuple{f_{1, 1, 1}}{f_{1, m_1, k_1^{m_1}}}}{\bftuple{f_{n, 1,
1}}{f_{n, m_n, k_n^{m_n}}}},	\\
D'	&=&	
(
\bftuple{\bo f_{1, 1, k_1^1} \of\cdots\of f_{1, 1, 1} \bc}{\bo
f_{1, m_1, k_1^{m_1}} \of\cdots\of f_{1, m_1, 1} \bc}
, \ldots,	\\
&&
\bftuple{\bo f_{n, 1, k_n^1} \of\cdots\of f_{n, 1, 1} 
\bc}{\bo f_{n, m_n, k_n^{m_n}} \of\cdots\of f_{n, m_n, 1} \bc} ),	\\
D''	&=&
\bftuple{\bftuple{f_{1, 1, 1}}{f_{1, 1, k_1^1}}}{\bftuple{f_{n, m_n,
1}}{f_{n, m_n, k_n^{m_n}}}}
\end{eqnarray*}
\item
identity: for any composable sequence \bftuple{f_1}{f_n} of 1-cells, the
diagrams
\[
\begin{diagram}[scriptlabels]
\bo f_n \of\cdots\of f_1 \bc	
&\rTo^{\bo \iota_{f_n} *\cdots* \iota_{f_1} \bc}
&\bo \bo f_n \bc \of\cdots\of \bo f_1 \bc \bc	\\
&\rdTo<{1}	
&\dTo>{\gamma_{\bftuple{(f_1)}{(f_n)}}}	\\
&& \bo f_n \of\cdots\of f_1 \bc	\\
\end{diagram}
\mbox{\hspace{3em}}
\begin{diagram}[scriptlabels]
\bo \bo f_n \of\cdots\of f_1 \bc \bc
&\lTo^{\iota_{\bo f_n \of\cdots\of f_1 \bc}}
&\bo f_n \of\cdots\of f_1 \bc	\\
\dTo<{\gamma_{(\bftuple{f_1}{f_n})}}
&\ldTo>{1}	& \\
\bo f_n \of\cdots\of f_1 \bc & & \\
\end{diagram}
\]
commute.
\end{itemize}
\end{defn}

\begin{remarks}{rmks:u-bicat}
\item
The associativity axiom is less fearsome than it might appear. It says that any
two ways of removing brackets are equivalent, just as the associativity
axiom does for a monad such as `free group' on \Set. If we allow different
styles of brackets then it says, for instance, that
\begin{diagram}[leftflush=2.5em]
 & &\{[(h \of g) \of (f \of e)] \of [(d \of c) \of (b \of a)]\} & & \\
 &\ldTo<{(\gamma * \gamma)} & &\rdTo>{\gamma} & \\
\{[h \of g \of f \of e] \of [d \of c \of b \of a]\}
& & & &
\{(h \of g) \of (f \of e) \of (d \of c) \of (b \of a)\} \\
 &\rdTo<{\gamma} & &\ldTo>{\gamma} & \\
 & &\{h \of g \of f \of e \of d \of c \of b \of a \} & & \\
\end{diagram}
commutes.

\item	\label{rmk:axioms-obvious}
The coherence axioms for an unbiased bicategory are rather obvious, in
contrast to the situation for classical bicategories: they look just like the
associativity and unit axioms for a monoid. 

\item
An \emph{unbiased monoidal category} may be defined as an unbiased bicategory
with precisely one object; we would then write $\otimes$ in place of both
$\of$ and $*$.

\item	\label{rmk:lax-bicat}
If we drop the condition that $\gamma$ and $\iota$ are invertible, then we
obtain what might be called a \emph{lax} or \emph{relaxed bicategory}. (Or
perhaps `colax' would be more appropriate.)  A one-object lax bicategory is
then a relaxed monoidal category in the sense of \cite[4.4]{GECM}. In the
other direction, let us define an \emph{unbiased 2-category} as an unbiased
bicategory in which the components of $\gamma$ and $\iota$ are all
identities. (Clearly unbiased 2-categories are in one-to-one correspondence
with ordinary 2-categories.)  So we have three classes of structures:
\[
\!\!\!\!\!
\{ \mbox{unbiased 2-categories} \}
\sub
\{ \mbox{unbiased bicategories} \}
\sub
\{ \mbox{lax bicategories} \}.
\]
For the moment this is just a statement about classes (large sets), but soon
we will define maps between these structures and thus be able to compare the
\emph{categories} they form. 

\item	\label{rmk:Cat-Gph}
We have given a very explicit definition of unbiased bicategory, but a more
abstract version is possible. There is a 2-category
\label{p:Cat-Gph}%
$\Cat\hyph\Gph$, an object of which is a set $\Bee_0$ together with an
indexed family
\[
(\homset{\Bee}{B}{B'})_{B, B' \in \Bee_0}
\] 
of categories (a `\Cat-graph'). An arrow $F: \Bee \go \cat{C}$ consists of a
function $F_0: \Bee_0 \go \cat{C}_0$ and a functor
\[
F_{B,B'}: \homset{\Bee}{B}{B'} \go \homset{\cat{C}}{F_0 B}{F_0 B'}
\]
for each $B, B' \in \Bee_0$. There is only a 2-cell
\[
\Bee\ctwo{F}{G}{}\cat{C}
\]
if $F_0 = G_0$, and in this case such a 2-cell $\alpha$ is a family of
natural transformations $\alpha_{B,B'}: F_{B,B'} \go G_{B,B'}$. Now, there is
a 2-monad `free 2-category' on $\Cat\hyph\Gph$, and a (small) unbiased
bicategory is, in a suitable sense, a weak algebra for this 2-monad. The
definition of relaxed monoidal category in \cite[4.4]{GECM} implicitly uses
this approach, but with lax algebras rather than weak algebras. For more on
this point of view, see~\cite{KS} and~\cite{Pow}. We also use this approach
in Appendix~\ref{app:unbiased}. 

\item
The notation $\bo f_n \of\cdots\of f_1 \bc$ for the composite of a diagram
\[
A_0 \goby{f_1} A_1 \goby{f_2} \cdots \goby{f_n} A_n
\]
is sometimes inadequate in the case $n=0$. When $n=0$ the data to be
composed is just a single object $A_0$, and we might prefer to write
$1_{A_0}$ rather than the standard notation, $\bo \bc$.
\end{remarks}

\begin{defn}
Let \Bee\ and \Beep\ be unbiased bicategories. An \emph{unbiased lax functor}
$\pr{F}{\phi}: \Bee \go \Beep$ consists of
\begin{itemize}
\item
a function $F_0: \Bee_0 \go \Bee'_0$ (usually just written $F$)
\item
for each $A, B \in \Bee_0$, a functor $F_{A,B}: \homset{\Bee}{A}{B} \go
\homset{\Beep}{F_0 A}{F_0 B}$ (again, usually just written $F$)
\item
for each composable sequence \bftuple{f_1}{f_n} of 1-cells, a 2-cell
\[
\phi_{\bftuple{f_1}{f_n}}:
\bo Ff_n \of\cdots\of Ff_1 \bc
\go
F \bo f_n \of\cdots\of f_1 \bc,
\]
\end{itemize}
with the properties that
\begin{itemize}
\item
$\phi_{\bftuple{f_1}{f_n}}$ is natural in each $f_i$
\item
for each double sequence
\bftuple{\bftuple{f_1^1}{f_1^{k_1}}}{\bftuple{f_n^1}{f_n^{k_n}}} of 1-cells
such that the following composites make sense, the diagram
\[
\begin{diagram}
\bo
   \bo Ff_n^{k_n} \of\cdots\of Ff_n^1 \bc
   \of\cdots\of
   \bo Ff_1^{k_1} \of\cdots\of Ff_1^1 \bc
\bc
&\rTo^{\gamma'_{\bftuple{\bftuple{Ff_1^1}{Ff_1^{k_1}}}{\bftuple{Ff_n^1}{Ff_n^{k_n}}}}}
&\bo Ff_n^{k_n} \of\cdots\of Ff_1^1 \bc	\\
\dTo>{\bo \phi_{\bftuple{f_n^1}{f_n^{k_n}}} *\cdots*
          \phi_{\bftuple{f_1^1}{f_1^{k_1}}} \bc}
& & \\
\bo
    F \bo f_n^{k_n} \of\cdots\of f_n^1 \bc
    \of\cdots\of
    F \bo f_1^{k_1} \of\cdots\of f_1^1 \bc
\bc
&&\dTo<{\phi_{\bftuple{f_1^1}{f_n^{k_n}}}}	\\
\dTo>{\phi_{\bftuple{\bo f_1^{k_1} \of\cdots\of f_1^1 \bc}{\bo f_n^{k_n}
\of\cdots\of f_n^1 \bc}}} 
& & \\
F \bo
     \bo f_n^{k_n} \of\cdots\of f_n^1 \bc
     \of\cdots\of
     \bo f_1^{k_1} \of\cdots\of f_1^1 \bc
  \bc
&\rTo_{F\gamma_{\bftuple{\bftuple{f_1^1}{f_1^{k_1}}}{\bftuple{f_n^1}{f_n^{k_n}}}}}
&F \bo f_n^{k_n} \of\cdots\of f_1^1 \bc \\
\end{diagram}
\]
commutes
\item
for each 1-cell $f$, the diagram
\begin{diagram}
Ff	&\rTo^{\iota'_{Ff}}	&\bo Ff \bc	\\
\dEquals&			&\dTo>{\phi_{\tuplebts{f}}}	\\
Ff	&\rTo_{F\iota_f}	&F\bo f \bc	\\
\end{diagram}
commutes.
\end{itemize}
An \emph{unbiased weak functor} is an unbiased lax functor \pr{F}{\phi} for
which each $\phi_{\bftuple{f_1}{f_n}}$ is invertible.  An \emph{unbiased
strict functor} is an unbiased lax functor \pr{F}{\phi} for which each
$\phi_{\bftuple{f_1}{f_n}}$ is the identity (so that $F$ preserves composites
and identities strictly). 
\end{defn}

We noted in Remark~\bref{rmk:axioms-obvious} that the coherence axioms for an
unbiased bicategory were rather obvious, having the shape of the axioms for a
monoid or monad. Perhaps the coherence axioms for an unbiased lax functor are
a little less obvious; however, they are the same shape as the axioms for a
monad functor given in Street's paper \cite{StrFTM}, and in any case seem to
be quite canonical in some vague sense.

Naturally, we would like to be able to compose lax functors. Given unbiased
lax functors
\[
\Bee \goby{\pr{F}{\phi}} \Beep \goby{\pr{F'}{\phi'}} \cat{B''},
\]
define the composite \pr{G}{\psi} by $G_0 = F'_0 \of F_0$, $G_{A,B} =
F'_{FA,FB} \of F_{A,B}$, and by taking $\psi_{\bftuple{f_1}{f_n}}$ to be the
composite of
\[
\!\!\!\!
\bo GFf_n \of\cdots\of GFf_1 \bc
\goby{\phi'_{\bftuple{Ff_1}{Ff_n}}}
G \bo Ff_n \of\cdots\of Ff_1 \bc
\goby{G\phi_{\bftuple{f_1}{f_n}}}
GF \bo f_n \of\cdots\of f_1 \bc.
\]
Also define the identity unbiased lax functor \pr{G}{\psi} on an
unbiased bicategory \Bee\ by $G_0=\id$, $G_{A,B}=\id$, and
$\psi_{\bftuple{f_1}{f_n}} = \id$. It is straightforward to check that
composition is associative and that the identity functors live up to their
name. We therefore obtain a category \UBilax\ of unbiased bicategories and
unbiased lax functors. Evidently there are subcategories
\[
\UBistr \sub \UBiwk \sub \UBilax,
\]
with the same objects and with arrows which are, respectively, unbiased
strict functors and unbiased weak functors. 

In fact, the definitions of unbiased lax functor and of their composites and
identities work just as well for \emph{lax}
bicategories~(\ref{rmks:u-bicat}\bref{rmk:lax-bicat}). So there are $3\times
3 = 9$ possible categories we might consider: for both the objects and the
arrows, we choose one of `strict', `weak' or `lax'. With what I hope is
obvious notation, the inclusions are as follows:
\begin{diagram}[height=1.2em]
\fcat{LBicat}_\mr{str}	&\sub	&\fcat{LBicat}_\mr{wk}	&\sub
&\fcat{LBicat}_\mr{lax} \\
\rotsub	&	&\rotsub	&	&\rotsub	\\	
\UBistr	&\sub	&\UBiwk		&\sub	&\UBilax	\\
\rotsub	&	&\rotsub	&	&\rotsub	\\	
\fcat{U2\hyph Cat}_\mr{str}	&\sub	&\fcat{U2\hyph Cat}_\mr{wk}	&\sub
&\fcat{U2\hyph Cat}_\mr{lax}. \\
\end{diagram}
Of these nine, we might consider the three on the diagonal (bottom-left to
top-right) to be the most conceptually natural. We will not actually need to
discuss anything except for the middle row in the rest of this work. However,
these remarks demonstrate the cleanliness of the unbiased theory when
compared to the biased (classical) theory. In the latter, the top row is
obscured---that is, there is no very satisfactory way to weaken the classical
definition of bicategory to get a lax version. Admittedly, one can drop the
condition that the classical associativity and unit maps are isomorphisms (as
in \cite{BorxI}, after Definition 7.7.1); but somehow this does not seem
quite right.

Another advertisement for the unbiased theory follows. To give it we need
some preliminary basic constructions. Firstly, for any bicategory \Bee\
(biased or unbiased), there is an opposite bicategory $\Bee^\op$, obtained by
reversing the 1-cells only: thus to each 2-cell
\[
A\ctwo{f}{g}{\alpha}B
\]
in \Bee\ there corresponds a 2-cell
\[
A\ctwoop{f}{g}{\alpha}B
\]
in $\Bee^\op$. Secondly, one may form the product $\cat{A} \times \Bee$ of
any two (biased or unbiased) bicategories in the obvious way (and this is the
categorical product in each of the lax, weak and strict contexts). Thirdly,
there is a 2-category \Cat\ of all (small) categories, functors and natural
transformations, and there is a corresponding unbiased 2-category \Cat.

Now, we would like to form a functor
\[
\begin{array}{rrcl}
\Hom:	&\Bee^\op \times \Bee	&\go	&\Cat,	\\
	&\pr{A}{B}		&\goesto&\homset{\Bee}{A}{B}\\
\end{array}
\]
for each \Bee\ (ignoring questions of size). In the biased case this is not
possible without making an arbitrary choice. For if $A'\goby{f}A$ and
$B\goby{g}B'$ in \Bee\ then applying \Hom\ should give us a function
\[
\homset{\Bee}{A}{B} \go \homset{\Bee}{A'}{B'},
\]
and this might reasonably be either $p \goesto (g\of p)\of f$ or $p \goesto
g\of (p\of f)$. Although we could, say, consistently choose the first option
and thereby get a weak functor \Hom, neither choice is `canonical'. However,
in the unbiased case one has a ternary composite $\bo g\of p\of f \bc$,
giving a \emph{canonical} weak functor
\[
\Hom: \Bee^\op\times\Bee \go \Cat.
\]

\section{Biased \emph{vs.}\ Unbiased}	\label{sec:versus}

In this section we define a forgetful functor $V: \UBilax \go \Bilax$, which
turns out to be full, faithful and surjective on objects. (Proofs
are deferred to Appendix~\ref{app:unbiased}.) Thus the categories of biased
and unbiased bicategories, with lax functors as maps, are equivalent; and the
same in fact goes for weak functors, although not strict ones. So we will
more or less be able to ignore the biased-unbiased distinction.

The primary reason for setting out the theory of unbiased bicategories in
this thesis is that in Chapter~\ref{ch:defn} we give a definition of weak
$n$-category, and a weak 2-category is exactly an unbiased bicategory. We
therefore want to know that unbiased and classical bicategories are
essentially the same, as a test of the reasonability of our proposed
definition. 

This somewhat practical motivation provides an answer to a question which the
reader may have been asking: where are the unbiased transformations and
modifications? Quite simply, we don't mention them because we don't need
them: the unbiased and classical theories can be compared without going
above the level of functors.

An equally important answer is that transformations and modifications between
unbiased bicategories are not defined because there seems to be no properly
`unbiased' way to do it. Of course, we can `cheat' by transporting the
definitions from \Bilax\ along the functor
\[
V: \UBilax \go \Bilax.
\]
This immediately gives a coherence theorem: every unbiased bicategory is
biequivalent to an unbiased 2-category. More honest coherence results, of the
form `every diagram commutes', appear in Appendix~\ref{app:unbiased}. 

Note also that the equivalence $\UBilax \eqv \Bilax$ is two levels better
than we might have expected: if \Bee\ and \Beep\ are two unbiased
bicategories with $V(\Bee) = V(\Beep) \in \Bilax$, then \Bee\ and \Beep\ are
not just biequivalent in \UBilax, or even just equivalent: they are actually
isomorphic. Put another way, we have a comparison which takes place at the
1-dimensional level, without having to resort to 2- or 3-dimensional
structures. 

To business: let us define the forgetful functor $V$. Given an unbiased
bicategory \Bee, attempt to define a biased bicategory $\cat{C} = V(\Bee)$
by:
\begin{itemize}
\item
$\cat{C}_0 = \Bee_0$
\item
$\homset{\cat{C}}{A}{B} = \homset{\Bee}{A}{B}$
\item
composition 
\[
\homset{\cat{C}}{B}{C} \times \homset{\cat{C}}{A}{B}
\go \homset{\cat{C}}{A}{C}
\]
in \cat{C} is
\[
\comp_{\triple{A}{B}{C}}: 
\homset{\Bee}{B}{C} \times \homset{\Bee}{A}{B}
\go \homset{\Bee}{A}{C}
\]
\item
the identity in \cat{C} on an object $A$ is (the image of) 
\[
\comp_{\tuplebts{A}}: \wun \go \homset{\Bee}{A}{A}
\]
\item
the associativity isomorphism $(h\of g)\of f \go h\of (g\of f)$ is the
composite of the 2-cells
\[
\!\!\!\!\!\!\!\!
\!\!\!\!\!\!\!\!
\!\!\!\!\!\!\!\!
\!\!\!\!\!\!\!\!
\bo \bo h \of g \bc \of f \bc
\goby{\bo 1 * \iota_f \bc}
\bo \bo h \of g \bc \of \bo f \bc \bc
\goby{\gamma_{\pr{\tuplebts{f}}{\pr{g}{h}}}}
\bo h \of g \of f \bc
\goby{\gamma^{-1}_{\pr{\pr{f}{g}}{\tuplebts{h}}}}
\bo \bo h \bc \of \bo g \of f \bc \bc
\goby{\bo \iota^{-1}_h * 1 \bc}
\bo h \of \bo g \of f \bc \bc
\]
\item
the left unit isomorphism $1\of f \go f$ is the composite of the 2-cells
\[
\bo \bo \bc \of f \bc
\goby{\bo 1 * \iota_f \bc}
\bo \bo \bc \of \bo f \bc \bc
\goby{\gamma_{\pr{\tuplebts{f}}{\tuplebts{}}}}
\bo f \bc
\goby{\iota^{-1}_f}
f,
\]
and dually for the right unit.
\end{itemize}
Given an unbiased lax functor $\pr{F}{\phi}: \Bee \go \Beep$, attempt to
define a lax functor $\pr{G}{\psi} = V\pr{F}{\phi}: V(\Bee) \go V(\Beep)$ by
\[
G_0=F_0,
\diagspace
G_{A,B} = F_{A,B},
\diagspace
\psi_{f,g} = \phi_{\pr{f}{g}},
\diagspace
\psi_A = \phi_{()}.
\]
Here the symbol $\phi_{()}$ denotes $\phi_{\bftuple{f_1}{f_n}}$ in the case
$n=0$, where
\[
A = A_0 \goby{f_1} A_1 \goby{f_2} \cdots \goby{f_n} A_n.
\]

In Appendix~\ref{app:unbiased} we prove:
\begin{thm}	\label{thm:biased-comparison}
With these definitions, 
\begin{enumerate}
\item $V(\Bee)$ is a bicategory and $V\pr{F}{\phi}$ is a lax functor
\item $V$ preserves composition and identities, so forms a functor 
\[
\UBilax \go \Bilax
\] 
\item $V$ is full, faithful and surjective on objects.
\end{enumerate}
\end{thm}

If \pr{F}{\phi} is a weak (respectively, strict) functor then
$V\pr{F}{\phi}$ is one too, so $V$ restricts to give functors
\[
V_\mr{wk}: \UBiwk \go \Biwk, \diagspace
V_\mr{str}: \UBistr \go \Bistr.
\] 
In the appendix we prove:
\begin{cor}	\label{cor:wk-biased-comparison}
The restricted functor $V_\mr{wk}:\UBiwk \go \Biwk$ is also full, faithful
and surjective on objects.
\end{cor}
Thus $\UBilax \eqv \Bilax$ and $\UBiwk \eqv \Biwk$.

Finally, what about the strict case---is $V_\mr{str}$ an equivalence of
categories? Certainly $V_\mr{str}$ is surjective on objects and faithful
(since the same is true of $V$), so the only question is whether it is
full. It is not. For let \cat{C} be any bicategory, and construct from
\cat{C} an unbiased bicategory \cat{L} with $V(\cat{L}) = \cat{C}$, defining
composition in \cat{L} by associating to the left: e.g.\ the composite
$(f_4\of f_3\of f_2\of f_1)$ in \cat{L} is the composite $((f_4 \of f_3) \of
f_2) \of f_1$ in \cat{C}. (Appendix~\ref{app:unbiased} shows that this
construction is possible.) Dually, define an unbiased bicategory \cat{R} with
$V(\cat{R}) = \cat{C}$ by associating to the right. If $F: \cat{L} \go
\cat{R}$ is an unbiased strict functor with $V(F) = 1_{\cat{C}}$ then $F$
must be the identity (since the data for an unbiased strict functor is just a
graph map), and so $\cat{L} = \cat{R}$. But we can choose a bicategory
\cat{C} in which $(h\of g)\of f \neq h\of (g\of f)$ for some 1-cells $f$,
$g$, $h$, so that $\cat{L} \neq \cat{R}$. Hence the identity on \cat{C} does
not lift to a strict functor $\cat{L} \go \cat{R}$, and therefore
$V_\mr{str}$ is not full.

\chapter{Operads and Multicategories}
\label{ch:mtione}

In this chapter we introduce the language of operads and multicategories to
be used in the rest of the thesis. The simplest kind of operad---a
\emph{plain operad}---consists of a sequence $C(0)$, $C(1)$, \ldots
of sets together with an `identity' element of $C(1)$ and `composition'
functions
\[
C(n)\times C(k_{1})\times\cdots\times C(k_{n})
\go C(k_{1}+\cdots +k_{n}),
\]
obeying associativity and identity laws. (In the original definition,
\cite{MayGIL}, the $C(n)$'s were not just sets but spaces with symmetric
group action. Our operads never have symmetric group actions.) The simplest
kind of multicategory---a \emph{plain multicategory}---consists of a
collection $C_{0}$ of objects, and arrows
\[
a_{1}, \ldots, a_{n} \goby{\theta} a
\]
($a_{1}, \ldots, a_{n}, a \elt C_{0}$), together with composition functions
and identity elements obeying associativity and unit laws. (See
\cite[p.~103]{Lam} for the details.) A plain operad is therefore a one-object
plain multicategory.

The general idea now is that there's nothing special about \emph{sequences}
of objects: the domain of an arrow might form another shape instead, such as
a tree of objects or just a single object (as in a normal category). Indeed,
the objects do not even need to form a set. Maybe a graph or a category would
do just as well. Together, what these generalizations amount to is the
replacement of the free-monoid monad on \Set\ with some other monad on some
other category.

This generalization is put into practice as follows. The graph structure of a
plain multicategory is a diagram
\[
\begin{slopeydiag}
	&	&C_1	&	&	\\
	&\ldTo<{\dom}	&	&\rdTo>{\cod}	&	\\
TC_0	&	&	&	&C_0	\\
\end{slopeydiag}
\]
in \Set, where $T$ is the free-monoid monad. Now, just as a (small) 
category can be described as a diagram 
\[
\begin{slopeydiag}
		&	&D_1		&	&	\\
		&\ldTo	&		&\rdTo	&	\\
D_0		&	&		&	&D_0	\\
\end{slopeydiag}
\]
in \Set\ together with identity and composition functions
\[
D_{0}\go D_{1},
\diagspace
D_{1}\times_{D_{0}}D_{1}\go D_{1}
\]
satisfying some axioms, so we may describe the multicategory structure on
\gph{C_1}{C_0} 
by manipulation of certain diagrams in \Set. In general, we take a category
\Eee\ and a monad $T$ on \Eee\ satisfying some simple conditions, and
define `\Cartpr-multicategory'. Thus a category is a
\pr{\Set}{\id}-multicategory.

Section~\ref{sec:cart-mon} describes the simple conditions on \Eee\ and $T$
required in order that everything that follows will work. Many examples are
given. Section~\ref{sec:multi} explains what \Cartpr-multicategories are,
and what \Cartpr-operads are---namely, one-object
\Cartpr-multicategories. Section~\ref{sec:alg} defines and explains
\emph{algebras} for multicategories, which are a generalization
of \Set-valued functors on a category. If an operad is thought of as a kind
of algebraic theory (in which the elements of $C(n)$ are $n$-ary operations)
then an algebra for an operad is a model of that theory.

\section{Cartesian Monads}	\label{sec:cart-mon}

In this section we introduce the conditions required of a monad \Mnd\ on a
category \Eee, in order that we may (in~\ref{sec:multi}) define the notions
of \Cartpr-multicategory and \Cartpr-operad. The conditions are that the
category and the monad are both cartesian, as defined now.

\begin{defn}
A category is called \emph{cartesian} if it has all finite limits.
\end{defn}

\begin{defn}		%\label{defn-cart}
A monad \Mnd\ on a category \Eee\ is called \emph{cartesian} if
\begin{enumerate}
\item $\eta$ and $\mu$ are cartesian natural transformations, i.e.\ for any
$X \goby{f} Y$ in \Eee\ the naturality squares
\[\left.
\begin{diagram}
X	&\rTo^{\eta_{X}}&TX		\\
\dTo<{f}&		&\dTo>{Tf}	\\
Y	&\rTo^{\eta_{Y}}&TY		\\
\end{diagram}
\right.
\diagspace
\left.
\begin{diagram}
T^2 X		&\rTo^{\mu_{X}}	&TX		\\
\dTo<{T^2 f}	&		&\dTo>{Tf}	\\
T^2 Y		&\rTo^{\mu_{Y}}	&TY		\\
\end{diagram}
\right.\]
are pullbacks, and %\label{nts-cart}
\item $T$ preserves pullbacks. %\label{pb-pres}
\end{enumerate}
\end{defn}

We often write $T$ to denote the whole monad \Mnd, as is customary.

It would perhaps be more consistent to call a category cartesian just if it
has pullbacks, and indeed this is all that is necessary in order to make the
theory of general multicategories work. However, all of our examples have a
terminal object too (and therefore all finite limits), and it is convenient
to assume that this is always the case. For instance, the definition of
\Cartpr-operad only makes sense when \Eee\ has a terminal object.

\begin{eg}{eg:cart-mnds}

\item The identity monad on any category is clearly cartesian.

\item 	\label{eg:mon-monoids}
Let $\Eee = \Set$ and let $T$ be the monoid monad, i.e. the
monad arising from the adjunction
\begin{diagram}
\fcat{Monoid}	&\pile{\rTo \\ \ \ \top \\ \lTo}	&\Set.
\end{diagram}
Certainly \Eee\ is cartesian. It is easy to calculate that $T$, too, is
cartesian (\cite[1.4(ii)]{GOM}), although the theory explained in
Example~\bref{eg:mon-CJ} below renders this unnecessary.

\item 	\label{eg:comm-not-cart}
A non-example. Let $\Eee=\Set$ and let \Mnd\ be the free commutative monoid
monad. This is not cartesian: e.g.\ the naturality square for $\mu$ at $2\go
1$ is not a pullback. See also Example~\ref{egs:multicats}\bref{eg:sym-ops}
for some related thoughts.

\item 	\label{eg:mon-CJ}
Let $\Eee = \Set$. Any finitary algebraic theory gives a monad
on \Eee; which are cartesian? Without answering this question
completely, we indicate a certain class of theories which do give
cartesian monads. An equation (made up of variables and finitary
operations) is said to be \emph{strongly regular} if the same variables
appear in the same order, without repetition, on each side. Thus
\[\begin{array}{ccc}
(x.y).z=x.(y.z) & \mr{and} & (x\uparrow y)\uparrow z = x\uparrow(y.z),\\
\end{array}\]
but not
\[\begin{array}{cccc}
x+(y+(-y))=x, & x.y=y.x & \mr{or} & (x.x).y=x.(x.y),\\
\end{array}\]
qualify. A theory is called \emph{strongly regular} if it can be presented by
operations and strongly regular equations. In
Example~(\ref{eg:mon-monoids}), the only property of the theory of
monoids that we actually needed was its strong regularity: for in general,
the monad yielded by any strongly regular theory is cartesian.

This last result, and the notion of strong regularity, are due to Carboni and
Johnstone. They show in \cite{CJ} (Proposition 3.2 via Theorem 2.6) that a
theory is strongly regular if and only if $\eta$ and $\mu$ are cartesian
natural transformations and $T$ preserves wide pullbacks. A \emph{wide
pullback} is by definition a limit of shape
\begin{diagram}
\cdot	&\cdot	&\cdot		&	&\cdots	&		&\cdot	\\
	&\rdTo(4,2)&\rdTo(3,2)	&\rdTo(2,2)&	&\ldTo(2,2)	&	\\
	&	&		&	&\cdot	&		&,	\\
\end{diagram}
where the top row is a set of any size (perhaps infinite). When the set is of
size 2 this is an ordinary pullback, so the monad from a strongly regular
theory is indeed cartesian. (Examples~\bref{eg:mon-exceptions},
\bref{eg:mon-action} and \bref{eg:mon-tree} can also be found in \cite{CJ}.)%

\item 	\label{eg:mon-exceptions}
Let $\Eee = \Set$, let $E$ be a fixed set, and let + denote binary
coproduct: then the endofunctor $\dashbk +E$ on \Eee\ has a natural monad
structure. This monad is cartesian, corresponding to the algebraic theory
consisting only of one constant for each member of $E$. In particular, if
$E=1$ then this is the theory of pointed sets.

\item 	\label{eg:mon-action}
Let $\Eee=\Set$ and let $M$ be a monoid: then the endofunctor
$M\times\dashbk$ on \Eee\ has a natural monad structure. This monad is
cartesian, corresponding to an algebraic theory consisting only of unary
operations. 

\item 	\label{eg:mon-tree}
Let $\Eee = \Set$, and consider the finitary algebraic theory
on $\Eee$ generated by one $n$-ary operation for each $n\elt\nat$ and
no equations. This theory is strongly regular, so the induced monad
\Mnd\ on $\Eee$ is cartesian.

If $X$ is any set then $TX$ can be described inductively by
\begin{itemize}
	\item if $x\elt X$ then $x\elt TX$
	\item if $\range{t_1}{t_n} \elt TX$ then
$\abftuple{t_1}{t_n}\elt TX$.
\end{itemize}
We can draw any element of $TX$ as a tree with leaves labelled by
elements of $X$:
\begin{itemize}
	\item $x\elt X$ is drawn as \nl{x}
	\item if \range{t_1}{t_n} are drawn as \range{T_1}{T_n} then
\abftuple{t_1}{t_n} is drawn as
$\left.
\begin{tree}
\raisebox{1ex}{\ensuremath{T_1}}&	&\raisebox{1ex}{\ensuremath{T_2}}&	&	&\cdots	& & &\raisebox{1ex}{\ensuremath{T_n}}\\
	&\rt{4}	&	&\rt{2}	&	&	& &\lt{4} &	\\
	&	&	&	&\node	&	& & &		\\
\end{tree}
\right.$~%
, or if $n=0$, as $\left.\EmptyOne\right.$.
\end{itemize}
Thus the element $\atuplebts{\atuplebts{x_{1},x_{2},\atuplebts{}}, x_{3},
\atuplebts{x_{4},x_{5}}}$ of $TX$ is drawn as
\[\left.\Oak{x_1}{x_2}{x_3}{x_4}{x_5}\right. .\]

The unit $X\go TX$ is $x\goesto\nl{x}$, and multiplication $T^2 X
\go TX$ takes a $TX$-labelled tree (e.g.\ 
\[
\left.\Pear{t_1}{t_2}{}\right.\ ,
\]
with
\[
t_{1}=\left.\Orange{x_1}{x_2}{}\right.
\mr{\ and\ }
t_{2}=\left.\Apple{x_3}{x_4}{}\right. \ ) 
\]
and gives an $X$-labelled tree by substituting at the leaves (here,
\[
\left.\MixedFruit{x_1}{x_2}{x_3}{x_4}{}\right. \ ).
\]

\item
On the category \Cat\ of small categories and functors, there is the free
strict monoidal category monad. Both \Cat\ and the monad are cartesian.

\item		\label{eg:glob-mnd}
In Chapter~\ref{ch:defn} we will examine the free strict $\omega$-category
monad on the category of globular sets. Both category and monad are
cartesian.

\item
A \emph{double category} may be defined as a category object in
\fcat{Cat}. More descriptively, the graph structure consists of collections
of
\begin{itemize}
	\item 0-cells $A$
	\item horizontal 1-cells $f$
	\item vertical 1-cells $p$
	\item 2-cells $\alpha$
\end{itemize}
and various source and target functions, as illustrated by the picture
\[
\begin{diagram}
A_1		&\rTo^{f_1}		&A_2		\\
\dTo<{p_1}	&\Downarrow \alpha	&\dTo>{p_2}	\\
A_3		&\rTo_{f_2}		&A_4		\\
\end{diagram}
\ .
\]
The category structure consists of identities and composition functions for
2-cells and both kinds of 1-cell, obeying strict associativity, identity and
interchange laws; see \cite{KS} for more details.

More generally, let us define \emph{$n$-cubical set} for any $n\elt\nat$; the
intention is that a 2-cubical set will be the underlying graph of a double
category. Let \scat{H} be the category $(1 \parpair{\sigma}{\tau} 0)$, so
that a functor $\scat{H} \go \Set$ is a directed graph, and define an
$n$-cubical set as a functor $\scat{H}^n \go \Set$. Then, for instance, a
functor $X: \scat{H}^2 \go \Set$ becomes a two-dimensional graph of the type
just described, via
\begin{itemize}
\item $X(0,0) = \{0\mbox{-cells}\}$
\item $X(1,0) = \{\mbox{horizontal }1\mbox{-cells}\}$
\item $X(0,1) = \{\mbox{vertical }1\mbox{-cells}\}$
\item $X(1,1) = \{2\mbox{-cells}\}$,
\end{itemize}
and the map
\[
(\sigma,1): (1,1) \go (0,1)
\]
in $\scat{H}^2$ induces the map 
\[
\{2\mbox{-cells}\} \go \{\mbox{vertical }1\mbox{-cells}\}
\]
which sends $\alpha$ to $p_1$ in the diagram above.

We may now define a \emph{(strict) $n$-tuple category} to be an $n$-cubical
set together with various compositions and identities, as for double
categories, all obeying strict laws. The category of $n$-cubical sets has on
it the free strict $n$-tuple category monad; both category and monad are
cartesian. Since we will not need to use cubical sets or $n$-tuple
categories, this construction is not made precise and no proof is offered
that the monad is cartesian.

\end{eg}

\section{Multicategories}	\label{sec:multi}

We now describe what an \Cartpr-multicategory is, where $T$ is a cartesian
monad on a cartesian category \Eee. As mentioned in the introduction to this
chapter, this is a generalization of the well-known description of a small
category as a monad object in the bicategory of spans.

We will use the phrase `\Cartpr\ is cartesian' to mean that \Eee\ is a
cartesian category and \Mnd\ is a cartesian monad on \Eee.

\begin{constn} \label{constn:bicat} \end{constn}
Let \Cartpr\ be cartesian. We construct a bicategory $\Span\Cartpr$ from
\Cartpr, which in the case $T=\id$ is the usual bicategory of spans in
\Eee. Hermida calls $\Span\Cartpr$ the `Kleisli bicategory of spans' in
\cite{HerRM}; the formal similarity between the definition of $\Span\Cartpr$
and the usual construction of a Kleisli category is evident. 

\begin{description}
\item[0-cell:] Object $S$ of \Eee.
\item[1-cell $R \go S$:] Diagram
\begin{slopeydiag}
	&	&M	&	&	\\
	&\ldTo	&	&\rdTo	&	\\
TR	&	&	&	&S	\\
\end{slopeydiag}
in \Eee.
\item[2-cell $M \go M'$:] Commutative diagram
\begin{slopeydiag}
	&	&M	&	&	\\
	&\ldTo	&	&\rdTo	&	\\
TR	&	&\dTo	&	&S	\\
	&\luTo	&	&\ruTo	&	\\
	&	&M'	&	&	\\
\end{slopeydiag}
in \Eee.
\item[1-cell composition:] To define this we need to choose particular
pullbacks in \Eee, and in everything that follows we assume this has
been done. Take 
\[\left.
\begin{slopeydiag}
	&	&M	&	&	\\
	&\ldTo<{d}&	&\rdTo>{c}&	\\
TR	&	&	&	&S	\\
\end{slopeydiag}
\right.
\mr{\ and\ }
\left.
\begin{slopeydiag}
	&	&N	&	&	\\
	&\ldTo<{q}&	&\rdTo>{p}&	\\
TS	&	&	&	&Q	\\
\end{slopeydiag}
\right.
;
\]
then their composite is given by the diagram
\begin{slopeydiag}
   &       &   &       &   &       &N\of M\Spbk&  &   &       &   \\
   &       &   &       &   &\ldTo  &      &\rdTo  &   &       &   \\
   &       &   &       &TM &       &      &       &N  &       &   \\
   &       &   &\ldTo<{Td}&&\rdTo>{Tc}&   &\ldTo<{q}& &\rdTo>{p}& \\
   &       &T^2 R&     &   &       &TS    &       &   &       &Q  \\
   &\ldTo<{\mu_R}&&    &   &       &      &       &   &       &   \\
TR &       &   &       &   &       &      &       &   &       &   \\
\end{slopeydiag}
where the right-angle mark in the top square indicates that the square
is a pullback.

\item[1-cell identities:] The identity on $S$ is
\begin{slopeydiag}
	&	&S	&	&	&\\
	&\ldTo<{\eta_S}&&\rdTo>{1}&	&\\
TS	&	&	&	&S	&.\\
\end{slopeydiag}

\item[2-cell identities and compositions:] Identities and vertical
composition are as in \Eee. Horizontal composition is defined in an
obvious way.

\end{description}

Because the choice of pullbacks is arbitrary, 1-cell composition does
not obey strict associative and identity laws. That it obeys them up
to invertible 2-cells is a consequence of the fact that \Mnd\ is
cartesian. \done

\begin{defn}		\label{defn:multicat}
Let \Cartpr\ be cartesian. Then an 
\emph{\Cartpr-multicategory} is a monad in $\Span\Cartpr$.
\end{defn}

An \Cartpr-multicategory therefore consists of a diagram
\[
\spaan{C_1}{TC_{0}}{C_0}{d}{c} 
\]
in \Eee\ and maps 
\[
C_{0} \goby{\ids} C_{1},
\diagspace
C_{1}\of C_{1} \goby{\comp} C_{1}
\]
satisfying associative and identity laws. Think of $C_0$ as `objects', $C_1$
as `arrows', $d$ as `domain' and $c$ as `codomain'.  Such a multicategory
will be called an \Cartpr-multicategory \emph{on $C_0$}, and a
\Cartpr-multicategory on the terminal object $1$ will be called an
\emph{\Cartpr-operad}.

(Plain multicategories are often called `coloured operads' in the literature,
where the `colours' are the objects of the multicategory: thus an operad is a
single-coloured operad. A two-object plain multicategory would be called an
`operad of two colours', typically black and white.  Baez and Dolan, in
\cite{BaDoHDA3}, use `operad' or `typed operad' for the same kind of purpose
as we use `multicategory', and `untyped operad' where we use `operad'.)

It is inherent that everything is small: when $\Eee=\Set$, for instance, the
objects and arrows form sets, not classes. For plain multicategories, at
least, there seems to be no practical difficulty in using large versions too.

In order to say what maps between \Cartpr-multicategories are, we first
introduce the notion of an \Cartpr-graph.

\begin{defn}
Let \Cartpr\ be cartesian. An \emph{\Cartpr-graph} (on an object $C_0$) is a
diagram \gph{C_1}{C_0} in \Eee. A \emph{map of \Cartpr-graphs}
\[\left.
\begin{slopeydiag}
	&	&C_1	&	&	\\
	&\ldTo	&	&\rdTo	&	\\	
TC_{0}	&	&	&	&C_{0}	\\
\end{slopeydiag}
\right.
\go
\left.
\begin{slopeydiag}
		&	&\twid{C_1}&	&	\\	
		&\ldTo	&	&\rdTo	&	\\
T\twid{C_0}	&	&	&	&\twid{C_0}\\
\end{slopeydiag}
\right.\]
is a pair \pr{C_{0} \goby{f_0} \twid{C_0}}{C_{1}\goby{f_1} \twid{C_1}} of maps
in
\Eee\ such that
\begin{slopeydiag}
	&	&C_1	&	&	\\
	&\ldTo	&	&\rdTo	&	\\
TC_{0}	&	&\dTo>{f_1}&	&C_0	\\
\dTo<{Tf_{0}}&	&\twid{C_1}&	&\dTo>{f_0}\\
	&\ldTo	&	&\rdTo	&	\\
T\twid{C_0}&	&	&	&\twid{C_0}\\
\end{slopeydiag}
commutes.
\end{defn}

This definition uses two different notions of a map between objects of \Eee:
on the one hand, genuine maps in \Eee, and on the other, spans (i.e.\ 1-cells
of $\Span\Cartpr$). A possible approach to formalizing this situation is
via the `equipments' of~\cite{CKVW}. But this is not our approach: as
explained in~\ref{sec:fcm} and~\ref{sec:bim}, \fc-multicategories are the
structures that capture exactly what we want.

Any \Cartpr-multicategory has an underlying \Cartpr-graph, enabling the
following definition to be made.

\begin{defn}	\label{defn:multifunctor}
A \emph{map of \Cartpr-multicategories} $C \go \twid{C}$ is a map $f$ of
their underlying graphs such that the diagrams
\[\left.
\begin{diagram}
C_0		&\rTo^{f_0}		&\twid{C_0}		\\	
\dTo<{\ids}	&			&\dTo>{\twid{\ids}}	\\
C_1		&\rTo^{f_1}		&\twid{C_1}		\\
\end{diagram}
\right.
\ \ \ 
\left.
\begin{diagram}
&C_{1}\of C_1	&\rTo^{f_{1}*f_1}	&\twid{C_1}\of\twid{C_1}	\\
&\dTo<{\comp}	&			&\dTo>{\twid{\comp}}	\\
&C_1		&\rTo^{f_1}		&\twid{C_1}		\\	
\end{diagram}
\right.\]
commute. (Here $f_1 * f_1$ is the evident map induced by two copies
of $C_1 \goby{f_1} C_1$.)
\end{defn}

With these definitions we obtain categories
\[
\Cartpr\hyph\Graph, \diagspace \Cartpr\hyph\Multicat,
\]
and a forgetful functor from the second to the first. Wherever possible we
drop the `\Eee' and refer simply to $T$-multicategories, $T$-operads,
$T\hyph\Graph$, etc. 

It is also possible to define modules (profunctors) and natural
transformations for $T$-multicategories, which we eventually do
in~\ref{egs:bim}\bref{eg:Bim-Cartpr}.

\begin{remarks}{rmks:maps}
\item 	\label{rmk:map-fixed-obj}
Fix $S\elt\Eee$. Then we may consider the category of
$T$-graphs on $S$, whose morphisms $f = \pr{S \goby{f_0} S}
{C_{1} \goby{f_1} \twid{C_1}}$ all have $f_{0}=1$. This is just the slice
category $\frac{\Eee}{TS\times S}$. It is also the full sub-bicategory of
$\Span\Cartpr$ whose only object is $S$, and is therefore a monoidal
category. The category of $T$-multicategories on $S$ is then the
category $\Mon(\frac{\Eee}{TS\times S})$ of monoids in 
$\frac{\Eee}{TS\times S}$. In particular, $\Eee/T1$ is a monoidal category,
and a monoid therein is a $T$-operad; this is a style of definition of plain
operad sometimes found in the literature. 

\item A choice of pullbacks in \Eee\ was made; changing that choice gives an
isomorphic category of \Cartpr-multicategories.

\item
If \pr{\Eee'}{T'} is also cartesian then a cartesian monad functor from
\Cartpr\ to \pr{\Eee'}{T'} induces a functor 
\[
\Cartpr\hyph\Multicat \go
\pr{\Eee'}{T'}\hyph\Multicat, 
\]
and the same is true of monad opfunctors. See~\ref{sec:change} for an
explanation.
\end{remarks}

\begin{eg}{egs:multicats}

\item
\begin{sloppypar}
Let $\Cartpr=\pr{\Set}{\id}$. Then $\Span\Cartpr$ is the usual `bicategory of
spans', and a monad in $\Span\Cartpr$ is just a (small) category. Thus
categories are \pr{\Set}{\id}-multicategories.  Functors are maps of
such. More generally, if \Eee\ is any cartesian category then
\pr{\Eee}{\id}-multicategories are internal categories in \Eee, and
similarly, \id-operads are monoids. 

\item
Let $\Cartpr=\pr{\Set}{\mr{free\ monoid}}$. Specifying a $T$-graph
\[
\graph{C_1}{C_0}{d}{c} 
\]
is equivalent to specifying a set $C_0$ (`of objects') together with a set
\multihom{C}{\range{a_1}{a_n}}{a} for each $n\geq 0$ and
\range{a_1}{a_n,a}\elt$C_0$. An element
$\theta\elt\multihom{C}{\range{a_1}{a_n}}{a}$ is illustrated by
\[
\range{a_1}{a_n} \goby{\theta} a
\]
or
\[
\ctransistor{\theta}{a_1}{a_2}{a_n}{a}
\]
or
\[
\begin{opetope}
	&	&	&\cnr	&\ldots	&	&	&	\\
	&\cnr	&\ruLine(2,1)<{a_2}&&	&	&\cnr	&	\\
\ruLine(1,2)<{a_1}&&	&	&\Downarrow \theta&&	&\rdLine(1,2)>{a_n}\\
\cnr	&	&	&\rLine_{a}&	&	&	&\cnr	\\
\end{opetope}
\ \ \ \ .
\]
When $n=0$, the first version looks like
\[
\cdot \goby{\theta} a,
\]
the second has no legs on the left-hand (`input') side, and the third is
drawn as
\[
\topezs{a}{\Downarrow\!\theta}\ .
\]
\end{sloppypar}

In $\Span\Cartpr$, the identity 1-cell \graph{C_0}{C_0}{\eta_{C_0}}{1} on
$C_0$ has
\[
\multihom{C_0}{\range{a_1}{a_n}}{a}=
\left\{
\begin{array}{ll}
1		&\mr{if\ }n=1\mr{\ and\ }a_{1}=a	\\
\emptyset	&\mr{otherwise.}
\end{array}
\right.
\]
The composite 1-cell $C_{1}\of C_1$ in $\Span\Cartpr$ is
\[
\{\pr{\bftuple{\theta_1}{\theta_n}}{\theta} \such 
d\theta=\bftuple{c\theta_1}{c\theta_n}\},
\]
i.e.\ is the set of diagrams
\begin{equation} 	\label{diag:multi-comp}
\comptrans{\theta_1}{\theta_2}{\theta_n}{\theta}.
\end{equation}

If $C$ is a $T$-multicategory then we have a function \ids\ assigning to each
$a\elt C_0$ a member of \multihom{C}{a}{a}, and a function \comp\ composing
diagrams like~(\ref{diag:multi-comp}). These are required to obey associative
and identity laws. Thus a \pr{\Set}{\mr{free\ monoid}}-multicategory is just
a plain multicategory and a \pr{\Set}{\mr{free\ monoid}}-operad is a plain
operad.

\item 	\label{eg:sym-ops}
Suppose we want to realise symmetric operads as $T$-operads for some $T$. By
a \emph{symmetric operad} I mean a plain operad $C$ with an action of the
$n$th symmetric group $S_n$ on $C(n)$ for each $n$, satisfying certain
axioms: in other words, an operad in the usual sense of topologists
(e.g.~\cite{MayDOA}), except that the $C(n)$'s are sets rather than spaces or
graded modules etc.

A first attempt might be to take the free commutative monoid monad $T$
on $\Set$. But this is both misguided and doomed to failure: misguided
because the maps 
\[
\dashbk\cdot\sigma: C(n) \go C(n)
\]
coming from permutations $\sigma\in S_n$ are only isomorphisms, not
identities; and doomed because $T$ is not
cartesian~(\ref{eg:cart-mnds}\bref{eg:comm-not-cart}). 

A more promising approach is to take $T$ to be the free symmetric strict
monoidal category monad on \Cat, and to try to identify the symmetric operads
as certain special $T$-operads. I have not investigated how well this works,
but this idea seems to be related to the structures called
`symmetric operads' at the beginning of~\cite{BaDoHDA3} and explored further
in~\cite{CheROM} and~\cite{CheEAT}.
 
\item
Let $\Eee=\Set$, and consider the monad $\dashbk+1$
of~\ref{eg:cart-mnds}(\ref{eg:mon-exceptions}). A
\pr{\Set}{\dashbk+1}-graph is a diagram \spaan{C_1}{C_0 +1}{C_0}{d}{c} of
sets; this is like an ordinary \pr{\Set}{\id}-graph, except
that some arrows have domain 0---an extra element not in $C_0$. (Thus
$1=\{0\}$ here.) If we put
\[
Y(a) = \{y\elt C_1 \such dy=0 \mbox{ and } cy=a\}
\]
for each $a\elt C_0$, then a multicategory structure on the graph provides a
function
\[
\begin{array}{rcl}
Y(a)	&\go	&Y(a')	\\
y	&\goesto&\theta\of y
\end{array}
\cdiscfibpic
\]
for each $\theta\elt C_1$ with $d(\theta)=a\elt C_0$ and $c(\theta)=a'$. It
also provides a category structure on \spaan{D_1}{D_0}{D_0}{d}{c}, where $D_0
= C_0$ and $D_1 = \{\theta\elt C_1\such d\theta\elt C_0\}$. Thus a
\pr{\Set}{\dashbk+1}-multicategory turns out to be a (small) category
$D$ together with a functor $Y: D\go\Set$. Similarly, a
\pr{\Set}{\dashbk+E}-multicategory is a category $D$ together with an
$E$-indexed family of functors $D\go\Set$.

To put it another way, an \Cartpr-multicategory is a discrete
opfibration. More exactly, the category of \Cartpr-multicategories is
equivalent to the category whose objects are discrete opfibrations between
small categories and whose morphisms are commutative squares.

\item		\label{eg:M-times-mti}
Let $M$ be a monoid and $\Cartpr = \pr{\Set}{M\times\dashbk}$. Then a
$T$-mul\-ti\-cat\-e\-gory consists of a category $C$ together with a functor
$C\go M$, and in fact $T\hyph\Multicat \iso \Cat/M$.

\item		\label{eg:tree-mti}
Let $\Cartpr = \pr{\Set}{\mr{tree\ monad}}$, as
in~\ref{eg:cart-mnds}(\ref{eg:mon-tree}). A $T$-multicategory consists of a
set $C_0$ of objects, and hom-sets like
\[
C\left(\Pear{a_1}{a_2}{a}\right)
\]
($a_{1}, a_{2}, a \elt C_0$), together with a unit element of each
$C(\nlal{a}{a})$ and composition functions like
\begin{eqnarray*}
C\left( \Pear{a_1}{a_2}{a}\right)
\times
\left\{	C\left( \Orange{b_1}{b_2}{a_1}\right)
	\times
	C\left( \Apple{b_3}{b_4}{a_2}\right)
\right\}					
\\
\go
C\left( \MixedFruit{b_1}{b_2}{b_3}{b_4}{a}\right)
\end{eqnarray*}
($b_{1}, b_{2}, b_{3}, b_{4} \elt C_0$). These are to satisfy associativity
and identity laws.

When $C_{0}=1$, so that we're considering $T$-operads, the graph
structure is comprised of sets like 
\[
C\left(\Pear{}{}{}\right).
\]

\begin{sloppypar}
The $T$-multicategories are a simpler version of Soibelman's pseudo-monoidal
categories (\cite{Soi}) or Borcherds' relaxed multilinear categories
(\cite{Borh}, \cite{SnyEBG}, \cite{SnyRMS}); they omit the aspect of maps
between trees. See the end of~\ref{sec:enrich} for comments on the
unsimplified version.
\end{sloppypar}

\item		%\label{eg:soib-op}
When $\Eee=\fcat{Cat}$ and $T$ is the free strict monoidal category monad, a
$T$-operad is what Soibelman calls a strict monoidal
2-operad~(\cite[2.1]{Soi}). Such a structure might also be thought of as a
plain operad enriched in \Cat, in a sense not made precise here but explained
in detail in \cite{GECM} and outlined in~\ref{sec:enrich} below.

\item		%\label{eg:glob-ops}
Let $\Cartpr=\pr{\mr{globular\ sets}}{\mr{free\ strict\
}\omega\hyph\mr{category}}$, as in
Example \ref{eg:cart-mnds}\bref{eg:glob-mnd}. A $T$-operad is exactly a
globular operad in the sense of Batanin: see Chapter~\ref{ch:defn}.

\item
Operads for $\Cartpr = \pr{n\hyph\mr{cubical\ sets}}{\mr{free\ strict\
}n\hyph\mr{tuple\ category}}$ can be understood in much the same way as
Batanin's globular operads (again, see Chapter~\ref{ch:defn}), with cubical
rather than globular shapes. For instance, a cell in the free strict
$n$-tuple category on the terminal $n$-cubical set can be represented as a
cuboid whose edge-lengths are natural numbers; a $T$-operad associates a set
(`of operations') to each such cuboid, and has composition functions
according to ways of combining cuboids. (I will not take this example any
further.) 

\item	\label{eg:multi-alg}
Let \Cartpr\ be cartesian, let $X\in\Eee$, and let $TX \goby{h} X$ be a
map. Then the $T$-graph $(\graph{TX}{X}{1}{h})$ can be given the structure of
a $T$-multicategory in at most one way, and this is possible if and only if
$TX \goby{h} X$ is an algebra for the monad $T$. (If it \emph{is} possible
then $\ids=\eta$ and $\comp=\mu$.) Maps between $T$-multicategories of this
form are, similarly, just $T$-algebra maps. So we have a full and faithful
functor
\[
M: \Eee^T \go T\hyph\Multicat
\]
turning algebras into multicategories.

\end{eg}

\section{Algebras}	\label{sec:alg}

The motivating idea in the definition of a (plain) operad is that it is some
kind of algebraic theory, with the $n$th component $C(n)$ of an operad $C$
being the set of $n$-ary operations. One therefore defines an \emph{algebra}
for an operad $C$ to be a set $X$ together with a suitable family of
functions 
\[
C(n) \times X^n \go X,
\]
one for each $n\in\nat$. More generally, a plain multicategory can be
regarded as a many-sorted theory, and in an algebra for a multicategory one
has not just a single set $X$, but one set $X(a)$ for each object $a$ of
$C$. Thus if \Set\ denotes the (large, plain) multicategory whose
objects are sets and in which a map
\[
S_1, \ldots, S_n \go S
\]
is a function 
\[
S_1 \times\cdots\times  S_n \go S,
\]
\label{p:multifunctor}%
then an algebra for a plain multicategory $C$ can be defined as a map $C\go
\Set$ of multicategories. 

In this section we generalize these ideas to arbitrary \Cartpr. That is, we
define a category $\Alg(C)$ of algebras for any \Cartpr-multicategory $C$.

\begin{constn} \label{constn:ind-mnd} \end{constn}
Let \Cartpr\ be cartesian: then any \Cartpr-multicategory $C$ gives rise to a
monad \triple{T_C}{\unit}{\mult} on $\Eee/C_0$. In what follows, I will
write
$T_C (X \goby{p} C_0)$ as $X' \goby{p'} C_0$. 
\begin{itemize}
\item Given $\slob{X}{p}{C_0} \in \Eee/C_0$, we define
\slob{X'}{p'}{C_0} to be the right-hand diagonal of the diagram
% %
% \begin{equation}	\label{eq:alg-defn}
\begin{slopeydiag}
	&	&\cdot\Spbk&	&	&	&	\\
	&\ldTo	&	&\rdTo	&	&	&	\\
TX	&	&	&	&C_1	&	&	\\
	&\rdTo<{Tp}&	&\ldTo>{d}&	&\rdTo>{c}&	\\
	&	&TC_0	&	&	&	&C_0	\\
\end{slopeydiag}
% \end{equation}
% %
(recalling that the right-angle mark denotes a pullback square).

\item If
$\left.
\begin{triangdiag}
X	&	&\rTo^{f}&	&Y	\\
	&\rdTo<{p}&	&\ldTo>{q}&	\\
	&	&C_0	&	&	\\
\end{triangdiag}
\right.$
is a map in $\Eee/C_0$ then there is a unique map $f': X' \go Y'$ making 
\begin{slopeydiag}
   &   &X'\Spbk&   &   &
   &   &   &   &   & 
   &   &   &   &   \\
   &\ldTo&   &\rdTo\rdGet(10,2)^{f'}&   &
   &   &   &   &   &
   &   &   &   &   \\
TX &   &   &   &C_1&
   &   &   &   &   &
   &   &Y'\Spbk&   &   \\
   &\rdTo<{Tp}\rdTo(10,2)^{Tf}&   &\ldTo>{d}&   &
\rdTo(10,2)^{1}&   &   &   &   &
   &\ldTo&   &\rdTo&   \\
   &   &TC_0&   &   &
   &   &   &   &   &
TY &   &   &   &C_1\\
   &   &   &\rdTo(10,2)^{1}&     &
   &   &   &   &   &
   &\rdTo<{Tq}&   &\ldTo>{d}&   \\
   &   &   &   &   &
   &   &   &   &   &
   &   &TC_0&  &   \\
\end{slopeydiag}
commute, and we define 
\[
T_C(f) = f': \bktdvslob{X'}{p'}{C_0} \go \bktdvslob{Y'}{q'}{C_0}.
\]

\item The unit at \slob{X}{p}{C_0} is given by
\[
\begin{slopeydiag}
	&			&X		&	&	\\
	&\ldTo(2,5)<{\eta_X}	&		&\rdTo>{p}&	\\
	&			&\dGet~{\unit_p}	&	&C_0	\\
	&			&X'\Spbk	&	&\dTo>{\ids}\\
	&\ldTo			&		&\rdTo	&	\\
TX	&			&		&	&C_1	\\
	&\rdTo<{Tp}		&		&\ldTo>{d}&	\\
	&			&TC_0.		&	&	\\
\end{slopeydiag}
\]

\item For multiplication, we have a commutative diagram
\begin{slopeydiag}
   &	&	&	&	&	&X''	&	&    &	&  \\
   &	&	&	&	&\ldTo	&	&\rdTo(4,4)& &	&  \\
   &	&	&	&TX'	&	&	&	&    &	&  \\
   &	&	&\ldTo	&	&\rdTo	&	&	&    &	&  \\
   &	&T^2 X	&	&	&	&TC_1	&	&    &	&C_1\\
 &\ldTo<{\mu_X}&&\rdTo<{T^2 p}& &\ldTo>{Td}&	&\rdTo<{Tc}&&\ldTo>{d}&\\
TX &	&	&	&T^2 C_0&	&	&	&TC_0 &  &  \\
\end{slopeydiag}
and a pullback square
\[
\begin{slopeydiag}
	&	&C_1\of C_1\Spbk&&	\\
	&\ldTo	&	&\rdTo	&	\\
TC_1	&	&	&	&C_1	\\
	&\rdTo<{Tc}&	&\ldTo>{d}&	\\
	&	&TC_0.	&	&	\\
\end{slopeydiag}
\]
From these we deduce that there are maps
\[\begin{slopeydiag}
	&	&X''	&	&		\\
	&\ldTo	&	&\rdTo	&		\\
TX	&	&	&	&C_1 \of C_1,	\\
\end{slopeydiag}
\]
and the multiplication at \slob{X}{p}{S} is then given by
\[	%\left.
\begin{slopeydiag}
	&			&X''		&	&	\\
	&\ldTo(2,5)		&		&\rdTo	&	\\
	&			&\dGet~{\mult_p}&	&C_1 \of C_1\\
	&			&X'\Spbk	&	&\dTo>{\comp}\\
	&\ldTo			&		&\rdTo	&	\\
TX	&			&		&	&C_1	\\
	&\rdTo<{Tp}		&		&\ldTo>{d}&	\\
	&			&TC_0.		&	&	\\
\end{slopeydiag}
\]
\end{itemize}

It is now straightforward, though tedious, to check that
\triple{T_C}{\unit}{\mult} forms a monad on $\Eee/C_0$.
\done

\begin{defn}
Let \Cartpr\ be cartesian and let $C$ be a $T$-multicategory. Then the
category $\Alg(C)$ of \emph{algebras} for $C$ is the category of algebras for
the monad $T_C$ on $\Eee/C_0$. 
\end{defn}
We sometimes say \emph{$C$-algebra} instead of `algebra for $C$'.

\begin{eg}{eg:algs}

\item 
When $\Cartpr=\pr{\Set}{\id}$, so that an \Cartpr-multicategory is an
ordinary (small) category $C$, we have $\Alg(C) \eqv \ftrcat{C}{\Set}$.

\item	\label{eg:alg-plain}
When $\Cartpr=\pr{\Set}{\mr{free\ monoid}}$, so that an \Cartpr-multicategory
is a plain multicategory, we already have an idea of what an algebra for $C$
should be: a map $C \go \Set$ of multicategories
(p.~\pageref{p:multifunctor}). That is, an algebra for $C$ should consist of:
\begin{itemize}
\item for each $a\elt C_0$, a set $X(a)$
\item for each $\range{a_1}{a_n}\goby{\theta}a$ in $C$, a function
$X(a_{1})\times\cdots\times X(a_{n})$\linebreak$ \go X(a)$,   
\end{itemize}
preserving identities and composition. This is the same as the
definition of algebra just given. To see this, let $(X \goby{p} C_0)$ be an
object of $\Eee/C_0$: then, writing $X(a) = p^{-1} \{ a \}$ for $a\in C_0$,
and similarly $X'(a) = (p')^{-1} \{ a \}$, we have 
\begin{eqnarray*}
\lefteqn{X'(a) =}							\\
	& &	\{\pr{\bftuple{x_1}{x_n}}{\theta} \such
			x_i \in X, \theta \in C_1,
			d\theta=\bftuple{px_1}{px_n}, c\theta = a\}	\\
	&=&	\{X(a_{1})\times\cdots\times X(a_{n}) \times
			\multihom{C}{\range{a_1}{a_n}}{a} \such
			\range{a_1}{a_{n}}\elt C_0\}.
\end{eqnarray*}
An algebra structure on $(X \goby{p} C_0)$ therefore consists of a
function
\[
X(a_{1})\times\cdots\times X(a_{n}) \goby{\bar{\theta}} X(a)
\]
for each 
\[
a_1, \ldots, a_n \goby{\theta} a
\]
in $C$, with the assignation $\theta \goesto \bar{\theta}$ subject to certain
rules. These turn out to say exactly that we have a multicategory map $C \go
\Set$.

\item
When $\Cartpr=\pr{\Set}{\dashbk+1}$, a $T$-multicategory is an ordinary
category $D$ together with a functor $D\goby{Y}\Set$. A
\pr{D}{Y}-algebra is then a functor $D \goby{X} \Set$ together with a
natural transformation
\[
\begin{diagram}
D	&\pile{\rTo^Y \\ \Downarrow \\ \rTo_X}	&\Set	\\
\end{diagram}
\mr{\ .}
\] 
In terms of fibrations, a $T$-multicategory is a discrete opfibration
$Y$ over a small category $D$, and an algebra for
$Y$ consists of another discrete opfibration $X$ over $D$ together with a
map from $Y$ to $X$ (of opfibrations over $D$).

\item 		\label{eg:M-times-alg}
Let $M$ be a monoid and let $\Cartpr = \pr{\Set}{M\times\dashbk}$, so that a
$T$-mul\-ti\-cat\-e\-gory is a category $C$ together with a functor $C
\goby{\pi} M$. Then the category of algebras for \pr{C}{\pi} is simply
\ftrcat{C}{\Set}, regardless of what $\pi$ is.

\item
Let \Cartpr\ be the tree monad on \Set. For simplicity, let us just consider
algebras for $T$-\emph{operads} $C$---thus the object-set $C_0$ is $1$. An
algebra for $C$ consists of a set $X$ together with a function $X'\go X$
satisfying some axioms. One can calculate that an element of $X'$ consists of
an $X$-labelling of the leaves of a tree $\tau$ together with a member of
$C(\tau)$.  An $X$-labelling of an $n$-leafed tree $\tau$ is just a member of
$X^n$, so one can view the algebra structure $X'\go X$ on $X$ as: for each
number $n$, each $n$-leafed tree $\tau$, and each element of $C(\tau)$, a
function $X^{n}\go{X}$.  These functions are required to be compatible with
composition and identities in $C$.

\item
For $\Cartpr=\pr{\mr{globular\ sets}}{\mr{free\ strict\ }
\omega\hyph\mr{category}}$, we will consider in Chapter~\ref{ch:defn} a
certain operad $L$, the initial `operad-with-contraction'. A weak
$\omega$-category is then defined to be an $L$-algebra.

\item		\label{eg:alg-terminal}
The graph \graph{T1}{1}{1}{!} is terminal amongst all \Cartpr-graphs. It
carries a unique multicategory structure, since a terminal object in a
monoidal category always carries a unique monoid structure. It then becomes
the terminal \Cartpr-multicategory. The induced monad on $\Eee/1$ is just
\Mnd, and so an algebra for the terminal multicategory is just a $T$-algebra.

This can aid recognition of when a theory of operads or multicategories fits
into our scheme. For instance, if we were to read Batanin's paper and learn
that, in his terminology, an algebra for the terminal operad is a strict
$\omega$-category (\cite[\S 7, example 3]{Bat}), then we might suspect that
his operads were \Cartpr-operads for the free strict $\omega$-category monad
$T$ on an appropriate category \Eee---as indeed they are.

\item		\label{eg:alg-to-multi}
If $T$ is a monad on a category \Eee, and $h= (TX \goby{h} X)$ is a
$T$-algebra, then there is a monad $T/h$ on $\Eee/X$ whose functor part acts
on objects by
\[
\vslob{Y}{p}{X}
\diagspace \goesto \diagspace
\begin{diagram}[height=1.5em]
TY	\\ \dTo>{Tp}	\\ TX	\\ \dTo>{h}	\\ X.
\end{diagram}
\]
Writing \Alg\ for the category of algebras of a monad, we then have
\[
\Alg(T/h) \iso \Alg(T)/h,
\]
where the right-hand side is $\Alg(T)$ sliced over $h$. 

Now recall from Example~\ref{egs:multicats}\bref{eg:multi-alg} that when
\Cartpr\ is cartesian, the algebra $h$ defines a $T$-multicategory
\[
\graph{TX}{X}{1}{h}.
\]
Naturally enough, it turns out that the monad on $\Eee/X$ induced by this
multicategory is $T/h$: so the category of algebras for this multicategory is
$\Alg(T)/h$. (Example~\bref{eg:alg-terminal} above is a special case.)

\end{eg}

We have seen how to associate to each \Cartpr-multicategory $C$ a category
$\Alg(C)$, and we would expect some kind of functoriality. When
$\Cartpr = \pr{\Set}{\id}$, a functor $C \go C'$ induces a functor
\[
\Alg(C) = \ftrcat{C}{\Set} \og \ftrcat{C'}{\Set} = \Alg(C'),
\]
and it is obvious that the same phenomenon holds for $\Cartpr =
(\Set,\mr{free\ monoid})$ if we view $C$-algebras as
multicategory maps $C\go \Set$~(\ref{eg:algs}\bref{eg:alg-plain}).

In general, given a map $f: C \go C'$ of \Cartpr-multicategories, we obtain a
functor $\Alg(C) \og \Alg(C')$ as follows. First of all, we have the functor
\[
\Eee/C_0 \ogby{f_0^*} \Eee/C'_0
\]
defined by pullback along $f_0: C_0 \go C'_0$. Then, as it turns out, there
is a naturally-arising natural transformation
\[
\begin{ntdiag}
\Eee/C_0	&	&\lTo^{f_0^*}	&	&\Eee/C'_0	\\
		&	&		&	&		\\
\dTo<{T_C}	&	&\rdTo>{\phi}	&	&\dTo>{T_{C'}}	\\
		&	&		&	&		\\
\Eee/C_0	&	&\lTo_{f_0^*}	&	&\Eee/C'_0,	\\
\end{ntdiag}
\]
and this satisfies the axioms for a monad functor from $T_{C'}$ to
$T_C$. (Monad functors are defined in~\ref{sec:change}, and details of this
construction are left to the reader.) Hence there is an induced functor from
the category of $T_{C'}$-algebras to the category of $T_C$-algebras---that
is, from $\Alg(C')$ to $\Alg(C)$. 

Because these induced functors are defined by pullback, the map
\[
\Alg: (\Cartpr\hyph\Multicat)^{\op} \go \fcat{CAT}
\]
inevitably preserves composition and identities only up to canonical
isomorphism; in other words, it is a weak functor or pseudo-functor. In fact,
there is a notion of natural transformation for $T$-multicategories, so that
$\Cartpr\hyph\Multicat$ is a 2-category; and $\Alg$ is then a weak functor
between 2-categories.  Transformations for $T$-multicategories are discussed
in Example~\ref{egs:bim}\bref{eg:Bim-Cartpr}, where we see that the natural
structure formed by $T$-multicategories is not a 2-category but something
richer: an \fc-multicategory.

\chapter{More on Operads and Multicategories}
\label{ch:mtitwo}

This chapter is an assortment of further topics in the general theory of
multicategories. Some will be used in the discussion of weak $n$-categories
in the final chapter. Others are not used there, but answer naturally-arising
questions or have applications outside this thesis. One of the
sections~(\ref{sec:enrich}, Enrichment) is an introduction to a topic too
large to include in full. 

\begin{sloppypar}
The contents of the sections are as follows.
\begin{description}
\item[\ref{sec:struc} Structured Categories] We look at $T$-structured
categories, which are to $T$-multicategories as strict monoidal categories
are to plain multicategories.
\item[\ref{sec:change} Change of Base] Here we ask whether the
passage from $\Cartpr$ to $\Cartpr\hyph\Multicat$ is functorial. It turns
out that it is, in not one but three different ways. (The `base' is \Cartpr.) 
\item[\ref{sec:free-multicats} Free Multicategories] This section concerns
when and how the free $T$-multicategory on a $T$-graph can be
formed. Details are deferred to Appendix~\ref{app:free}.
\end{description}
\end{sloppypar}

The next two sections each give an alternative (but equivalent) definition of
algebra for a $T$-multicategory.
\begin{description}
\item[\ref{sec:fib} Algebras via Fibrations] In ordinary category theory
there is a correspondence between \Set-valued functors and discrete
fibrations. We extend this to $T$-multicategories, giving an alternative
definition of an algebra.
\item[\ref{sec:endo} Algebras via Endomorphisms] We give a second alternative
definition of an algebra, generalizing the definition of algebra for an operad
often used by topologists.
\end{description}

The final sections are on \fc-multicategories: what they are, and two
familiar categorical ideas for which they provide generalized contexts.
\begin{description}
\item[\ref{sec:fcm} \fc-Multicategories] This really belongs as an example in
the previous chapter, and would be there but for its length. We examine
$T$-mul\-ti\-cat\-e\-go\-ries in the case when $T$ is the free category monad
on the category of directed graphs.
\item[\ref{sec:bim} The Bimodules Construction] We show how an
\fc-multicategory $V$ gives rise to a new \fc-multicategory \Bim{V}, by
taking bimodules ($=$ modules, $=$ profunctors, $=$ distributors) in $V$.
\item[\ref{sec:enrich} Enrichment] There is an interesting theory of
enrichment for $T$-mul\-ti\-cat\-e\-gories. Applied to the most basic case,
categories, it provides a theory of categories enriched in an
\fc-multicategory. All of this is explained properly in~\cite{GECM}; here we
sketch the ideas.
\end{description}

None of these sections is necessary in order to read the bulk of
Chapter~\ref{ch:defn}. The last part of Chapter~\ref{ch:defn}, on weak
$n$-categories, does rely on the material of~\ref{sec:change} (Change of
Base). It also contains inessential references to~\ref{sec:endo}
and~\ref{sec:fcm}.  Appendix~\ref{app:initial}, which supports
Chapter~\ref{ch:defn}, uses free multicategories~(\ref{sec:free-multicats}).

\section{Structured Categories}		\label{sec:struc}

The observation from which this section takes off is that any strict monoidal
category has an underlying multicategory. (For the time being, all monoidal
categories and maps between them are strict, and `multicategory' means plain
multicategory.)  Explicitly, if \pr{D}{\otimes} is a monoidal category, then
the underlying multicategory $C$ has the same objects as $D$ and has hom-sets
defined by
\[
\multihom{C}{\range{a_1}{a_n}}{a} = 
\homset{D}{a_{1}\otimes\cdots\otimes a_{n}}{a}
\]
for objects \range{a_1}{a_n,a}. Composition and identities in $C$ are easily
defined.

There is a converse process: given any multicategory $C$, there is a `free'
monoidal category $D$ on it.  An object (respectively, arrow) of $D$ is a
sequence of objects (respectively, arrows) of $C$.  Thus the objects of $D$
are of the form \bftuple{a_1}{a_n} ($a_i \elt C_0$), and a typical arrow
\[
\tuplebts{a_1,a_2,a_3,a_4,a_5} \go \tuplebts{a'_1,a'_2,a'_3}
\]
is a sequence
\tuplebts{\theta_1,\theta_2,\theta_3} of maps in $C$ with domains and
codomains as illustrated:
\begin{equation}	\label{diag:arrows-in-mon-cat} 
\freemoncatpic
\end{equation}
The tensor in $D$ is juxtaposition.

For example, the terminal multicategory \fcat{1} has one object and, for each
$n\elt\nat$, one arrow of the form
\[
n \left\{\rule{0ex}{5.3ex}\right.\left.
\ctransistor{}{}{}{}{}
\right.
\ ;
\]
diagram~\bref{diag:arrows-in-mon-cat} suggests that the `free' monoidal
category on the multicategory \fcat{1} is $\Delta$, the category of finite
ordinals (including 0), with addition as~$\otimes$.

The name `free' is justified: that is, there is an adjunction
\begin{diagram}[height=2em]
\mbox{(monoidal categories)}	\\
\uTo\dashv\dTo			\\
\mbox{(multicategories)}	\\
\end{diagram}
where the two functors are those described above, and (monoidal categories)
denotes the category of strict monoidal categories and strict monoidal
functors. Moreover, this adjunction is monadic. 

(Note that the forgetful functor does not provide a \emph{full} embedding of
(monoidal categories) into (multicategories). For example,
there is a multicategory map $\fcat{1} \go \Delta$ sending the unique object
of \fcat{1} to the object 1 of $\Delta$, and this map does not preserve the
monoidal structure. If $D$ and $D'$ are strict monoidal categories then a map
$UD \go UD'$ of their underlying multicategories is actually the same as a
\emph{lax} monoidal functor $D\go D'$.)

Naturally, we would like to generalize from $\Cartpr=\pr{\Set}{\mr{free\
monoid}}$ to any cartesian \Cartpr. To do this, we need a notion of
`\Cartpr-structured category' which in the case \pr{\Set}{\mr{free\
monoid}} means monoidal category. A monoidal category is a category
object in \fcat{Monoid}, so it is reasonable to define a
\emph{$T$-structured category} to be an
\pr{\Eee^T}{\id}-multicategory---that is, an internal category in the
category $\Eee^T$ of algebras for the monad $T$ on \Eee. We write
\[
T\hyph\Struc = \pr{\Eee^T}{\id}\hyph\Multicat.
\]
The fact that $\Eee$ is cartesian guarantees that $\Eee^T$ is too.

(In this section, $\pr{\cat{D}}{S}\hyph\Multicat$ is treated as a mere
(1-)category, for any cartesian \cat{D} and $S$.)

It is now possible to describe a monadic adjunction
\begin{diagram}[height=2em]
T\hyph\Struc			\\
\uTo<{F}\dashv\dTo>{U}		\\
T\hyph\Multicat			\\
\end{diagram}
generalizing that above. The effect of the functors $U$ and $F$ on objects is
as outlined now. Given a $T$-structured category $D$, with algebraic
structure $TD_{0} \goby{h_0} D_{0}$ and $TD_{1} \goby{h_1} D_{1}$, the graph
$(\gph{C_1}{D_0})$ of $UD$ is given by
\[
\begin{slopeydiag}
	&		&C_{1}\Spbk&	&	&	&	\\
	&\ldTo		&	&\rdTo	&	&	&	\\
TD_{0}	&		&	&	&D_{1}	&	&	\\
	&\rdTo<{h_0}	&	&\ldTo	&	&\rdTo	&	\\
	&		&D_0	&	&	&	&D_0.	\\
\end{slopeydiag}
\]
Given a $T$-multicategory $C$, the category $FC$ has graph
\[
\begin{slopeydiag}
	&	&	&	&TC_{1}	&	&	&	&	\\
	&	&	&\ldTo	&	&\rdTo(4,4)&	&	&	\\
	&	&T^2 C_0&	&	&	&	&	&	\\
	&\ldTo<{\mu_{C_0}}&&	&	&	&	&	&	\\
TC_{0}	&	&	&	&	&	&	&	&TC_{0}	\\
\end{slopeydiag}
\]
and the algebraic structures $T^2 C_0 \goby{h_0} T C_{0}$ and $T^2 C_{1}
\goby{h_1} TC_{1}$ are components of $\mu$.

For an example of $U$ in action, take a $T$-algebra $(TX \goby{h} X)$. The
diagram
\[
X \ogby{1} X \goby{1} X
\]
determines a (discrete) internal category in $\Eee^T$; that is, a
$T$-structured category, 
\label{p:discrete-struc}%
$D(X,h)$. Then $U(D(X,h))$ is a $T$-multicategory with graph isomorphic to 
\[
TX \ogby{1} TX \goby{h} X.
\]
So $U(D(X,h))$ is isomorphic to the $T$-multicategory $M(X,h)$ of
Example~\ref{egs:multicats}\bref{eg:multi-alg}, and we have a triangle of
functors 
\[
\begin{slopeydiag}
\Eee^T	&		&\rTo^{D}	&		&T\hyph\Struc	\\
	&\rdTo<{M} 	&		&\ldTo>{U}	&		\\
	&		&T\hyph\Multicat&		&		\\
\end{slopeydiag}
\]
which commutes up to natural isomorphism.

For an example of $F$, take \Eee\ to be \Set\ and $T=\blank^*$ to be the free
monoid monad. Take the terminal plain multicategory $1$, which has graph
\[
\nat \ogby{1} \nat \goby{!} 1.
\]
Then $F(1)$ is a strict monoidal category with graph
\[
\nat \ogby{+} \nat^* \goby{!^*} \nat.
\]
In other words, the objects of $F(1)$ are the natural numbers, and an arrow
$(m \go n)$ in $F(1)$ is a sequence \bftuple{m_1}{m_n} of natural numbers
such that $m_1 +\cdots+ m_n = m$. That is, the objects are the finite
ordinals and the arrows are the order-preserving functions. So we find that
$F(1) \iso \Delta$, as claimed above.

As the reader may have noticed, a monoidal category does not have to be
strict in order to have an underlying plain multicategory: any monoidal
category will do. If $D$ is the monoidal category then we can define a plain
multicategory $C$ with the same objects as $D$ and with
\[
\multihom{C}{\range{a_1}{a_n}}{a} =
\homset{D}{a_1 \otimes\cdots\otimes a_n}{a}.
\]
In order for this to make sense, $D$ must have $n$-fold tensor products for
all $n$, not just for $n=0$ and $n=2$. There are various attitudes we can
take to this. One is to abandon the usual definition of monoidal category,
and work instead with unbiased monoidal categories, as defined in
Chapter~\ref{ch:bicats}. Another is to use the traditional definition, but to
derive $n$-fold tensors by, for instance, `associating to the left' (as in
Appendix~\ref{app:unbiased}); but this is really just a roundabout version of
the first attitude.

A third is more sophisticated. Take an $n$-leafed tree $\tau$ in which all
nodes have either 0 or 2 outgoing edges: in the language introduced on
page~\pageref{p:classical-trees}, a `classical tree'. This gives a method of
tensoring together $n$ objects in a classical monoidal category, which will
be written
\[
\bftuple{a_1}{a_n} \goesto \ovln{\tau}\bftuple{a_1}{a_n}.
\]
For instance, if $n=2$ then $\tau$ might be the first tree illustrated in
Example~\ref{egs:multicats}\bref{eg:tree-mti} (without its labels), in which
case
\[
\ovln{\tau}\pr{a_1}{a_2} = (a_1 \otimes I) \otimes a_2.
\]
If $\tau$ and $\tau'$ are two $n$-leafed classical trees then there is a
canonical isomorphism 
\[
\omega_{\tau, \tau'}: \ovln{\tau}\bftuple{a_1}{a_n} \goiso
\ovln{\tau'}\bftuple{a_1}{a_n}.
\]
Now, start with a monoidal category $D$, and define from it a plain
multicategory $C$ with the same objects as $D$ and in which a map
$\range{a_1}{a_n} \go a$ is a family 
\[
(f_\tau: \ovln{\tau}\bftuple{a_1}{a_n} \go a)
\]
of maps in $C$ indexed by $n$-leafed classical trees $\tau$, such
that $f_{\tau'} \of \omega_{\tau,\tau'} = f_\tau$ for all $\tau$ and
$\tau'$. Since all of the $f_\tau$'s are determined by any single one of
them, the multicategory $C$ is isomorphic to the one obtained by
associating to the left; however, our new construction does not have the
element of arbitrary choice.

Choosing one version or another of this process, we can compose with the
functor $F$ above to obtain a functor from (non-strict) monoidal categories
to strict monoidal categories. Let $D$ be a monoidal category and $E$ the
resulting strict monoidal category. Then an object of $E$ is a sequence of
objects of $D$, and an arrow 
\[
\bftuple{a_1}{a_m} \go \bftuple{b_1}{b_n}
\]
in $E$ consists of a sequence of arrows
\[
\begin{array}{ccc}
a_1^1 \otimes\cdots\otimes a_1^{k_1} 	&\go 	&b_1,		\\
\ldots					&	&\ldots		\\
a_n^1 \otimes\cdots\otimes a_n^{k_n} 	&\go 	&b_n
\end{array}
\]
in $D$ (with $n$-fold tensors interpreted in the chosen way), such that
the sequence
\[
a_1^1, \ldots, a_1^{k_1}, \ldots, a_n^1, \ldots, a_n^{k_n}
\]
is equal to \range{a_1}{a_m}. Tensor of both objects and arrows in $E$ is by
juxtaposition, and composition comes from the composition in $D$. (This $E$
is not to be confused with the strict monoidal category $\mb{st}(D)$ defined
in \cite[\S 1]{JS}, which is monoidally equivalent to $D$.)

It does not seem straightforward to generalize the notion of (non-strict)
monoidal category to give a notion of weak $T$-structured category, so
for now these observations must be confined to the context of monoidal
categories.

\section{Change of Base}	\label{sec:change}

So far we have only discussed \Cartpr-multicategories for a \emph{fixed}
\Cartpr. In this section we look at what happens when \Cartpr\ varies: in
other words, at how the \Multicat\ construction is functorial. We also
examine functoriality of the structured categories construction, \Struc.

Throughout this section $\Cartpr\hyph\Multicat$ will be regarded as
a (1-)category. Section~\ref{sec:bim} contains an outline of a more
advanced treatment of this material, in which $\Cartpr\hyph\Multicat$ is
treated as an \fc-mul\-ti\-cat\-e\-gory---a categorical structure containing
much more information than a mere category.

First we need to `recall' some definitions from Street's paper~\cite{StrFTM}.

Let $T$ and $T'$ be monads on respective categories \Eee\ and
\Eeep\ (not necessarily cartesian). A \emph{monad functor} 
$\Cartpr \goby{\pr{Q}{\psi}} \Cartprp$ consists of a functor
$\Eee \goby{Q} \Eee'$ together with a natural transformation
\begin{ntdiag}
\Eee		&	&\rTo^{T}	&	&\Eee		\\
		&	&		&\	&		\\
\dTo<{Q}	&	&\ruTo>{\psi}	&	&\dTo>{Q}	\\
		&\	&		&	&		\\
\Eee		&	&\rTo_{T'}	&	&\Eeep		\\
\end{ntdiag}
making the diagrams 
\[
\begin{diagram}[height=2em]
T'^2 Q		&\rTo^{T'\psi}	&T'QT		&\rTo^{\psi T}	&QT^2	\\
\dTo<{\mu' Q}	&		&		&		&\dTo>{Q\mu}\\
T'Q		&		&\rTo_{\psi}	&		&QT	\\
\end{diagram}
\diagspace
\begin{diagram}[height=2em]
Q		&\rEquals	&Q		\\
\dTo<{\eta' Q}	&		&\dTo>{Q \eta}	\\
T'Q		&\rTo_{\psi}	&QT		\\
\end{diagram}
\]
commute. If $\Cartpr \goby{\pr{R}{\chi}} \Cartprp$ is another monad functor
then a \emph{monad functor transformation} $\pr{Q}{\psi} \go \pr{R}{\chi}$ is
a natural transformation $Q \goby{\alpha} R$ such that $(\alpha
T)\of\psi = \chi\of(T'\alpha)$. There is consequently a 2-category \fcat{Mnd}
whose 0-cells are pairs \pr{\Eee}{T}, whose 1-cells are monad functors, and
whose 2-cells are monad functor transformations.

(In fact, \cite{StrFTM} concerns monads and monad functors etc.\ in an
arbitrary 2-category \cat{V}. We are only interested in the case $\cat{V} =
\Cat$.) 

A crucial property of monad functors is that they induce maps between
(Eilenberg-Moore) categories of algebras: thus if \pr{Q}{\psi} is a monad
functor as above then there is an induced functor $\ovln{Q}: \Eee^T \go
\Eee'^{T'}$.

Dually, there is a notion of a \emph{monad opfunctor}, which is just like a
monad functor except that $\psi$ travels in the opposite direction;
similarly, \emph{monad opfunctor transformations}. This gives another
2-category, \fcat{Mnd'}. A monad opfunctor $\Cartpr \go \Cartprp$ induces a
functor $\Eee_T \go \Eee'_{T'}$ between Kleisli categories.

We will also need a third 2-category, $\fcat{Mnd}^{\dashv}$. Again, an object
is a category \Eee\ equipped with a monad $T$. A 1-cell from \Cartpr\ to
\Cartprp\ consists of a monad functor $\pr{Q}{\psi}: \Cartpr \go \Cartprp$, a
monad opfunctor $\pr{P}{\phi}: \Cartprp \go \Cartpr$, and an adjunction
$P\ladj Q$ compatible with the two monads. (Explicitly, this compatibility
means that if $\gamma$ and $\delta$ are the unit and counit of the adjunction
then the diagrams
\[
\begin{diagram}[height=2em]
T'		&\rTo^{T'\gamma}	&T'QP		\\
\dTo<{\gamma T'}&			&\dTo>{\psi P}	\\
QPT'		&\rTo_{Q\phi}		&QTP		\\
\end{diagram}
\diagspace
\begin{diagram}[height=2em]
PT'Q		&\rTo^{\phi Q}		&TPQ		\\
\dTo<{P \psi}	&			&\dTo>{T\delta}	\\
PQT		&\rTo_{\delta T}	&T		\\	
\end{diagram}
\]
commute.) A 2-cell in $\fcat{Mnd}^{\dashv}$ consists of a monad functor
transformation and a monad opfunctor transformation obeying further
compatibility laws. Composition and identities in $\fcat{Mnd}^{\dashv}$ are
defined in the evident way. 

(Incidentally, if we are given a monad opfunctor $\pr{P}{\phi}: \Cartprp
\go\Cartpr$ and a functor $Q$ right adjoint to $P$, then $Q$
naturally becomes a monad functor \pr{Q}{\psi} by taking $\psi$ to be the
mate of $\phi$. The two compatibility squares then commute, so we get a
1-cell of $\fcat{Mnd}^{\dashv}$. This fact is used in the proof of
Proposition~\ref{propn:three-adjns}\bref{propn:n-(n-1)}.)

A monad functor \pr{Q}{\psi} will be called \emph{cartesian} if $Q$ preserves
pullbacks; then cartesian pairs \Cartpr, cartesian monad functors, and all
monad functor transformations form a sub-2-category $\fcat{CartMnd}$ of
$\fcat{Mnd}$. A monad opfunctor \pr{P}{\phi} will be called \emph{cartesian}
if $P$ preserves pullbacks \emph{and $\phi$ is a cartesian natural
transformation}; then cartesian pairs \Cartpr, cartesian monad opfunctors,
and all monad opfunctor transformations form a sub-2-category
$\fcat{CartMnd'}$ of $\fcat{Mnd'}$. Finally, we get a sub-2-category
$\fcat{CartMnd}^{\dashv}$ of $\fcat{Mnd}^{\dashv}$ by allowing only cartesian
pairs \Cartpr\ as objects, 1-cells $(P,\phi,Q,\psi,\gamma,\delta)$ in which
\pr{P}{\phi} is a cartesian monad opfunctor and \pr{Q}{\psi} a cartesian
monad functor, and all 2-cells.

These definitions are rather haphazard: natural transformations are
apparently required to be cartesian (or not) at random. The only
justification I can give is that they seem to be necessary in order to make
the constructions in the rest of this section work. Pulling in the other
direction, if we want the \Struc\ example (diagram~\bref{eq:Struc})
to work then we cannot modify the definition of cartesian monad functor to
include the condition that the natural transformation part $\psi$ is
cartesian---for in that case, it isn't.

We have now collected together the definitions we need, and are ready to see
the three different ways in which the \Multicat\ construction is
functorial. Only an outline of each construction is presented; the details
are easily filled in.

Firstly, let \Cartpr\ and \Cartprp\ be cartesian and let $\pr{Q}{\psi}:
\Cartpr \go\Cartprp$ be a cartesian monad functor. Then there is
an induced functor
\[
\ovln{Q}: \Cartpr\hyph\Multicat \go \Cartprp\hyph\Multicat
\]
defined by pullback. That is, if $C$ is a $T$-multicategory then
$\ovln{Q}(C)$ is a $T'$-multicategory on $QC_0$ whose underlying graph is
given by 
\[
\begin{slopeydiag}
	&	&(\ovln{Q}C)_1\Spbk&	&	&	&	\\
	&\ldTo	&		&\rdTo	&	&	&	\\
T'Q C_0 &	&		&	&QC_1	&	&	\\
	&\rdTo<{\psi_{C_0}}&	&\ldTo>{Qd}&	&\rdTo>{Qc}&	\\
	&	&QT C_0		&	&	&	&QC_0.	\\
\end{slopeydiag}
\]

Dually, let $\pr{P}{\phi}: \Cartprp \go \Cartpr$ be a cartesian monad
opfunctor. Then there is an induced functor
\[
\ovln{P}: \Cartprp\hyph\Multicat \go \Cartpr\hyph\Multicat
\]
defined by composition. That is, if $C'$ is a $T'$-multicategory then
$\ovln{P}(C')$ is a $T$-multicategory on $PC'_0$ whose underlying graph is
given by 
\[
\begin{slopeydiag}
	&	&	&	&PC'_1	&	&	&	&	\\
	&	&	&\ldTo<{Pd}&	&\rdTo(4,4)>{Pc}&&	&	\\
	&	&PT' C'_0&	&	&	&	&	&	\\
	&\ldTo<{\phi_{C'_0}}&&	&	&	&	&	&	\\
TP C'_0	&	&	&	&	&	&	&	&PC'_0.	\\
\end{slopeydiag}
\]

After filling in all the details we get two maps of 2-categories:
\[
\fcat{CartMnd} \go \fcat{CAT}, 
\diagspace
\fcat{CartMnd'} \go \fcat{CAT}.
\]
The first is defined using pullbacks, so is only a weak functor
(pseudo-functor); the second is a strict functor.  On $0$-cells, both
functors send \Cartpr\ to the (large) category $\Cartpr\hyph\Multicat$. At
the `intersection' of \fcat{CartMnd} and \fcat{CartMnd'} is the 2-category
whose 1-cells are what might be called cartesian weak maps of monads: that
is, cartesian monad functors---or equivalently opfunctors---whose natural
transformation part is an isomorphism. Our two functors agree, up to
isomorphism, on these 1-cells.

For the third construction, take a 1-cell in $\fcat{CartMnd}^{\dashv}$ as
shown:
\begin{diagram}[height=2em]
\Cartpr			\\
\uTo<{\pr{P}{\phi}}	
\ladj			
\dTo>{\pr{Q}{\psi}}	\\
\Cartprp		\\
\end{diagram}
Let $\ovln{P}$ and $\ovln{Q}$ be the induced functors just described. Then
there naturally arises an adjunction
\begin{diagram}[height=2em]
\Cartpr\hyph\Multicat	\\
\uTo<{\ovln{P}}		
\ladj			
\dTo>{\ovln{Q}}		\\
\Cartprp\hyph\Multicat.	\\
\end{diagram}
This construction gives a weak functor from $\fcat{CartMnd}^{\dashv}$ to a
suitable 2-category of categories and adjunctions.

As an application of this third construction, take any cartesian
\Cartpr. Then there is a 1-cell
\begin{equation}		\label{eq:Struc}
\begin{diagram}[height=2em]
\pr{\Eee^T}{\id}	\\
\uTo<{\pr{F}{\nu}}	
\ladj			
\dTo>{\pr{U}{\epsln}}	\\
\Cartpr			\\
\end{diagram}
\end{equation}
in $\fcat{CartMnd}^{\dashv}$, in which $F$ and $U$ are the free and forgetful
$T$-algebra functors, and $\nu$ and $\epsln$ are certain canonical natural
transformations which the reader may easily identify. Applying the
construction gives exactly the adjunction
\begin{diagram}[height=2em]
\Cartpr\hyph\Struc	\\
\uTo			
\ladj			
\dTo			\\
\Cartpr\hyph\Multicat	\\
\end{diagram}
of section~\ref{sec:struc}. 

Let us now look at change of base for structured categories. Let \Cartpr\ and
\Cartprp\ be cartesian, and let
\[
\pr{Q}{\psi}: \Cartpr \go \Cartprp
\]
be a cartesian monad functor. This induces a pullback-preserving functor
$\ovln{Q}: \Eee^T \go \Eee'^{T'}$. In turn, \emph{this} induces a functor
from the internal categories in $\Eee^T$ to those in $\Eee'^{T'}$,
\[
\ovln{Q}: \Cartpr\hyph\Struc \go \Cartprp\hyph\Struc.
\]
(The same induced functor results if instead of thinking in terms of internal
categories, we think of the monad functor or opfunctor
\begin{equation}	\label{eq:triv-mon-ftr}
\pr{\ovln{Q}}{1}: \pr{\Eee^T}{\id} \go \pr{\Eee'^{T'}}{\id}
\end{equation}
and use change of base for multicategories. This point of view will be useful
later on.)

Note that this is compatible with the construction of a $T$-structured
category from a $T$-algebra (the functor $D$ on
page~\pageref{p:discrete-struc}), in the sense that the square
\[
\begin{diagram}[height=2em,width=5em]
\Eee^T		&\rTo^{\ovln{Q}}	&\Eee'^{T'}		\\
\dIncl<{D}	&			&\dIncl>{D}		\\
T\hyph\Struc	&\rTo_{\ovln{Q}}	&T'\hyph\Struc		\\
\end{diagram}
\]
commutes. Moreover, change of base for multicategories extends change of base
for structured categories, in the sense that the square
\[
\begin{diagram}[height=2em,width=5em]
T\hyph\Struc	&\rTo^{\ovln{Q}}	&T'\hyph\Struc		\\
\dTo<{\ovln{U}}	&			&\dTo>{\ovln{U}}	\\
T\hyph\Multicat	&\rTo_{\ovln{Q}}	&T'\hyph\Multicat	\\
\end{diagram}
\]
commutes up to canonical isomorphism. Here both $\ovln{U}$'s are the functors
denoted $U$ in~\ref{sec:struc}, and \pr{Q}{\psi} is a cartesian monad functor
(as above). To see that the square commutes, take the monad functor
\pr{\ovln{Q}}{1} of~\bref{eq:triv-mon-ftr} above, and consider the square of
monad functors 
\begin{diagram}[height=2em,width=5em]
\pr{\Eee^T}{\id}	&\rTo^{\pr{\ovln{Q}}{1}}&\pr{\Eee'^{T'}}{\id}	\\
\dTo<{\pr{U}{\epsln}}	&			&\dTo>{\pr{U}{\epsln}}	\\
\Cartpr			&\rTo_{\pr{Q}{\psi}}	&\Cartprp.		\\
\end{diagram}
This square commutes, so by (weak) functoriality the previous square commutes
up to isomorphism.

One might expect a dual to all this, involving monad opfunctors and Kleisli
categories. I do not know what this might be.

\section{Free Multicategories}	\label{sec:free-multicats}

Just as one can form the free category on a directed graph, one can form the
free \Cartpr-multicategory on an \Cartpr-graph, assuming that \Eee\ and $T$
are suitably pleasant. In Appendix~\ref{app:free} we define what it means for
\Cartpr\ to be \emph{suitable} (which is stronger than being cartesian), and
prove the following result:
\begin{thm}	\label{thm:free-main}
Let \Cartpr\ be suitable. Then the forgetful functor
\[
\Cartpr\hyph\Multicat \go \Eeep = \Cartpr\hyph\Graph
\]
has a left adjoint, the adjunction is monadic, and if $T'$ is the resulting
monad on \Eeep\ then \Cartprp\ is suitable.
\end{thm}

When one takes the free category on an ordinary directed graph, the
collection of objects (vertices) is unchanged, and the corresponding fact
for multicategories is expressed in a variant of the theorem. If $S$ is an
object of \Eee\ then we write $\Cartpr\hyph\Multicat_S$ for the subcategory
of $\Cartpr\hyph\Multicat$ whose objects $C$ have $C_0=S$, and whose
morphisms $f$ have $f_0=1_S$; similarly, we write $\Eee'_S$ for the category
of \Cartpr-graphs on $S$ (see~\ref{rmks:maps}\bref{rmk:map-fixed-obj}). 
\begin{thm}	\label{thm:free-fixed}
Let \Cartpr\ be suitable and let $S\in\Eee$. Then the forgetful functor
\[
\Cartpr\hyph\Multicat_S \go \Eee'_S
\]
has a left adjoint, the adjunction is monadic, and if $T'_S$ is the resulting
monad on $\Eee'_S$ then \pr{\Eee'_S}{T'_S} is suitable. Moreover, if \Eee\
has filtered colimits and $T$ preserves them, then the same is true of
$\Eee'_S$ and $T'_S$.
\end{thm}

Most of the time we will only need the weaker conclusions that \Cartprp\ and
\pr{\Eee'_S}{T'_S} are cartesian (rather than suitable); the full recursive
power of the two theorems is only brought into play in a couple of passing
comments (pages~\pageref{p:hierarchy} and~\pageref{p:enrich-hierarchy}). A
point not mentioned elsewhere is that repeated application of
Theorem~\ref{thm:free-fixed} gives an instant definition of a sequence of
sets $(S_n)_{n\in\nat}$ looking very much like the $n$-dimensional opetopes
or multitopes. See \cite[Ch.~IV]{SHDCT} or \cite[4.1]{GOM} for this
construction, and \cite{BaDoHDA3}, \cite{HMP}, \cite{CheROM}
and~\cite{CheEAT} for background.

Our two theorems so far are useless without some instances of suitable
\Cartpr's:
\begin{thm}	\label{thm:free-gen}
Let \Eee\ be a category equivalent to a functor category
\ftrcat{\scat{E}}{\Set}, where \scat{E} is small, and let $T$ be a finitary
cartesian monad on \Eee. Then \Cartpr\ is suitable. 
\end{thm}
Almost all of the specific examples of \Eee\ in this thesis are of the form
\ftrcat{\scat{E}}{\Set}, and all of the monads $T$ are finitary.

The proofs of all three theorems are confined to the appendix. To give
the rough idea, here is a description of the free plain multicategory
construction. Let $\blank^*$ denote the free monoid functor on \Set, and let
\[
X_0^* \og X_1 \go X_0
\]
be a $\blank^*$-graph. This, then, is a set $X_0$ together with a set
$\multihom{X}{\range{x_1}{x_n}}{x}$ for each $x_1, \ldots, x_n, x \in
X_0$. The free plain multicategory $F(X)$ on $X$ has graph
\[
X_0^* \og A \go X_0
\]
where $A$ is defined recursively as follows:
\begin{itemize}
\item if $x\in X_0$ then $\multihom{A}{x}{x}$ has an element $I_x$
\item if $\theta \in \multihom{X}{\range{x_1}{x_n}}{x}$ and 
\[
\alpha_1 \in \multihom{A}{\range{x_1^1}{x_1^{k_1}}}{x_1},\ 
\ldots,\ 
\alpha_n \in \multihom{A}{\range{x_n^1}{x_n^{k_n}}}{x_n}
\] 
then $\multihom{A}{\range{x_1^1}{x_n^{k_n}}}{x}$ has an element
$\theta\langle\range{\alpha_1}{\alpha_n}\rangle$.
\end{itemize}
Here $I_x$ and $\theta\langle\range{\alpha_1}{\alpha_n}\rangle$ are `formal
symbols', and identities and composition in $F(X)$ are defined in ways
suggested by these symbols. 

What this means is that an arrow in $F(X)$ is a tree of arrows in $X$. So for
instance, if $\theta_1, \theta_2, \theta_3, \theta_4$ are arrows in $X$ with
appropriately-matching domains and codomains, then
\[
\theta_1 \langle \theta_2 \langle I_x, \theta_3 \rangle, I_y, \theta_4
\rangle
\]
is an arrow in $F(X)$. Similarly, if $\theta: \range{x_1}{x_n}
\go x$ is any arrow in $X$ then $A$ has an element
\[
\theta \langle \range{I_{x_1}}{I_{x_n}} \rangle. 
\]
This provides the map $X_1\go A$ that determines the unit of the adjunction
at $X$. It also explains why we did not specify that any element of $X$ was
an element of $A$ in the recursive definition---this comes about
automatically. 
 
As a special case, if $X$ is the terminal $\blank^*$-graph (that is, the
terminal object of $\Set/\nat$) then $F(X)$ is the operad \tr\ of
(unlabelled) trees, as described in section~\ref{sec:coherence}. For more on
trees see Example \ref{eg:cart-mnds}\bref{eg:mon-tree}, where labels get
attached to leaves rather than internal nodes.

\section{Algebras via Fibrations}	\label{sec:fib}

It is well-known that for a small category $C$, the functor category
\ftrcat{C}{\Set} is equivalent to the category of discrete opfibrations over
$C$. In this section we extend the notion of discrete opfibration from
categories to general $T$-multicategories, and show that the category of
discrete opfibrations over a given $T$-multicategory is equivalent to its
category of algebras.

By definition, a functor $g: D \go C$ between ordinary categories is a
discrete opfibration if and only if, for any object $b$ of $D$ and arrow
$g(b) \goby{\theta} a$ in $C$, there is a unique arrow $b \goby{\chi} b'$ in
$D$ such that $g(\chi) = \theta$. Another way of saying this is that in the
diagram
\[
\begin{slopeydiag}
	&		&D_1		&		&	\\
	&\ldTo<{d}	&		&\rdTo>{c}	&	\\
D_0	&		&\dTo>{g_1}	&		&D_0	\\
\dTo<{g_0}&		&C_1		&		&\dTo>{g_0}\\
	&\ldTo<{d}	&		&\rdTo>{c}	&	\\
C_0	&		&		&		&C_0	\\
\end{slopeydiag}
\]
depicting $g$, the left-hand `square' is a pullback. 

Generalizing to all cartesian \Cartpr's, let us say that a map $D\goby{g}C$
of $T$-multicategories is a \emph{discrete opfibration} if the square
\begin{diagram}
TD_{0}		&\lTo^{d}	&D_1		\\
\dTo<{Tg_0}	&		&\dTo>{g_1}	\\
TC_{0}		&\lTo^{d}	&C_1		\\
\end{diagram}
is a pullback. We obtain, for any $T$-multicategory $C$, the
category $\fcat{DOpfib}(C)$ of discrete opfibrations over $C$, in
which an object is a discrete opfibration with codomain $C$ and an arrow from
$(D \goby{g} C)$ to $(D' \goby{g'} C)$ is a $T$-multicategory map $D \goby{f}
D'$ such that $g' \of f = g$. (This $f$ is automatically a discrete
opfibration too, by a standard lemma on pasting of pullback squares.)

Notice, incidentally, that being a discrete opfibration is really a property
of maps between $T$-graphs rather than $T$-multicategories. In this sense,
the notion of a discrete opfibration between categories exists at a more
primitive level than the full notion of opfibration.

\begin{thm}	\label{thm:alt-alg}
Let \Cartpr\ be cartesian and let $C$ be a $T$-multicategory. Then there is
an equivalence of categories
\[
% \Alg(C) \eqv \fcat{DOpfib}(C).
\fcat{DOpfib}(C) \eqv \Alg(C).
\]
\end{thm}

\paragraph*{Remark} 
\begin{sloppypar}
A more precise statement is that the forgetful functor
from $\fcat{DOpfib}(C)$ to $\Eee/C_0$ (sending $g$ to $g_0$) is monadic, and
that the induced monad is isomorphic to $T_C$.
\end{sloppypar}

\paragraph*{Proof}
Recall from~\ref{constn:ind-mnd} that a $C$-algebra is an algebra for the
monad $T_C$ on $\Eee/C_0$. The effect of $T_C$ on an object $(X \goby{p}
C_0)$ of $\Eee/C_0$ is given by the pullback diagram
\[
\begin{slopeydiag}
	&	&X'\Spbk&	&	&	&	\\
	&\ldTo<{\phi_X}&&\rdTo>{\pi_X}&	&	&	\\
TX	&	&	&	&C_1&	&	\\
	&\rdTo<{Tp}&	&\ldTo<{d}&	&\rdTo>{c}&	\\
	&	&TC_0	&	&	&	&C_0	\\
\end{slopeydiag}
\]
and the formula $T_C(X \goby{p} C_0) = (X' \goby{c \pi_X} C_0)$. So a
$C$-algebra consists of $(X \goby{p} C_0)$ together with a map $h: X' \go X$
satisfying axioms. 

Given a $C$-algebra \pr{X \goby{p} C_0}{h}, then, we get a commutative
diagram
\[
\begin{slopeydiag}
	&		&X'		&		&	\\
	&\ldTo<{\phi_X}	&		&\rdTo>{h}	&	\\
TX	&		&\dTo>{\pi_X}	&		&X	\\
\dTo<{Tp}&		&C_1		&		&\dTo>{p}\\
	&\ldTo<{d}	&		&\rdTo>{c}	&	\\
TC_0	&		&		&		&C_0.	\\
\end{slopeydiag}
\]
The top part of this diagram defines a $T$-graph $D$, and there is a map $g:
D \go C$ defined by $g_0=p$ and $g_1=\pi_X$. With some calculation we see
that $D$ is naturally a $T$-multicategory and $g$ a map of
$T$-multicategories.  (Composition in $D$ is defined using composition in
$C$, and similarly identities. In the case $\Cartpr=\pr{\Set}{\id}$, we are
dealing with the familiar Grothendieck opfibration.) Moreover, the left-hand
half of the diagram is a pullback, so we have constructed from the
$C$-algebra a discrete opfibration over $C$.

We thus arrive at a functor from $\Alg(C)$ to $\fcat{DOpfib}(C)$, which is
easily checked to be full, faithful and essentially surjective on objects.
\done

\paragraph*{}
Let us take a closer look at the $T$-multicategory $D$ corresponding to a
$C$-algebra $h=\pr{X \goby{p} C_0}{h}$. We could call $D$ the multicategory
of elements or the Grothendieck opfibration of $h$; for reasons soon to
be apparent, I will write $D=C/h$.

A natural question to ask is: given a multicategory $C$ and an algebra $h$
for $C$, what are the algebras for $C/h$? To answer it we recall the process
of slicing a monad by an algebra, as in
Example~\ref{eg:algs}\bref{eg:alg-to-multi}: for any monad $S$ on a category
\cat{D} and any $S$-algebra $k$, there is a monad $S/k$ on \cat{D} with the
property that 
\[
\Alg(S/k) \iso \Alg(S)/k. 
\]
(Here and below, $\Alg$ means the category of algebras for either a monad or
a multicategory. So for instance, $\Alg(C) = \Alg(T_C)$.)

The following two results answer the question. Both proofs are easy.
\begin{propn}	\label{propn:slice-multicat}
Let \Cartpr\ be cartesian, let $C$ be a $T$-multicategory, and let $h$ be a
$C$-algebra. Then there is an isomorphism of monads $T_{C/h} \iso T_C /h$. 
\done
\end{propn}
\begin{cor}	\label{cor:slice-multicat}
In the situation of the proposition, there is an isomorphism of categories
$\Alg(C/h) \iso \Alg(C)/h$.
\done
\end{cor}

The corollary generalizes the familiar fact that when $C$ is a category and
$C/h$ is the category of elements of a functor $h: C\go \Set$,
\[
\ftrcat{C/h}{\Set} \iso \ftrcat{C}{\Set}/h.
\]

In addition to the corollary, we have:
\begin{propn}	\label{propn:slice-multicat-fib}
Let \Cartpr\ be cartesian, let $C$ be a $T$-multicategory, and let $h$ be a
$C$-algebra. Then there is an isomorphism of categories
\[
\fcat{DOpfib}(C/h) \iso \fcat{DOpfib}(C)/h.
\]
\end{propn}
\begin{proof} \sloppy
This follows from standard results on the pasting of pullback squares.\\
\mbox{ }\done
\end{proof}

It is very nearly possible to deduce either one of~\ref{cor:slice-multicat}
or~\ref{propn:slice-multicat-fib} from the other. The only obstacle is that
both results assert the \emph{isomorphism} of a pair of categories, whereas
$\Alg(D)$ and $\fcat{DOpfib}(D)$ are only \emph{equivalent}, for
$T$-multicategories $D$.

As an example, let $C$ be the terminal $T$-multicategory $1$. We
have $T_1 \iso T$ and so $\Alg(1) \iso \Alg(T)$ (Example
\ref{eg:algs}\bref{eg:alg-terminal}). Given a $T$-algebra $h$, we therefore
obtain a $T$-multicategory $1/h$; plausibly enough, this is the
$T$-multicategory of Example~\ref{egs:multicats}\bref{eg:multi-alg}, with
graph
\[
TX \ogby{1} TX \goby{h} X.
\]
The results above tell us that $T_{1/h} \iso T/h$ and $\Alg(1/h) \iso
\Alg(T)/h$, as we also saw in Example~\ref{eg:algs}\bref{eg:alg-to-multi}.

As another application, let us construct the \emph{slice multicategory} $C^+$
of a $T$-multicategory $C$, which will have the property that 
\[
\Alg(C^+) \eqv T\hyph\Multicat/C.
\]
In detail, let \Cartpr\ be suitable, let $\Eeep=T\hyph\Graph$, and let
$T'$ be the free $T$-multicategory monad, as
in~\ref{sec:free-multicats}. Then $C$ is an algebra for the terminal
$T'$-multicategory $1$ (that is, a $T'$-algebra), so  
\[
\Alg(1/C) \iso \Alg(T')/C \eqv T\hyph\Multicat/C.
\]
We therefore define $C^+ = 1/C$, and this has the required property. Notice
that we have moved up a level: whereas $C$ was a $T$-multicategory, $C^+$ is
a $T'$-multicategory.

The slice multicategory construction was first proposed by Baez and Dolan for
their definition~\cite{BaDoHDA3} of weak $n$-category, where it plays a
central part. Their construction takes place in a different and more
specialized context than ours, but there is an evident similarity between the
two. See also~\cite{CheROM} and~\cite{CheEAT} for an elucidation of
Baez-Dolan slicing, and~\cite[IV.4]{SHDCT} for further thoughts on our
version.

\section{Algebras via Endomorphisms}	\label{sec:endo}

The prototypical example of a plain operad arises from substitution. That is,
if $X$ is a set then there is a plain operad $\END(X)$ with
\[
(\END(X))(n) = \homset{\Set}{X^n}{X},
\]
with the identity element of $(\END(X))(1)$ provided by the
identity function on $X$, and with composition in the operad defined by
\[
\theta \of \bftuple{\theta_1}{\theta_n} =
\theta \of (\theta_1 \times\cdots\times \theta_n).
\]

For any plain operad $C$ and set $X$, there is a one-to-one correspondence
between $C$-algebra structures on $X$ and operad maps $C \go
\END(X)$. Indeed, this is often used to \emph{define} what an algebra
for an operad is: for instance, in many accounts of the classical theory of
operads, and in Batanin's account~\cite{Bat} of his globular operads. (In the
classical case a symmetric group action is usually involved too, but we
ignore this elaboration.) In this short section we show that for a large
class of cartesian \Cartpr, this alternative definition of algebra is also
possible.

As motivation, let's consider what the appropriate definition of
$\END$ is for plain \emph{multicategories}. An algebra for a
plain multicategory $C$ consists of a family $(X(a))_{a\in C_0}$ of sets
together with a function
\[
\multihom{C}{\range{a_1}{a_n}}{a} \times
X(a_1) \times\cdots\times X(a_n) 
\go X(a)
\]
for each $a_1, \ldots, a_n, a \in C_0$, satisfying certain axioms. In other
words, a $C$-algebra consists of an object $X \go C_0$ of $\Set/C_0$ together
with a function
\[
\multihom{C}{\range{a_1}{a_n}}{a} 
\go
\homset{\Set}{X(a_1) \times\cdots\times X(a_n)}{X(a)}
\]
for each $a_1, \ldots, a_n, a$, again satisfying axioms. With some work we
see that given any object $X \goby{p} C_0$ of $\Set/C_0$, there is a plain
multicategory $\END(X)$ with object-set $C_0$, with hom-sets
\begin{equation}	\label{eq:endo-plain}
\multihom{(\END(X))}{\range{a_1}{a_n}}{a}
=
\homset{\Set}{X(a_1) \times\cdots\times X(a_n)}{X(a)},
\end{equation}
and with composition and identities given by substitution and identities of
functions; we also see that a $C$-algebra structure on $X$ is exactly a
multicategory map $C \go \END(X)$ which is the identity on
objects.

Analysing this further, let $T$ be the free monoid functor and, given $X
\goby{p} C_0$, consider the following two $T$-graphs on $C_0$:
\[
\begin{diagram}[size=1.5em]
   &          &   &          &TX \times C_0 &            &   &   &    \\
   &          &   &\ldTo<{\mr{pr}_1}&       &\rdTo(4,4)>{\mr{pr}_2}&&&\\
   &          &TX &          &              &            &   &   &    \\
   &\ldTo<{Tp}&   &          &              &            &   &   &    \\
TC_0&         &   &          &              &            &   &   &C_0 \\
\end{diagram}
\diagspace
\begin{diagram}[size=1.5em]
   &       &   &          &TC_0 \times X &            &   &       &   \\
   &       &   &\ldTo(4,4)<{\mr{pr}_1}&  &\rdTo>{\mr{pr}_2}&&     &   \\
   &       &   &          &              &            &X  &       &   \\
   &       &   &          &              &            &   &\rdTo>{p}& \\
TC_0&      &   &          &              &            &   &      &C_0.\\
\end{diagram}
\]
Call these graphs $G_1(X)$ and $G_2(X)$ respectively. In $G_1(X)$, the set of
arrows (that is, elements of $TX \times C_0$) with domain \bftuple{a_1}{a_n}
and codomain $a$ is $X(a_1) \times\cdots\times X(a_n)$; in $G_2(X)$, the set
of arrows with this domain and codomain is $X(a)$. Let $[\ ,\ ]$ denote
exponential in the category $\Set/(TC_0 \times C_0)$ of $T$-graphs on
$C_0$. Then in the $T$-graph $[G_1(X), G_2(X)]$, the set of arrows with the
aforementioned domain and codomain is the right-hand side
of~\bref{eq:endo-plain}. Hence $[G_1(X), G_2(X)]$ is the underlying $T$-graph
of the endomorphism multicategory $\END(X)$ described above.

It is now easy to move to the general case. Let \Cartpr\ be cartesian, and
suppose that each slice $\Eee/Z$ of \Eee\ is cartesian closed. (This happens
if $\Eee \eqv \ftrcat{\scat{E}}{\Set}$ for some small category \scat{E}, as
in almost all of our examples.) Let $S\in\Eee$ and let $X\goby{p} S$ be an
object of $\Eee/S$. Define $T$-graphs $G_1(X)$ and $G_2(X)$ on $S$ by the
same diagrams as above, replacing $C_0$ by $S$ throughout, and define a
$T$-graph 
\[
\END(X) = [G_1(X), G_2(X)],
\]
where $[\ ,\ ]$ is exponential in the category $\Eee/(TS \times S)$ of
$T$-graphs on $S$. Then $\END(X)$ carries a natural $T$-multicategory
structure, as may be verified. Moreover, if $C$ is any $T$-multicategory with
$C_0=S$ then $T$-algebra structures on $X$ correspond one-to-one with those
$T$-multicategory maps $C \goby{h} \END(X)$ which are the identity on objects
(that is, $h_0=1$, in the terminology of~\ref{defn:multifunctor}). Put
another way, an algebra for $C$ is an object $X$ over $C_0$ together with a
map $C \go \END(X)$ of multicategories on $C_0$.

To discuss maps between $C$-algebras (for a fixed $T$-multicategory $C$) we
define
\[
\homset{\HOM}{X}{Y} = [G_1(X), G_2(Y)]
\]
for $T$-graphs $X$ and $Y$ on $S=C_0$. Since both $G_1$ and $G_2$ are
functors, so too is $\HOM$. If \pr{X}{h} and \pr{Y}{k} are
$C$-algebras, then an algebra map $\pr{X}{h} \go \pr{Y}{k}$ is exactly a map
$X \goby{f} Y$ in $\Eee/C_0$ such that the diagram
\begin{diagram}
C	&\rTo^{h}		&\homset{\HOM}{X}{X}	\\
\dTo<{k}&			&\dTo>{\homset{\HOM}{1}{f}}	\\
\homset{\HOM}{Y}{Y}		&
\rTo_{\homset{\HOM}{f}{1}}	&
\homset{\HOM}{X}{Y}						\\
\end{diagram}
commutes. Put formally, we have just given an alternative definition of the
category of algebras for a $T$-multicategory $C$, and this alternative
category is isomorphic to the official category $\Alg(C)$.

\section{\fc-Multicategories}	\label{sec:fcm}

In this section we take a close look at $T$-multicategories in the case where
$T$ is the free category monad, \fc, on the category of directed graphs,
$\pr{\Set}{\id}\hyph\Graph$. This case is interesting for a variety of
reasons. First of all, it arises naturally as soon as one thinks about
categories and the fact that $\Cat$ is monadic over
$\pr{\Set}{\id}\hyph\Graph$. It is therefore the first step in an infinite
hierarchy: that is, if we define
\begin{eqnarray*}	\label{p:hierarchy}
\pr{\Eee^{(0)}}{T^{(0)}}&=	&\pr{\Set}{\id},			\\
\Eee^{(n+1)}		&=	&T^{(n)}\hyph\Graph,			\\
T^{(n+1)}		&=	&\mbox{free } T^{(n)} \mbox{-multicategory}
\end{eqnarray*}
then a $T^{(1)}$-multicategory is an \fc-multicategory. (The validity of
these definitions is guaranteed by Theorems~\ref{thm:free-main}
and~\ref{thm:free-gen}; in particular, they say that the monad \fc\ is
cartesian.) We will only consider this first step here; more can be found in
\cite[3.4]{GECM}.

Secondly, \fc-multicategories encompass many familiar `two-dimensional'
categorical structures, including bicategories, double categories, monoidal
categories and plain multicategories. They also include structures we will
call \emph{weak double categories}, in which composition of horizontal
1-cells only obeys associativity and unit laws up to coherent isomorphism,
and include structures resembling the 2-opetopic sets of Baez and Dolan.

Thirdly, there are a couple of well-known categorical ideas for which
\fc-multicategories provide a more general context than is traditional:
the bimodules construction (usually performed on bicategories), and
the enrichment of categories (usually done in monoidal categories, or
occasionally bicategories). These subjects are treated in,
respectively, sections~\ref{sec:bim} and~\ref{sec:enrich}.

Let us begin by finding out what an \fc-multicategory is in explicit
terms. An \fc-graph $V$ is a diagram
\begin{diagram}[noPS]
	&	&V_1=(V_{11}&\pile{\rTo\\ \rTo}	&V_{10})&	&	\\
	&	&\ldTo	&			&\rdTo	&	&	\\
\fc(V_0)=(V_{01}^*&\pile{\rTo\\ \rTo}&V_{00})&		&
V_0=(V_{01} &\pile{\rTo\\ \rTo}&V_{00}),\\
\end{diagram}
where $V_1$ and $V_0$ are directed graphs, the $V_{ij}$ are sets, $V_{01}^*$
is the set of paths in $V_0$, the horizontal arrows are set maps, and the
diagonal arrows are maps of directed graphs. Think of elements of $V_{00}$ as
\emph{objects} or \emph{0-cells}, elements of $V_{01}$ as \emph{horizontal
1-cells}, elements of $V_{10}$ as \emph{vertical 1-cells}, and elements of
$V_{11}$ as \emph{2-cells}, as in the picture
\begin{equation}	\label{eq:two-cell}
\begin{diagram}[height=2em]
x_0	&\rTo^{m_1}	&x_1	&\rTo^{m_2}	&\ 	&\cdots	
&\ 	&\rTo^{m_n}	&x_n	\\
\dTo<{f}&		&	&		&\Downarrow\,\theta&
&	&		&\dTo>{f'}\\
x	&		&	&		&\rTo_{m}	&	
&	&		&x'	\\
\end{diagram}
\end{equation}
($n\geq 0$, $x_i, x, x' \in V_{00}$, $m_i, m \in V_{01}$, $f, f' \in V_{10}$,
$\theta\in V_{11}$). An \fc-multicategory structure on the \fc-graph $V$
firstly makes 
\[
(\spn{V_{10}}{V_{00}}{V_{00}}), 
\]
the objects and vertical 1-cells, into a category. It also gives a
composition function for 2-cells,
\begin{equation}	\label{eq:pasted-two-cells}
\begin{diagram}[width=.5em,height=1em]
\blob&\rTo^{m_1^1}&\cdots&\rTo^{m_1^{k_1}}&
\blob&\rTo^{m_2^1}&\cdots&\rTo^{m_2^{k_2}}&\blob&
\ &\cdots&\ &
\blob&\rTo^{m_n^1}&\cdots&\rTo^{m_n^{k_n}}&\blob\\
\dTo<{f_0}&&\Downarrow\theta_1&&
\dTo&&\Downarrow\theta_2&&\dTo&
\ &\cdots&\ &
\dTo&&\Downarrow\theta_n&&\dTo>{f_n}\\
\blob&&\rTo_{m_1}&&
\blob&&\rTo_{m_2}&&\blob&
\ &\cdots&\ &
\blob&&\rTo_{m_n}&&\blob\\
\dTo<{f}&&&&&&&&\Downarrow\theta &&&&&&&&\dTo>{f'}\\
\blob&&&&&&&&\rTo_{m}&&&&&&&&\blob\\
\end{diagram}
\end{equation}
$\goesto$
\begin{diagram}[width=.5em,height=1em]
\blob&\rTo^{m_1^1}&\ &&
&&&&\cdots&
&&&
&&\ &\rTo^{m_n^{k_n}}&\blob\\
&&&&&&&&&&&&&&&&\\
\dTo<{f\of f_0}&&&&&&&&\Downarrow\theta\of\tuple{\theta_1}{\theta_2}{\theta_n}
&&&&&&&&\dTo>{f'\of f_n}\\
&&&&&&&&&&&&&&&&\\
\blob&&&&&&&&\rTo_{m}&&&&&&&&\blob\\
\end{diagram}
($n\geq 0, k_i\geq 0$, with \blob's representing objects), and an identity function
\[
\begin{diagram}[width=1em,height=2em]
x&\rTo^{m}&x'\\
\end{diagram}
\diagspace\goesto\diagspace
\begin{diagram}[width=1em,height=2em]
x&\rTo^{m}&x'\\
\dTo<{1_x}&\Downarrow 1_{m}&\dTo>{1_{x'}}\\
x&\rTo_{m}&x'.\\
\end{diagram}
\]
The composition and identities obey associativity and identity laws, which
ensure that any 2-cell diagram with a rectangular boundary has a well-defined
composite. 

The pictures in the nullary case are worth a short comment.
When $n=0$, the 2-cell of diagram~\bref{eq:two-cell} is drawn as
\begin{diagram}[height=2em]
x_0		&\rEquals		&x_0		\\
\dTo<{f}	&\Downarrow\,\theta	&\dTo>{f'}	\\
x		&\rTo_{m}		&x',		\\
\end{diagram}
and the diagram of pasted-together 2-cells in the domain
of~\bref{eq:pasted-two-cells} is drawn as
\begin{diagram}[size=1.5em]
w_0		&\rEquals		&w_0		\\
\dTo<{f_0}	&=			&\dTo>{f_0}	\\
x_0		&\rEquals		&x_0		\\
\dTo<{f}	&\Downarrow\,\theta	&\dTo>{f'}	\\
x		&\rTo_{m}		&x'.		\\
\end{diagram}
The composite%
\label{p:null-notation}
of this last diagram will be written as $\theta\of f_0$.

As such, \fc-multicategories are not familiar, but various degenerate cases
are. These are explained in the following examples, and summarized in
Figure~\ref{fig:degens}.
\begin{figure}\protect\small
\begin{tabular}{r|lll}	
			&Not `representable'	&`Representable'	
&`Uniformly 	\\
			&			&			
&representable'	\\
\hline
No degeneracy		&\fc-multicategory	&Weak double
&Double category\\
			&			&\ category
&		\\
All vertical 1-cells	&Vertically discrete	&Bicategory
&2-category	\\
are identities		&\ \fc-multicategory	&
&		\\
Only one object and	&Plain multicategory	&Monoidal category	
&Strict monoidal \\
one vertical 1-cell	&			&
&\ category
\end{tabular}
\caption{Some of the possible degeneracies of an \fc-multicategory. The
left-hand column refers to degeneracies in the category formed by the objects
and vertical 1-cells. The top row refers to whether the \fc-multicategory
structure arises from a composition rule for horizontal 1-cells. See
Examples~\ref{egs:fcms}.}
\label{fig:degens}
\end{figure}

\begin{eg}{egs:fcms}

\item \label{eg:strict-double} 
Any double category gives an \fc-multicategory, in which a 2-cell as
in~\bref{eq:two-cell} is a 2-cell
\begin{diagram}[height=2em]
x_0	&\rTo^{m_n \of \cdots \of m_1}	&x_n\\
\dTo<{f}&\Downarrow			&\dTo>{f'}\\
x	&\rTo_{m}			&x'\\
\end{diagram}
in the double category.

\item	\label{eg:weak-double}
In fact,~\bref{eg:strict-double} works even when the double category is
`horizontally weak'. A typical example of such a structure---a \emph{weak
double category}---has rings (not necessarily commutative) as its 0-cells,
bimodules as its horizontal 1-cells, ring homomorphisms as its vertical
1-cells, and `homomorphisms of bimodules with respect to the vertical changes
of base' as 2-cells. In other words, a 2-cell looks like
\begin{diagram}[height=2em]
R	&\rTo^{M}		&R'	\\
\dTo<{f}&\Downarrow\theta	&\dTo>{f'}\\
S	&\rTo_{N}		&S',	\\
\end{diagram}
where $R$, $R'$, $S$, $S'$ are rings, $M$ is an \pr{R'}{R}-bimodule (i.e.\
simultaneously a left $R'$-module and a right $R$-module) and $N$ similarly,
$f$ and $f'$ are ring homomorphisms, and $\theta: M \go N$ is an abelian group
homomorphism such that
\[
\theta(r' \cdot m \cdot r) = f'(r') \cdot \theta(m) \cdot f(r).
\]
Composition of
horizontal 1-cells is tensor, composition of vertical 1-cells is the usual
composition of ring homomorphisms, and composition of 2-cells is defined in
an evident way. The essential point is that although the 0-cells and vertical
1-cells form a category, the same cannot be said of the horizontal structure:
tensor only obeys the associative and unit laws up to coherent
isomorphism. 

I will not write down the full definition of weak double category, since it
is just an easy extension of the definition of a bicategory. It is most
convenient to extend the definition of \emph{unbiased} bicategory, since in
order to have a 1-cell `$m_n \of\cdots\of m_1$', as in the diagram
of~\bref{eg:strict-double}, we need $n$-fold composition.

Another example has small categories as 0-cells, profunctors (bimodules) as
horizontal 1-cells, functors as vertical 1-cells, and `morphisms of
profunctors with respect to the vertical functors' as 2-cells.  We will
explore both of these examples further in section~\ref{sec:bim}.

\item	\label{eg:vdfcm}
Suppose that all vertical 1-cells are identities, that is, $V_{10}=V_{00}$
and
\[
(\spn{V_{10}}{V_{00}}{V_{00}}) = (\spaan{V_{00}}{V_{00}}{V_{00}}{1}{1}). 
\]
The category formed by the objects and vertical 1-cells is discrete, so we
may call the \fc-multicategory $V$ \emph{vertically discrete}. In this case,
an alternative way of drawing the underlying \fc-graph of $V$ is as
\[
\begin{slopeydiag}
	&	&V_{11}	&	&	\\
	&\ldTo	&	&\rdTo	&	\\
V_{01}^*&	&	&	&V_{01}	\\
\dTo	&\rdTo(4,2)&	&\ldTo(4,2)&\dTo\\
V_{00}	&	&	&	&V_{00}.\\
\end{slopeydiag}
\]
Thus a vertically discrete \fc-multicategory consists of
some objects $x, x', \ldots$, some 1-cells $m, m',\ldots$, and some
2-cells looking like
\begin{diagram}[size=1.5em,tight]
		&	&	&x_2	&\ldots	&	&	&	\\
		&x_1	&\ruEdge(2,1)^{m_2}&&	&	&x_{n-1}&	\\
\ruEdge(1,2)<{m_1}&	&	&	&\Downarrow \theta&&	
&\rdEdge(1,2)>{m_{n}}							\\
x_0		&	&	&\rEdge_{m}&	&	&	&x_n,	\\
\end{diagram}
together with a composition function
\[
\!\!\!\!\!\!\!\!
\!\!\!\!\!\!\!\!
\!
\piccy{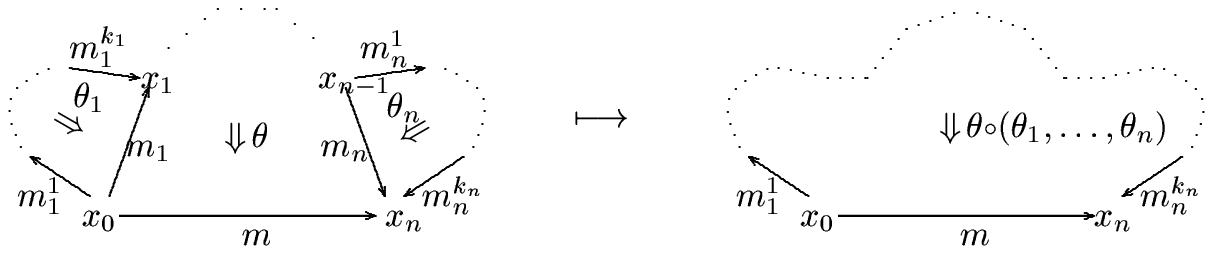}
\]
and an identity function
\[
\piccy{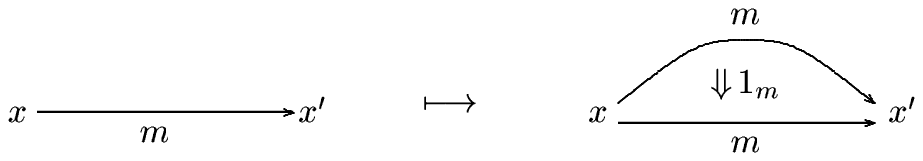}
\]
obeying the inevitable associativity and identity laws. A vertically discrete
\fc-graph bears a strong resemblance to a 2-opetopic set in the sense
of~\cite{BaDoHDA3}, or a 2-truncated multitopic set in the sense
of~\cite{HMP}; see also \cite{CheROM}, \cite{CheEAT} and
\cite[Ch.~IV]{SHDCT}.

\item	\label{eg:bicat-fcm}
Any bicategory gives rise to a vertically discrete \fc-multicategory, in
which a 2-cell as at~\bref{eq:two-cell} is a 2-cell
\[
x_0
\ctwo{m_n \of\cdots\of m_1}{m}{}
x_n
\]
in the bicategory (with $x_0=x$ and $x_n=x'$). This is a special case
of~\bref{eg:weak-double}. 

\item	\label{eg:mon-cat-fcm}
Any monoidal category $M$ gives an \fc-multicategory 
in which there is one
object and one vertical 1-cell, and a 2-cell
\begin{equation}	\label{eq:vt-two-cell}
\begin{diagram}[size=2em,abut]
\bullet	&\rLine^{m_1}	&\bullet	&\rLine^{m_2}	&\bullet&\cdots 
&\bullet	&\rLine^{m_n}	&\bullet	\\
\dLine<1&		&	&		&\Downarrow&
&	&		&\dLine>1\\
\bullet	&		&	&		&\rLine_{m}&
&	&		&\bullet\\
\end{diagram}
\end{equation}
is a morphism $m_n\otimes\cdots\otimes m_1 \go m$ in $M$. This, in turn, is a
special case of~\bref{eg:bicat-fcm}.

\item
Similarly, any plain multicategory $M$ gives an \fc-multicategory:
there is one object, one vertical 1-cell, and a 2-cell \bref{eq:vt-two-cell}
is a map 
\[
\range{m_1}{m_n}\go m
\]
in $M$. In fact, a plain multicategory is exactly an \fc-multicategory in
which the category formed by the objects and vertical 1-cells is \One, the
terminal category.

\item 	\label{eg:fcm-span}
Let \Cartpr\ be cartesian, and define an \fc-multicategory $V$ as
follows. The objects are the objects of \Eee, and the horizontal 1-cells are
the same as the 1-cells of the bicategory $\Span\Cartpr$ defined
in~\ref{constn:bicat}. A vertical 1-cell is a morphism in \Eee, and a 2-cell
\begin{diagram}[height=2em]
X_0	&\rTo^{M_1}	&X_1	&\rTo^{M_2}	&\ 	&\cdots	
&\ 	&\rTo^{M_n}	&X_n	\\
\dTo<{f}&		&	&		&\Downarrow\	&
&	&		&\dTo>{f'}\\
X	&		&	&		&\rTo_{M}	&	
&	&		&X'	\\
\end{diagram}
is a function $\theta$ making
\begin{diagram}[width=2em,height=1em]
	&	&M_n\of\cdots\of M_1	&	&	\\
	&\ldTo	&			&\rdTo	&	\\
TX_0	&	&			&	&X_n	\\
	&	&\dTo>{\theta}		&	&	\\
\dTo<{Tf}&	&			&	&\dTo>{f'}\\
	&	&M			&	&	\\
	&\ldTo	&			&\rdTo	&	\\
TX	&	&			&	&X'	\\
\end{diagram}
commute, where $M_n\of\cdots\of M_1$ is the composite in
$\Span\Cartpr$. Composition and identities in $V$ are defined in the obvious
way. (This is actually not just an \fc-multicategory, but a weak double
category. Strictly speaking, $\Eee$ should be small; but having given an
elementary description of what an \fc-multicategory is, I will feel free to
ignore this restriction.) 

Given any \fc-multicategory, we obtain a vertically discrete
\fc-mul\-ti\-cat\-e\-gory simply by discarding all non-identity vertical
1-cells. Applying this process to $V$ gives the same vertically discrete
\fc-mul\-ti\-cat\-e\-gory as arises from the bicategory $\Span\Cartpr$ by the
method of~\bref{eg:bicat-fcm}. For this reason we also write $\Span\Cartpr$
for the \fc-multicategory $V$. In the next section we will see that it is
useful---and perhaps more natural---to regard $\Span\Cartpr$ as an
\fc-multicategory rather than a bicategory.

\end{eg}

So far all of our examples of \fc-multicategories have been degenerate in
some way: either weak double categories or vertically discrete. The next
section provides some non-degenerate examples.

\section{The Bimodules Construction}	\label{sec:bim}

Bimodules have traditionally been discussed in the context of
bicategories. Thus given a bicategory \cat{B}, one constructs a new bicategory
\Bim{\cat{B}} whose 1-cells are bimodules in \cat{B} (see \cite{CKW} or
\cite{Kos}). The drawback is that this is only possible when \cat{B} has
certain properties concerning the existence and behaviour of local reflexive
coequalizers. 

Here we extend the \fcat{Bim} construction from bicategories to
\fc-mul\-ti\-cat\-e\-go\-ries, which allows us to drop the technical
assumptions. In other words, we will construct an honest functor
\[
\fcat{Bim}: \fc\hyph\Multicat \go \fc\hyph\Multicat. 
\]
This provides lots of new examples of \fc-multicategories.

I would like to be able to, but at present cannot, place the \fcat{Bim}
construction in a more abstract setting: as it stands it is somewhat
\emph{ad hoc}. Possibly there is some connection with the contractions of
Chapter~\ref{ch:defn}.

Let $V$ be an \fc-multicategory. The \fc-multicategory \Bim{V} is defined as
follows.
\begin{description}
\item[0-cells]
\begin{sloppypar}
A 0-cell of \Bim{V} is a multicategory map $1\go V$. That is, it is a 0-cell
$x$ of $V$ together with a horizontal 1-cell $x\goby{t}x$ and 2-cells
\[
\begin{diagram}[height=2em]
x	&\rTo^{t}	&x			&\rTo^{t}	&x	\\
\dTo<{1}&		&\Downarrow\,\mu	&		&\dTo>{1}\\
x	&		&\rTo_{t}		&		&x	\\
\end{diagram}
\diagspace
\begin{diagram}[height=2em]
x	&\rEquals		&x	\\
\dTo<{1}&\Downarrow\,\eta	&\dTo>{1}\\
x	&\rTo_{t}		&x	\\
\end{diagram}
\]
satisfying the usual axioms for a monad, $\mu\of\pr{\mu}{1_t} =
\mu\of\pr{1_t}{\mu}$ and $\mu\of\pr{\eta}{1_t} = 1_t = \mu\of\pr{1_t}{\eta}$.
\end{sloppypar}

\item[Horizontal 1-cells]
A horizontal 1-cell $(x,t,\eta,\mu) \rTo (x',t',\eta',\mu')$ consists of a
horizontal 1-cell $x\goby{m}x'$ in $V$ together with 2-cells
\[
\begin{diagram}[height=2em]
x	&\rTo^{t}	&x			&\rTo^{m}	&x'	\\
\dTo<{1}&		&\Downarrow\,\theta	&		&\dTo>{1}\\
x	&		&\rTo_{m}		&		&x'	\\
\end{diagram}
\diagspace
\begin{diagram}[height=2em]
x	&\rTo^{m}	&x'			&\rTo^{t'}	&x'	\\
\dTo<{1}&		&\Downarrow\,\theta'	&		&\dTo>{1}\\
x	&		&\rTo_{m}		&		&x'	\\
\end{diagram}
\]
satisfying the usual module axioms $\theta\of\pr{\eta}{1_m}=1_m$,
$\theta\of\pr{\mu}{1_m} = \theta\of\pr{1_t}{\theta}$, and dually for
$\theta'$, and the `commuting actions' axiom
$\theta'\of\pr{\theta}{1_{t'}} = \theta\of\pr{1_t}{\theta'}$.

\item[Vertical 1-cells]
A vertical 1-cell
\vslob{(x,t,\eta,\mu)}{}{(\hat{x},\hat{t},\hat{\eta},\hat{\mu})} in \Bim{V}
is a vertical 1-cell \vslob{x}{f}{\hat{x}} in $V$ together with a 2-cell
\begin{diagram}[height=2em]
x		&\rTo^{t}		&x		\\
\dTo<{f}	&\Downarrow\,\omega	&\dTo>{f}	\\
\hat{x}		&\rTo_{\hat{t}}		&\hat{x}	\\
\end{diagram}
such that $\omega\of\mu = \hat{\mu}\of\pr{\omega}{\omega}$ and $\omega\of\eta
= \hat{\eta}\of f$. (The notation on the right-hand side of the second
equation is explained on page~\pageref{p:null-notation}.)

\item[2-cells]
A 2-cell 
\begin{diagram}[height=2em]
t_0	&\rTo^{m_1}	&t_1	&\rTo^{m_2}	&\ 	&\cdots	
&\ 	&\rTo^{m_{n}}	&t_n	\\
\dTo<{f}&		&	&		&\Downarrow	&	
&	&		&\dTo>{f'}\\
t	&		&	&		&\rTo_{m}	&	
&	&		&t'	\\
\end{diagram}
in \Bim{V}, where $t$ stands for $(x,t,\eta,\mu)$, $m$ for
\triple{m}{\theta}{\theta'}, $f$ for \pr{f}{\omega}, and so on, consists of a
2-cell 
\begin{diagram}[height=2em]
x_0	&\rTo^{m_1}	&x_1	&\rTo^{m_2}	&\ 	&\cdots	
&\ 	&\rTo^{m_{n}}	&x_n	\\
\dTo<{f}&		&	&		&\Downarrow\,\alpha&	
&	&		&\dTo>{f'}\\
x	&		&	&		&\rTo_{m}	&	
&	&		&x'	\\
\end{diagram}
in $V$, satisfying the `external equivariance' axioms
\begin{eqnarray*}
\alpha\of(\theta_1,\range{1_{m_2}}{1_{m_n}}) 		&=&
\theta\of\pr{\omega}{\alpha}				\\
\alpha\of(\range{1_{m_1}}{1_{m_{n-1}}},\theta'_n)	&=&
\theta'\of\pr{\alpha}{\omega'}
\end{eqnarray*}
and the `internal equivariance' axioms
\begin{eqnarray*}
\lefteqn{\alpha\of(\range{1_{m_1}}{1_{m_{i-2}}}, \theta'_{i-1}, 1_{m_{i}},
\range{1_{m_{i+1}}}{1_{m_n}})		=}				\\
&&\alpha\of(\range{1_{m_1}}{1_{m_{i-2}}}, 1_{m_{i-1}}, \theta_i,
\range{1_{m_{i+1}}}{1_{m_n}}) 
\end{eqnarray*}
for $2\leq i\leq n$.

\item[Composition and identities]
For both 2-cells and vertical 1-cells in \Bim{V}, composition is
defined directly from the composition in $V$, and similarly identities.

\end{description}

We have now defined an \fc-multicategory \Bim{V} for each
\fc-mul\-ti\-cat\-e\-gory $V$, and it is clear how to do the same thing for
maps of \fc-multicategories, so that we have a functor
\[
\fcat{Bim}: \fc\hyph\Multicat \go \fc\hyph\Multicat.
\]
We could go further and treat $\fc\hyph\Multicat$ as a 2-category (cf.\
the remarks at the end of Chapter~\ref{ch:mtione}). Further still, it is
really more natural to regard $\fc\hyph\Multicat$ as a (large)
\fc-multicategory itself, as we shall see very shortly. Such extensions are
left to the consideration of the reader.

\pagebreak

\begin{eg}{egs:bim}
\item
Let \Bee\ be a bicategory satisfying the conditions on local reflexive
coequalizers mentioned in the first paragraph of this section, so that it is
possible to construct a bicategory \Bim{\Bee} in the traditional way. Let
$V$ be the \fc-multicategory coming from \Bee. Then a 0-cell of \Bim{V} is a
monad in \Bee, a horizontal 1-cell $t\go t'$ is a \pr{t'}{t}-bimodule, and a
2-cell of the form
\begin{diagram}[height=2em]
t_0	&\rTo^{m_1}	&t_1	&\rTo^{m_2}	&\ 	&\cdots	
&\ 	&\rTo^{m_{n}}	&t_n	\\
\dTo<{1}&		&	&		&\Downarrow	&	
&	&		&\dTo>{1}\\
t_0	&		&	&		&\rTo_{m}	&	
&	&		&t_n	\\
\end{diagram}
is a map 
\[
m_n \otimes_{t_{n-1}}\cdots\otimes_{t_1} m_1 \go m
\] 
of \pr{t_n}{t_0}-bimodules, i.e.\ a 2-cell in \Bim{\Bee}. So if we discard
the non-identity 1-cells of \Bim{V} to get a vertically discrete
\fc-multicategory, then this is precisely the \fc-multicategory associated
with the bicategory \Bim{\Bee}.

\item	\label{eg:fcm-Cat}
Let $V$ be the \fc-multicategory \Span\pr{\Set}{\id}, as defined
in~\ref{egs:fcms}\bref{eg:fcm-span}. Then \Bim{V} has
\begin{description}
\item[objects:]
small categories
\item[vertical 1-cells:]
functors
\item[horizontal 1-cells:]
profunctors (that is, a horizontal 1-cell $C\go C'$ is a functor
$C^{\op} \times C' \go \Set$) 
\item[2-cells:]
a 2-cell
\begin{equation}		\label{eq:cat-2-cell}
\begin{diagram}[height=2em]
C_0	&\rTo^{M_1}	&C_1	&\rTo^{M_2}	&\ 	&\cdots	
&\ 	&\rTo^{M_n}	&C_n	\\
\dTo<{F}&		&	&		&\Downarrow	&
&	&		&\dTo>{F'}\\
C	&		&	&		&\rTo_{M}	&	
&	&		&C'	\\
\end{diagram}
\end{equation}
is a family of functions
\[
M_n(a_{n-1},a_n) \times\cdots\times M_1(a_0,a_1) \go M(Fa_0,F'a_n),
\]
one for each $a_0 \in C_0$, \ldots, $a_n \in C_n$, natural in the
$a_i$'s. 
\end{description}

The 2-cells can be described another way. Firstly, there is a profunctor
$M': C_0 \go C_n$ defined by $M'(a_0,a_n) = M(Fa_0,F'a_n)$. Secondly, we can
tensor together (compose) the profunctors $M_i$ to obtain the profunctor $M_n
\otimes\cdots\otimes M_1: C_0 \go C_n$. A 2-cell as shown above is then a
morphism $M_n \otimes\cdots\otimes M_1 \go M'$ of profunctors, in the usual
sense. 

In particular, our \fc-multicategory (which could reasonably be called \Cat)\
incorporates natural transformations. For let $D$ and $C$ be categories and
$F,F': D \go C$ functors; write $I_D$ and $I_C$ for the identity
profunctors on $D$ and $C$, i.e. $I_D = \Hom_D$ and $I_C = \Hom_C$.
Then by a simple Yoneda argument, a 2-cell
\begin{equation}	\label{eq:transf-square}
\begin{diagram}[height=2em]
D		&\rTo^{I_D}	&D		\\
\dTo<{F}	&\Downarrow	&\dTo>{F'}	\\
C		&\rTo_{I_C}	&C		\\
\end{diagram}
\end{equation}
in \Bim{\Span\pr{\Set}{\id}} is just a natural transformation $F\go F'$.

\item		\label{eg:Bim-Cartpr} 
More generally, consider \Bim{\Span\Cartpr} for any cartesian
\Cartpr. As we would expect, an object is a $T$-multicategory and a vertical
1-cell is a map of $T$-multicategories. A horizontal 1-cell $C\go C'$ is a
\emph{profunctor} or \emph{(bi)module} between $T$-multicategories: that is,
a span
\[
TC_0 \og M \go C'_0
\]
together with maps (`actions') $M\of C_1 \go M$ and $C'_1 \of M \go
M$ obeying the usual rules for a bimodule. Here and in what
follows, `$\of$' indicates composition of 1-cells in the bicategory
$\Span\Cartpr$. A 2-cell as pictured in~\bref{eq:cat-2-cell} is a map
$\theta$ in \Eee\ making the diagram
\begin{diagram}[width=2em,height=1em]
	&	&M_n \of\cdots\of M_1	&	&	\\
	&\ldTo	&			&\rdTo	&	\\
TC_0	&	&			&	&C_n	\\
	&	&\dTo>{\theta}		&	&	\\
\dTo<{TF}&	&			&	&\dTo>{F'}\\
	&	&M			&	&	\\
	&\ldTo	&			&\rdTo	&	\\
TC	&	&			&	&C'	\\
\end{diagram}
commute and satisfying compatibility axioms for the actions by the $C_i$'s,
$C$ and $C'$. 

\begin{sloppypar}
This provides a family of examples of \fc-multicategories which are not
degenerate in any of the ways described in~\ref{egs:fcms}. In other words,
\Bim{\Span\Cartpr} does not usually form a weak double category. For recall
that in order to form the tensor of ordinary profunctors (as in the previous
example), one needs to use a certain coend, which is effectively a reflexive
coequalizer in the category of sets. Similarly, in order to form a tensor of
profunctors in the \Cartpr\ setting we need \Eee\ to possess certain
reflexive coequalizers, and in order for tensor to obey (weak) associative
and unit laws we need $T$ to preserve such coequalizers.  In general \Eee\
and $T$ will not have these properties.
\end{sloppypar}

(A rather self-referential example is provided by $T=\fc$, which does not
preserve all reflexive coequalizers. This is all at a pragmatic level; I have
not actually got a proven counterexample to the claim that for all cartesian
\Cartpr, the \fc-multicategory $\Bim{\Span\Cartpr}$ comes from a weak double
category.)

So, for a fixed \Cartpr, the natural structure formed by
$T$-mul\-ti\-cat\-e\-go\-ries is an \fc-multicategory. As for ordinary
categories, this incorporates a sensible notion of natural
transformation. Formally, if $C$ is a $T$-multicategory then let $I_C$ denote
the profunctor $C\go C$ consisting of the span
\[
TC_0 \ogby{d} C_1 \goby{c} C_0
\]
with left and right $C$-actions defined by composition in $C$. Let $D$ and
$C$ be $T$-multicategories and $F, F': D \go C$ maps of $T$-multicategories:
then a \emph{transformation} $F\go F'$ is a 2-cell in
\Bim{\Span\Cartpr} as shown in diagram~\bref{eq:transf-square}. An elementary
definition of transformation is given in \cite[1.1.1]{GECM}. For plain
multicategories, a transformation $\alpha: F\go F'$ consists of an arrow
$\alpha_d: Fd \go F'd$ for each $d\in D_0$, such that
\[
\alpha_d \of (Fg) = F'g \of \bftuple{\alpha_{d_1}}{\alpha_{d_n}}
\]
for all arrows $g: \range{d_1}{d_n} \go d$ in $D$. 

\item		\label{eg:Bim-Ab}
For a less taxing example, let $V$ be the \fc-multicategory coming from the
monoidal category \triple{\Ab}{\otimes}{\integers} (as
in~\ref{egs:fcms}\bref{eg:mon-cat-fcm}). Then \Bim{V} has
\begin{description}
\item[objects:]
rings
\item[vertical 1-cells:]
ring homomorphisms
\item[horizontal 1-cells $R\go R'$:]
\pr{R'}{R}-bimodules
\item[2-cells:]
A 2-cell
\begin{diagram}[height=2em]
R_0	&\rTo^{M_1}	&R_1	&\rTo^{M_2}	&\ 	&\cdots	
&\ 	&\rTo^{M_n}	&R_n	\\
\dTo<{f}&		&	&		&\Downarrow\,\theta&
&	&		&\dTo>{f'}\\
R	&		&	&		&\rTo_{M}	&	
&	&		&R'	\\
\end{diagram}
is a function $M_n\times\cdots\times M_1 \goby{\theta} M$ which preserves
addition in each coordinate (is `multi-additive'), and satisfies
\begin{eqnarray*}
\theta(r_n \cdot m_n, m_{n-1}, \ldots)	&=&
f'(r_n) \cdot \theta(m_n, m_{n-1}, \ldots)	\\
\theta(m_n \cdot r_{n-1}, m_{n-1}, \ldots)	&=&
\theta(m_n, r_{n-1} \cdot m_{n-1}, \ldots)
\end{eqnarray*}
etc.
\end{description}
This, then, is the \fc-multicategory arising from the weak double category
of~\ref{egs:fcms}\bref{eg:weak-double}. As in the last example, it is only a
weak double category because certain reflexive coequalizers exist and behave
well in the monoidal category \triple{\Ab}{\otimes}{\integers}. 

\item
\begin{sloppypar}
The previous example can be repeated with the monoidal category
\triple{\Set}{\times}{1}, with obviously analogous results. 
\end{sloppypar}

\item 
Let $W$ be a 2-category. We construct from $W$ an \fc-multicategory $V$,
different from the vertically discrete \fc-multicategory
of~\ref{egs:fcms}\bref{eg:bicat-fcm}. The objects of $V$ are the objects of
$W$, the vertical and horizontal 1-cells of $V$ are both just the 1-cells of
$W$, and a 2-cell
\begin{diagram}[height=2em]
x_0	&\rTo^{g_1}	&x_1	&\rTo^{g_2}	&\ 	&\cdots	
&\ 	&\rTo^{g_n}	&x_n	\\
\dTo<{f}&		&	&		&\Downarrow	&
&	&		&\dTo>{f'}\\
x	&		&	&		&\rTo_{g}	&	
&	&		&x'	\\
\end{diagram}
in $V$ is a 2-cell
\[
x_0%
\ctwo{g\of f}{f'\of g_n \of\cdots\of g_1}{}%
x_n
\]
in $W$. Composition and identities are defined by pasting of 2-cells in
$W$. (Effectively we are associating a (strict) double category to $W$,
and obtaining from that an \fc-multicategory as
in~\ref{egs:fcms}\bref{eg:strict-double}.)

Now consider the \fc-multicategory \Bim{V}. The category formed by the
objects and vertical 1-cells is the category of monads and monad functors in
$W$, in the sense of~\ref{sec:change} and~\cite{StrFTM}. A horizontal 1-cell
is what might be called a (bi)module between monads. (For an application of
such modules to `hard-nosed mathematics'---homotopy theory, in fact---see
\cite[9.4]{MayGIL}.) The description of a general 2-cell is omitted.

Dually, we can reverse direction of the 2-cells in $V$ to obtain another
\fc-multicategory $V'$ from $W$, and then the objects and vertical 1-cells of
\Bim{V'} form the category of monads and monad opfunctors in $W$.

\end{eg}

Example~\bref{eg:Bim-Cartpr} points the way towards a truly uncompromising,
but logically superior, approach to writing up the general theory of
multicategories. This approach would start with an elementary definition
of (possibly large) \fc-multicategory, which would look like the description
at the beginning of~\ref{sec:fcm}. It would continue with the definition of
the \fc-multicategory $\Span\Cartpr$, for any cartesian \Cartpr, and a
definition of the bimodules construction. By applying the latter to the
former it would arrive at the \fc-multicategory \Cartpr\hyph\Multicat. A
\Cartpr-multicategory would be, by definition, an object of this
\fc-multicategory, and similarly maps, modules, etc. So in this approach, it
would not be necessary to treat $\Span\Cartpr$ as a bicategory at all.

The `change of base' discussed in section~\ref{sec:change} can also be
explained using the bimodules construction. Given cartesian \Cartpr\ and
\Cartprp, a cartesian monad functor from \Cartpr\ to \Cartprp\ gives rise to
a map
\[
\Span\Cartpr \go \Span\Cartprp
\]
of \fc-multicategories. (This is a direct and explicit construction, of
which no further explanation is offered.) The same is true for a cartesian
monad \emph{op}functor, using a different construction. Applying \fcat{Bim}
then gives an \fc-multicategory map
\[
\Cartpr\hyph\Multicat \go \Cartprp\hyph\Multicat.
\]
So a cartesian monad (op)functor tells us not just how to turn a
$T$-mul\-ti\-cat\-e\-gory into a $T'$-mul\-ti\-cat\-e\-gory, and a functor
between $T$-multicategories into a functor between $T'$-multicategories (as
in~\ref{sec:change}), but also works on profunctors, transformations, and the
general 2-cells described in Example~\bref{eg:Bim-Cartpr}.

\section{Enrichment}	\label{sec:enrich}

There is a quite surprising theory of enrichment for general
multicategories. In this section I will give a short outline of the shape of
the theory, referring the reader to~\cite{GECM} for a more full account.

The main surprise is what one enriches in. Given a category \Eee\ and a monad
$T$ on \Eee, which are `suitable' in the sense of~\ref{sec:free-multicats},
define
\begin{eqnarray*}
\Eeep	&=	&\Cartpr\hyph\Graph,	\\
T'	&=	&\mbox{free } \Cartpr\mbox{-multicategory}.
\end{eqnarray*}
Fix a $T'$-multicategory $V$ (which makes sense as \Cartpr\ is
suitable). Then we will talk about `$T$-multicategories enriched in $V$'. In
other words, we enrich $T$-multicategories in $T'$-multicategories. This
means that we can take, for instance, the hierarchy
\label{p:enrich-hierarchy}%
$\pr{\Eee^{(n)}}{T^{(n)}}$ of monads defined at the beginning
of~\ref{sec:fcm}, and consider $T^{(n)}$-multicategories enriched in
$T^{(n+1)}$-multicategories; thus a structure of one type gets enriched in a
structure of a more complicated type.

(In general there appears to be no such thing as the `underlying'
$T$-mul\-ti\-cat\-e\-gory of a $V$-enriched $T$-multicategory, in contrast
with the familiar situation for categories enriched in a monoidal category.)

The definition itself is very simple. Given an object $C_0$ of \Eee, we can
form $I(C_0)$ (with $I$ for indiscrete), the unique $T$-multicategory with
graph
\[
TC_0 \ogby{\mr{pr}_1} TC_0 \times C_0 \goby{\mr{pr}_2} C_0.
\]
Then $I(C_0)$ is a $T'$-algebra, say $h: T'(I(C_0)) \go I(C_0)$. By
Example~\ref{egs:multicats}\bref{eg:multi-alg}, we get from this $MI(C_0)$,
the unique $T'$-multicategory with graph
\[
T'(I(C_0)) \ogby{1} T'(I(C_0)) \goby{h} I(C_0).
\]
For a $T'$-multicategory $V$, a \emph{$V$-enriched $T$-multicategory} is an
object $C_0$ of \Eee\ together with a map $MI(C_0) \go V$ of
$T'$-multicategories. Maps between $V$-enriched $T$-multicategories are
also defined in a simple way (see \cite{GECM}), thus giving a category.

The simplest case is $\Cartpr=\pr{\Set}{\id}$. Then $\Eee'$ is the category
of directed graphs, $T'=\fc$, and we have a theory of categories enriched in
an \fc-multicategory. This extends the usual theory of categories enriched in
a monoidal category, as well as the less popular theory of categories
enriched in a bicategory (\cite{BCSW}, \cite{CKW}, \cite{Wal}) and the
evident but hardly-written-up theory of categories enriched in a plain
multicategory. Categories enriched in an \fc-multicategory are examined in
each of~\cite{FCM}, \cite{GECM} and~\cite{GEC}.

The theory of bimodules interacts with the theory of enrichment in an
\fc-multicategory in the following way. Write $V\hyph\Cat$ for the category
of categories enriched in an \fc-multicategory $V$. We then have some facts:
\begin{enumerate}
\item given a map $V_1 \go V_2$ of \fc-multicategories, there is an induced
functor $V_1 \hyph\Cat \go V_2 \hyph\Cat$
\item 	\label{fact:forgetful}
there is a forgetful map $\Bim{V} \go V$, for any $V$
\item 	\label{fact:iso}
the forgetful map $\Bim{MI(C_0)} \go MI(C_0)$ is an isomorphism for
any set $C_0$ (which takes a little thought)
\item by~\bref{fact:iso}, a $V$-enriched category $(MI(C_0) \goby{\gamma} V)$
gives rise to a \Bim{V}-enriched category 
\[
MI(C_0) \goiso \Bim{MI(C_0)} \goby{\Bim{\gamma}} \Bim{V}
\]
\item the same goes for maps, so there is a functor
\[
V\hyph\Cat \go \Bim{V}\hyph\Cat.
\]
\end{enumerate}
(As it happens, this functor is right adjoint to the functor induced by the
forgetful map of~\bref{fact:forgetful}.) 

For instance, a category $C$ enriched in the monoidal category \Ab\ gives
rise to a category enriched in the \fc-multicategory \Bim{\Ab}
of~\ref{egs:bim}\bref{eg:Bim-Ab}. In concrete terms, this works because the
abelian group \homset{C}{a}{a} is naturally a ring, and the abelian group
\homset{C}{a}{b} is naturally a left \homset{C}{b}{b}-module and a right
\homset{C}{a}{a}-module, for any $a,b \in C_0$. Further explanation is
in~\cite{GECM} and~\cite{FCM}.

This piece of theory again illustrates the advantages of working with
\fc-multicategories instead of bicategories. Let $\Bee$ be a bicategory
satisfying the usual conditions on local reflexive coequalizers, so that
there is a bicategory \Bim{\Bee}. Then, just as above, any $\Bee$-enriched
category gives rise to a \Bim{\Bee}-enriched category. However, this
construction is not functorial: a map between $\Bee$-enriched categories does
not give rise to a map between the associated \Bim{\Bee}-enriched
categories. Essentially, the problem is that the definition of a map between
$\cat{C}$-enriched categories (for a bicategory $\cat{C}$, which in this case
is \Bim{\Bee}) is too restrictive; and in turn, this restrictive definition
is forced because bicategories do not have any vertical 1-cells. Once again,
the reader is referred elsewhere for elucidation of cryptic remarks: see
\cite{GECM} or \cite{GEC}.

The second-most simple case of enrichment for general multicategories is when
$\Eee=\Set$ and $T$ is the free monoid monad. This has an interesting
application, concerning the structures called `pseudo-monoidal categories' by
Soibelman and `relaxed multicategories' by Borcherds. (See~\cite{Soi}
and~\cite{Borh}, and~\cite{SnyEBG} and~\cite{SnyRMS} for further
explanation. Borcherds actually used relaxed multi\emph{linear} categories,
where the hom-sets are not just sets but vector spaces.)
In~\cite[Ch.~4]{GECM} it is shown that, for a certain naturally-arising
$T'$-multicategory $V$, relaxed multicategories are exactly plain
multicategories enriched in $V$.

\chapter{A Definition of Weak $\omega$-Category}	\label{ch:defn}

In this chapter we present a definition of weak $\omega$-category, a
variation on that given by Batanin in~\cite{Bat}. We start~(\ref{sec:formal})
by giving the definition in purely formal terms, which can be done very
quickly. However, it is the explanation of why it is a \emph{reasonable}
definition that occupies most of the chapter
(\ref{sec:pds}--\ref{sec:examples}).

It turns out that there are (at least) two natural ways to use our definition
of weak $\omega$-category to give a definition of weak $n$-category. We show
that these two definitions are equivalent in a strong
sense~(\ref{sec:weak-n}). Moreover, we show~(\ref{sec:weak-2}) that weak
2-categories are the same as unbiased bicategories.

In order to make the definition of weak $\omega$-category we need to rely on
certain technical results, which are confined to
Appendices~\ref{app:strict-omega} and~\ref{app:initial}.

Our definition of weak $\omega$-category is very close to Batanin's, although
not the same. Both definitions involve two main ideas: operads and
contractions. The operads he uses \emph{are} the same as the \Cartpr-operads
here (for the particular choice of \Cartpr\ that we will make), and part of
the purpose of this chapter is to explain in elementary language and pictures
what these \Cartpr-operads are, so that the knowledgeable reader may
understand that the two kinds of operad are the same. (More precisely, our
\Cartpr-operads are what Batanin calls `$\omega$-operads in \textit{Span}'.)
However, Batanin's notion of contraction is different from the one here. The
difference between the two definitions is explained further at the end
of~\ref{sec:the-defn}.

\section{Formal Account}	\label{sec:formal}

Let \scat{G} be the category whose objects are the natural numbers
$0,1,\ldots$, and whose arrows are generated by
\[
\sigma_n, \tau_n: n \go n-1
\]
for each $n\geq 1$, subject to equations
\[
\sigma_{n-1} \of \sigma_n = \sigma_{n-1} \of \tau_n,
\diagspace
\tau_{n-1} \of \sigma_n = \tau_{n-1} \of \tau_n
\]
($n\geq 2$). A functor $X: \scat{G} \go \Set$ is called a \emph{globular
set}; I will write $s$ instead of $X(\sigma_n)$, and $t$ instead of
$X(\tau_n)$. 

Any (small) strict $\omega$-category has an underlying globular set $X$, in
which $X(n)$ is the set of $n$-cells and $s$ and $t$ are the source and
target maps. Moreover, a strict $\omega$-functor induces a map of underlying
globular sets, so there is a forgetful functor from the category
$\omega\hyph\Cat$ (of strict $\omega$-categories and strict
$\omega$-functors) to the category \ftrcat{\scat{G}}{\Set} of globular sets.
In Appendix~\ref{app:strict-omega} we put this into exact terms and
establish:
\begin{propn}	\label{propn:formal-strict-omega}
The forgetful functor $\omega\hyph\Cat \go \ftrcat{\scat{G}}{\Set}$ has a
left adjoint, and the induced monad \triple{T}{\eta}{\mu} on
\ftrcat{\scat{G}}{\Set} is cartesian.
\end{propn}

This proposition means that it makes sense to talk about $T$-operads. Let $C$
be a $T$-operad. The underlying $T$-graph of $C$ is a diagram $(C \goby{d}
T1)$ in \ftrcat{\scat{G}}{\Set}; if $\nu\in (T1)(n)$, write 
\[
C(\nu) = \{ \theta \in C(n) \such d(\theta) = \nu \}.
\]

For $n\geq 2$ and $\pi\in (T1)(n)$, define
\[
P_{\pi}(C) = 
\{ \pr{\theta_0}{\theta_1} \in C(s(\pi)) \times C(t(\pi)) \such
s(\theta_0) = s(\theta_1) \mbox{ and } t(\theta_0) = t(\theta_1)\},
\]
and for $\pi\in (T1)(1)$, define
\[
P_{\pi}(C) = C(s(\pi)) \times C(t(\pi)).
\]
% \label{p:defn}%
A \emph{contraction} $\kappa$ on $C$ is a family of functions
\[
(\kappa_\pi: P_{\pi}(C) \go C(\pi))_{n\geq 1, \pi\in (T1)(n)} ,
\]
satisfying
\[
s(\kappa_\pi \pr{\theta_0}{\theta_1} ) = \theta_0,
\diagspace
t(\kappa_\pi \pr{\theta_0}{\theta_1} ) = \theta_1
\]
for every $n\geq 1$, $\pi\in (T1)(n)$ and $\pr{\theta_0}{\theta_1} \in
P_{\pi}(C)$. 

An \emph{operad-with-contraction} is a pair \pr{C}{\kappa} in which $C$ is a
$T$-operad and $\kappa$ is a contraction on $C$. Let \fcat{OWC} be the
category whose objects are operads-with-contraction, and in which a map
$\pr{C}{\kappa} \go \pr{C'}{\kappa'}$ is a map $F: C\go C'$ of $T$-operads
such that for all $n\geq 1$, $\pi\in (T1)(n)$ and $\pr{\theta_0}{\theta_1}
\in P_{\pi}(C)$,
\[
F(\kappa_{\pi} \pr{\theta_0}{\theta_1}) =
\kappa'_{\pi} \pr{F(\theta_0)}{F(\theta_1)}.
\]
(It is easy to verify that the right-hand side makes sense, i.e.\
that $\pr{F(\theta_0)}{F(\theta_1)} \in P_{\pi}(C')$.)

In Appendix~\ref{app:initial} we prove the following:
\begin{propn}
\fcat{OWC} has an initial object.
\end{propn}
Write \pr{L}{\lambda} for the initial object. This determines the $T$-operad
$L$ up to isomorphism; and since the algebras construction is functorial, the
category $\Alg(L)$ is determined up to isomorphism.  

\begin{defn}
A \emph{weak $\omega$-category} is an $L$-algebra. 
\end{defn}

It is not meant to be obvious why this is a reasonable definition of weak
$\omega$-category, and the next few sections are devoted to an explanation.%
\label{p:end-of-defn}

\section{Pasting Diagrams}	\label{sec:pds}

Before understanding weak $\omega$-categories, we must first understand
strict ones, and in particular we need to know about the free strict
$\omega$-category monad on the category of globular sets. In
Appendix~\ref{app:strict-omega} we prove the existence and relevant
properties of this monad, and that the category of strict $\omega$-categories
is monadic over the category of globular sets.
Here we give pictorial descriptions.

First let us contemplate the globular set $T(\One)$, where
\[
\One = (\cdots \pile{\rTo \\ \rTo} 1 \pile{\rTo \\ \rTo} \cdots
\pile{\rTo \\ \rTo} 1)	\\
\] 
is the terminal globular set. The free strict $\omega$-category functor takes
a globular set $X$ and creates formally all possible composites in it, to
make $TX$. Thus a typical element of $(T\One)(2)$ looks like
\begin{equation}
\gfst{}%
\gfour{}{}{}{}{}{}{}%
\grgt{}%
\gone{}%
\glft{}%
\gtwo{}{}{}%
\glst{},
\label{pic:four-one-three-glob}
\end{equation}
where each $k$-cell drawn represents the unique member of $\One(k)$. Note
that because of identities (which we think of throughout as nullary
composites), this diagram might be thought of as representing an element of
$(T\One)(n)$ for any given $n\geq 2$. Let us call an element of $(T\One)(n)$
(or the picture representing it) an \emph{$n$-pasting diagram}, and define
$\pd=T\One$. (The sets $\pd(m)$ and $\pd(n)$ are considered disjoint, when
$m\neq n$.) This 2-pasting diagram~(\ref{pic:four-one-three-glob}) has a
source and a target, both of which are the 1-pasting diagram
\[
\gfst{}%
\gone{}%
\gblw{}%
\gone{}%
\gblw{}%
\gone{}%
\glst{}.
\]
Since all cells in \One\ have the same source and target---are
`endomorphisms'---it is inevitable that the same should be true in $T\One =
\pd$.

It is not hard to give a concrete description of the globular set
$\pd$. Write $\blank^*$ for the free monoid functor on \Set: then $\pd(0)=1$
and $\pd(n+1) = \pd(n)^*$. In other words, an $(n+1)$-pasting diagram is a
sequence of $n$-pasting diagrams. For example, the $2$-pasting diagram
depicted in~\bref{pic:four-one-three-glob} is the sequence
\[
(\gfst{}\gone{}\gblw{}\gone{}\gblw{}\gone{}\glst{},
\gzero{},
\gfst{}\gone{}\glst{})
\]
of $1$-pasting diagrams, so if $\pd(0) = \{ \blob \}$
then~\bref{pic:four-one-three-glob} is the double sequence
\[
((\blob, \blob, \blob), (), (\blob)) 
\in \pd(2).
\]
The source and target maps $s,t: \pd(n+1) \go \pd(n)$ are equal, and we will
write both as $\bdry$ (`boundary'); $\bdry$ is defined inductively by
\[
(\pd(n+1) \goby{\bdry} \pd(n)) = (\pd(n) \goby{\bdry} \pd(n-1))^*.  
\]
The correctness of this description of $T\One$ follows from the results of
Appendix~\ref{app:strict-omega}.

Having described \pd\ as a globular set, we next turn to its strict
$\omega$-category structure: in other words, how pasting diagrams may be
composed.

Typical binary compositions are illustrated by
\[
\piccy{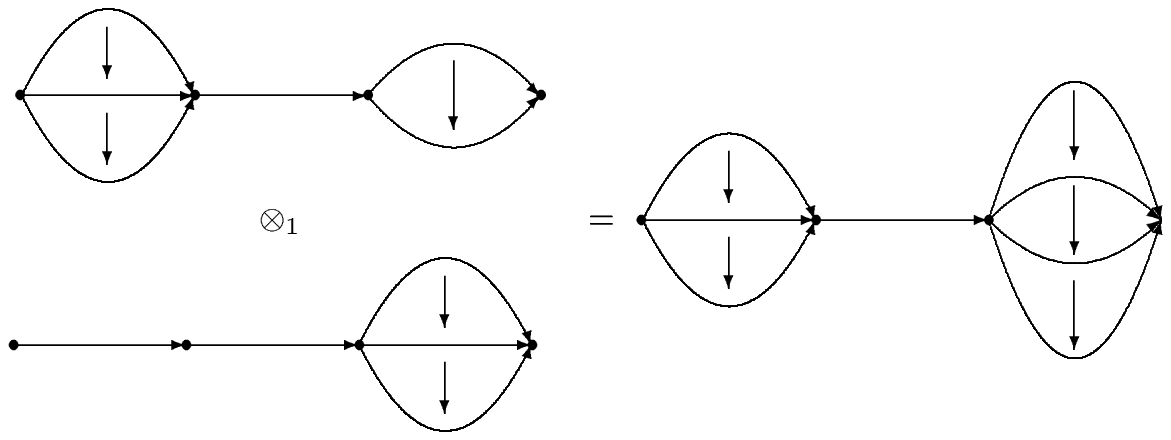}%\label{pic:bin-comps}
\]
and
\[
\piccy{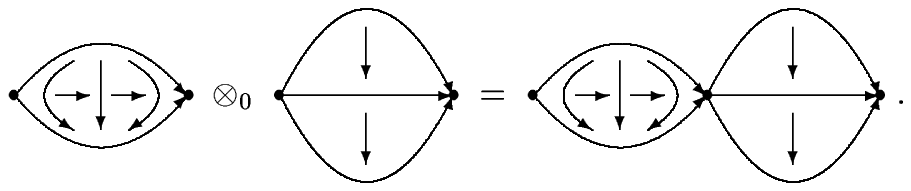}
\]
These compositions are possible because the sources/targets match
appropriately: e.g.\ in the first calculation, where we are gluing along
1-cells (indicated by $\otimes_1$), the 1-dimensional parts of the two
arguments are the same. A typical 
nullary composition---identity---is
\[
\!\!\!\!\!\!\!\!
\!\!
\begin{array}{ccc}	
\gfst{}\gone{}\gblw{}\gone{}\gblw{}\gone{}\glst{}	&
\goesto	&
\gfst{}\gone{}\gblw{}\gone{}\gblw{}\gone{}\glst{}	\\
\elt \pd(1)	&	&\elt \pd(2).
\end{array}
\]

It is helpful to ponder not just binary and nullary composition in \pd, but
composition indexed by arbitrary shapes, in the sense now explained.  We may
think of the first binary composition above as indexed
by 
\[
\gfst{}\gthree{}{}{}{}{}\glst{} \ \elt \ \pd(2),
\]
because we were composing one 2-cell with another by joining along their
bounding 1-cells.  The composition can be represented as
\begin{equation}	%\label{pic:first-bin-comp}
\piccy{pain35b.ps}.
\end{equation}
In general, the ways of composing pasting diagrams are indexed by pasting
diagrams themselves. For instance,
\begin{equation}	\label{pic:gen-comps}
\piccy{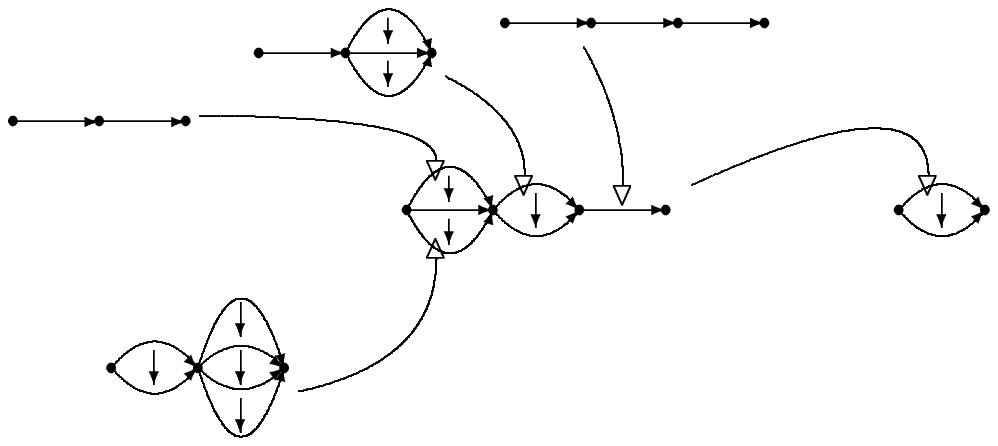}
\end{equation}
represents the composition
\[
\piccy{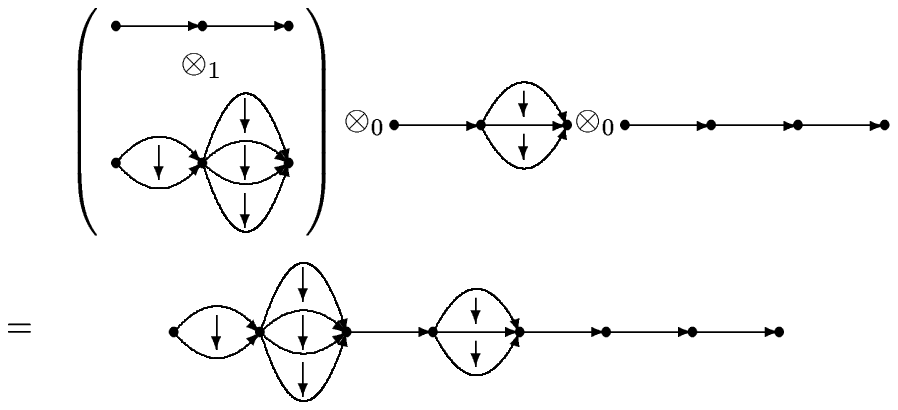}.
\]

We have now described the strict $\omega$-category $\pd=T\One$. More
generally, what does $TX$ look like for an arbitrary globular set $X$? A
globular set is a diagram 
\[
\cdots \pile{\rTo^s \\ \rTo_t} X(n+1) 
\pile{\rTo^s \\ \rTo_t} X(n)
\pile{\rTo^s \\ \rTo_t} \cdots
\pile{\rTo^s \\ \rTo_t} X(0)
\]
of sets, in which $s$ and $t$ obey the `globularity' relations given
in~\ref{sec:formal}; elements of $X(k)$ are called \emph{$k$-cells} of
$X$.  An element of $(TX)(n)$ is an $n$-pasting diagram labelled by elements
of $X$: for example, a typical element of $(TX)(2)$ is a diagram
\[
\piccy{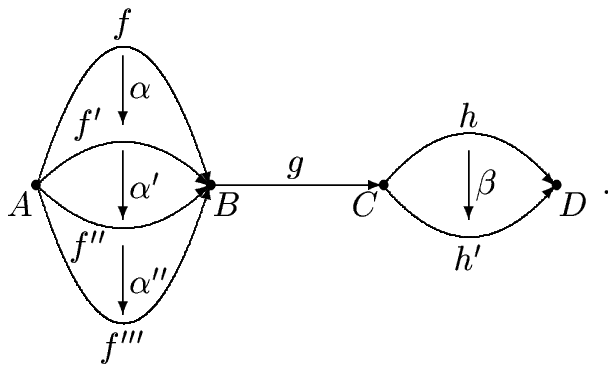}
\]
where $\range{A}{D}\elt X(0)$, $\range{f}{h'}\elt X(1)$,
$\range{\alpha}{\beta}\elt X(2)$, and $s(\alpha)=f$, $t(\alpha)=f'$,
etc. 

To state this more precisely, we first associate to each pasting diagram
$\pi$ a globular set $\rep{\pi}$---the globular set `looking like $\pi$'. For
instance, if $\pi$ is the $2$-pasting diagram~\bref{pic:four-one-three-glob}
then
\[
|\rep{\pi}(k)| = 
\left\{
\begin{array}{ll}
4	&\mbox{if } k=0	\\
7	&\mbox{if } k=1 \\
4	&\mbox{if } k=2 \\
0	&\mbox{if } k\geq 3 
\end{array}
\right.
\]
since (the picture of) $\pi$ has $4$ $0$-cells, $7$ $1$-cells, and so on. We
construct $\rep{\pi}$, for $\pi\in\pd(n)$, recursively on $n$. If
$\pi$ is the unique element of $\pd(0)$ then define
\[
\rep{\pi} = (\cdots \pile{\rTo \\ \rTo} \emptyset \pile{\rTo \\ \rTo} 
		\emptyset \pile{\rTo \\ \rTo} 1).
\]
If $n\geq 0$ and $\pi\in\pd(n+1)$ then $\pi = \bftuple{\pi_1}{\pi_r}$
for some $\range{\pi_1}{\pi_r} \in \pd(n)$, and define
\begin{equation}	\label{eqn:rep-constr}
\rep{\pi} = (\cdots \pile{\rTo \\ \rTo} \coprod_{i=1}^{r} \rep{\pi_i}(1) 
		\pile{\rTo \\ \rTo} \coprod_{i=1}^{r} \rep{\pi_i}(0)
		\pile{\rTo \\ \rTo} \{0, 1, \ldots, r\} ).
\end{equation}
The source and target maps in all but the bottom dimension are the evident
disjoint unions; in the bottom dimension, they are defined by
\[
s(x) = i-1, \diagspace t(x) = i, 
\diagspace
\mbox{for } x\in\rep{\pi_i}(0).
\]

Having defined for each pasting diagram $\pi$ its `representation'
$\rep{\pi}$, we can formalize our guess as to what an element of $(TX)(n)$
is. A `labelling of $\pi$ by elements of $X$' is a map $\rep{\pi} \go X$, so
we are guessing that
\[
\label{p:pasting-rep}%
(TX)(n) \iso \coprod_{\pi\in\pd(n)} 
\homset{\ftrcat{\scat{G}}{\Set}}{\rep{\pi}}{X}.
\]
This is indeed the case, as is shown in section~\ref{sec:pasting-rep}. (Note
that by taking $X=\One$ we recover the fact that $(T\One)(n) \iso \pd(n)$.)

With a little more effort we could define the source and target inclusions
$s, t: \rep{\bdry\pi} \go \rep{\pi}$, to give a concrete description of
the source and target maps in $TX$, and hence of the functor $T$. 
\label{p:effort}%
With an appreciable amount of effort, we could do the same thing for the
monad structure on $T$; but we do not, as the constructions involved for
multiplication are rather complex and not especially illuminating.

There is an alternative way to represent elements of $(T\One)(n)$, used by
Batanin in his paper \cite{Bat}: as trees. (These trees differ slightly from
those which occur elsewhere in this dissertation, and serve a different
conceptual purpose.) For example, we translate the pasting diagram
\[
\gfst{}%
\gfour{}{}{}{}{}{}{}%
\grgt{}%
\gone{}%
\glft{}%
\gtwo{}{}{}%
\glst{}
\]
into the tree 
\[
\treedc.
\]
The thinking here is that the pasting diagram is 3 1-cells
long, so we start the tree as $\treec$; then the first column is 3 2-cells
high, the second 0, and the third 1, so the tree becomes 
\[
\treedc; 
\]
finally, there are no 3-cells so the tree stops there.

Formally, let us define an
\emph{$n$-stage tree} ($n\elt\nat$) to be a diagram
\[
\tau(n)\go\tau(n-1)\go\cdots\go\tau(1)\go\tau(0)=1
\]
in the category $\Delta$ of all finite ordinals, and write $\fcat{Bt}(n)$ for
the set of all $n$-stage trees (with \fcat{Bt} for `Batanin trees'). The
element of $\fcat{Bt}(2)$ just drawn corresponds to a certain diagram $4\go
3\go 1$ in $\Delta$, for example; note that if $\tau$ is an $n$-stage tree
with $\tau(n)=0$ then the height of the picture of $\tau$ will be less than
$n$. The source/target $\bdry\tau$ of an $n$-stage tree $\tau$ is the
$(n-1)$-stage tree obtained by removing all the nodes at height $n$, or
formally, truncating
\[
\tau(n)\go\tau(n-1)\go\cdots\go\tau(1)\go\tau(0)
\]
to
\[
\tau(n-1)\go\cdots\go\tau(1)\go\tau(0).
\]
We thus have a diagram
\begin{equation}	\label{eq: Bt}
\cdots \go \fcat{Bt}(n+1) \goby{\bdry} \fcat{Bt}(n) \go \cdots 
\goby{\bdry} \fcat{Bt}(0)
\end{equation}
in \Set, and so a globular set \fcat{Bt} whose source and target maps are
equal. This is isomorphic to $T\One$, by the following result.

\begin{propn}
The diagram~\bref{eq: Bt} in \Set\ is isomorphic to
\[
\cdots \go \pd(n+1) \goby{\bdry} \pd(n) \go \cdots \goby{\bdry} \pd(0).
\]
\end{propn}
\begin{proof}
$\pd(0)$ and $\fcat{Bt}(0)$ are both 1-element sets, hence isomorphic in a
unique way. Suppose inductively that $n\geq 0$ and that we have constructed a
commuting diagram 
\begin{diagram}
\pd(n)	&\rTo^{\bdry}	&\pd(n-1)	&\rTo^{\bdry}	&\cdots	&\rTo^{\bdry}
&\pd(0)	\\
\dTo<{\alpha}&		&\dTo<{\alpha}	&		&	&
&\dTo<{\alpha}	\\
\fcat{Bt}(n)&\rTo^{\bdry}&\fcat{Bt}(n-1)&\rTo^{\bdry}	&\cdots	&\rTo^{\bdry}
&\fcat{Bt}(0).\\	
\end{diagram}
If $\pi\in\pd(n+1)$ then $\pi = \bftuple{\pi_1}{\pi_r}$ for some $r\in\nat$
and $\pi_i \in\pd(n)$; then define $\alpha(\pi)$ to be
\[
\sum_{i=1}^{r} (\alpha(\pi_i))(n) \go \cdots \go
\sum_{i=1}^{r} (\alpha(\pi_i))(0) \go 1.
\]
It is easy to check that the map $\alpha: \pd(n+1) \go \fcat{Bt}(n+1)$ thus
defined is a bijection and commutes with the $\bdry$'s. 
\done
\end{proof}

Composition and identities in the strict $\omega$-category $\fcat{Bt}$
($\iso T\One$) can also be expressed in the pictorial language of trees,
in a simple and compelling way; for that the reader is referred to \cite{Bat}
or \cite[Ch.\ II]{SHDCT}.

\section{Globular Operads and Algebras}		\label{sec:glob-ops}

Let $T$ be the free strict $\omega$-category monad on the category
\ftrcat{\scat{G}}{\Set} of globular sets. This section is an attempt at an
elementary explanation of $T$-operads and their algebras.

A \emph{collection} is a $T$-graph on $1$: that is, it is a globular set $C$
together with a map $C \go \pd$. Put another way, a collection consists of a
set $C(\pi)$ for each $n$-pasting diagram $\pi$, together with a pair of
functions 
$
\begin{diagram}
C(\pi)	&\pile{\rTo^{s} \\ \rTo_{t}}	&C(\bdry\pi)	\\
\end{diagram}
$
(when $n\geq 1$), satisfying the usual globularity equations $ss=st$ and
$ts=tt$.

A $T$-operad is a collection $C \goby{d} \pd$ together with identities and
compositions satisfying suitable axioms. The elements of $C(\pi)$ are to be
thought of as the operations of shape or arity $\pi$: in other words, as the
functions
\begin{equation}	\label{eq:alg-struc}
\homset{\ftrcat{\scat{G}}{\Set}}{\rep{\pi}}{X} \go X(n)
\end{equation}
which exist as part of the structure of a $C$-algebra $X$. (Recall
that \homset{\ftrcat{\scat{G}}{\Set}}{\rep{\pi}}{X} is the set of
`labellings of $\pi$ by elements of $X$'.)

The identities consist of an element of $C(\iota_n)$ for each $n$, where
$\iota_{n} \elt \pd(n)$ is the $n$-pasting diagram looking like a single
$n$-cell: formally, $\iota_0$ is the unique element of $\pd(0)$ and
\[
\iota_{n+1} = (\iota_n) \in (\pd(n))^* = \pd(n+1).
\]

Composition is a map $C\of C \go C$ over \pd, where the collection
$C\of C \goby{\twid{d}} \pd$ is the left-hand diagonal of the diagram
\[
\begin{slopeydiag}
   &       &   &       &   &       &C\of C\Spbk&  &   \\
   &       &   &       &   &\ldTo  &      &\rdTo  &   \\
   &       &   &       &T(C)&      &      &       &C  \\
   &       &   &\ldTo<{T(d)}&&\rdTo<{T(!)}&&\ldTo>{d}&\\
   &       &T(\pd)&    &   &       &\pd   &       &   \\
   &\ldTo<{\mu_\One}&& &   &       &      &       &   \\
\pd&       &   &       &   &       &      &       &   \\
\end{slopeydiag}
\ \ \ \ .
\]
A typical element of $(C\of C)(2)$ is depicted in the following diagram:
\begin{equation}	\label{pic:compn-in-operad}
\piccy{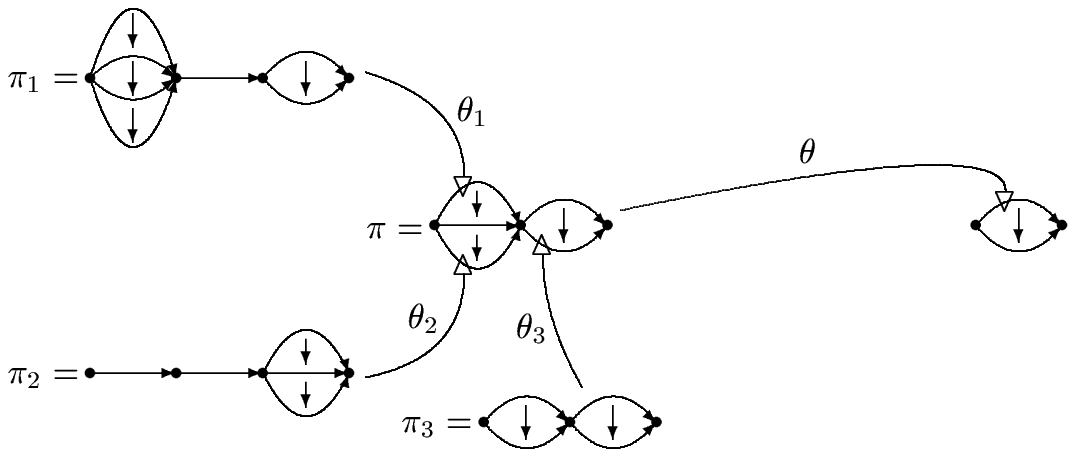}.
\end{equation}
This diagram is meant to indicate that $\theta_1 \in C(\pi_1)$, $\theta_2 \in
C(\pi_2)$, $\theta_3 \in C(\pi_3)$, $\theta \in C(\pi)$, and that $\theta_1$,
$\theta_2$, $\theta_3$ match suitably on their sources and targets (e.g.\
$t(\theta_1) = s(\theta_2)$). The left-hand half of the diagram is an element
of the fibre over $\pi$ in the map $T(C)\goby{T(!)}\pd$; the right-hand half
is an element of $C(\pi)$ (which is the fibre over $\pi$ in the map $C
\goby{d} \pd$); hence the whole diagram is an element of $(C\of C)(2)$. The
map $C\of C \goby{\twid{d}} \pd$ sends this element to the 2-pasting diagram
\[
\pi\of(\pi_{1},\pi_{2},\pi_{3}) = 
\gzersu\gfoursu\gzersu\gonesu\gzersu%
\gfoursu\gzersu\gtwosu\gzersu\gtwosu\gzersu
\]
(which is the composite of $\pi$ with $\pi_1$, $\pi_2$, $\pi_3$ in the
$\omega$-category \pd; cf.\ diagram~\bref{pic:gen-comps}).  So, composition
sends the data assembled in~\bref{pic:compn-in-operad} to an element of
$C(\pi\of(\pi_{1},\pi_{2},\pi_{3}))$, which may be drawn as
\[
\piccy{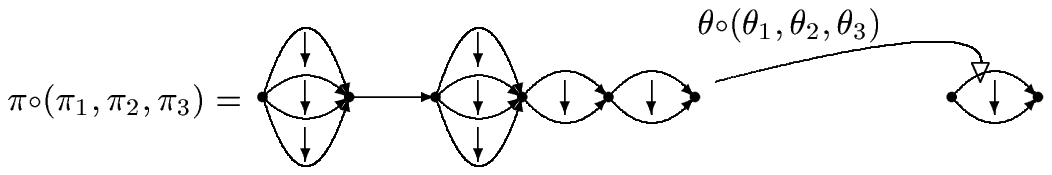}.
\]

(The `linear' notation $\pi\of(\pi_{1},\pi_{2},\pi_{3})$ and
$\theta\of(\theta_{1},\theta_{2},\theta_{3})$ should not be taken too
seriously. There is evidently no natural order in which to put the $\pi_i$'s
and $\theta_i$'s; the notation is just being used temporarily for
convenience.) 

The composition and identities in a $T$-operad $C$ are required to commute
with the source and target maps and, of course, to obey associativity and
identity laws. For example, if we have a diagram
\begin{diagram}[height=1.5em]
\cdot	&	&	&	&	&	&	\\
	&\rdTo>{\theta_{11}}&&	&	&	&	\\
\cdot	&\rTo^{\theta_{12}}&\cdot&&	&	&	\\
	&\ruTo>{\theta_{13}}&&\rdTo(2,3)>{\theta_1}&&&	\\
\cdot	&	&	&	&	&	&	\\
	&	&	&	&\cdot	&\rTo^{\theta}&\cdot\\
	&	&	&\ruTo>{\theta_2}&&	&	\\
	&	&\cdot	&	&	&	&	\\
\end{diagram}
of the same kind as~(\ref{pic:compn-in-operad}), then
\[
\theta\of(\theta_{1}\of(\theta_{11},\theta_{12},\theta_{13}),\theta_{2}) =
(\theta\of(\theta_{1},\theta_{2})) \of (\theta_{11},\theta_{12},\theta_{13},1).
\]

We have now seen that an operad consists of a set $C(\pi)$ for each pasting
diagram $\pi$, with source and target functions, and compositions between the
$C(\pi)$'s according to the pasting-together of pasting diagrams. We have
already argued (equation~\bref{eq:alg-struc}) that an algebra for $C$ `ought'
to consist of a globular set $X$ together with a function
\[
\ovln{\theta}:
\homset{\ftrcat{\scat{G}}{\Set}}{\rep{\pi}}{X} \go X(n)
\]
for each $\theta\in C(\pi)$ ($\pi\in\pd(n)$, $n\in\nat$), obeying suitable
axioms. So for instance, suppose that
\[
\pi = 
\gfst{}\gfour{}{}{}{}{}{}{}\grgt{}\gone{}\glst{},
\]
that $\theta\in C(\pi)$, and that
\[
\gfst{A}%
\gfour{p}{p'}{p''}{p'''}{\alpha}{\alpha'}{\alpha''}%
\grgt{B}%
\gone{q}%
\glst{C}
\]
is a diagram of cells in $X$; then $\ovln{\theta}$ assigns to this picture a
2-cell in $X$.

This is indeed what the general theory says a $C$-algebra is. For an algebra
structure on $X$ is a map $T_C X \goby{h} X$ obeying suitable laws, where
$T_C (X)$ is the pullback
\[
\begin{slopeydiag}
	&		&T_C X\Spbk	&		&	\\
	&\ldTo		&		&\rdTo		&	\\
TX	&		&		&		&C	\\
	&\rdTo<{T(!)}	&		&\ldTo>{d}	&	\\
	&		&\pd		&		&	\\
\end{slopeydiag}
\ \ \ \ ,
\]
and this means that 
\[
(T_C X)(n)=\coprod_{\pi\elt\pd(n)}
C(\pi) \times \homset{\ftrcat{\scat{G}}{\Set}}{\rep{\pi}}{X}.
\]
Hence $h$ consists of a function
\[
C(\pi) \times \homset{\ftrcat{\scat{G}}{\Set}}{\rep{\pi}}{X}
\go
X(n)
\]
for each number $n$ and $n$-pasting diagram $\pi$. Writing $h(\theta,
\dashbk)$ as $\ovln{\theta}$, we see that this is just the description
above. 

We have now discussed what operads and their algebras look like, and it is
time to come to the main point of the chapter.%

\section{Contractions}		\label{sec:contractions}

We start with an informal description of what a weak $\omega$-category
`should' be, centred around the idea of contraction, and then see how this is
expressed by the formal definition of contraction.

\begin{sloppypar}
The graph structure of an $\omega$-category consists of 0-cells \gzero{A},
1-cells \gfst{A}\gone{f}\glst{B}, 2-cells
\gfst{A}\gtwo{f}{g}{\alpha}\glst{B}, \ldots. There are then various ways of
composing these cells; just how many ways and how they interact depends on
whether we are dealing with strict or weak $\omega$-categories, or something
in between. In a strict $\omega$-category, there will be precisely one way of
composing a diagram like
\begin{equation}	\label{pic:labelled-2-pd}
\gfst{A}%
\gthree{f}{f'}{f''}{\alpha}{\alpha'}%
\grgt{B}%
\gone{g}%
\glft{C}%
\gtwo{h}{h'}{\gamma}%
\glst{D}
\end{equation}
to obtain a 2-cell: that is, any two different ways of doing it (such as
`compose $\alpha'$ with $\alpha$, and $\gamma$ with $g$, then the two of
these together') give exactly the same resulting 2-cell. In a weak
$\omega$-category there will be many ways, but the resulting 2-cells will all
be equivalent in a suitably weak sense.
\end{sloppypar}

Our method of describing what ways of composing are available in a weak
$\omega$-category depends on one simple principle, the \emph{contraction
principle}.  Take, for example, the diagram~(\ref{pic:labelled-2-pd}) above.
%\label{p:star}
Suppose we have already constructed two ways of composing a generic diagram
\[
\gfst{}%
\gone{p}%
\gblw{}%
\gone{q}%
\gblw{}%
\gone{r}%
\glst{}
\]
of 1-cells, namely $(rq)p$ and $r(qp)$. Then the contraction principle says
that there is a way of composing diagram~(\ref{pic:labelled-2-pd}) to get a
2-cell of the form \gfst{A}\gtwo{(hg)f}{h'(gf'')}{}\glst{D}. As another
example of the principle, this time in one higher dimension, take a diagram
\begin{equation}	\label{pic:labelled-3-pd}
\gfst{A}%
\gspecialone{f}{f'}{f''}{\alpha}{\alpha'}{\alpha''}{\beta}{x}{y}%
\grgt{B}%
\gonew{g}%
\glft{C}%
\gthreecell{h}{h'}{\gamma}{\gamma'}{z}%
\glst{D}.
\end{equation}
Suppose we have constructed two ways of composing a generic diagram of the
shape of~(\ref{pic:labelled-2-pd}) to a 2-cell, each of which invokes the
same way of composing the 1-cells
\[
\gfst{}%
\gone{}%
\gblw{}%
\gone{}%
\gblw{}%
\gone{}%
\glst{}
\]
along the top and bottom. Say, for instance, that the first way of
composing~(\ref{pic:labelled-2-pd}) results in a 2-cell
\gfst{A}\gtwo{(hg)f}{h'(gf'')}{\delta}\glst{D} and the second way in a 2-cell
\gfst{A}\gtwo{(hg)f}{h'(gf'')}{\delta'}\glst{D}. Then the contraction
principle says that there is a way of
composing~(\ref{pic:labelled-3-pd}) to get a 3-cell of the form
\[
\gfst{A}\gthreecell{(hg)f}{h'(gf'')}{\delta}{\delta'}{}\glst{D}.
\]

In general, the contraction principle can be stated as follows. Suppose we
are given an $n$-dimensional diagram and two ways of composing the
$(n-1)$-dimensional diagram at its source/target, such that these two ways
match on the $(n-2)$-dimensional source and target. Then there's a way of
composing the $n$-dimensional diagram, inducing the first way on its source
and the second way on its target. (In our first example, we implicitly used
the fact that the two ways of composing
\[
\gfst{}%
\gone{p}%
\gblw{}%
\gone{q}%
\gblw{}%
\gone{r}%
\glst{},
\]
$(rq)p$ and $r(qp)$, do the same thing to the bounding 0-cells: nothing at
all.) The ways of composing in a weak $\omega$-category are to be generated
by this principle, and this principle alone.

How does this idea of contraction compare to the definition given
in~\ref{sec:formal}? The structure encoding `ways of composing' is, of
course, a $T$-operad $C$. For $\pi\in\pd(n)$ ($n\geq 2$), we defined
\[
P_{\pi}(C) = 
\{ \pr{\theta_0}{\theta_1} \in C(\bdry\pi)^2 \such
s(\theta_0) = s(\theta_1) \mbox{ and } t(\theta_0) = t(\theta_1)\},
\]
and for $\pi\in\pd(1)$,
\[
P_{\pi}(C) = C(\bdry\pi)^2.
\]
Thus an element of $P_\pi(C)$ can be thought of as a way $\theta_0$ of
composing the $(n-1)$-dimensional source of an $n$-dimensional diagram of
shape $\pi$, together with a way $\theta_1$ of composing its target, such
that these two ways match on the $(n-2)$-dimensional part. A contraction
$\kappa$ on $C$ was defined as a function
\[
\kappa_\pi: P_\pi (C) \go C(\pi)
\]
for each $\pi$, such that
\[
s(\kappa_\pi \pr{\theta_0}{\theta_1} ) = \theta_0,
\diagspace
t(\kappa_\pi \pr{\theta_0}{\theta_1} ) = \theta_1.
\]
In other words, it extends $\theta_0$ and $\theta_1$ to a way $\kappa_\pi
\pr{\theta_0}{\theta_1}$ of composing a whole $\pi$-shaped diagram. This is
exactly the effect of the informal contraction principle.

Notice, incidentally, that if $\kappa$ is a contraction on a $T$-operad $C$
then the functions $\kappa_\pi$ are \emph{not} required to be compatible with
the operad structure on $C$ in any way. So the natural entity on which to
\label{p:entity}%
define a contraction is not in fact a $T$-operad but a collection (i.e.\ a $T$-graph on $1$).

An important feature of the contraction idea is what happens with degenerate
pasting diagrams. There is not only a 2-pasting diagram $\sigma$ shaped like
diagram~\bref{pic:labelled-2-pd}, but also a (degenerate) 3-pasting diagram
$\pi$ shaped like it too: thus $\bdry\pi = \sigma$. Now, suppose that
$\theta_0, \theta_1 \in C(\sigma)$ with $s(\theta_0) = s(\theta_1)$ and
$t(\theta_0) = t(\theta_1)$. Then there is an element $\theta = \kappa_\pi
\pr{\theta_0}{\theta_1}$ of $C(\pi)$ with $s(\theta) = \theta_0$ and
$t(\theta) = \theta_1$. This means that $\theta$ assigns to the data
in~\bref{pic:labelled-2-pd} a 3-cell
\[
\gfst{A}\gthreecell{}{}{\delta_0}{\delta_1}{}\glst{D},
\]
in which $\delta_0$ and $\delta_1$ are respectively the results of applying
$\theta_0$ and $\theta_1$ to~\bref{pic:labelled-2-pd}. This is the kind of
argument we would use to prove that any two composites of a given diagram
are, in a suitable sense, equivalent.

\section{The Definition}	\label{sec:the-defn}

A weak $\omega$-category is defined to be an $L$-algebra, where
\pr{L}{\lambda} is the initial operad-with-contraction. We have seen what an
operad-with-contraction is, and what an algebra for one is; now we have just
to see why the \emph{initial} one gives us what we want.

Another way of saying that \pr{L}{\lambda} is initial in \fcat{OWC} is that
\pr{L}{\lambda} is the operad-with-contraction freely generated by the empty
collection $(\emptyset \go \pd)$. That is, we start with the empty
collection and freely add in just enough to make it into a $T$-operad $L$ with
a contraction $\lambda$ on it. 

So, for a start there is an identity element $1\in L(\blob)$, where
$\blob\in\pd(0)$. Next, take the 1-pasting diagram
\[
\pi_n = 
(\gfst{}%
\gone{}%
\gblw{}%
\gone{}%
\diagspace%
\cdots%
\diagspace%
\gone{}%
\glst{})
\]
of length $n$. The contraction gives us an element $\psi_n =
\lambda_{\pi_n}\pr{1}{1}$ of $L(\pi_n)$. Thus in an $L$-algebra, $\psi_n$
provides a way of composing a diagram
\[
\gfst{A_0}%
\gone{f_1}%
\gblw{A_1}%
\gone{f_2}%
\diagspace%
\cdots%
\diagspace%
\gone{f_n}%
\glst{A_n}
\]
to give a 1-cell \gfst{A_0}\gone{\psi_n\bftuple{f_1}{f_n}}\glst{A_n}; let us
write
\[
(f_n \of\cdots\of f_1) = \psi_n\bftuple{f_1}{f_n} 
\]
($n\geq 1$), and $1 = \psi_0 ()$. Next, the operad structure on $L$ gives us
1-dimensional elements of $L$ such as
\[
\psi_3 \of (\psi_3, \psi_0, \psi_2) \in L(\pi_5),
\]
which is interpreted in an $L$-algebra as the function
\begin{eqnarray*}
&
\gfst{}%
\gone{f_1}%
\gblw{}%
\gone{f_2}%
\gblw{}%
\gone{f_3}%
\gblw{}%
\gone{f_4}%
\gblw{}%
\gone{f_5}%
\glst{}\\
\goesto	&
\gfst{}%
\gone{((f_{5}\of f_{4})\of 1\of (f_{3}\of f_{2}\of f_{1}))}%
\glst{}.
\end{eqnarray*}
This analysis might lead us to suspect that $L(\pi_n)$ is the set $\tr(n)$ of
$n$-leafed trees (described in~\ref{sec:free-multicats}
and~\ref{sec:coherence}), which in fact it is.

Moving now to the 2-dimensional level, if $\pi$ is the 2-pasting diagram
shaped like diagram~\bref{pic:labelled-2-pd} then any pair
\pr{\theta_0}{\theta_1} of elements of $L(\pi_3)$ gives rise to an element
$\theta = \lambda_{\pi}\pr{\theta_0}{\theta_1}$ of $L(\pi)$. (Since
$L(\blob)$ has only one element, there is no need to worry about $\theta_0$
and $\theta_1$ matching at the 0-dimensional level.) Generally, the
2-dimensional part of $L$ contains elements obtained by contraction from the
1-dimensional parts, together with all the elements obtained by pasting them
together (using the operad structure of $L$). To take a reasonably manageable
example, let
\[
\begin{array}{cc}
\pi = \gfst{}\gthree{}{}{}{}{}\glst{},
\diagspace
\pi' = \gfst{}\gtwo{}{}{}\gfbw{}\gtwo{}{}{}\glst{},		\\
\pi'' = \ \gfst{}\gthree{}{}{}{}{}\gfbw{}\gthree{}{}{}{}{}\glst{}.
\end{array}
\]
Then:
\begin{itemize}
\item $L(\pi)$ has an element $\psi = \lambda_{\pi}\pr{1}{1}$ (where $1\in
L(\pi_1)$) 
\item $L(\pi')$ has an element $\psi' = \lambda_{\pi'}\pr{\psi_2}{\psi_2}$
\item $L(\pi'')$ has an element $\lambda_{\pi''}\pr{\psi_2}{\psi_2}$ (`compose
all four cells at once') 
\item $L(\pi'')$ also has an element which might reasonably be
denoted $\psi \of \pr{\psi'}{\psi'}$ (`first compose horizontally,
then compose vertically')
\item $L(\pi'')$ has a third element $\psi' \of \pr{\psi}{\psi}$ (`first
compose vertically, then compose horizontally').
\end{itemize}
These elements $\psi \of \pr{\psi'}{\psi'}$ and $\psi' \of
\pr{\psi}{\psi}$ of $L(\pi)$ are familiar from the interchange law in the
definition of 2-category. Of course, the three elements of $L(\pi)$ we have
mentioned are not its \emph{only} elements; there are infinitely many, since
$\tr(2)$ is an infinite set. 

This concludes our explanation of why weak $\omega$-categories can reasonably
be defined as objects of $\Alg(L)$. 

Notice, however, that weak $\omega$-functors are \emph{not} defined as
morphisms in $\Alg(L)$.  On the contrary, a morphism in $\Alg(L)$ preserves
the $L$-algebra structure strictly, so should be thought of as a
\label{p:strict-omega-ftrs}%
strict map of weak $\omega$-categories.

Here is a sketch of how `weak $\omega$-functor' might be defined. This is
only speculation, and no proper definition is attempted here. As in the
definition of weak $\omega$-category, the idea is to take a theory of strict
things and a notion of contraction to create a theory of weak things.

So, there is a $T$-multicategory \fcat{Map} such that a \fcat{Map}-algebra is
a pair \pr{X}{Y} of strict $\omega$-categories together with a strict
$\omega$-functor $X\go Y$. (The objects-object $\fcat{Map}_0$ of \fcat{Map}
is the coproduct $\One + \One$ of two copies of the terminal globular set.)
There is also a notion of what a contraction on a map of $T$-multicategories
is. Hence there is a category of $T$-multicategories with contraction over
\fcat{Map}, in which an object consists of a $T$-multicategory $D$, a map $d:
D \go \fcat{Map}$ of $T$-multicategories, and a contraction $\delta$ on
$d$. This category has an initial object \pr{M \goby{m} \fcat{Map}}{\nu}, and
a weak $\omega$-functor is defined as an $M$-algebra.

The notion of a contraction on a map of $T$-multicategories has the property
that for $T$-operads $C$, a contraction on the unique map from $C$ to the
terminal $T$-operad is precisely a contraction on $C$ in the sense of the
rest of this chapter. This means that the two inclusions $1 \pile{\rTo\\
\rTo} \fcat{Map}$ induce another pair of maps $L \pile{\rTo\\ \rTo} M$, and
hence a pair of functors $\Alg(M) \pile{\rTo\\ \rTo}\Alg(L)$. 
These are the functors assigning to a weak $\omega$-functor its domain and
codomain.

Batanin's paper~\cite{Bat} contains a definition (\S 8) of weak
$\omega$-functor, which unfortunately I have not been able to
understand. However, I think I can explain how Batanin's definition of weak
$\omega$-category differs from the present one, as follows.

Let $C$ be a $T$-operad. Firstly, a \emph{system of compositions} on $C$
consists of a chosen element $\theta_\pi$ of $C(\pi)$ for each pasting
diagram $\pi$ that represents a binary composition: for instance, $\pi$
might be one of
\[
\gzersu\gonesu\gzersu\gonesu\gzersu\ ,
\diagspace
\gzersu\gtwosu\gzersu\gtwosu\gzersu\ ,
\diagspace
\gzersu\gthreesu\gzersu
\ .
\]
These chosen elements are required to be consistent with one another: e.g.\
if $\pi_1$ and $\pi_2$ are the first and second of these three diagrams, then
\[
s(\theta_{\pi_2}) = \theta_{\pi_1} = t(\theta_{\pi_2}).
\]
Secondly, a \emph{contraction} $\kappa$ on $C$ is a family $(\kappa_\pi)$ of
functions of a certain kind, exactly as in our definition, except that now
$\pi$ only ranges over those $n$-pasting diagrams satisfying
$\rep{\pi}(n) = \emptyset$. The latter condition means that $\pi$ is
`degenerate', as discussed earlier in the chapter.

Now consider the full subcategory \cat{Q} of $T\hyph\fcat{Operad}$ whose
objects are those $T$-operads on which there exists both a system of
compositions and a contraction. Batanin constructs a certain weakly initial
object $K$ of \cat{Q}, and defines a weak $\omega$-category to be a
$K$-algebra.

`Weakly initial' means that there is at least one map from $K$ to any other
object of \cat{Q}. So $K$ is not determined by its weak initiality,
and this means that if we want to know what a Batanin weak
$\omega$-category is then we actually need the details of the construction of
$K$ in~\cite{Bat}. It might be the case that if we take the category \cat{Q'}
of $T$-operads \emph{equipped} with a system of compositions and a
contraction, then $K$ (together with its system of compositions and
contraction) is initial in \cat{Q'}, and of course this \emph{would}
determine $K$. A remark in~\cite{Bat} (just before Definition 8.6) suggests
that this is true.

The idea behind the Batanin definition appears to be that the theory of weak
$\omega$-categories---that is, the operad $K$ for which they are
algebras---is generated by two things: operations and equations. The
operations are binary compositions of various dimensions, and these are
provided by the system of compositions. The `equations' should really be
called `equivalences', and are provided by the contraction: compare the use
of degenerate pasting diagrams at the end of~\ref{sec:contractions} above.
In our approach these two ingredients are merged into one: the more
comprehensive notion of contraction.

I do not know if the present definition of weak $\omega$-category is in any
sense equivalent to Batanin's. I would certainly imagine so, but there is
little chance of providing a comparison before weak $\omega$-functors are
understood.

\section{Examples}	\label{sec:examples}

At this point it would be nice to give a fully worked-out non-trivial example
of a weak $\omega$-category. Unfortunately I do not yet have one for which
all the details have been settled. However, the following remarks may provide
partial satisfaction.

Recall from~\ref{sec:alg} that a map between $T$-operads induces a map in
the opposite direction between their categories of algebras, and that an
algebra for the terminal $T$-operad is just a $T$-algebra. Hence the unique
$T$-operad map $L\go 1$ induces a functor
\[
(\mbox{strict } \omega\mbox{-categories}) = \Alg(1)
\go
\Alg(L) = (\mbox{weak } \omega\mbox{-categories}).
\]
That is, `every strict $\omega$-category is a weak
$\omega$-category'. Incidentally, the terminal $T$-operad $1$ carries a
unique contraction, and is then the terminal operad-with-contraction: so
algebras for the terminal operad-with-contraction are strict
$\omega$-categories, and algebras for the initial operad-with-contraction are
weak $\omega$-categories.

More generally, for any operad-with-contraction \pr{C}{\kappa} there
is a unique contraction-preserving operad map $L\go C$, and this induces a
functor 
\[
\Alg(C) \go \Alg(L).
\]
This provides a means of finding examples of
weak $\omega$-categories. For instance, suppose we wanted to define a weak
$\omega$-category $\Pi_\omega (S)$ for every topological space $S$, its
`fundamental $\omega$-groupoid'. It is clear what the globular set
$\Pi_\omega (S)$ should be, and our strategy might then be to find a
$T$-operad $C$ such that
\begin{itemize}
\item $\Pi_\omega (S)$ is naturally a $C$-algebra for every space $S$, and
\item there is a contraction on $C$.
\end{itemize}
Any way of doing this will give the globular set $\Pi_\omega (S)$ the
structure of a weak $\omega$-category. (The rough idea is that $C(\pi)$ is
the set of continuous maps from the closed $n$-ball to the contractible
space which looks like the usual picture of $\pi$ (for
$\pi\in\pd(n)$), subject to conditions on boundary-preservation. Something
like this is done in \cite[\S 9]{Bat}.)

In the next section, weak $n$-categories will be defined as weak
$\omega$-cat\-e\-go\-ries of a special kind. We will subsequently show that
weak $2$-categories are essentially the same as bicategories. Thus any
bicategory provides a (degenerate) example of a weak $\omega$-category.

\section{Weak $n$-Categories}	\label{sec:weak-n}

Our definition of weak $\omega$-category suggests not just one, but two
plausible definitions of weak $n$-category. In this section we present both
of these definitions and show that the two different categories of weak
$n$-categories (with strict $n$-functors as morphisms) are equivalent. 

Let us say that a globular set $X$ is \emph{$n$-dimensional} (for $n\in\nat$)
if for all $m\geq n$,
\[
s=t: X(m+1) \go X(m)
\]
and this map is an isomorphism. 
\begin{defn}
A \emph{weak $n$-category} is a weak $\omega$-category whose underlying
globular set is $n$-dimensional.
\end{defn}
This formalizes the idea that an $n$-category is an $\omega$-category in
which the only cells of dimension greater than $n$ are identities. Let us
write $\fcat{Wk}\hyph n\hyph\Cat$ for the full subcategory of $\Alg(L)$ whose
objects are weak $n$-categories. 

The alternative approach does not use the definition of weak
$\omega$-category directly, but instead imitates it. Write $\cat{G} =
\ftrcat{\scat{G}}{\Set}$ for the category of globular sets. Let $\scat{G}_n$
be the full subcategory of $\scat{G}$ with objects $0, \ldots, n$, let
$\cat{G}_n = \ftrcat{\scat{G}_n}{\Set}$, call objects of $\cat{G}_n$
\emph{$n$-globular sets}, and let $T_n$ be the free strict $n$-category monad
on $\cat{G}_n$. Theorem~\ref{thm:strict-n} tells us that $T_n$ is a cartesian
monad on $\cat{G}_n$, so we can discuss $T_n$-operads.

Let $C$ be a $T_n$-operad. If $1\leq k\leq n$ and $\pi\in\pd(k)$, we may
define the set $P_\pi(C)$ just as in section~\ref{sec:formal}. A
\emph{precontraction} on $C$ is a family of functions
\[
(\kappa_\pi: P_{\pi}(C) \go C(\pi))_{1\leq k\leq n, \pi\in \pd(k)}
\]
satisfying the same equations as in~\ref{sec:formal}.  If $C$ has the
property that for all $\pi\in\pd(n)$ and $\theta_0, \theta_1 \in \pd(\pi)$,
\[
s(\theta_0) = s(\theta_1) 
\mbox{ and }
t(\theta_0) = t(\theta_1) 
\diagspace\Longrightarrow\diagspace
\theta_0 = \theta_1
\]
then any precontraction on $C$ is called a \emph{contraction}. (There is then
no choice about what the contraction does in the top dimension.)  We
therefore obtain a category $\fcat{OWC}_n$, in which an object is a
$T_n$-operad equipped with a contraction and a map is a map of operads
preserving contractions, defined analogously to $\fcat{OWC}$
in~\ref{sec:formal}. 

Later we will show that $\fcat{OWC}_n$ has an initial object
\pr{L_n}{\lambda_n}. The alternative definition of weak $n$-category is as an
$L_n$-algebra. As in the case of $\omega$-categories, the morphisms in
$\Alg(L_n)$ should be interpreted as \emph{strict} maps. 

The aim of the rest of this section is to show that these two definitions are
equivalent, in the following strong sense. (`Strong', because we do not have
to resort to weak $n$-functors in order to be able to compare the objects of
the two categories.)
\begin{thm}	\label{thm:n-eqv}
There is an equivalence of categories
\[
\fcat{Wk}\hyph n\hyph\Cat  \eqv  \Alg(L_n).
\]
\end{thm}

The proof is in two parts: first we express the initial object
\pr{L_n}{\lambda_n} of $\fcat{OWC}_n$ in terms of \pr{L}{\lambda}, and then
we are in a position to compare algebras for $L_n$ and for $L$.

So, the inclusions $\scat{G}_{n-1} \rIncl \scat{G}_n$ and $\scat{G}_k \rIncl
\scat{G}$ induce `restriction' functors
\[
R^n_{n-1}: \cat{G}_n \go \cat{G}_{n-1}, 
\diagspace
R^{\omega}_{k}: \cat{G} \go \cat{G}_k,
\]
for any $n\geq 1$ and $k\geq 0$. We then have:
\begin{propn}	\label{propn:three-adjns}
\begin{enumerate}
\item 		\label{propn:n-(n-1)}
For any $n\geq 1$, the functor $R^n_{n-1}$ has a right adjoint
$S^n_{n-1}$, and there is an induced adjunction
\begin{diagram}[height=2em]
\fcat{OWP}_n				\\
\dTo<{R^n_{n-1}} \ladj \uTo>{S^n_{n-1}}	\\
\fcat{OWP}_{n-1}
\end{diagram}
(abusing notation by reusing the symbols $R^n_{n-1}$ and $S^n_{n-1}$)
\item  		\label{propn:n-(n-1)-eqv}
This adjunction restricts to an equivalence of categories
\begin{diagram}[height=2em]
\fcat{OWC}_n				\\
\dTo<{R^n_{n-1}} \eqv \uTo>{S^n_{n-1}}	\\
\fcat{OWP}_{n-1}
\end{diagram}
\item  		\label{propn:omega-k}
For any $k\geq 0$, the functor $R^{\omega}_k$ has a right adjoint
$S^{\omega}_k$, and there is an induced adjunction
\begin{diagram}[height=2em]
\fcat{OWC}				\\
\dTo<{R^{\omega}_k} \ladj \uTo>{S^{\omega}_k}	\\
\fcat{OWP}_k.
\end{diagram}
\end{enumerate}
\end{propn}
\begin{proof}
\begin{enumerate}
\item That $R^n_{n-1}: \cat{G}_n \go \cat{G}_{n-1}$ has a right adjoint is
immediate: it is the right Kan extension of the inclusion $\scat{G}_{n-1}
\rIncl \scat{G}_n$. However, it will be useful to have the following explicit
description of $S^n_{n-1}$: if $X\in\cat{G}_{n-1}$ then 
\[
\!\!\!\!
\begin{array}{l}
(S^n_{n-1} X)(k) = X(k) \diagspace\mbox{for }0\leq k\leq n-1,\\
(S^n_{n-1} X)(n) =
\{ \pr{x_0}{x_1}\in (X(n-1))^2 \such s(x_0)=s(x_1), t(x_0)=t(x_1) \}.
\end{array}
\]
(When $n=1$ the second line does not make sense, and we instead define
$(S^1_0 X)(1)$ as $X(0)\times X(0)$; essentially we are `taking $X(-1)=1$'.)
The source and target maps are the obvious ones.

\begin{sloppypar}
As is shown in Appendix~\ref{app:strict-omega}, $R^n_{n-1}$ is naturally a
monad opfunctor $\pr{\cat{G}_n}{T_n} \go \pr{\cat{G}_{n-1}}{T_{n-1}}$, whose
natural transformation part
\[
R^n_{n-1} \of T_n  \go  T_{n-1} \of R^n_{n-1}
\]
is an isomorphism. Under the adjunction $R^n_{n-1} \ladj S^n_{n-1}$, the mate
of this isomorphism is a natural transformation
\[
T_n \of S^n_{n-1}  \go  S^n_{n-1} \of T_{n-1},
\]
and this gives $S^n_{n-1}$ the structure of a monad functor
$\pr{\cat{G}_{n-1}}{T_{n-1}} \go $\linebreak$\pr{\cat{G}_n}{T_n}$. Further
checks reveal that the conditions of section~\ref{sec:change} are satisfied,
so that there is an induced adjunction between categories of multicategories;
moreover, $R^n_{n-1}$ and $S^n_{n-1}$ each preserve terminal objects, so this
restricts to an adjunction
\begin{equation}			\label{eq:n-oper-adjn}
\begin{diagram}[height=2em]
T_n\hyph\fcat{Operad}			\\
\dTo<{R^n_{n-1}} \ladj \uTo>{S^n_{n-1}}	\\
T_{n-1}\hyph\fcat{Operad}.		\\
\end{diagram}
\end{equation}
$R^n_{n-1}$ has the obvious restriction effect on $T_n$-operads; in the other
direction, if $D$ is a $T_{n-1}$-operad, $0\leq k\leq n$ and $\pi\in\pd(k)$,
then
\[
(S^n_{n-1} D)(\pi) = 
\left\{
\begin{array}{ll}
D(\pi)		&\mbox{for }0\leq k\leq n-1	\\
P_\pi (D)	&\mbox{for }k=n.		
\end{array}
\right.
\]
\end{sloppypar}

Next we bring in precontractions. Any precontraction on a $T_n$-operad $C$
evidently gives rise to a precontraction on $R^n_{n-1} C$; conversely, any
precontraction on a $T_{n-1}$-operad $D$ extends uniquely to a precontraction
on $S^n_{n-1} D$. The precontractions produced by these two constructions are
preserved by the unit and counit maps of the
adjunction~\bref{eq:n-oper-adjn}, so we obtain an adjunction
\begin{diagram}[height=2em]
\fcat{OWP}_n				\\
\dTo<{R^n_{n-1}} \ladj \uTo>{S^n_{n-1}}	\\
\fcat{OWP}_{n-1}
\end{diagram}
as required.

\item Any adjunction $F\ladj G:\cat{D} \go \cat{C}$ restricts to an
equivalence between $\cat{C'}$ and $\cat{D'}$, where $\cat{C'}$ is the full
subcategory of \cat{C} whose objects are those at which the unit of the
adjunction is an isomorphism, and similarly \cat{D'} with the counit. In the
present case we have $R^n_{n-1} \of S^n_{n-1} = 1$, and the counit of the
adjunction is the identity transformation. On the other hand, let
\pr{C}{\kappa} be a $T_n$-operad with precontraction and consider the unit
map
\[
\pr{C}{\kappa} \go S^n_{n-1} R^n_{n-1} \pr{C}{\kappa}.
\]
This is the identity in dimensions less than $n$, and in dimension $n$ it
consists of the maps
\[
\pr{s}{t}: C(\pi) \go P_{\pi} (C)
\]
($\pi\in\pd(n)$). This is always surjective as $C$ carries a precontraction,
and is injective precisely when $C$ satisfies the condition for
precontractions on it to be called contractions. So the unit at
\pr{C}{\kappa} is an isomorphism if and only if \pr{C}{\kappa} is an object of
$\fcat{OWC}_n$. 

\item The proof is just like that of part~\bref{propn:n-(n-1)}. Again it will
be useful to have an explicit description of the right adjoint $S^\omega_k$
of $R^\omega_k$: it is given by 
\[
\!\!\!\!
\begin{array}{l}
(S^\omega_k X)(m) =					\\
\left\{
\begin{array}{ll}
X(m)	&\mbox{for }0\leq m\leq k		\\
\{ \pr{x_0}{x_1}\in (X(k))^2 \such
	s(x_0)=s(x_1), t(x_0)=t(x_1) \}
	&\mbox{for }m\geq k+1.
\end{array}
\right.
\end{array}
\]
The source and target maps in dimensions $\leq k$ are as in $X$; from
dimension $k+1$ to dimension $k$ they are first and second projection; and in
dimensions above $k+1$, they are identities. 
\done
\end{enumerate}
\end{proof}

\begin{figure}
\piccy{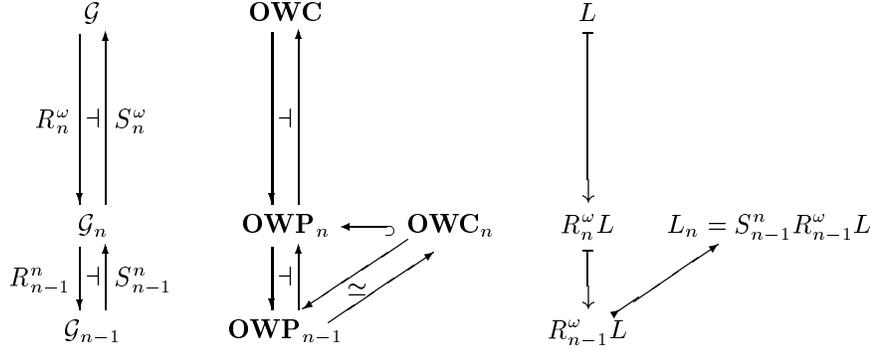}
\caption{Relating $L$ and $L_n$}
\label{fig:strategy}
\end{figure}

From this we deduce the following corollary, which shows incidentally that
$\fcat{OWC}_n$ does have an initial object. The overall strategy is depicted
in Figure~\ref{fig:strategy}.
\begin{cor}
$S^n_{n-1} R^\omega_{n-1} \pr{L}{\lambda}$ is an initial object of
$\fcat{OWC}_n$.
\end{cor}
\begin{proof}
The functor $R^\omega_{n-1}: \fcat{OWC} \go \fcat{OWP}_{n-1}$ constructed in
part~\bref{propn:omega-k} of the proposition has a right adjoint, so
$R^\omega_{n-1} \pr{L}{\lambda}$ is initial in $\fcat{OWP}_{n-1}$. The
functor $S^n_{n-1}: \fcat{OWP}_{n-1} \go \fcat{OWC}_n$ constructed in
part~\bref{propn:n-(n-1)-eqv} is an equivalence, so $S^n_{n-1}
(R^\omega_{n-1} \pr{L}{\lambda})$ is initial in $\fcat{OWC}_n$.
\done
\end{proof}

We write \pr{L_n}{\lambda_n} for the initial object of $\fcat{OWC}_n$: that
is,
\[
\pr{L_n}{\lambda_n} = S^n_{n-1} R^\omega_{n-1} \pr{L}{\lambda}.
\]

Before moving to the second half of the proof of Theorem~\ref{thm:n-eqv}, let
us recall some notation. Fix $n\in\nat$. To any $m$-pasting diagram $\pi$
there is associated the globular set $\rep{\pi}$, and we may turn $\rep{\pi}$
into an $n$-globular set by restriction (truncation). If $m\leq n$ then this
only amounts to ignoring some $\emptyset$'s, since $\rep{\pi}(k) = \emptyset$
for $k>m$. In Appendix~\ref{app:strict-omega} we show that if $X$ is a
globular set and $m\leq n$ then
\[
(T_n X)(m) = 
\coprod_{\pi\in\pd(m)} \homset{\cat{G}_n}{\rep{\pi}}{X}.
\]
Given a $T_n$-operad $C$, a $C$-algebra structure on $X$ consists of a map
\[
h_{\pi}: C(\pi) \times \homset{\cat{G}_n}{\rep{\pi}}{X}
\go X(m)
\]
for each $m\leq n$ and $\pi\in\pd(m)$, subject to various axioms. For
$\theta\in C(\pi)$, we write
\[
\ovln{\theta} = h_\pi \pr{\theta}{\dashbk}: 
\homset{\cat{G}_n}{\rep{\pi}}{X} \go X(m).
\]

Now, any weak $n$-category is isomorphic to a `strictly' $n$-dimensional weak
$\omega$-category: that is, to one whose underlying globular set is of the
form 
\begin{equation}	\label{eq:extended-glob-set}
\cdots \pile{\rTo^1 \\ \rTo_1} X(n) 
\pile{\rTo^1 \\ \rTo_1} X(n)
\pile{\rTo^s \\ \rTo_t} \cdots
\pile{\rTo^s \\ \rTo_t} X(0).
\end{equation}
So to prove Theorem~\ref{thm:n-eqv} it is enough to prove that the category
of strictly $n$-dimensional weak $\omega$-categories is equivalent to
$\Alg(L_n)$; indeed, we will prove that these two categories are isomorphic. 

Let $X$ be an $n$-globular set. An $L$-algebra structure on the globular
set~\bref{eq:extended-glob-set} consists precisely of an $(R^\omega_n
L)$-algebra structure on $X$ together with a dotted arrow
\begin{diagram}[width=3.5em]
L(\sigma) \times \homset{\cat{G}_n}{\rep{\bdry\sigma}}{X}
&\rGet
&X(n)							\\
\dTo<{s\times 1} \dTo>{t\times 1}
&
&\dTo<1 \dTo>1						\\
L(\bdry\sigma) \times \homset{\cat{G}_n}{\rep{\bdry\sigma}}{X}
&\rTo_{h_{\bdry\sigma}}
&X(n)							\\
\end{diagram}
making the diagram commute serially, for each $\sigma\in\pd(n+1)$. To see
this, note first that a map $\rep{\bdry\sigma} \go X$ is the same as a map
from $\rep{\sigma}$ to the globular set~\bref{eq:extended-glob-set}, so an
algebra structure on~\bref{eq:extended-glob-set} yields such a dotted arrow
for each $\sigma$. Conversely, given such arrows, all the $L$-algebra
structure in higher dimensions is uniquely determined, and the algebra axioms
are automatically satisfied. There is at most one way of choosing the dotted
arrows, and such a way exists if and only if
\begin{quote}
for all $\sigma\in\pd(n+1)$ and $\theta\in L(\sigma)$, 
\[
\ovln{s\theta} = \ovln{t\theta}:
\homset{\cat{G}_n}{\rep{\bdry\sigma}}{X} \go X(n).
\]
\end{quote}
Since $L$ admits a
contraction, and for each $\pi\in\pd(n)$ there exists $\sigma\in\pd(n+1)$
with $\bdry\sigma = \pi$, this condition is equivalent to: 
\begin{quote}
for all $\pi\in\pd(n)$ and $\pr{\theta_0}{\theta_1}\in P_\pi (L)$, 
\begin{equation}	\label{eq:tidy}
\ovln{\theta_0} = \ovln{\theta_1}:
\homset{\cat{G}_n}{\rep{\pi}}{X} \go X(n).
\end{equation}
\end{quote}
So a strictly $n$-dimensional weak $\omega$-category consists precisely of an
$(R^\omega_n L)$-algebra $X$ satisfying condition~\bref{eq:tidy}. 

Working from the other end, let 
\[
C = R^\omega_n L \in T_n\hyph\fcat{Operad},
\]
and let $u$ be the unit map $C\go S^n_{n-1} R^n_{n-1} C$ coming from the
adjunction $R^n_{n-1} \ladj S^n_{n-1}$. By the description of this adjunction
in the proofs of Proposition~\ref{propn:three-adjns}\bref{propn:n-(n-1)}
and~\bref{propn:n-(n-1)-eqv}, $u_\pi(\theta_0) = u_\pi(\theta_1)$ whenever
$\pi\in\pd(n)$ and $\pr{\theta_0}{\theta_1} \in P_\pi (L)$; since $C$ admits
a precontraction, $u$ is (surjective and therefore) the \emph{universal} map
out of $C$ with this property. It follows that an algebra for $S^n_{n-1}
R^n_{n-1} C$ is exactly an algebra $X$ for $C$ satisfying the
condition~\bref{eq:tidy}. (The details of this step are omitted; the idea is
perhaps most naturally expressed in terms of endomorphism
operads~(\ref{sec:endo}).) So we have
\begin{eqnarray*}
\lefteqn{\mbox{(strictly $n$-dimensional weak $\omega$-categories)}}	\\
&\iso	&((R^\omega_n L)\mbox{-algebras $X$ satisfying~\bref{eq:tidy}})	\\
&\iso	&\Alg(S^n_{n-1} R^n_{n-1} C)					\\
&= 	&\Alg(S^n_{n-1} R^\omega_{n-1} L)				\\
&\iso	&\Alg(L_n).
\end{eqnarray*}
We have only discussed the objects of these categories, and not their
morphisms; but everything works as it should since in each case the morphisms
are the maps strictly preserving all the structure. This proves
Theorem~\ref{thm:n-eqv}.

\section{Weak $2$-Categories}	\label{sec:weak-2}

A polite person proposing a definition of weak $n$-category should explain
what happens when $n=2$. With our definition, the category $\fcat{Wk}\hyph
2\hyph\Cat$ of weak 2-categories turns out to be equivalent to \UBistr, the
category of small unbiased bicategories and unbiased strict functors. This is
the main result of this section.

Note that because the morphisms in $\fcat{Wk}\hyph 2\hyph\Cat$ are
\emph{strict} maps (as noted on page~\pageref{p:strict-omega-ftrs}), we
obtain an equivalence with \UBistr, not \UBiwk\ or \UBilax; and unlike the
weak and lax versions, \UBistr\ is not equivalent to the corresponding
category of classical bicategories (at least, the obvious functor is not an
equivalence). So we cannot conclude that $\fcat{Wk}\hyph 2\hyph\Cat$ is
equivalent to $\Bistr$. Nevertheless, the results of Chapter~\ref{ch:bicats}
mean that it is fair to regard classical bicategories as `essentially the
same as' unbiased bicategories, and therefore, by the results below,
`essentially the same as' weak 2-categories. If the definition of weak
functor between $n$-categories were in place, we would expect there to be a
genuine equivalence between \Biwk\ and the category of weak 2-categories and
weak 2-functors.

Before embarking on the analysis of $n=2$, let us check that things are as
they should be for $n=0$ and $n=1$. In all cases, we will analyse $\Alg(L_n)$
rather than the equivalent category $\fcat{Wk}\hyph n\hyph\Cat$,  where
\pr{L_n}{\lambda_n} is the initial $T_n$-operad with contraction. 

\begin{thm}
$\fcat{Wk}\hyph 0\hyph\Cat \eqv \Set$.
\end{thm}
\paragraph*{Proof}
$T_0$ is the identity monad on the category $\cat{G}_0$ of sets, so a
$T_0$-operad is a monoid. Any $T_0$-operad carries a unique contraction
(vacuously), so $\fcat{OWC}_0$ is the category of monoids; the initial object
of $\fcat{OWC}_0$ is the monoid $1$. An algebra for the terminal $T_0$-operad
is just a $T_0$-algebra (see~\ref{eg:algs}\bref{eg:alg-terminal}), so 
$
\Alg(L_0) \iso \cat{G}_0^{T_0} 
$\linebreak$
\iso \Set.  
$ 
\done

\begin{thm}
$\fcat{Wk}\hyph 1\hyph\Cat \eqv \Cat$.
\end{thm}
\paragraph*{Proof}
$T_1$ is the free category monad \fc\ on the category $\cat{G}_1$ of directed
graphs, so a $T_1$-operad is an \fc-operad (see~\ref{sec:fcm}). A
$T_1$-operad $C$ admits at most one contraction, and does admit one if and
only if the function
\[
\pr{s}{t}: C(\pi) \go C(\zeropd) \times C(\zeropd)
\]
is a bijection for each 1-pasting diagram $\pi$ (where
$\zeropd\in\pd(0)$). It follows that the terminal $T_1$-operad is the initial
object, $L_1$, of $\fcat{OWC}_1$. Hence
$
\Alg(L_1) \iso \cat{G}_1^{T_1} 
$\linebreak$
\iso \Cat.
$
\done

\paragraph*{} The full proof that $\fcat{Wk}\hyph 2\hyph\Cat \eqv \UBistr$
involves rather more detailed manipulation than the reader would probably
like to see. To keep the presentation light, I will use the coherence results
of Appendix~\ref{app:unbiased} for unbiased bicategories in the inexact form
`all diagrams commute'. In the same spirit, I will use the following
formulations of the notions of functor and natural transformation:
\begin{itemize}
\item A functor%
\label{p:ftrs-and-nts}
$F: \cat{A}_1 \times\cdots\times \cat{A}_n \go \cat{A}$
consists of
\begin{itemize}
\item a function $F_0: \ob\cat{A}_1 \times\cdots\times \ob\cat{A}_n \go
\ob\cat{A}$ 
\item a function assigning to each array of maps
\begin{equation}	\label{eq:array}
\begin{array}{cl}
a_1^0 \goby{\alpha_1^1} \cdots \goby{\alpha_1^{k_1}} a_1^{k_1}	&
\mbox{in }\cat{A}_1,						\\
\ldots\diagspace\ldots	&					\\
a_n^0 \goby{\alpha_n^1} \cdots \goby{\alpha_n^{k_n}} a_n^{k_n}	&
\mbox{in }\cat{A}_n					
\end{array}
\end{equation}
a map
\[
F_0 \bftuple{a_1^0}{a_n^0} \go F_0 \bftuple{a_1^{k_1}}{a_n^{k_n}}
\]
in \cat{A}, 
\end{itemize}
obeying `all reasonable coherence axioms'.
\item A natural transformation
\[
\cat{A}_1 \times\cdots\times \cat{A}_n%
\ctwo{F}{F'}{\phi}%
\cat{A}
\]
consists of a function assigning to each array of maps~\bref{eq:array} a map
\[
F_0 \bftuple{a_1^0}{a_n^0} \go F'_0 \bftuple{a_1^{k_1}}{a_n^{k_n}},
\]
in such a way that `all reasonable coherence axioms' hold.  
\end{itemize}
In all parts of the proof where such sweeping language is used, the diligent
reader should not find it difficult to fill in the details.

It will also be useful to have some notation for $m$-pasting diagrams when
$m\leq 2$. The unique $0$-pasting diagram will be denoted \zeropd. We have
$\pd(1) \iso (\pd(0))^*$\linebreak$\iso \nat$, and the element of $\pd(1)$
corresponding to $n\in\nat$ will be denoted $\pi_n$; so $\pi_n$ is usually
drawn as
\[
\gfst{}\gone{}\diagspace\cdots\diagspace\gone{}\glst{}
\]
($n$ arrows). Similarly, $\pd(2) \iso (\pd(1))^* \iso \nat^*$, and we write
$\pi_{k_1, \ldots, k_n}$ for the $2$-pasting diagram corresponding to
$\bftuple{k_1}{k_n} \in \nat^*$, which is usually drawn as a diagram
\[
\gfst{}\gdots{}{}{}{}{}{}\gfbw{}%
\diagspace\cdots\cdots\diagspace%
\gfbw{}\gdots{}{}{}{}{}{}\glst{}%
\]
with $n$ columns and $k_i$ 2-cells in the $i$th column. 

\pagebreak

\begin{thm}
$\fcat{Wk}\hyph 2\hyph\Cat \eqv \UBistr$.
\end{thm}
\begin{proof}
\begin{sloppypar}
First we identify the initial object \pr{L_2}{\lambda_2} of
$\fcat{OWC}_2$; and since $\fcat{OWC}_2 \eqv \fcat{OWP}_1$, this means
examining $\fcat{OWP}_1$. A precontraction on a $T_1$-operad $C$ consists of
a function
\[
\kappa_{\pi_n}: C(\zeropd) \times C(\zeropd) \go C(\pi_n)
\]
for each $n\in\nat$, such that
\[
s(\kappa_{\pi_n} \pr{\theta_0}{\theta_1} ) = \theta_0,
\diagspace
t(\kappa_{\pi_n} \pr{\theta_0}{\theta_1} ) = \theta_1
\]
for all $\theta_0$, $\theta_1$. A $T_1$-operad $C$ with $C(\zeropd)=1$ is
merely a plain operad---call it $\twid{C}$---and a precontraction on $C$
consists of a distinguished element of $\twid{C}(n)$ for each $n\in\nat$. The
operad $\tr$ described in~\ref{sec:free-multicats} and~\ref{sec:coherence},
together with the element $\nu_n = \bftuple{\utree}{\utree}$ of $\tr(n)$ for
each $n$, therefore defines a $T_1$-operad with precontraction. Using the
fact that $\tr$ is the free plain operad on the terminal (free~monoid)-graph,
it is easy to see that this is the initial object of $\fcat{OWP}_1$. By
Proposition~\ref{propn:three-adjns}, \pr{L_2}{\lambda_2} is $S^2_1$ applied
to this initial object: that is,
\begin{eqnarray*}
L_2 (\zeropd)			&=	&1,		\\
L_2 (\pi_n)			&=	&\tr(n),	\\
L_2 (\pi_{k_1,\ldots,k_n})	&=	&\tr(n) \times \tr(n)
\end{eqnarray*}
($n,k_i \in\nat$). In dimension $1$, the $T_2$-operad structure is as in the
plain operad \tr. Given that the source and target functions
\[
L_2 (\pi_{k_1,\ldots,k_n}) \parpair{}{} L_2 (\pi_n)
\]
are first and second projection, the $T_2$-operad structure in dimension $2$
is uniquely determined. 
\end{sloppypar}

This fully describes $L_2$. An algebra for $L_2$ is, therefore:
\begin{itemize}
\item a 2-globular set $X(2) \parpair{}{} X(1) \parpair{}{} X(0)$
\item for each $n\in\nat$ and $\tau\in\tr(n)$, a function
\[
\ovln{\tau}: \homset{\cat{G}_2}{\rep{\pi_n}}{X} \go X(1)
\]
\item for each $n, k_1, \ldots, k_n \in\nat$ and $\tau, \tau' \in \tr(n)$, a
function
\[
\ovln{\pr{\tau}{\tau'}}:  
\homset{\cat{G}_2}{\rep{\pi_{k_1,\ldots,k_n}}}{X} \go X(2)
\]
\end{itemize}
satisfying axioms concerning the source and target of
$\ovln{\pr{\tau}{\tau'}}$ in terms of $\ovln{\tau}$ and $\ovln{\tau'}$,
together with the axioms for an algebra (which we regard as `all reasonable
coherence axioms'). 

Rephrasing this a little, an algebra for $L_2$ consists of
\begin{itemize}
\item a set $\scat{B}_0$ (which is the $X(0)$ of the previous paragraph)
\item for each $A, B \in \scat{B}_0$, a directed graph
\[
\homset{\scat{B}}{A}{B} = 
( \homset{\scat{B}}{A}{B}_1 \parpair{}{} \homset{\scat{B}}{A}{B}_0 )
\]
\item for each $\tau\in\tr(n)$ and $A_0, \ldots, A_n \in \scat{B}_0$, a
function 
\[
\homset{\scat{B}}{A_0}{A_1}_0 \times\cdots\times 
\homset{\scat{B}}{A_{n-1}}{A_n}_0
\goby{\ovln{\tau}}
\homset{\scat{B}}{A_0}{A_n}_0
\]
\item for each $\tau, \tau' \in \tr(n)$, each $A_0, \ldots, A_n \in
\scat{B}_0$, and each array of arrows
\[
\begin{array}{cl}
f_1^0 \goby{\alpha_1^1} \cdots \goby{\alpha_1^{k_1}} f_1^{k_1}	&
\mbox{in }\homset{\scat{B}}{A_0}{A_1},				\\
\ldots\diagspace\ldots	&					\\
f_n^0 \goby{\alpha_n^1} \cdots \goby{\alpha_n^{k_n}} f_n^{k_n}	&
\mbox{in }\homset{\scat{B}}{A_{n-1}}{A_n},				
\end{array}
\]
an arrow
\[
\ovln{\tau}\bftuple{f^0_1}{f^0_n} \go 
\ovln{\tau'}\bftuple{f^{k_1}_1}{f^{k_n}_n}
\]
in \homset{\scat{B}}{A_0}{A_n}, 
\end{itemize}
satisfying `all reasonable coherence axioms'. These axioms imply that if
$\tau=\utree \in\tr(1)$ then the function
\[
\ovln{\tau}: \homset{\scat{B}}{A_0}{A_1}_0 \go \homset{\scat{B}}{A_0}{A_1}_0
\]
is the identity. Now taking $n=1$ and $\tau=\tau'=\utree$ in the fourth
item, we have a function which assigns to each string of arrows
\[
f^0 \goby{\alpha^1} \cdots \goby{\alpha^k} f^k
\]
in \homset{\scat{B}}{A}{B} an arrow $\ovln{\utree}(f^0) \go
\ovln{\utree}(f^k)$, that is, $f^0 \go f^k$. This gives the directed graph
\homset{\scat{B}}{A}{B} the structure of a category. By the preliminary
comments on functors and natural transformations
(page~\pageref{p:ftrs-and-nts}), an $L_2$-algebra therefore consists of 
\begin{itemize}
\item a set $\scat{B}_0$
\item for each $A,B\in\scat{B}_0$, a category \homset{\scat{B}}{A}{B}
\item for each $\tau\in\tr(n)$ and $A_0, \ldots, A_n \in\scat{B}_0$, a
functor
\[
\ovln{\tau}: 
\homset{\scat{B}}{A_0}{A_1} \times\cdots\times
\homset{\scat{B}}{A_{n-1}}{A_n}
\go 
\homset{\scat{B}}{A_0}{A_n}
\]
\item for each $\tau,\tau' \in\tr(n)$ and $A_0, \ldots, A_n \in \scat{B}_0$,
a natural transformation 
\[
\ovln{\pr{\tau}{\tau'}}: \ovln{\tau} \go \ovln{\tau'}
\]
\end{itemize}
satisfying `all reasonable coherence axioms'. Writing $\ovln{\tau}$ as
$\comp_\tau$ and $\ovln{\pr{\tau}{\tau'}}$ as $\omega_{\tau,\tau'}$, we
see that this is just the description of a (small) unbiased bicategory given
by Theorem~\ref{thm:coh-u-bicats} and the comments thereafter.  \done
\end{proof}

We have proved that in the cases $n=0,1,2$, the category $\fcat{Wk}\hyph
n\hyph\Cat$ is equivalent to, respectively, \Set, \Cat\ and \UBistr. In fact,
we have proved that $\Alg(L_0)$ is \emph{isomorphic} to \Set, and similarly
$\Alg(L_1)$ to \Cat. The analogous property for $n=2$ does not quite hold,
because an unbiased bicategory is defined to be a structure on a `graph of
directed graphs' (that is, a set $\scat{B}_0$ together with a directed graph
\homset{\scat{B}}{A}{B} for each $A,B\in\scat{B}_0$) whereas an $L_2$-algebra
is a structure on a 2-globular set, and the category $\cat{G}_2$ of
2-globular sets is merely \emph{equivalent} to the category of graphs of
directed graphs. However, the proof reveals that this difference is the only
obstacle to the equivalence $\Alg(L_2) \eqv \UBistr$ becoming an
isomorphism. 

This concludes the material on weak $\omega$- and $n$-categories, and indeed
the main body of the thesis. From the explanation of the formal definition of
weak $\omega$-category, and the analysis of the case $n=2$, I hope that the
reader is persuaded that the proposed definition is a reasonable one.
Nonetheless, we have clearly only touched the beginning of a theory of weak
higher-dimensional categories.

\appendix

\chapter{Biased \emph{vs.}\/ Unbiased Bicategories}
\label{app:unbiased}

In this appendix we prove the following results from Chapter~\ref{ch:bicats},
concerning the forgetful functor $V: \UBilax \go \Bilax$:
\begin{trivlist} \item 
\textbf{Theorem \ref{thm:biased-comparison}}\ \itshape
With the definitions given in section~\ref{sec:versus}, 
\begin{enumerate}
\item $V(\Bee)$ is a bicategory and $V\pr{F}{\phi}$ is a lax functor
\item $V$ preserves composition and identities, so forms a functor 
\[
\UBilax \go \Bilax
\] 
\item $V$ is full, faithful and surjective on objects.
\end{enumerate}
\end{trivlist}
\begin{trivlist} \item 
\textbf{Corollary \ref{cor:wk-biased-comparison}}\ \itshape
The restricted functor $V_\mr{wk}:\UBiwk \go \Biwk$ is also full, faithful
and surjective on objects.
\end{trivlist}
It is possible to do the proofs in a thoroughly explicit way, as a very long
sequence of calculations. At the other extreme, it is possible to state and
prove a very general result, as follows. In the classical definition of
bicategory, there is one nullary and one (horizontal) binary composition
operation. In the unbiased definition, there is one $n$-ary
operation for each $n\in\nat$. Given a sequence $(\Omega_n)_{n\in\nat}$ of
sets, there is a notion of `bicategory' in which there is one $n$-ary
operation for each member of $\Omega_n$, and corresponding notions of lax and
weak functors. So the classical case has $\Omega_n=1$ for $n=0,2$ and
$\Omega_n=\emptyset$ otherwise, and the unbiased case has $\Omega_n=1$ for
all $n$. As long as $\Omega_0 \neq \emptyset$ and $\Omega_n \neq \emptyset$
for some $n\geq 2$, this gives a category of `bicategories' and lax functors
which is equivalent to \Bilax. This is the method employed for monoidal
categories in~\cite{WMC}.

To keep things short, we shun both extremes and follow a Third Way. The
strategy is to start by proving some coherence results for unbiased
bicategories and functors, of the form `every diagram commutes', and to
recall similar coherence results for classical bicategories and
functors. (All of this so far would be necessary even in the abstract
approach outlined above.) We can then use these results as an aid to
calculation when proving that $V$ is well-defined and an equivalence; indeed,
they are so powerful that detailed calculations can almost entirely be
avoided.

Incidentally, the proofs of the coherence results for the unbiased theory are
all absolutely straightforward. Just a little care is needed to keep track of
the subscripts, but the proofs call for none of the ingenuity required in
proving coherence for classical bicategories (see e.g.~\cite[1.1]{JS}).

The issue of large \emph{vs.}\ small structures is not addressed here; it is
left as a matter of conscience to the reader.

\section{Coherence}		\label{sec:coherence}

\subsection*{Preliminaries}

To state our results we need some new language. 

First recall the 2-category $\Cat\hyph\Gph$ from page~\pageref{p:Cat-Gph}
(Remark~\bref{rmk:Cat-Gph}). There is some extra structure on $\Cat\hyph\Gph$:
if \Bee, \Beep\ are \Cat-graphs with $\Bee_0 = \Bee'_0 = S$, say, then there
is a \Cat-graph $\Bee\otimes\Beep$ defined by
\[
(\Bee\otimes\Beep)_0 = S,
\diagspace
(\Bee\otimes\Beep)(s_1, s_2) =
\coprod_{s\in S} \Bee(s_1, s) \times \Beep(s, s_2).
\]
There is also an object $\cat{I}_S$ of $\Cat\hyph\Gph$ defined by
\begin{eqnarray*}
(\cat{I}_S)_0 	&= 	&S,	\\
(\cat{I}_S)(s_1, s_2)	&=	&
\left\{
\begin{array}{ll}
\One	&\mbox{if }s_1 = s_2	\\
\emptyset	&\mbox{otherwise.}
\end{array}
\right.
\end{eqnarray*}
This defines a monoidal category structure on $\Cat^{S\times S}$ for each set
$S$. 

Furthermore, if $\Bee \goby{F} \cat{C}$ and $\Beep \goby{F'} \cat{C'}$ are
maps in $\Cat\hyph\Gph$ with $\Bee=\Beep=S$, say, $\cat{C}_0 = \cat{C}'_0$,
and $F_0 = F'_0$, then there is a map $F\otimes F': \Bee\otimes\Beep \go
\cat{C}\otimes\cat{C'}$ in $\Cat\hyph\Gph$ defined by
\[
\begin{array}{c}
(F\otimes F')_0 = F_0 = F'_0,	\\
(F\otimes F')_{s_1, s_2} (p,p') = (F_{s_1,s}(p), F'_{s,s_2}(p'))
\end{array}
\]
for $s_1,s_2,s\in S$, $p\in\Bee(s_1,s)$ and $p'\in\Beep(s,s_2)$. In
particular, if \Bee\ is a \Cat-graph then there is a \Cat-graph
$\Bee^{\otimes n}$ for each $n\in\nat$, and if $F: \Bee \go \cat{C}$ is a map
of \Cat-graphs then there is a map $F^{\otimes n}: \Bee^{\otimes n} \go
\cat{C}^{\otimes n}$. (So, for instance, the free 2-category functor on
$\Cat\hyph\Gph$ is given by $\Bee \goesto \coprod_{n\in\nat} \Bee^{\otimes
n}$.)

I will not attempt to describe exactly what structure is formed by
$\Cat\hyph\Gph$ together with these tensor operations, although we will
implicitly use some of its fairly obvious properties (such as functoriality
of tensor). If we were discussing monoidal categories rather than
bicategories, then the place of $\Cat\hyph\Gph$ would be taken by the
monoidal category \triple{\Cat}{\times}{\One}.

The definitions of unbiased bicategory and unbiased lax/weak functor can now
be recast as follows. An unbiased bicategory consists of a \Cat-graph \Bee\
together with a functor $\comp_n: \Bee^{\otimes n} \go \Bee$ for each
$n\in\nat$ and natural isomorphisms
\[
\begin{array}{c}
\begin{diagram}[size=1.5em]
\Bee^{\otimes (k_1 +\cdots+ k_n)}	&	&
\rTo^{\comp_{k_1} \otimes\cdots\otimes \comp_{k_n}}	&	&
\Bee^{\otimes n}	\\
	&\rdTo(4,4)<{\comp_{k_1 +\cdots+ k_n}}	
&	&\swarrow \gamma_{k_1,\ldots,k_n}&	\\
	&		&		&		&\dTo>{\comp_n}	\\
	& 		&		&		&		\\
	&		&		&		&\Bee		\\
\end{diagram}
\\
\Bee^{\otimes 1}%
\ctwo{\diso}{\comp_1}{\iota}%
\Bee
\end{array}
\]
satisfying associativity and identity axioms. An unbiased lax functor
$\pr{F}{\phi}: \Bee \go \Beep$ consists of a map $F: \Bee\go\Beep$ of
\Cat-graphs together with a natural transformation
\[
\begin{ntdiag}
\Bee^{\otimes n}	&&\rTo^{F^{\otimes n}}	&&\Beep^{\otimes n}	\\
		&		&		&\		&	\\
\dTo<{\comp_n}		&&\ldTo>{\phi_n}	&&\dTo>{\comp_n}	\\
		&\		&		&		&	\\
\Bee			&&\rTo_{F}		&&\Beep
\end{ntdiag}
\]
for each $n$, satisfying axioms. (So unbiased bicategories are weak algebras,
and unbiased lax functors are lax maps of weak algebras, for the free
2-category 2-monad on $\Cat\hyph\Gph$.)

We will also need the language of trees. By definition, \tr\ is the free
(non-symmetric) operad (of sets) on the terminal object of $\Set^\nat$, as
explained more fully in~\ref{sec:free-multicats}. Explicitly, we can define
for each $n\in\nat$ a set $\tr(n)$ of \emph{$n$-leafed trees} by the
following recursive clauses:
\begin{itemize}	\label{p:tree-defn}
\item $\tr(1)$ has an element \utree\ (a formal symbol)
\item if $n\in\nat$ and $\tau_1\in\tr(k_1), \ldots, \tau_n\in\tr(k_n)$, then
$\tr(k_1 +\cdots+ k_n)$ has an element \bftuple{\tau_1}{\tau_n}.
\end{itemize}
(See Example~\ref{eg:cart-mnds}\bref{eg:mon-tree} for why the word `tree' is
used.) We call \utree\ the \emph{unit} tree, and define for each
$\tau\in\tr(n)$ and $\tau_1\in\tr(k_1), \ldots, \tau_n\in\tr(k_n)$ a
\emph{composite} tree $\tau\of\bftuple{\tau_1}{\tau_n}$ as follows.
\pagebreak
\begin{itemize}
\item If $\tau=\utree$ then $\tau\of(\tau_1) = \tau_1$
\item Suppose $\tau = \bftuple{\sigma_1}{\sigma_r}$ with $\sigma_i \in
\tr(n_i)$ and $n_1 +\cdots+ n_r = n$: then we may write the sequence $\tau_1,
\ldots, \tau_n$ as $\tau_1^1, \ldots, \tau_1^{n_1}, \ldots, \tau_r^1, \ldots,
\tau_r^{n_r}$ and define 
\[
\tau\of\bftuple{\tau_1}{\tau_n} =
\bftuple{\sigma_1 \of \bftuple{\tau_1^1}{\tau_1^{n_1}}}{\sigma_r \of
\bftuple{\tau_r^1}{\tau_r^{n_r}}}.
\]
\end{itemize}
Composition and unit obey associativity and unit laws: in other words, \tr\
forms a non-symmetric operad. Note also that if $\nu_n$ is the $n$-leafed tree
\bftuple{\utree}{\utree} then $\bftuple{\tau_1}{\tau_n} = \nu_n \of
\bftuple{\tau_1}{\tau_n}$.

\subsection*{Coherence for unbiased bicategories}

Fix an unbiased bicategory \Bee. Define for each $n\in\nat$ and
$\tau\in\tr(n)$ a functor $\comp_\tau: \Bee^{\otimes n} \go \Bee$, as
follows:
\begin{itemize}
\item $\comp_\utree$ is the identity on \Bee
\item if $\tau_1\in\tr(k_1), \ldots, \tau_n\in\tr(k_n)$ then
$\comp_\bftuple{\tau_1}{\tau_n}$ is the composite
\begin{diagram}
\Bee^{\otimes (k_1 +\cdots+ k_n)}
&\rTo^{\comp_{\tau_1} \otimes\cdots\otimes \comp_{\tau_n}}
&\Bee^{\otimes n} 
&\rTo^{\comp_n} 
&\Bee.
\end{diagram}
\end{itemize}
(More accurately, $\comp_\utree$ is not the identity but the canonical
isomorphism $\Bee^{\otimes 1} \go \Bee$. I will ignore such distinctions.)

\begin{propn}
\begin{enumerate}
\item 	\label{propn:comp-a}
If $\tau\in\tr(n), \tau_1\in\tr(k_1), \ldots, \tau_n\in\tr(k_n)$ then
\[
\comp_{\tau\of\bftuple{\tau_1}{\tau_n}} = 
\comp_{\tau} \of (\comp_{\tau_1} \otimes\cdots\otimes \comp_{\tau_n})
\]
\item  	\label{propn:comp-b}
$\comp_\utree = \id$
\item  	\label{propn:comp-c}
$\comp_{\nu_n} = \comp_n$.
\end{enumerate}
\end{propn}
\begin{proof}
Part~\bref{propn:comp-a} is a straightforward induction on the
structure of $\tau$. Part \bref{propn:comp-b} is just the definition of
$\comp_\utree$. Part~\bref{propn:comp-c} is also straightforward. 
\done
\end{proof}

Next define for each tree $\tau\in\tr(n)$ a natural isomorphism $\omega_\tau:
\comp_\tau \go $\linebreak$\comp_n$, by 
\begin{itemize}
\item $\omega_\utree = \iota: \id \go \comp_1$ 
\item if $\tau_1\in\tr(k_1), \ldots, \tau_n\in\tr(k_n)$ then
$\omega_{\bftuple{\tau_1}{\tau_n}}$ is the composite
\begin{eqnarray*}
\comp_{\bftuple{\tau_1}{\tau_n}}	&=&
\comp_n \of (\comp_{\tau_1} \otimes\cdots\otimes \comp_{\tau_n})	\\
	&\goby{1* (\omega_{\tau_1} \otimes\cdots\otimes \omega_{\tau_n})}
	&\comp_n \of (\comp_{k_1} \otimes\cdots\otimes \comp_{k_n})	\\
	&\goby{\gamma_{k_1,\ldots,k_n}}
	&\comp_{k_1 +\cdots+ k_n}.					\\
\end{eqnarray*}
\end{itemize}
The $\omega_\tau$'s fit together coherently, as expressed by the following
result.
\begin{propn}	\label{propn:single-omega}
\begin{enumerate}
\item		\label{propn:single-omega-a}
If $\tau\in\tr(n), \tau_1\in\tr(k_1), \ldots, \tau_n\in\tr(k_n)$ then 
\begin{diagram}
\comp_{\tau\of\bftuple{\tau_1}{\tau_n}}	&\rEquals
&\comp_\tau \of (\comp_{\tau_1} \otimes\cdots\otimes \comp_{\tau_n})	\\
&&\dTo>{\omega_\tau * (\omega_{\tau_1} \otimes\cdots\otimes\omega_{\tau_n})}
\\
\dTo<{\omega_{\tau\of\bftuple{\tau_1}{\tau_n}}}	&	&
\comp_n \of (\comp_{k_1} \otimes\cdots\otimes \comp_{k_n})	\\
&&\dTo>{\gamma_{k_1,\ldots,k_n}}					\\
\comp_{k_1 +\cdots+ k_n} 	&\rEquals	&\comp_{k_1 +\cdots+ k_n}\\
\end{diagram}
commutes
\item 	\label{propn:single-omega-b}
The diagram
\begin{diagram}
\comp_{\utree}		&\rEquals	&\id		\\
\dTo<{\omega_\utree}	&		&\dTo>{\iota}	\\
\comp_1			&\rEquals	&\comp_1	
\end{diagram}
commutes
\item 	\label{propn:single-omega-c}
$\omega_{\nu_n} = 1$, 
$\omega_{\nu_n \of \bftuple{\nu_{k_1}}{\nu_{k_n}}} =
\gamma_{k_1, \ldots, k_n}$, 
and $\omega_{\utree} = \iota$.
\end{enumerate}
\end{propn}
\begin{proof}
As in the previous proof,~\bref{propn:single-omega-a} is by induction on
$\tau$, \bref{propn:single-omega-b} is immediate,
and~\bref{propn:single-omega-c} is straightforward.
\done
\end{proof}

Everything so far works for \emph{lax} bicategories, but the next part does
not. For each $\tau, \tau' \in\tr(n)$, define a natural isomorphism
\[
\omega_{\tau, \tau'} = 
(\comp_\tau \goby{\omega_\tau} \comp_n \goby{\omega_{\tau'}^{-1}}
\comp_{\tau'}). 
\]
The $\omega_{\tau,\tau'}$'s also fit together coherently:
\pagebreak
\begin{thm}	\label{thm:coh-u-bicats}
\begin{enumerate}
\item  	\label{propn:double-omega-a}
If $\tau,\tau' \in\tr(n), \tau_1,\tau'_1 \in\tr(k_1), \ldots,
\tau_n,\tau'_n \in\tr(k_n)$ then
\begin{diagram}[scriptlabels]
\comp_{\tau\of\bftuple{\tau_1}{\tau_n}}
&\rEquals
&\comp_\tau \of (\comp_{\tau_1} \otimes\cdots\otimes \comp_{\tau_n})	\\
\dTo<{\omega_{\tau\of\bftuple{\tau_1}{\tau_n},
\tau'\of\bftuple{\tau'_1}{\tau'_n}}} 	&	&
\dTo>{\omega_{\tau,\tau'} * (\omega_{\tau_1,\tau'_1} \otimes\cdots\otimes
\omega_{\tau_n,\tau'_n})}						\\
\comp_{\tau'\of\bftuple{\tau'_1}{\tau'_n}}
&\rEquals
&\comp_{\tau'} \of (\comp_{\tau'_1} \otimes\cdots\otimes \comp_{\tau'_n})	\\
\end{diagram}
commutes
\item 	\label{propn:double-omega-b}
$\omega_{\tau',\tau''} \of \omega_{\tau,\tau'} = \omega_{\tau,\tau''}$ and
$\omega_{\tau, \tau} = 1$
\item 	\label{propn:double-omega-c}
$\omega_{\nu_n \of \bftuple{\nu_{k_1}}{\nu_{k_n}}, \nu_{k_1 +\cdots+ k_n}}
= \gamma_{k_1, \ldots, k_n}$
and $\omega_{\utree,\nu_1} = \iota$
\end{enumerate}
\end{thm}
\begin{proof}
\bref{propn:double-omega-b} is immediate, and~\bref{propn:double-omega-a}
and~\bref{propn:double-omega-c} follow from
Proposition~\ref{propn:single-omega}. 
\done
\end{proof}

The theorem says that for any pair of $n$-leafed trees $\tau$ and $\tau'$,
there is precisely one map $\comp_\tau \go \comp_{\tau'}$ which can be built up
from $\gamma$ and $\iota$. In short, there is a single canonical isomorphism
between $\comp_\tau$ and $\comp_{\tau'}$: `coherence for an unbiased
bicategory'. 

This is a little different from the usual formulation of bicategorical
coherence, in that we have not directly discussed graph maps $\Bee^{\otimes
n} \go \Bee^{\otimes m}$ (or transformations between them) built up from the
bicategory operations, except in the case $m=1$. This is a feature of the
tree-based (operadic) approach; it seems cleaner and, in any case, what we
have done is enough for our present purpose.

\subsection*{Coherence for unbiased lax functors}

Fix an unbiased lax functor $\pr{F}{\phi}: \Bee\go\Beep$. I will use the same
notation $\gamma$, $\iota$, $\comp$ and $\omega$ in both \Bee\ and \Beep;
confusion should not arise.

Define for each $n\in\nat$ and $\tau\in\tr(n)$ a natural transformation
\[
\phi_\tau: \comp_\tau \of F^{\otimes n} \go F \of \comp_\tau
\]
by
\begin{itemize}
\item $\phi_\utree$ is the identity (or again, more precisely, the canonical
isomorphism)  
\item if $\tau_1\in\tr(k_1), \ldots, \tau_n\in\tr(k_n)$ then
$\phi_{\bftuple{\tau_1}{\tau_n}}$ is the composite
\begin{eqnarray*}
\lefteqn{\comp_{\bftuple{\tau_1}{\tau_n}} \of F^{\otimes (k_1 +\cdots+ k_n)}}
\\
&=&\comp_n \of ((\comp_{\tau_1} \of F^{\otimes k_1}) \otimes\cdots\otimes
(\comp_{\tau_n} \of F^{\otimes k_n}))			\\
&\goby{1*(\phi_{\tau_1} \otimes\cdots\otimes \phi_{\tau_n})}
&\comp_n \of ((F\of \comp_{\tau_1}) \otimes\cdots\otimes 
(F\of \comp_{\tau_n}))					\\
&=
&\comp_n \of F^{\otimes n} \of (\comp_{\tau_1}  \otimes\cdots\otimes
\comp_{\tau_n})						\\
&\goby{\phi_n * 1}
&F\of \comp_n \of (\comp_{\tau_1}  \otimes\cdots\otimes \comp_{\tau_n})	\\
&=
&F\of \comp_{\bftuple{\tau_1}{\tau_n}}.
\end{eqnarray*}
\end{itemize}
Once again we have a coherence result.
\begin{propn}
\begin{enumerate}
\item 	\label{propn:single-phi-a}
If $\tau\in\tr(n), \tau_1\in\tr(k_1), \ldots, \tau_n\in\tr(k_n)$ then
\[
\begin{diagram}
\comp_{\tau\of\bftuple{\tau_1}{\tau_n}} \of F^{\otimes (k_1 +\cdots+ k_n)}
&\rEquals
&\comp_{\tau} \of ((\comp_{\tau_1} \of F^{\otimes k_1}) \otimes\cdots\otimes
(\comp_{\tau_n} \of F^{\otimes k_n})) 	\\
&&\dTo>{1*(\phi_{\tau_1} \otimes\cdots\otimes \phi_{\tau_n})}	\\
\dTo<{\phi_{\tau\of\bftuple{\tau_1}{\tau_n}}}
&&\comp_\tau \of F^{\otimes n} \of 
(\comp_{\tau_1}	\otimes\cdots\otimes \comp_{\tau_n})		\\
&&\dTo>{\phi_\tau * 1}						\\
F\of \comp_{\tau\of\bftuple{\tau_1}{\tau_n}}
&\rEquals
&F\of \comp_\tau \of (\comp_{\tau_1} \otimes\cdots\otimes \comp_{\tau_n})
\end{diagram}
\]
commutes
\item 	\label{propn:single-phi-b}
The diagram
\begin{diagram}
\comp_{\utree} \of F^{\otimes 1}	&\rEquals	&F		\\
\dTo<{\phi_\utree}			&		&\dTo>{1}	\\
F\of \comp_\utree			&\rEquals	&F		\\
\end{diagram}
commutes
\item 	\label{propn:single-phi-c}
$\phi_{\nu_n} = \phi_n$.
\end{enumerate}
\end{propn}
\begin{proof}
\bref{propn:single-phi-a} is by induction on $\tau$;
\bref{propn:single-phi-b} and~\bref{propn:single-phi-c} are immediate. 
\done
\end{proof}

At this point, we have for each $\tau\in\tr(n)$ a canonical map
\[
\phi_\tau: \comp_\tau \of F^{\otimes n} \go F\of \comp_\tau
\]
built up from $\phi_n$'s only. Next we bring in the coherence isomorphisms
$\omega$ of \Bee\ and \Beep.

\begin{propn}	\label{propn:phi-square}
If $\tau, \tau' \in \tr(n)$ then 
\begin{equation}	\label{eqn:omega-phi}
\begin{diagram}
\comp_\tau \of F^{\otimes n}	&\rTo^{\phi_\tau}
&F\of \comp_\tau					\\
\dTo<{\omega_{\tau,\tau'} * 1}	&
&\dTo>{1*\omega_{\tau,\tau'}}				\\
\comp_{\tau'} \of F^{\otimes n}	&\rTo^{\phi_{\tau'}}
&F\of \comp_{\tau'}					\\
\end{diagram}
\end{equation}
commutes.
\end{propn}
\begin{proof}
It is enough to prove this when $\tau' = \nu_n$, in which case
$\omega_{\tau,\tau'} = \omega_\tau$. The proof is then another easy induction
on $\tau$.
\done
\end{proof}

For $\tau, \tau' \in \tr(n)$, define
\[
\phi_{\tau, \tau'}: \comp_\tau \of F^{\otimes n} \go F\of \comp_{\tau'}
\]
as the diagonal of~\bref{eqn:omega-phi}. We then have:
\begin{thm}
\begin{enumerate}
\item		\label{thm:double-phi-a}
If $\tau, \tau', \sigma, \sigma' \in \tr(n)$ then the diagrams
\[
\begin{diagram}
\comp_\tau \of F^{\otimes n}	&\rTo^{\phi_{\tau, \tau'}}
&F\of \comp_{\tau'}						\\
				&\rdTo<{\phi_{\tau,\sigma'}}	
&\dTo>{1 * \omega_{\tau', \sigma'}}				\\
		&		&F\of \comp_{\sigma'}		\\
\end{diagram}
\diagspace
\begin{diagram}
\comp_\tau \of F^{\otimes n}	&\rTo^{\omega_{\tau, \sigma} * 1}
&\comp_{\sigma} \of F^{\otimes n}				\\
				&\rdTo<{\phi_{\tau,\sigma'}}	
&\dTo>{\phi_{\sigma, \sigma'}}					\\
		&		&F\of \comp_{\sigma'}		\\
\end{diagram}
\]
commute
\item		\label{thm:double-phi-b} 
If $\tau, \tau'\in\tr(n), \tau_1, \tau'_1 \in\tr(k_1), \ldots, \tau_n,
\tau'_n \in\tr(k_n)$, then
\[
\!\!\!\!\!\!\!\!
\!\!\!\!\!\!\!\!
\!\!\!\!\!\!\!\!
\!\!\!\!\!\!\!\!
\!\!\!\!\!\!\!\!
\!\!\!\!\!\!
\begin{diagram}[scriptlabels]
\comp_{\tau\of\bftuple{\tau_1}{\tau_n}} \of F^{\otimes (k_1 +\cdots+ k_n)}
&\rEquals
&\comp_{\tau} \of ((\comp_{\tau_1} \of F^{\otimes k_1}) \otimes\cdots\otimes
(\comp_{\tau_n} \of F^{\otimes k_n})) 	\\
&&\dTo>{1*(\phi_{\tau_1, \tau'_1} \otimes\cdots\otimes \phi_{\tau_n,
\tau'_n})}	\\ 
\dTo<{\phi_{\tau\of\bftuple{\tau_1}{\tau_n},
\tau'\of\bftuple{\tau'_1}{\tau'_n}}} 
&&\comp_\tau \of F^{\otimes n} \of 
(\comp_{\tau'_1} \otimes\cdots\otimes \comp_{\tau'_n})		\\
&&\dTo>{\phi_{\tau, \tau'} * 1}					\\
F\of \comp_{\tau'\of\bftuple{\tau'_1}{\tau'_n}}
&\rEquals
&F\of \comp_{\tau'} \of (\comp_{\tau'_1} \otimes\cdots\otimes \comp_{\tau'_n})
\end{diagram}
\]
commutes
\item		\label{thm:double-phi-c}
$\phi_{\utree,\utree}: \comp_\utree \of F^{\otimes 1} \go F\of \comp_\utree$
is the identity
\item		\label{thm:double-phi-d}
$\phi_{\nu_n, \nu_n} = \phi_n$.
\end{enumerate}
\end{thm}
\begin{proof}
These all follow from the last two propositions.
\done
\end{proof}

This theorem is `coherence for an unbiased lax functor' \pr{F}{\phi}: there
is precisely one map
\[
\comp_{\tau} \of F^{\otimes n} \go F \of \comp_{\tau'}
\]
built up from $\phi$ and the coherence cells $\gamma$ and $\iota$ of \Bee\
and \Beep.

A warning is due here. We have shown that, for instance, any two maps 
\[
((Ff_4 \of Ff_3) \of 1 \of (Ff_2 \of Ff_1))
\parpair{}{}
F(f_4 \of (f_3 \of f_2 \of f_1))
\]
built up from coherence cells are equal. The form of the codomain is
important, being $F$ applied to a composite of 1-cells in \Bee. In contrast,
a counterexample in the introduction to~\cite{Lew} shows that there can be
two distinct maps
\[
F1 \parpair{}{} F1 \of F1
\]
built up from coherence cells. (The counterexample is stated in the context
of classical bicategories---in fact, monoidal categories---but translates
easily to the unbiased context.)

\subsection*{Summary}

We have articulated the following coherence principles for the unbiased
theory:
\begin{description}
\item[\UB] In an unbiased bicategory \Bee, there is a unique natural
isomorphism 
\[
\comp_\tau \go \comp_{\tau'}
\]
built up from $\gamma$ and $\iota$, for any pair $\tau, \tau'$ of trees with
the same number of leaves
\item[\UF] For an unbiased lax functor $\pr{F}{\phi}: \Bee\go\Beep$, there is
a unique natural transformation
\[
\comp_{\tau} \of F^{\otimes n} \go F \of \comp_{\tau'}
\]
built up from $\phi$, $\gamma$ and $\iota$, for any pair $\tau$, $\tau'$ of
$n$-leafed trees. 
\end{description}

We will also need to use similar coherence principles for classical
bicategories. To state them, we define the set $\fcat{ctr}(k)$ of
\label{p:classical-trees}%
\emph{$k$-leafed classical trees} for each $k\in\nat$ by exactly the same
recursive clauses as we used in the definition of \tr\
(page~\pageref{p:tree-defn}), but only allowing $n\in\{0,2\}$ instead of
$n\in\nat$ in the second clause. As in the unbiased case, we can define for
each classical bicategory \cat{C}, each $n\in\nat$ and each
$\tau\in\fcat{ctr}(n)$, a functor $\comp_\tau: \cat{C}^{\otimes n} \go
\cat{C}$. We then have the following coherence principles for the classical
theory:
\begin{description}
\item[\CB] In a classical bicategory \cat{C}, there is a unique natural
isomorphism 
\[
\comp_\tau \go \comp_{\tau'}
\]
built up from the associativity and unit isomorphisms, for any pair $\tau,
\tau'$ of classical trees with the same number of leaves
\item[\CF] For a classical lax functor $\pr{G}{\psi}: \cat{C}\go\cat{C'}$,
there is a unique natural transformation
\[
\comp_{\tau} \of G^{\otimes n} \go G \of \comp_{\tau'}
\]
built up from $\psi$ and the associativity and unit isomorphisms, for any
pair $\tau$, $\tau'$ of $n$-leafed classical trees.
\end{description}
Principle \CB\ follows from the classical coherence theorem for bicategories,
in the form `every diagram commutes'. \CF\ comes from Lewis's paper
\cite{Lew}.

\section{The Proof}

We can now prove that $\UBilax \eqv \Bilax$ and $\UBiwk \eqv \Biwk$ with
almost no real work.

Recall from~\ref{sec:versus} that we attempted to construct a functor 
\[
V: \UBilax \go \Bilax;
\]
that is, we specified all the necessary \emph{data} for $V$ but did not check
any of the axioms. Here we must check these axioms, and must prove that $V$
is full, faithful and surjective on objects. The easiest way to deduce the
latter from our results so far is to construct a pseudo-inverse $L$ to $V$,
with $V\of L = 1$.

Explicitly, take a (classical) bicategory \cat{C}, and write its composition
and identity as \Cat-graph maps
\[
\cat{I}_{\cat{C}_0} \goby{\ids} \cat{C} \ogby{\comp} \cat{C} \otimes \cat{C}.
\]
Attempt to define an unbiased bicategory $\Bee = L(\cat{C})$ by setting \Bee\
equal to \cat{C} as a \Cat-graph, putting
\[
\begin{array}{c}
\comp_0	= \ids,						
\diagspace
\comp_1 = 1_{\Bee},
\\
\comp_{n+1} = (\Bee^{\otimes (n+1)} \iso 
\cat{C}^{\otimes n} \otimes \cat{C} \goby{\comp_n \otimes 1_{\cat{C}}}
\cat{C} \otimes \cat{C} \goby{\comp}
\cat{C} = \Bee)
\end{array}
\]
($n\geq 1$), and taking $\gamma$ and $\iota$ to be the canonical isomorphisms
(which exist by coherence principle \CB). (So this choice of a pseudo-inverse
is an arbitrary one; we have decided to `associate to the left'.) Given a
classical lax functor $\pr{G}{\psi}: \cat{C} \go \cat{C'}$, attempt to define
an unbiased lax functor
\[ 
\pr{F}{\phi} = L\pr{G}{\psi}: L(\cat{C}) \go L(\cat{C'})
\]
by setting $F=G$ and taking $\phi_{f_1, \ldots, f_n}$ to be the canonical map
\[
(Ff_n \of\cdots\of Ff_1) \go F(f_n \of\cdots\of f_1),
\]
which makes sense by coherence principle \CF. 

So far we have attempted to construct functors
\[
\UBilax \oppair{V}{L} \Bilax,
\]
and we will show that $VL=1$ and $LV \iso 1$. For the latter we attempt to
construct unbiased weak functors
\[
\Bee \oppair{\pr{\Theta_\Bee}{\theta_\Bee}}{\pr{\Xi_\Bee}{\xi_\Bee}}
LV(\Bee)
\]
for each unbiased bicategory \Bee. This is done by taking $\Theta_\Bee$ and
$\Xi_\Bee$ each to be the identity on \Bee\ (in $\Cat\hyph\Gph$), and by
taking $\theta_\Bee$ and $\xi_\Bee$ to be the canonical isomorphisms (which
exist by coherence principle \UB).

Theorem~\ref{thm:biased-comparison} now follows from:
\begin{propn}
With the definitions above, $V$ and $L$ are both functors, $VL=1$, and 
$1 \oppair{\pr{\Theta}{\theta}}{\pr{\Xi}{\xi}} LV$ are mutually inverse natural
transformations. 
\end{propn}
\begin{proof}
Essentially we have to check that our data satisfies a large collection of
axioms, but our coherence results cover almost all of these checks
automatically. Here is the list of the things to be checked and which
coherence result each one can be inferred from.
\begin{itemize}
\item $V$ is a functor $\UBilax \go \Bilax$. This means:
\begin{itemize}
\item $V(\Bee)$ is a bicategory for any \Bee: \UB
\item $V\pr{F}{\phi}$ is a lax morphism for any \pr{F}{\phi}: \UF
\item $V$ preserves identities: \UB
\item $V$ preserves composition: really we should deduce this from `coherence
for a composable pair of unbiased lax morphisms' (which we did not prove),
but a direct check is easy.
\end{itemize}
\item $L$ is a functor $\Bilax \go \UBilax$. This means:
\begin{itemize}
\item $L(\cat{C})$ is an unbiased bicategory for any \cat{C}: \CB
\item $L\pr{G}{\psi}$ is an unbiased lax functor for any \pr{G}{\psi}: \CF
\item $L$ preserves identities: \CB
\item $L$ preserves composition: as for $V$ above.
\end{itemize}
\item $VL=1$. This means:
\begin{itemize}
\item $VL(\cat{C}) = \cat{C}$ for any \cat{C}: by construction, $VL(\cat{C})$
and \cat{C} are the same in all respects except perhaps their associativity
and unit isomorphisms; and these too are equal by \CB
\item $VL\pr{G}{\psi} = \pr{G}{\psi}$ for any \pr{G}{\psi}: \CF.
\end{itemize}
\item $1 \oppair{\pr{\Theta}{\theta}}{\pr{\Xi}{\xi}} LV$ are natural
transformations. This means:
\begin{itemize}
\item $\Bee \oppair{\pr{\Theta_\Bee}{\theta_\Bee}}{\pr{\Xi_\Bee}{\xi_\Bee}}
LV(\Bee)$ are unbiased lax functors for any \Bee: \UB
\item \pr{\Theta_\Bee}{\theta_\Bee} and \pr{\Xi_\Bee}{\xi_\Bee} are natural
in \Bee: \UF.
\end{itemize}
\item $\pr{\Theta_\Bee}{\theta_\Bee} \of \pr{\Xi_\Bee}{\xi_\Bee} = 1$ and 
$\pr{\Xi_\Bee}{\xi_\Bee} \of \pr{\Theta_\Bee}{\theta_\Bee} = 1$ for any \Bee:
\UB. 
\done
\end{itemize}
\end{proof}

Evidently $L$ sends weak functors to unbiased weak functors, and so restricts
to a functor $L_\mr{wk}: \Biwk \go \UBiwk$. Moreover, both
\pr{\Theta_\Bee}{\theta_\Bee} and \pr{\Xi_\Bee}{\xi_\Bee} are unbiased
\emph{weak} functors, for any unbiased bicategory \Bee. Hence:
\begin{cor}
The functors $\UBiwk \oppair{V_\mr{wk}}{L_\mr{wk}} \Biwk$ satisfy
$V_\mr{wk} L_\mr{wk}=1$ and $L_\mr{wk} V_\mr{wk} \iso 1$.
\done
\end{cor}

Corollary~\ref{cor:wk-biased-comparison} follows immediately.

\chapter{The Free Multicategory Construction}
\label{app:free}

In this appendix we define `suitability' and sketch proofs of Theorems
\ref{thm:free-main}, \ref{thm:free-fixed} and \ref{thm:free-gen}. First we
need some terminology.

Let \Eee\ be a category with pullbacks, \scat{I} a small category, $D:
\scat{I}\go\Eee$ a functor for which a colimit exists, and $(D(I)\go
Z)_{I\in \scat{I}}$ a colimit cone. We say that the colimit is
\emph{stable under pullback} if for any map $Z'\go Z$ in \Eee, the cone
$(D'(I)\go Z')_{I\in \scat{I}}$ is a colimit cone; here $D'$ and the new cone
are obtained by pullback, so that
\begin{diagram}[size=1.5em]
D'	&\rTo	&D	\\
\dTo	&	&\dTo	\\
Z'	&\rTo	&Z	\\
\end{diagram}
is a pullback square in the functor category \ftrcat{\scat{I}}{\Eee}.

The morphisms $k_I$ in a colimit cone $(D(I)\goby{k_I} Z)_{I\in \scat{I}}$
will be called the \emph{coprojections} of the colimit, and in particular we
say that the colimit of $D$ `has monic coprojections' to mean that each $k_I$
is monic.

A category will be said to have \emph{disjoint finite coproducts} if it has
finite coproducts, these coproducts have monic coprojections, and for any
pair $A, B$ of objects, the square
\begin{diagram}[size=1.5em]
0	&\rTo	&B	\\
\dTo	&	&\dTo	\\
A	&\rTo	&A+B	\\
\end{diagram}
is a pullback.

Let $\omega$ be the natural numbers with their usual ordering. A \emph{nested
sequence} in a category \Eee\ is a functor $\omega\go\Eee$ in which the image
of every morphism of $\omega$ is monic. In other words, it is a diagram
\[
A_0 \monic A_1 \monic \cdots
\]
in \Eee, where as usual \monic\ indicates a monic. Note that a functor which
preserves pullbacks also preserves monics, so it makes sense for such a
functor to `preserve colimits of nested sequences'. Similarly, it makes sense
to say that colimits of nested sequences commute with pullbacks, where
`commute' is used in the same sense as when we say that filtered colimits
commute with finite limits in \Set.

A category \Eee\ is \emph{suitable} if it satisfies
\begin{description}
\item[C1] \Eee\ is cartesian
\item[C2] \Eee\ has disjoint finite coproducts which are stable under
pullback
\item[C3] \Eee\ has colimits of nested sequences; these commute with
pullbacks and have monic coprojections.
\end{description}
A monad \Mnd\ is \emph{suitable} if it satisfies
\begin{description}
\item[M1] \Mnd\ is cartesian
\item[M2] $T$ preserves colimits of nested sequences.
\end{description}
We say that \Cartpr\ is \emph{suitable} when \Mnd\ is a suitable monad on a
suitable category \Eee.

We now sketch a proof of the main theorem,~\ref{thm:free-main}, on the
formation of free multicategories, which for convenience is re-stated here.
\begin{trivlist} \item 
\textbf{Theorem \ref{thm:free-main}}\ \itshape
Let \Cartpr\ be suitable. Then the forgetful functor
\[
\Cartpr\hyph\Multicat \goby{U} \Eeep = \Cartpr\hyph\Graph
\]
has a left adjoint, the adjunction is monadic, and if $T'$ is the resulting
monad on \Eeep\ then \Cartprp\ is suitable.
\end{trivlist}

\noindent\textbf{Sketch proof}
We proceed in four steps:
\begin{enumerate}
\item 	\label{ftr}
construct a functor $F: \Eeep \go \Cartpr\hyph\Multicat$
\item 	\label{adjn}
construct an adjunction between $F$ and $U$
\item 	\label{primesuit}
check that \Cartprp\ is suitable
\item 	\label{monadic}
check that the adjunction is monadic.
\end{enumerate}

Each step goes roughly as follows.
\begin{enumerate}
\item \emph{Construct a functor $F: \Eeep \go \Cartpr\hyph\Multicat$}\\ 
Let $X$ be a $T$-graph. Define for each $n$ a graph
\spaan{A_n}{TX_0}{X_0}{d_n}{c_n}, by
\begin{itemize}
\item $A_0=X_0$, $d_0=\eta_{X_0}$ and $c_0=1$
\item 
\begin{sloppypar}
$A_{n+1} = X_0 + X_1\of A_n$, where $X_1\of A_n$ is the 1-cell
composite in $\Cartpr\hyph\Span$, with the obvious choices of $d_{n+1}$ and
$c_{n+1}$.
\end{sloppypar}
\end{itemize}
Define for each $n$ a map $A_n \goby{i_n} A_{n+1}$, by
\begin{itemize}
\item $i_0: X_0 \go X_0 + X_1 \of X_0$ is first coprojection
\item $i_{n+1} = 1_{X_0} + (1_{X_1} * i_n)$.
\end{itemize}
Then the $i_n$'s are monic, and by taking $A$ to be the colimit of 
\[
A_0 \rMonic^{i_0} A_1 \rMonic^{i_1} \cdots
\]
we obtain a graph \spn{A}{TX_0}{X_0}. This graph naturally has the structure
of a multicategory: the identities map $X_0 \go A$ is just the colimit
coprojection $A_0 \monic A$, and composition comes from maps \linebreak$A_m
\of A_n \go A_{m+n}$ which piece together to give a map $A\of A \go A$. The
composition construction needs many of the suitability axioms.

We have now described what effect $F$ is to have on objects, and extension to
morphisms is straightforward.

(The colimit of the nested sequence of $A_n$'s appears, in light disguise, as
the recursive description of the free plain multicategory monad in
section~\ref{sec:free-multicats}: $A_n$ is the set of formal expressions
which can be obtained from the first clause and up to $n$ applications of the
second clause.)

\item \emph{Construct an adjunction between $F$ and $U$}\\ 
We do this by constructing unit and counit transformations and verifying the
triangle identities. Both transformations are the identity on the object of
objects, so we only need to define them on the object of arrows. For the unit
$\eta'$, if $X\in\Eeep$ then $\eta'_X: X_1 \go A$ is the composite
\[
X_1 \goiso X_1 \of X_0 \rMonic^{\mr{copr}_2} X_0 + X_1 \of X_0 = A_1
\rMonic A.
\]
For the counit $\epsln'$, let $C\in\Cartpr\hyph\Multicat$. Write $A$ and
$A_n$ for the objects used in the construction of the free multicategory on
$U(C)$, as if $X=U(C)$ in part~\bref{ftr}. Define for each $n$ a map
$\epsln'_{C,n}: A_n \go C_1$ by
\begin{itemize}
\item $\epsln'_{C,0} = (A_0 \goby{=} C_0 \goby{\ids} C_1)$
\item $\epsln'_{C,n+1} = 
(C_0 + C_1 \of A_n 
\goby{1+1*\epsln'_{C,n}} 
C_0 + C_1 \of C_1 
\goby{q} 
C_1)$, where $q$ is \ids\ on the first summand and \comp\ on the second,
\end{itemize}
and then there is a unique $\epsln'_C: A \go C_1$ such that 
\[
\epsln'_{C,n} = (A_n\monic A \goby{\epsln'_C} C_1)
\]
for all $n$. 

\item \emph{Check that \Cartprp\ is suitable}\\
This is quite routine.

\item \emph{Check that the adjunction is monadic}\\
We apply the Monadicity Theorem by checking that $U$ creates coequalizers for
$U$-absolute coequalizer pairs. This can be done quite separately from the
rest of the proof, and again is quite routine.
\done 
\end{enumerate}

We can now easily prove the fixed-object version,
Theorem~\ref{thm:free-fixed}. Recall that $\Eee'_S$ is the category of
$T$-graphs on $S$ (that is, $\Eee/(TS\times S)$) and
$\Cartpr\hyph\Multicat_S$ is the category of $T$-multicategories on $S$.
\begin{trivlist} \item 
\textbf{Theorem \ref{thm:free-fixed}}\ \itshape
Let \Cartpr\ be suitable and let $S\in\Eee$. Then the forgetful functor
\[
\Cartpr\hyph\Multicat_S \go \Eee'_S
\]
has a left adjoint, the adjunction is monadic, and if $T'_S$ is the resulting
monad on $\Eee'_S$ then \pr{\Eee'_S}{T'_S} is suitable. Moreover, if \Eee\
has filtered colimits and $T$ preserves them, then the same is true of
$\Eee'_S$ and $T'_S$. 
\end{trivlist}

\paragraph*{Sketch proof} 
\begin{sloppypar}
It is evident from the proof of~\ref{thm:free-main} that the adjunction
$(F,U,\eta',\mu')$ constructed there restricts to the subcategories $\Eee'_S$
and $\Cartpr\hyph\Multicat_S$, so we only have to check that
\pr{\Eee'_S}{T'_S} is suitable and the restricted adjunction is monadic. This
is again quite routine, and involves many of the same calculations. (The most
substantial difference between the two cases is that coproducts in \Eeep\ and
$\Eee'_S$ are calculated differently, i.e.\ the inclusion $\Eee'_S \rIncl
\Eeep$ does not preserve them.) `Moreover' is straightforward.  \done
\end{sloppypar}

\paragraph*{}

Finally, we have to prove Theorem~\ref{thm:free-gen}: that any functor
category \ftrcat{\scat{E}}{\Set}, and any finitary cartesian monad on it, is
suitable. Since the category $\omega$ is filtered, the monad part is
immediate. For the category part it is enough to see that $\Set$ is suitable,
and this follows straight away from standard results.

\chapter{Strict $\omega$-Categories}	\label{app:strict-omega}

In this appendix we prove:
\begin{thm}	\label{thm:strict-omega}
The forgetful functor $\omega\hyph\Cat \go \cat{G}_\omega$ has a left adjoint,
the adjunction is monadic, and the induced monad on $\cat{G}_\omega$ is
cartesian and finitary. 
\end{thm}
Here $\omega\hyph\Cat$ is the category of strict $\omega$-categories and
$\cat{G}_\omega$ is the category of globular sets. We need to know that the
left adjoint exists and that the induced monad $T$ is cartesian in order to
be able to talk about $T$-operads (as we do in Chapter~\ref{ch:defn}), we
need monadicity in order to understand the definition of weak
$\omega$-category (Chapter~\ref{ch:defn} again), and we need to know that $T$
is finitary in Appendix~\ref{app:initial}.

In section~\ref{sec:outline} we recall the basics of strict
$\omega$-categories and strict $n$-categories, and outline the strategy for
proving Theorem~\ref{thm:strict-omega}. Section~\ref{sec:proof} is devoted to
the proof itself. In section~\ref{sec:pasting-rep} we show that $T$ acts on
globular sets $X$ by the formula
\[
(TX)(n) \iso \coprod_{\pi\in\pd(n)} 
\homset{\ftrcat{\scat{G}}{\Set}}{\rep{\pi}}{X},
\]
as asserted in Chapter~\ref{ch:defn} (page~\pageref{p:pasting-rep}).

\section{Outline}	\label{sec:outline}

Let $\cat{V}$ be a category with finite products. Then there is a category
$\cat{V}\hyph\Cat$ of $\cat{V}$-enriched categories and $\cat{V}$-enriched
functors, which also has finite products. Moreover, if $F: \cat{V} \go
\cat{W}$ is a functor which preserves finite products then there is an
induced functor $ F_*: \cat{V}\hyph\Cat \go \cat{W}\hyph\Cat, $ which also
preserves finite products. Here, as everywhere in this chapter, the monoidal
structure on the categories we are enriching in is always the cartesian
product, and our enriched categories will always have just a \emph{set} of
objects---nothing larger. 

These observations allow us to make the following definitions. For
$n\in\nat$, define the category $n\hyph\Cat$ of \emph{strict $n$-categories}
and \emph{strict $n$-functors} by
\begin{eqnarray*}
0\hyph\Cat 	&= 	&\Set,			\\
(n+1)\hyph\Cat 	&= 	&(n\hyph\Cat)\hyph\Cat.
\end{eqnarray*}
Also define functors $S_n: (n+1)\hyph\Cat \go n\hyph\Cat$, by taking $S_0:
\Cat \go \Set$ to be the objects functor and $S_{n+1} = (S_n)_*$.

We thus have a diagram
\[
\cdots \go (n+1)\hyph\Cat \goby{S_n} n\hyph\Cat \goby{S_{n-1}} \cdots 
\goby{S_0} 0\hyph\Cat = \Set
\]
in \fcat{CAT}, and the category $\omega\hyph\Cat$ of \emph{strict
$\omega$-categories} and \emph{strict $\omega$-functors} is defined as the
limit of this diagram. (\fcat{CAT} is the category of all categories,
possibly large.) 

Now let $\scat{G}_\omega$ be the category denoted \scat{G} in
section~\ref{sec:formal}, and let $\cat{G}_\omega =
\ftrcat{\scat{G}_\omega}{\Set}$ (the category of globular sets). For
$n\in\nat$, let $\scat{G}_n$ be the full subcategory of $\scat{G}_\omega$
with objects $0, \ldots, n$, let $\cat{G}_n = \ftrcat{\scat{G}_n}{\Set}$, and
call objects of $\cat{G}_n$ \emph{$n$-globular sets}. The inclusions
$\scat{G}_n \rIncl \scat{G}_{n+1}$ give rise to a diagram
\[
\cdots \go \cat{G}_{n+1} \goby{R_n} \cat{G}_n \goby{R_{n-1}} \cdots 
\goby{R_0} \cat{G}_0 \iso \Set
\]
in \fcat{CAT}, of which $\cat{G}_\omega$ is the limit.

The next step is to see that there is a forgetful functor
\[
U_n: n\hyph\Cat \go \cat{G}_n
\]
for each $n$, expressing the idea that an $n$-globular set is the underlying
graph structure of an $n$-category. 

Formally, we first define for each category \cat{V} the category
$\cat{V}\hyph\Gph$, in which an object is a set $X_0$ together with an
indexed family $(X(x,x'))_{x,x'\in X_0}$ of objects of \cat{V}, and a map $f:
X \go Y$ consists of a function $f_0: X_0 \go Y_0$ together with a map
$f_{x,x'}: X(x,x') \go Y(f_0 x, f_0 x')$ in \cat{V} for each $x,x' \in
X_0$. Objects of $\cat{V}\hyph\Gph$ will be called \emph{\cat{V}-graphs} (not to be
confused with the $T$-graphs defined in~\ref{sec:multi}). Observe that:
\begin{itemize}
\item if \cat{V} has finite products then so does $\cat{V}\hyph\Gph$, and the
evident forgetful functor $\cat{V}\hyph\Cat \go \cat{V}\hyph\Gph$ preserves
finite products
\item if \cat{V} and \cat{W} have finite products and $\cat{V} \go \cat{W}$
is a functor preserving them, then the evident functor $\cat{V}\hyph\Gph \go
\cat{W}\hyph\Gph$ also preserves them
\item in the situation of the previous item, the diagram
\begin{diagram}
\cat{V}\hyph\Cat	&\rTo	&\cat{W}\hyph\Cat	\\
\dTo			&	&\dTo			\\
\cat{V}\hyph\Gph	&\rTo	&\cat{W}\hyph\Gph	\\
\end{diagram}
commutes, which means that there is an unambiguous functor
$\cat{V}\hyph\Cat\go $\linebreak$\cat{W}\hyph\Gph$ induced by the functor
$\cat{V} \go \cat{W}$.
\end{itemize}

To apply this to the current situation, note that $\cat{G}_{n+1} \eqv
\cat{G}_n \hyph\Gph$; then take $U_0: \Set \go \Set$ to be the identity and
define $U_{n+1}: (n+1)\hyph\Cat \go \cat{G}_{n+1}$ to be the
functor
\[
(n\hyph\Cat) \hyph\Cat \go \cat{G}_n \hyph\Gph
\]
induced by $U_n: n\hyph\Cat \go \cat{G}_n$. (All the conditions on finite
products go through.) These $U_n$'s commute with the restriction functors
$R_n$ and $S_n$, so we obtain a forgetful functor $U_\omega: \omega\hyph\Cat
\go \cat{G}_\omega$:
\[
% \ \ \  
\begin{diagram}
\omega\hyph\Cat	&	&\cdots	&\rTo	&(n+1)\hyph\Cat	&\rTo^{S_n}	
&n\hyph\Cat	&\rTo^{S_{n-1}}	&\cdots	&\rTo^{S_0}	&0\hyph\Cat	\\
\dTo>{U_\omega}&	&	&	&\dTo>{U_{n+1}}	&		
&\dTo>{U_n}	&		&	&		&\dTo>{U_0}	\\
\cat{G}_\omega	&	&\cdots	&\rTo	&\cat{G}_{n+1}	&\rTo_{R_n}	
&\cat{G}_n	&\rTo_{R_{n-1}}	&\cdots	&\rTo_{R_0}	&\cat{G}_0.	\\
\end{diagram}
\]

Having constructed $U_\omega$, we have given a precise meaning to
Theorem~\ref{thm:strict-omega}. (`Finitary' means `preserves filtered
colimits'.) In order to prove the Theorem, it is enough to prove:
\begin{thm}	\label{thm:strict-n}
Let $n\in\nat$. Then
\begin{enumerate}
\item 	\label{part:strict-n-objs}
the forgetful functor $U_n: n\hyph\Cat \go \cat{G}_n$ has a left
adjoint $F_n$, the adjunction is monadic, and the induced monad $T_n$ on
$\cat{G}_n$ is cartesian and finitary
\item \label{part:strict-n-arrs} $R_n$ is a weak map of monads
$\pr{\cat{G}_{n+1}}{T_{n+1}} \go \pr{\cat{G}_n}{T_n}$, and $S_n$ is the map
$\cat{G}_{n+1}^{T_{n+1}} \go \cat{G}_n^{T_n}$ induced by $R_n$.
\end{enumerate}
\end{thm}
By a \emph{weak map of monads} I mean a monad functor (or equivalently,
opfunctor) whose natural transformation part is an isomorphism: thus there is
an isomorphism between $T_n \of R_n$ and $R_n \of T_{n+1}$ which respects
unit and multiplication. Because it is a monad functor, there is an induced
functor $\cat{G}_{n+1}^{T_{n+1}} \go \cat{G}_n^{T_n}$, and therefore
$(n+1)\hyph\Cat \go n\hyph\Cat$. 

Theorem~\ref{thm:strict-omega} follows almost immediately from
Theorem~\ref{thm:strict-n}. The only sticking point is that the squares
\begin{equation}	\label{eqn:FRS}
\begin{diagram}
(n+1)\hyph\Cat	&\rTo^{S_n}	&n\hyph\Cat	\\
\uTo<{F_{n+1}}	&		&\uTo>{F_n}	\\
\cat{G}_{n+1}	&\rTo_{R_n}	&\cat{G}_n	\\
\end{diagram}
\end{equation}
do not \emph{a priori} commute strictly, only up to (canonical)
isomorphism. Since $\cat{G}_\omega$ and $\omega\hyph\Cat$ are strict (not
$2$-categorical) limits, this means that the functors $F_n$ do not
necessarily induce a functor $F_\omega: \cat{G}_\omega \go
\omega\hyph\Cat$. But we can, in fact, choose the left adjoints $F_n$ so that
each canonical isomorphism inside~\bref{eqn:FRS} is the identity, and the
situation is then rescued. The key is that the functors $S_n$ have the
following (easily proved) isomorphism-lifting property: if $C\in
(n+1)\hyph\Cat$ and $j: S_n(C) \goiso D$ is an isomorphism in $n\hyph\Cat$,
then there is an isomorphism $i: C \goiso C'$ in $(n+1)\hyph\Cat$ with $S_n
C' = D$ and $S_n i = j$. This allows us to choose left adjoints $F_0$, $F_1$,
\ldots\ successively so that everything is strictly commutative, which is just
what we need.

\section{The Proof}	\label{sec:proof}

In this section we prove Theorem~\ref{thm:strict-n}. The core of the argument
is contained in the following result:
\begin{propn}	\label{propn:core}
Let \scat{A} be a small category and $\cat{A} = \ftrcat{\scat{A}}{\Set}$. Let
\Mnd\ be a monad on \cat{A} such that $T$ preserves all coproducts. Then
\begin{enumerate}
\item 	\label{part:core-basic}
the forgetful functor $\cat{A}^T \hyph\Cat \go \cat{A}\hyph\Gph$ is
monadic and preserves all coproducts
\item  	\label{part:core-cart}
if \Mnd\ is cartesian then so is the induced monad
\triple{\twid{T}}{\twid{\eta}}{\twid{\mu}} on $\cat{A}\hyph\Gph$
\item  	\label{part:core-fin}
if $T$ is finitary then so is \twid{T}.
\end{enumerate}
\end{propn}
\paragraph*{Remarks} The `forgetful functor' in the first part is induced by
the forgetful functor $\cat{A}^T \go \cat{A}$, where $\cat{A}^T$ is the
category of $T$-algebras. Since \cat{A} has all limits
and a monadic functor creates limits, $\cat{A}^T$ has all limits---and in
particular pullbacks, so that it makes sense to discuss $\cat{A}^T
\hyph\Cat$.

Parts~\bref{part:core-cart} and~\bref{part:core-fin} make sense even if
$\cat{A}\hyph\Gph$ does not have all pullbacks or all filtered colimits. But
in fact, $\cat{A}\hyph\Gph$ has all limits and colimits. This follows from
the observation that $\cat{A}\hyph\Gph \eqv \ftrcat{\twid{\scat{A}}}{\Set}$,
where \twid{\scat{A}} is the category obtained from \scat{A} by adjoining a
new object $0$ and a pair of morphisms $A \parpair{\sigma_A}{\tau_A} 0$ for
each $A\in \scat{A}$, with $\sigma_A \of f = \sigma_{A'}$ and $\tau_A \of f
= \tau_{A'}$ for any morphism $f: A' \go A$ in \scat{A}.

It is not necessary to insist that \cat{A} is of the form
\ftrcat{\scat{A}}{\Set} in order to make the proof work. We could get by on
the assumption that \cat{A} has finite limits and all (small) colimits, and
that these interact in suitable ways: e.g.\ that $\times$ distributes over
coproduct. But we do not need such a precise result, and by working in
\ftrcat{\scat{A}}{\Set} we can manipulate limits and colimits as if we were
in \Set. 

\paragraph*{} Before proving the Proposition, let us apply it to prove
part~\bref{part:strict-n-objs} of Theorem~\ref{thm:strict-n}. The proof is by
induction on $n$, adding in the hypothesis that the functor $T_n$ preserves
all (small) coproducts. When $n=0$, the forgetful functor $U_0$ is an
isomorphism, and its inverse $F_0$ is a left adjoint; thus the induced monad
$T_0$ on $\cat{G}_0$ is the identity. For the inductive step we just take
$\scat{A} = \scat{G}_n$ and $T=T_n$ in Proposition~\ref{propn:core}, noting
that under the equivalences $\cat{A}^T \hyph\Cat \eqv (n+1)\hyph\Cat$ and
$\cat{A}\hyph\Gph \eqv \cat{G}_{n+1}$, the forgetful functor $\cat{A}^T
\hyph\Cat \go \cat{A}\hyph\Gph$ becomes $U_{n+1}: (n+1)\hyph\Cat \go
\cat{G}_{n+1}$. 

\paragraph*{Proof of Proposition~\ref{propn:core}} The strategy is to
construct two monads $P$ and $Q$ on $\cat{A}\hyph\Gph$ and a distributive
law $Q\of P \go P\of Q$ (in the sense of \cite[\S 6]{StrFTM}). This gives the
functor $\twid{T} = Q\of P$ the structure of a monad on
$\cat{A}\hyph\Gph$. We then show that $(\cat{A}\hyph\Gph)^{\twid{T}} \iso
\cat{A}^T \hyph\Cat$, and that the diagram
\begin{diagram}[height=1.5em]
(\cat{A}\hyph\Gph)^{\twid{T}}	&\iso		&\cat{A}^T \hyph\Cat	\\
				&\rdTo(1,2)\ldTo(1,2)	&		\\
				&\cat{A\hyph\Gph}	&		\\
\end{diagram}
(in which the two arrows are the forgetful functors)
commutes. Part~\bref{part:core-basic} follows, and by our construction of
\twid{T}, \bref{part:core-cart} and \bref{part:core-fin} are easy
consequences. The idea behind this strategy is that to form the free
$\cat{A}^T$-category on an \cat{A}-graph $X$, one first forms the free
$T$-algebra on each `hom-object' $X(x,x')$, then one forms the free
$\cat{A}^T$-category on the resulting $\cat{A}^T$-graph.

So, the functor $T: \cat{A} \go \cat{A}$ induces a functor $P:
\cat{A}\hyph\Gph \go \cat{A}\hyph\Gph$: explicitly, $(PX)_0 = X_0$
and $(PX)(x,x') = T(X(x,x'))$ for $x, x' \in X_0$. Similarly, the unit and
multiplication of $T$ give $P$ the structure of a monad on
$\cat{A}\hyph\Gph$.

A second monad $Q$ on $\cat{A}\hyph\Gph$ is given by the forgetful functor
\[
\cat{A}\hyph\Cat \go \cat{A}\hyph\Gph\]
and its left adjoint. Explicitly, if
$X\in \cat{A}\hyph\Gph$ then $(QX)_0 = X_0$ and 
\[
(QX)(x,x') = 
\coprod_{x=x_0, \ldots, x_r=x'}
X(x_0, x_1) \times\cdots\times X(x_{r-1}, x_r),
\]
where the coproduct is over all $r\in\nat$ and $x_0, \ldots, x_r \in X_0$
such that $x_0=x$ and $x_r = x'$.  Everything works in the familiar
way---that is, as for the free category monad on $\cat{G}_1$---because
\cat{A} is a functor category \ftrcat{\scat{A}}{\Set}.

A distributive law $\lambda: PQ \go QP$ is given as follows. If $X\in
\cat{A}\hyph\Gph$ then
\begin{eqnarray*}
(PQX)_0		&= 	&X_0,			\\
(PQX)(x,x')	&=	&T((QX)(x,x'))		\\
		&\iso	
	&\coprod_{x=x_0, \ldots, x_r=x'}
	T \{ X(x_0, x_1) \times\cdots\times X(x_{r-1}, x_r) \}
\end{eqnarray*}
(since $T$ preserves coproducts), and 
\begin{eqnarray*}
(QPX)_0		&= 	&X_0,			\\
(QPX)(x,x') 	&= 	&\coprod_{x=x_0, \ldots, x_r=x'}
		T(X(x_0, x_1)) \times\cdots\times T(X(x_{r-1}, x_r)),
\end{eqnarray*}
for $x,x' \in X_0$. So there is a map
\[
(PQX)(x,x') \go (QPX)(x,x')
\]
defined by projections, giving a map
\[
\lambda_X: PQX \go QPX
\]
of \cat{A}-graphs which is the identity on objects. The axioms for a
distributive law then hold.

$P$, $Q$ and $\lambda$ together define a monad
\triple{\twid{T}}{\twid{\eta}}{\twid{\mu}} on $\cat{A}\hyph\Gph$, where
$\twid{T} = Q\of P$ (again, see~\cite[\S 6]{StrFTM}). A \twid{T}-algebra is
an \cat{A}-graph $X$ equipped with a $P$-algebra structure $h$ and a
$Q$-algebra structure $k$ such that
\begin{diagram}[height=2em]
PQX		&\rTo^{Pk}	&PX	\\
\dTo<{\lambda_X}&		&	\\
QPX		&		&\dTo>{h}\\
\dTo<{Qh}	&		&	\\
QX		&\rTo_{k}	&X
\end{diagram}
commutes. In other words, it is an \cat{A}-graph $X$ together with a
$T$-algebra structure $h_{x,x'}$ on $X(x,x')$ for each $x,x' \in X_0$, and an
\cat{A}-category structure
\begin{diagram}
X(x_0, x_1) \times\cdots\times X(x_{r-1}, x_r)
&\rTo^{k_{x_0, \ldots, x_r}}
&X(x_0, x_r)
\end{diagram}
($x_i \in X_0$) on $X$, such that for all $x_0, \ldots, x_r \in X_0$,
\begin{diagram}[height=2em]
T \{ X(x_0, x_1) \times\cdots\times X(x_{r-1}, x_r) \}
&\rTo^{T(k_{x_0, \ldots, x_r})}	
&T(X(x_0, x_r))	\\
\dTo<{(T(\mr{pr}_1), \ldots, T(\mr{pr}_r))}
&		&		\\
T(X(x_0, x_1)) \times\cdots\times T(X(x_{r-1}, x_r))
&		&\dTo>{h_{x_0,x_r}}\\
\dTo<{h_{x_0,x_1} \times\cdots\times h_{x_{r-1},x_r}}	
&		&		\\
X(x_0, x_1) \times\cdots\times X(x_{r-1}, x_r)
&\rTo_{k_{x_0, \ldots, x_r}}	
&X(x_0, x_r)	\\
\end{diagram}
commutes. But the left-hand column of this diagram is the product in
$\cat{A}^T$ of the $T$-algebras $X(x_0, x_1), \ldots, X(x_{r-1}, x_r)$
(recalling the way in which a monadic functor creates limits): so a
\twid{T}-algebra is exactly a category enriched in $\cat{A}^T$, and
$(\cat{A}\hyph\Gph)^{\twid{T}} \iso \cat{A}^T \hyph\Cat$.

It is easy to see that the diagram of forgetful functors in the first
paragraph of the proof commutes, so the forgetful functor $\cat{A}^T \hyph\Cat
\go \cat{A}\hyph\Gph$ is monadic. Moreover, $P$ preserves coproducts since
$T$ does, and $Q$ evidently preserves coproducts, so the functor $\twid{T} =
Q\of P$ does too. This completes the proof of~\bref{part:core-basic}.

For~\bref{part:core-cart} and~\bref{part:core-fin}, note that $P$ is
cartesian (respectively, finitary) if $T$ is, and that $Q$ is cartesian and
finitary in any case. It only remains to prove that if the monad $T$ is
cartesian then the natural transformation $\lambda: PQ \go QP$ is also
cartesian, and this is straightforward.
\done
\paragraph*{}

The proof of Theorem~\ref{thm:strict-n}\bref{part:strict-n-objs} is now
done. For part~\bref{part:strict-n-arrs} we use the following result:
\begin{propn}	\label{propn:core-map}
Let $J: \scat{A'} \go \scat{A}$ be a functor between small categories, let
\triple{T'}{\eta'}{\mu'} be a monad on $\cat{A'} = \ftrcat{\scat{A'}}{\Set}$
such that $T'$ preserves all coproducts, and similarly \triple{T}{\eta}{\mu}
on $\cat{A} = \ftrcat{\scat{A}}{\Set}$. If $J^*$ is a
weak map of monads
\begin{equation}		\label{eqn:first-mm}
\pr{\cat{A}}{T} \go \pr{\cat{A'}}{T'},
\end{equation}
then the induced functor $J^*\hyph\Gph: \cat{A}\hyph\Gph \go
\cat{A'}\hyph\Gph$ also becomes a weak map of monads
\begin{equation}		\label{eqn:second-mm}
\pr{\cat{A}\hyph\Gph}{\twid{T}} \go \pr{\cat{A'}\hyph\Gph}{\twid{T'}},
\end{equation}
where \twid{T} and \twid{T'} are as in
Proposition~\ref{propn:core}. Moreover, the diagram
\begin{diagram}[height=1.5em]
\cat{A}^T \hyph\Cat		&\rTo	&\cat{A'}^{T'} \hyph\Cat	\\
\dTo<{\iso}			&	&\dTo>{\iso}			\\
(\cat{A}\hyph\Gph)^{\twid{T}}	&\rTo	&(\cat{A'}\hyph\Gph)^{\twid{T'}}\\
\end{diagram}
commutes, where the map along the top is induced by the monad
map~\bref{eqn:first-mm} and the map along the bottom by the monad
map~\bref{eqn:second-mm}. 
\end{propn}

\begin{proof}
Consider the diagram
\begin{diagram}[height=2em]
\cat{A}\hyph\Gph	&\rTo^{J^*\hyph\Gph}	&\cat{A'}\hyph\Gph	\\
\dTo<{P}		&			&\dTo>{P'}		\\
\cat{A}\hyph\Gph	&\rTo^{J^*\hyph\Gph}	&\cat{A'}\hyph\Gph	\\
\dTo<{Q}		&			&\dTo>{Q'}		\\
\cat{A}\hyph\Gph	&\rTo^{J^*\hyph\Gph}	&\cat{A'}\hyph\Gph,	\\
\end{diagram}
where $P$ and $Q$ are as in the proof of Proposition~\ref{propn:core}, and
similarly $P'$ and $Q'$. Applying $\blank\hyph\Gph$ to the isomorphism $T'\of
J^* \iso J^* \of T$ gives an isomorphism `inside' the upper square, making
$J^*\hyph\Gph$ into a weak map of monads
\[
\pr{\cat{A}\hyph\Gph}{P} \go \pr{\cat{A'}\hyph\Gph}{P'}.
\]
There is also a natural isomorphism inside the lower square, expressing the
fact that the free enriched category construction is natural in a suitable
sense, and this gives a weak map of monads
\[
\pr{\cat{A}\hyph\Gph}{Q} \go \pr{\cat{A'}\hyph\Gph}{Q'}.
\]
(The checks involved here use the fact that $J^*: \cat{A'} \go \cat{A}$ is
induced by $J: \scat{A} \go \scat{A'}$; again, this is an unnecessarily
strong hypothesis, but serves our purpose.) Gluing together these two weak
maps of monads gives a third weak map of monads,
\[
\pr{\cat{A}\hyph\Gph}{\twid{T}} \go \pr{\cat{A'}\hyph\Gph}{\twid{T'}},
\]
as required. One can easily check that the diagram in the last sentence of
the Proposition commutes. 
\done
\end{proof}

Theorem~\ref{thm:strict-n}\bref{part:strict-n-arrs} can now be proved by
induction on $n$. 

For the base step, take the monads $T_0$ on $\cat{G}_0 =
\ftrcat{\scat{G}_0}{\Set} \iso \Set$ and $T_1$ on $\cat{G}_1 =
\ftrcat{\scat{G}_1}{\Set}$, and the inclusion $J: \scat{G}_0 \go
\scat{G}_1$. Then $T_0$ is the identity monad, $T_1$ is the free category
monad, and $R_0=J^*: \cat{G}_1 \go \cat{G}_0$ assigns to a directed graph
its set of objects. Hence $R_0$ is naturally a weak map of monads. (With the
usual description of $T_1$, $R_0$ is in fact a \emph{strict} map of monads,
i.e.\ the isomorphism $T_0 \of R_0 \goiso R_0 \of T_1$ is the identity.)
Moreover, the map $\cat{G}_1^{T_1} \go \cat{G}_0^{T_0}$ induced by this monad
map is the objects functor $S_0$, once one has identified $\cat{G}_1^{T_1}$
with $1\hyph\Cat$ and $\cat{G}_0^{T_0}$ with $0\hyph\Cat$. 

For the inductive step, let $n\geq 1$ and apply
Proposition~\ref{propn:core-map} with $\scat{A} = \scat{G}_n$, $\scat{A'} =
\scat{G}_{n-1}$, the inclusion $J: \scat{G}_{n-1} \go \scat{G}_n$, the monad
$T=T_n$ on $\cat{A} = \cat{G}_n$, and the monad $T'=T_{n-1}$ on $\cat{A'} =
\cat{G}_{n-1}$. Then $J^*=R_{n-1}$, which by inductive hypothesis is a weak
map of monads. This makes $J^*\hyph\Gph$ into a weak map of monads
\[
\pr{\cat{A}\hyph\Gph}{\twid{T}} \go \pr{\cat{A'}\hyph\Gph}{\twid{T'}},
\]
and using the equivalences $\cat{G}_n\hyph\Gph \eqv \cat{G}_{n+1}$,
$\cat{G}_{n-1}\hyph\Gph \eqv \cat{G}_n$, this says that $R_n$ is a weak map
of monads
\[
\pr{\cat{G}_{n+1}}{T_{n+1}} \go \pr{\cat{G}_n}{T_n}.
\]
By the last part of Proposition~\ref{propn:core-map}, the functor from
$\cat{G}_{n+1}^{T_{n+1}}$ ($\eqv (n+1)\hyph\Cat$) to
$\cat{G}_n^{T_n}$ ($\eqv n\hyph\Cat$) induced by this map of monads is indeed
$S_n$.

\section{Representation by Pasting Diagrams}	\label{sec:pasting-rep}

We finish this appendix by showing that if $T=T_\omega$ is the free strict
$\omega$-category monad on $\cat{G}_\omega$, and $X$ a globular set, then
\[
(TX)(n) \iso \coprod_{\pi\in\pd(n)} 
\homset{\cat{G}_\omega}{\rep{\pi}}{X}
\]
for all $n\in\nat$. Really this is just the beginning of a longer story
which is not told here.  Having given concrete descriptions of the globular
set \pd\ and the globular sets $\rep{\pi}$, we could, as hinted on
page~\pageref{p:effort}, go on to specify further data which would determine
the whole monad structure \triple{T}{\eta}{\mu}. Such data would, for
instance, encode the process of composition in the strict $\omega$-category
\pd, i.e.\ the gluing together of pasting diagrams.

By analogy, the Carboni-Johnstone paper \cite{CJ} discusses how a family
$(\rep{\pi})_{\pi\in P}$ of sets gives rise to a cartesian endofunctor $T
= \coprod_{\pi\in P} \homset{\Set}{\rep{\pi}}{\dashbk}$ on \Set, and contains
the result that any cartesian endofunctor on \Set\ arises in this way. (To be
precise, the condition is that $T$ preserves \emph{wide} pullbacks, not just
ordinary pullbacks.) The paper also goes some of the way towards saying what,
in terms of the representing family $(\rep{\pi})_{\pi\in P}$, a cartesian
monad structure on such an endofunctor would be. 

What I envisage is that this theory extends from \Set\ to functor categories
\ftrcat{\scat{A}}{\Set}. This would mean that the free strict
$\omega$-category monad, purely on the grounds of being cartesian, is
familially representable in a suitable sense, and the theory should tell us
what the representing family is---namely, $(\rep{\pi})_{\pi\in\pd(n)}$ for
each $n$, together with the extra data alluded to above. This extended
theory seems to work perfectly well, but the details become so formidable
that an \emph{ad hoc} approach seems more sensible here.

Before proving our result we need some notation. If $X$ is a globular set
then denote by $X^{[n]}$ the $n$-globular set obtained by truncating $X$: in
other words, the image of $X$ under the limit-projection $\cat{G}_\omega \go
\cat{G}_n$. If $Y$ is an $(n+1)$-globular set and $y, y' \in Y(0)$ then there
is an $n$-globular set $Y(y,y')$ given by
\[
(Y(y,y'))(k) = \{ z\in Y(k+1) \such s^k(z)=y, t^k(z)=y' \}.
\]
The same definition can be made when $Y$ is an ($\omega$-)globular set, in
which case $Y(y,y')$ is also an ($\omega$-)globular set. 

Next, let $P_{n+1}: \cat{G}_{n+1} \go \cat{G}_{n+1}$ be the functor given by
\[
(P_{n+1}Y)(0) = Y(0), \diagspace
(P_{n+1}Y)(y,y') = T_n(Y(y,y'))
\]
($y, y' \in Y(0)$), and let $Q_{n+1}: \cat{G}_{n+1} \go \cat{G}_{n+1}$ be the
functor given by 
\begin{eqnarray*}
(Q_{n+1}Y)(0) 		&= 	&Y(0), 					\\
(Q_{n+1}Y)(y,y') 	&= 	&\coprod_{y=y_0, \ldots, y_r=y'} 
	Y(y_0,y_1) \times\cdots\times Y(y_{r-1},y_r).
\end{eqnarray*}
The arguments of the previous section established that $T_{n+1} \iso
Q_{n+1}\of P_{n+1}$.

The proof of the present result is by induction on $n$. First of all, if $X$
is a globular set then 
\[
(T_\omega X)(0) = (T_\omega X)^{[0]}(0) = (T_0 X^{[0]})(0) = X(0) \iso 
\coprod_{\pi\in\pd(0)} \homset{\cat{G}_\omega}{\rep{\pi}}{X},
\]
the second equality coming from the definition of $T_\omega$ as the limit of
the $T_n$'s. Now suppose that the theorem holds for some $n\geq 0$. We have
\begin{eqnarray*}
\lefteqn{(T_\omega X)(n+1)}					\\
&=	&(T_\omega X)^{[n+1]}(n+1)				\\
&=	&(T_{n+1} X^{[n+1]})(n+1)				\\
&\iso	&(Q_{n+1} P_{n+1} X^{[n+1]})(n+1)			\\
&\iso	&\coprod_{x,x' \in X(0)}
	((Q_{n+1} P_{n+1} X^{[n+1]})(x,x'))(n)			\\
&\iso	&\coprod_{x_0, \ldots, x_r \in X(0)}
	(T_n(X^{[n+1]}(x_0,x_1)))(n) \times\cdots\times 
	(T_n(X^{[n+1]}(x_{r-1},x_r)))(n)			\\
&=	&\coprod_{x_0, \ldots, x_r \in X(0)}
	(T_\omega(X(x_0,x_1)))(n) \times\cdots\times
	(T_\omega(X(x_{r-1},x_r)))(n)				\\
&\iso	&\coprod_{\begin{array}{c}\scriptstyle
			x_0, \ldots, x_r \in X(0),\\ \scriptstyle
			\pi_1, \ldots, \pi_r \in \pd(n)
			\end{array}}
	\homset{\cat{G}_\omega}{\rep{\pi_1}}{X(x_0,x_1)}
	\times\cdots\times
	\homset{\cat{G}_\omega}{\rep{\pi_r}}{X(x_{r-1},x_r)}	\\
&\iso	&\coprod_{\pi_1, \ldots, \pi_r \in \pd(n)}
	\homset{\cat{G}_\omega}{\rep{\bftuple{\pi_1}{\pi_r}}}{X}\\
&\iso	&\coprod_{\pi \in \pd(n+1)}
	\homset{\cat{G}_\omega}{\rep{\pi}}{X}
\end{eqnarray*}
where in the penultimate isomorphism we use the construction of
\rep{\bftuple{\pi_1}{\pi_r}} from $\rep{\pi_1}, \ldots, \rep{\pi_r}$
(equation~\bref{eqn:rep-constr}). This completes the induction.

\chapter{Existence of Initial Operad-with-Contraction}	
\label{app:initial}

Here we prove that the category \fcat{OWC} of operads-with-contraction has an
initial object, as required in~\ref{sec:formal}. 

\section{The Strategy}

The explanation in~\ref{sec:the-defn} suggests a way of constructing the
initial operad-with-contraction explicitly: ascend through the dimensions, at
each stage freely adding in elements got by contraction and then freely
adding in elements got by operadic composition. However, we do not take this
route here, instead relying on the following result from Kelly's
paper~\cite{KelUTT}:
\begin{thm}	\label{thm:comb-mon}
Let
\begin{diagram}[size=2em]
\cat{D}	&\rTo	&\cat{C}	\\
\dTo	&	&\dTo>Q		\\
\cat{B}&\rTo_P	&\cat{A}	\\
\end{diagram}
be a (strict) pullback diagram in \fcat{CAT}. If $\cat{A}$ is locally
finitely presentable and each of $P$ and $Q$ is finitary and monadic, then
the functor $\cat{D} \go \cat{A}$ is also monadic.
\done
\end{thm}

All we need to take from this is:
\begin{cor}
In the situation of Theorem~\ref{thm:comb-mon}, \cat{D} has an initial object.
\end{cor}
\begin{proof}
By definition, a locally finitely presentable category is cocomplete, so
\cat{A} has an initial object. The functor $\cat{D} \go \cat{A}$ has a left
adjoint (being monadic), which applied to the initial object of \cat{A} gives
an initial object of \cat{D}. 
\done
\end{proof}

We apply this corollary as follows. Let $T$ be the free strict
$\omega$-category monad on the category $\Eee = \ftrcat{\scat{G}}{\Set}$, as
in Chapter~\ref{ch:defn}. Write \fcat{Coll} for the category $\Eee/\pd$ of
collections (i.e.\ $T$-graphs on $1$: see~\ref{sec:glob-ops}). Write
\fcat{Oper} for the category of $T$-operads; then there is a forgetful
functor $\fcat{Oper} \go \fcat{Coll}$.  As observed on page~\pageref{p:entity},
the definition of a contraction on a $T$-operad is really a definition of a
contraction on a collection, which means that we have a category \fcat{CWC}
of collections-with-contraction and a (strict) pullback diagram
\begin{diagram}
\fcat{OWC}	&\rTo	&\fcat{Oper}	\\
\dTo		&	&\dTo		\\
\fcat{CWC}	&\rTo	&\fcat{Coll}	\\
\end{diagram}
in \fcat{CAT}.

All we need to do now is check that the hypotheses of
Theorem~\ref{thm:comb-mon} hold in this situation, and that is the content of
the next section.

\section{The Proof}

\paragraph{Hypothesis on \fcat{Coll}} We have first to check that \fcat{Coll}
is locally finitely presentable. Indeed, if $\mr{Gr}(\pd)$ is the Grothendieck
fibration (category of elements) of the functor $\pd: \scat{G} \go \Set$,
then
\[
\fcat{Coll} \iso \ftrcat{\scat{G}}{\Set}/\pd \eqv \ftrcat{\mr{Gr}(\pd)}{\Set},
\]
and any category of the form \ftrcat{\scat{A}}{\Set} (with \scat{A} small) is
locally finitely presentable: see \cite{BorxII}, Example 5.2.2(b).

\paragraph{Hypotheses on $U: \fcat{CWC} \protect\rTo \fcat{Coll}$} We have to
see that $U$ is finitary and monadic. It is straightforward to calculate that
$U$ creates filtered colimits; and since \fcat{Coll} possesses all filtered
colimits, $U$ preserves them too. It is also easy to calculate that $U$
creates coequalizers for $U$-split coequalizer pairs. Hence we have only to
show that $U$ has a left adjoint.

We construct a left adjoint $F$ explicitly. Let $C$ be a collection, and
define a new collection $FC$ and a map $\alpha_C: C \go FC$ inductively as
follows: 
\begin{itemize}
\item if $\pi\in\pd(0)$ then $(FC)(\pi) = C(\pi)$
\item if $\pi\in\pd(1)$ then
 $(FC)(\pi) = C(\pi) + (C(\bdry\pi) \times C(\bdry\pi))$ 
\item if $n\geq 2$ and $\pi\in\pd(n)$ then
\begin{eqnarray}	\label{eq:FC}
\lefteqn{(FC)(\pi) =} \\
&&C(\pi) + \{\pr{\psi_0}{\psi_1} \in (FC)(\bdry\pi)^2 \such 
s(\psi_0) = s(\psi_1) \mbox{ and } t(\psi_0) = t(\psi_1)\} \nonumber
\end{eqnarray}
\item $\alpha_{C,\pi}: C(\pi) \rIncl (FC)(\pi)$ is inclusion as the first
component, for all $\pi$
\item if $n\geq 1$ and $\pi\in\pd(n)$ then the source map $s: (FC)(\pi)
\go (FC)(\bdry\pi)$ is given by
\begin{itemize}
\item the composite $C(\pi) \goby{s} C(\bdry\pi) \goby{\alpha_{C,\bdry\pi}}
(FC)(\bdry\pi)$, on the first summand
\item first projection, on the second summand,
\end{itemize}
and the target map is defined similarly.
\end{itemize}

It is easy to check that the globularity relations in $FC$ are 
satisfied, so that $FC$ forms a collection, and that $\alpha_C: C \go FC$ is
a map of collections.

In the notation of \ref{sec:formal}, the set $\{\ldots\}$ in
equation~\bref{eq:FC} is $P_\pi(FC)$, so
\[
(FC)(\pi) = C(\pi) + P_\pi(FC)
\]
for any $n\geq 1$ and $\pi\in\pd(n)$. Thus we can define a contraction
$\kappa^C$ on $FC$ by taking $\kappa^C_\pi$ to be second inclusion
$P_\pi(FC) \rIncl (FC)(\pi)$. 

We have now associated to each collection $C$ a
collection-with-con\-trac\-tion \pr{FC}{\kappa^C} and a map $\alpha_C: C
\go FC$ of collections. Another easy check shows that $\alpha_C$ has the
appropriate universal property, so that $U$ has a left adjoint.

\paragraph{Hypotheses on $\fcat{Oper} \protect\rTo \fcat{Coll}$}  Again, we
have to see that this functor is finitary and monadic.

Monadicity will follow from Theorems~\ref{thm:free-fixed}
and~\ref{thm:free-gen} just as long as $T$ is finitary, which is true by
Theorem~\ref{thm:strict-omega}.

Let $T'_1$ be the monad on \fcat{Coll} induced by $\fcat{Oper} \go
\fcat{Coll}$ and its left adjoint. The fact that $T$ is finitary also implies
that $T'_1$ is finitary, by the `moreover' part of
Theorem~\ref{thm:free-fixed}. So our monadic adjunction is finitary, as
required.

\backmatter


\begin{thebibliography}{MMMM}		%\raggedright
\ucontents{chapter}{Bibliography} 

\small

\bibitem[$\!$]{ThisIsADummyValue}
\hspace{-16.5mm}\rule[-1mm]{1.3mm}{3.7mm}\hspace{15.2mm}%
\textbf{%
E-print numbers such as \url{math.CT/9810053} and \url{q-alg/9706008}
refer to the electronic mathematics archive, at
\url{http://arXiv.org}. For direct access to a paper, go
to the address of the form
\url{http://arXiv.org/abs/math.CT/9810053}.}


\bibitem[BD]{BaDoHDA3}
J. Baez, J. Dolan,
Higher-dimensional algebra III: $n$-categories and the algebra of opetopes
(1998). 
\emph{Advances in Mathematics} 135, pp.~145--206.
Also available via \url{http://math.ucr.edu/home/baez}. 

\bibitem[Bat]{Bat}
M. Batanin,
Monoidal globular categories as a natural environment for the theory of weak
$n$-categories (1997). 
\emph{Advances in Mathematics} 136, pp.~39--103. 

\bibitem[BJT]{BJT}
Hans-Joachim Baues, Mamuka Jibladze, Andy Tonks,
Cohomology of monoids in monoidal categories (1997).
In: \emph{Operads: Proceedings of Renaissance Conferences}, ed.\ Loday,
Stasheff, Voronov, Contemporary Mathematics 202, AMS. 

\bibitem[B\'{e}n]{Ben}
Jean B\'{e}nabou,
Introduction to bicategories (1967).
In: \emph{Reports of the Midwest Category Seminar}, ed.\ B\'{e}nabou et
al, Springer LNM 47. 

\bibitem[BCSW]{BCSW}
Renato Betti, Aurelio Carboni, Ross Street, Robert Walters,
Variation through enrichment (1983).
\emph{Journal of Pure and Applied Algebra} 29, pp.~109--127.

\bibitem[BKP]{BKP}
R. Blackwell, G. M. Kelly, A. J. Power,
Two-dimensional monad theory (1989).
\emph{Journal of Pure and Applied Algebra} 59, 1--41.

\bibitem[Borx1]{BorxI}
Francis Borceux,
\emph{Handbook of Categorical Algebra 1: Basic Category Theory}
(1994).
Cambridge University Press.

\bibitem[Borx2]{BorxII}
Francis Borceux,
\emph{Handbook of Categorical Algebra 2: Categories and Structures}
(1994).
Cambridge University Press.

\bibitem[Borh]{Borh}
R. Borcherds, 
Vertex algebras (1997). 
E-print \url{q-alg/9706008}. 

\bibitem[Bur]{Bur}
A. Burroni,
$T$-cat\'{e}gories (1971). 
\emph{Cahiers Top.\ Geom.\ Diff.}, Vol.~XII, No.~3, pp.~215--321. 

\bibitem[CJ]{CJ}
A. Carboni, P. Johnstone,
Connected limits, familial representability and Artin glueing (1995).
\emph{Mathematical Structures in Computer Science}, Vol.~5, pp.~441--459.

\bibitem[CKW]{CKW}
Aurelio Carboni, Stefano Kasangian, Robert Walters,
An axiomatics for bicategories of modules (1987).
\emph{Journal of Pure and Applied Algebra} 45, pp.~127--141.

\bibitem[CKVW]{CKVW}
A. Carboni, G. M. Kelly, D. Verity, R. J. Wood,
A 2-categorical approach to change of base and geometric morphisms II
(1998). 
\emph{Theory and Applications of Categories}, Vol.~4, No.~5, pp.~82--136. 

\bibitem[Che1]{CheROM}
Eugenia Cheng,
The relationship between the opetopic and multitopic approaches to weak
$n$-categories (2000).
Available via \url{http://www.dpmms.cam.ac.uk/$\sim$elgc2}.

\bibitem[Che2]{CheEAT}
Eugenia Cheng,
Equivalence between approaches to the theory of opetopes (2000).
Available via \url{http://www.dpmms.cam.ac.uk/$\sim$elgc2}.

\bibitem[Her1]{HerCSU}
Claudio Hermida,
From coherent structures to universal properties (2000).
E-print \url{math.CT/0006161}. 

\bibitem[Her2]{HerRM}
Claudio Hermida,
Representable multicategories (2000).
\emph{Advances in Mathematics} 151, pp.~164--225. Also available via
\texttt{http://www.cs.math.ist.utl.pt/s84.www/cs/claudio.html}.

\bibitem[HMP]{HMP}
C. Hermida, M. Makkai, J. Power,
On weak higher dimensional categories (1997).
Preprint.

\bibitem[JS]{JS}
Andr\'{e} Joyal, Ross Street,
Braided tensor categories (1993).
\emph{Advances in Mathematics} 102, pp.~20--78.

\bibitem[Kel1]{KelCD}			
G. M. Kelly,
On clubs and doctrines (1974).
In: \emph{Category Seminar}, Springer LNM 420, pp.~181--256.

\bibitem[Kel2]{KelUTT}
G. M. Kelly,
A unified treatment of transfinite constructions for free algebras, free
monoids, colimits, associated sheaves, and so on (1980).
\emph{Bulletin of the Australian Mathematical Society}, Vol.~22, pp.~1--83.

\bibitem[KS]{KS}			
G. M. Kelly, R. Street,
Review of the elements of 2-categories (1974).
In: \emph{Category Seminar}, Springer LNM 420, pp.~75--103.

\bibitem[Kos]{Kos}
J\"{u}rgen Koslowski,
Monads and interpolads in bicategories (1997).
\emph{Theory and Applications of Categories}, Vol.~3, No.~8, pp.~182--212.

\bibitem[Lam]{Lam}
Joachim Lambek,
Deductive systems and categories II: standard constructions and closed
categories (1969).
In: \emph{Category Theory, Homology Theory and their Applications I}, ed.\
P. Hilton, Springer LNM 86.

\bibitem[Lei1]{GOM}
Tom Leinster,
General operads and multicategories (1997).
E-print \url{math.CT/9810053}.

\bibitem[Lei2]{BB}
Tom Leinster,
Basic bicategories (1998).
E-print \url{math.CT/9810017}.

\bibitem[Lei3]{SHDCT}
Tom Leinster,
Structures in higher-dimensional category theory (1998).
Available via \url{http://www.dpmms.cam.ac.uk/$\sim$leinster}.

\bibitem[Lei4]{FCM}
Tom Leinster,
\fc-multicategories (1999).
E-print \url{math.CT/9903004}.

\bibitem[Lei5]{GECM}
Tom Leinster,
Generalized enrichment for categories and multicategories (1999).
E-print \url{math.CT/9901139}.

\bibitem[Lei6]{GEC}
Tom Leinster,
Generalized enrichment of categories (1999).
To appear in \emph{Journal of Pure and Applied Algebra}.

\bibitem[Lei7]{HAO}
Tom Leinster,
Homotopy algebras for operads (2000).
E-print \url{math.QA/0002180}.

\bibitem[Lei8]{WMC}
Tom Leinster,
What's a monoidal category?\ (2000).
Poster at CT2000, Como, Italy.

\bibitem[Lew]{Lew}
Geoffrey Lewis, 
Coherence for a closed functor (1972). 
In: \emph{Coherence in Categories}, ed. Kelly, Laplaza, Lewis, Mac Lane,
Springer LNM 281. 

\bibitem[May1]{MayGIL}
J. P. May,
\emph{The Geometry of Iterated Loop Spaces} (1972).
Springer LNM 271.

\bibitem[May2]{MayDOA}
J. P. May,
Definitions: operads, algebras and modules (1997). 
In: \emph{Operads: Proceedings of Renaissance Conferences}, ed.\ Loday,
Stasheff, Voronov, Contemporary Mathematics 202, AMS. 

\bibitem[Pow]{Pow}
A. J. Power,
A general coherence result (1989).
\emph{Journal of Pure and Applied Algebra} 57, no.~2, 165--173.

\bibitem[Sny1]{SnyEBG}
Craig T. Snydal,
Equivalence of Borcherds $G$-vertex algebras and axiomatic vertex algebras
(1999).
E-print \url{math.QA/9904104}.

\bibitem[Sny2]{SnyRMS}
Craig T. Snydal,
Relaxed multi category structure of a global category of rings and modules
(1999).
E-print \url{math.QA/9912075}.

\bibitem[Soi]{Soi}
Y. Soibelman,
Meromorphic tensor categories (1997).
E-print \url{q-alg/9709030}. 

\bibitem[Str1]{StrFTM}
Ross Street,
The formal theory of monads (1972).
\emph{Journal of Pure and Applied Algebra} 2, pp.~149--168.

\bibitem[Str2]{StrCS}
Ross Street,
Categorical structures (1995).
In: \emph{Handbook of Algebra,}\/ Vol.~1, ed.\ M. Hazewinkel, Elsevier
North-Holland. 

\bibitem[Str3]{StrRMB}
Ross Street,
The role of Michael Batanin's monoidal globular categories (1997).
Available via \url{http://www.mpce.mq.edu.au/$\sim$street}.

\bibitem[Wal]{Wal}
R. F. C. Walters,
Sheaves and Cauchy-complete categories (1981).
\emph{Cahiers Top.\ Geom.\ Diff.}, Vol.~XXII, No.~3. 



\end{thebibliography}
\end{document}